\documentclass[10pt,notitlepage]{amsart} 
\usepackage{amsmath}
\usepackage{amsthm} 
\usepackage{amsfonts}
\usepackage{scalerel}
\usepackage{amssymb,amscd,epsf,verbatim,mathtools,framed}
\usepackage{mathrsfs}
\usepackage{graphicx} 
\usepackage{latexsym}
\usepackage{lscape}
\usepackage[colorlinks=true]{hyperref}
\hypersetup{colorlinks, citecolor=blue, filecolor=black, linkcolor=purple, urlcolor=violet}
\usepackage{epstopdf}
\usepackage{pdflscape}
 
\usepackage{tikz}
\usepackage[all]{xy}
\usetikzlibrary{calc}
\usetikzlibrary{matrix,arrows,decorations.pathmorphing}
\usepackage{tabularx}
\usepackage{tikz-cd}
\usepackage{color}
\usepackage{multirow}  
\usepackage{scalefnt}
\usepackage{fancyhdr}
\usepackage[margin=1in]{geometry}
\usepackage[T1]{fontenc}

\usepackage{relsize} 
\usepackage[bbgreekl]{mathbbol} 
\DeclareSymbolFontAlphabet{\mathbb}{AMSb} 
\DeclareSymbolFontAlphabet{\mathbbl}{bbold}

\newcommand{\Z}{\mathbb{Z}}
\newcommand{\D}{\mathbb{D}}
\newcommand{\F}{\mathbb{F}}

\newcommand{\R}{\mathbb{R}}
\newcommand{\Q}{\mathbb{Q}}

\renewcommand{\L}{\mathbb{L}}

\newcommand{\C}{\mathbb{C}}
\newcommand{\T}{\mathbb{T}}

\newcommand{\mO}{\mathcal{O}}

\newcommand{\mX}{\mathcal{X}}

\newcommand{\fA}{\mathbf{A}}
\newcommand{\fB}{\mathbf{B}}
\newcommand{\fC}{\mathbf{C}}
\newcommand{\fN}{\mathbf{N}}

\newcommand{\tu}{\textup}
\newcommand{\cl}{\overline}

\renewcommand{\le}{\leqslant}
\renewcommand{\ge}{\geqslant}
\newcommand{\ra}{\rightarrow}

\newcommand{\sq}{\widetilde}

\DeclareMathOperator{\im}{im}
\DeclareMathOperator{\Hom}{Hom}

\DeclareMathOperator{\GL}{GL}

\DeclareMathOperator{\lra}{\: \longrightarrow \:} 
\DeclareMathOperator{\isom}{\;\xrightarrow{\: {}_{\sim} \:} \;} 
\DeclareMathOperator{\Gal}{Gal}

\newcommand{\gr}[1]{\langle {#1} \rangle} 
\DeclareMathOperator{\spec}{Spec}
\DeclareMathOperator{\spf}{Spf}
\DeclareMathOperator{\spa}{Spa}
\DeclareMathOperator{\Spa}{Spa}
 
\DeclareMathOperator{\ett}{\mathrm{\acute{e}t}}

\DeclareMathOperator{\proet}{\mathrm{pro\acute{e}t}}

\theoremstyle{theorem}
\newtheorem{theorem}{Theorem}[section] 
\newtheorem{mainthm}{Theorem}

\theoremstyle{definition}
\newtheorem{mainrmk}[mainthm]{Remark}
\newtheorem{example}[theorem]{Example}
\newtheorem{definition}[theorem]{Definition}
\newtheorem{proposition}[theorem]{Proposition}

\newtheorem{lemma}[theorem]{Lemma}
\newtheorem{corollary}[theorem]{Corollary}
\newtheorem{construction}[theorem]{Construction} 

\newtheorem{notation}[theorem]{Notation}
\newtheorem{remark}[theorem]{Remark}

\newtheorem{assumption}[theorem]{Assumption}

\newtheorem{setup}[theorem]{Set-up}

\makeatother
\title{Relative $(\varphi, \Gamma)$-modules and $p$-adic differential equations}

\author{Hansheng Diao, Yong Suk Moon, and Zijian Yao}

\date{}

\numberwithin{equation}{section}

\begin{document}

\maketitle
\thispagestyle{empty}

\begin{abstract} 
Let $X$ be an affinoid rigid analytic space over a $p$-adic field, equipped with a suitable \'etale map to a unit polydisk. We provide a formalism of \emph{imperfect relative period rings} over $X$, which lie inside the corresponding perfect relative period rings constructed by Kedlaya--Liu. Then we establish a relative version of the Fontaine--Cherbonnier--Colmez equivalence between $p$-adic local systems on $X$ and \'etale $(\varphi, \Gamma)$-modules over such imperfect relative period rings, generalizing previous works of Andreatta--Brinon and others. Using this equivalence, we construct $p$-adic differential equations attached to de Rham local systems, carrying both geometric and arithmetic differential operators. This generalizes the   work of Berger to the relative geometric setting. Along the way, we study a relative Fontaine--Sen theory on the decompletion of $\Gamma$-modules over the relative $\mathbf{B}_{\mathrm{dR}}^+$-period rings. 
\end{abstract}

\vspace{0.2in} 

\tableofcontents

\section{Introduction}\label{section: introduction}
 
\subsection{Main results} 

Let $K$ be a $p$-adic field, that is, a complete discrete valuation field of mixed characteristic $(0,p)$ with perfect residue field. Let $K_{\infty}$ be the cyclotomic extension of $K$. Let $\Gal_K$ denote the absolute Galois group of $K$ and write $\Gamma_K=\Gal(K_{\infty}/K)$.

A central achievement in classical $p$-adic Hodge theory is the conversion of $p$-adic Galois representations into (semi-)linear algebraic data of $(\varphi, \Gamma)$-modules over suitable period rings. The following are two major breakthroughs in this direction.
\begin{enumerate}
\item[(a)] Fontaine's construction \cite{Fontaine_90}, together with the overconvergence theorem of Cherbonnier--Colmez \cite{Cherbonnier_Colmez}, establishes an equivalence between the category of $p$-adic representations of $\Gal_K$ and that of \'etale $(\varphi, \Gamma_K)$-modules over certain \textit{overconvergent period rings}. 
\item[(b)] For a de Rham representation $V$ of $\Gal_K$, Berger \cite{Berger-differential} constructs a canonical modification $\mathbf{N}_{\mathrm{dR}}(V)$ of the aforementioned $(\varphi, \Gamma_K)$-module. It carries a differential operator induced from the infinitesimal action of $\Gamma_K$; in other words, $\mathbf{N}_{\mathrm{dR}}(V)$ is a \emph{$p$-adic differential equation}. This construction allows one to recover the crystalline and semistable periods of $V$ from the associated $p$-adic differential equation, and forms the cornerstone of the proof of the $p$-adic monodromy theorem \cite{Berger-differential, Andre_p-adic_monodormy, Kedlaya_p-adic_monodormy, Mebkhout_p-adic_monodormy}. 
\end{enumerate}

The goal of this article is to generalize both (a) and (b) to the relative setting by constructing functors
\begin{align*}
    \big\{p\text{-adic local systems} \big\} & \isom    \big\{\text{relative \'etale }  (\varphi, \Gamma)\text{-modules} \big\} \\
     \big\{\text{de Rham local systems} \big\} & \lra    \big\{\text{relative }  p\text{-adic differential equations} \big\}.
\end{align*}
This provides a crucial missing component in relative $p$-adic Hodge theory, and plays an essential role in the recent work of the authors on a relative $p$-adic monodromy theorem for smooth projective curves (see \cite{p_adic_monodromy_curve}).  

We now explain the precise statements. Let $\mathbb D_K^{d}$ denote the $d$-dimensional closed unit polydisk over $K$ with its coordinates denoted by $T_1, \ldots, T_d$. Let $X$ be a smooth affinoid rigid analytic space over $K$ equipped with an \'etale morphism $h\colon X \ra \mathbb D_K^{d}$ that is a composition of rational localizations and finite \'etale maps such that $h$ is away from the locus $T_1\cdots T_d = 0$; such an $h$ is referred to as a \emph{toric chart}. Consider the standard pro-finite Galois tower $ \mathbb D_{K, \infty}^{d} \ra \mathbb D_K^{d}$ obtained by extracting all $p$-power roots of unity and all $p$-power roots of $T_1, \ldots, T_d$. The Galois group $\Gamma$ of this tower can be viewed as the relative analogue of $\Gamma_K$ in the classical setting. Taking base change along $h$, we arrive at a pro-finite Galois tower $X_\infty \ra X$, still with Galois group $\Gamma$. 

\vspace{0.05in}

Our first task is to describe \'etale $\Z_p$-local systems on $X$ in terms of relative \'etale $(\varphi, \Gamma)$-modules over certain \textit{imperfect relative period rings}. To state this result, let $\mathbf{B}_{\Lambda}^{\dagger}$ and $\mathbf{C}_{\Lambda}^{\dagger}$ be the imperfect relative period rings associated with a decompleting datum $\Lambda$ for $(X,h)$. We shall return to the notion of ``imperfect relative period rings'' and ``decompleting data'' momentarily. Let us just remark that $\mathbf{B}_{\Lambda}^{\dagger}$ contains the classical period ring $\mathbf{B}_{K}^{\dagger}$ studied in \cite{Cherbonnier_Colmez, Berger-differential} as a subring. Similarly, $\mathbf{C}_{\Lambda}^{\dagger}$ contains the classical Robba ring $\mathbf{C}^{\dagger}_K$, which is denoted by $\mathbf{B}^{\dagger}_{\mathrm{rig}, K}$ in \cite{Berger-differential}. The following theorem is a relative version of the Fontaine--Cherbonnier--Colmez equivalence; see Theorem~\ref{thm: equivalence-isogeny-local-system-overconvergent-etale-(phi, Gamma)-modules-C} for the full statement. 

\begin{mainthm} \label{intro:theorem_equiv_of_categories}
There are natural $\otimes$-equivalences 
\[
D^{\dagger}\colon \tu{Loc}_{\Z_p} (X)_{\Q} \isom \mathrm{Mod}^{(\varphi, \Gamma), \ett}_{/\mathbf{B}_{\Lambda}^{\dagger}} 
\]
\[
D^{\dagger}_{\mathrm{rig}}\colon \tu{Loc}_{\Z_p} (X)_{\Q} \isom \mathrm{Mod}^{(\varphi, \Gamma), \ett}_{/\mathbf{C}_{\Lambda}^{\dagger}} 
\]
from the isogeny category of \'etale $\Z_p$-local systems on $X$ to the category of finite projective \'etale $(\varphi, \Gamma)$-modules over $\mathbf{B}_{\Lambda}^{\dagger}$ (resp. $\mathbf{C}_{\Lambda}^{\dagger}$). 
\end{mainthm}

\vspace{0.05in}

When the local system $\L$ is \textit{de Rham} in the sense of \cite{Scholze_p_adic_Hodge}, we construct another relative $(\varphi, \Gamma)$-module $\mathbf{N}_{\mathrm{dR}}(\L)$ from $D^{\dagger}_{\mathrm{rig}}(\L)$ by a modification process, carrying not only the arithmetic differential operator arising from the cyclotomic direction, but also differential operators from the geometric directions. This appears to be the first construction of a relative version of 
Berger's $p$-adic differential equation $\mathbf{N}_{\mathrm{dR}}(V)$ in the literature. To state our result, we assume that $\Lambda$ is in addition \emph{strongly decompleting} --- a notion we shall explain later in this introduction --- and let $\mathbf{C}^\dagger_{\Lambda}$ be the associated imperfect relative period ring.

\begin{mainthm} \label{intro:theorem_relative_Berger}
There is a natural $\otimes$-functor  
\[
\mathbf{N}_{\mathrm{dR}}\colon \mathrm{Loc}^{\mathrm{dR}}_{\Z_p}(X)_{\Q} \lra \mathrm{Mod}^{(\varphi, \Gamma)}_{/\mathbf{C}^{\dagger}_{\Lambda}}
\]
from the isogeny category of de Rham $\Z_p$-local systems to the category of finite projective $(\varphi, \Gamma)$-modules over $\mathbf{C}^{\dagger}_{\Lambda}$, such that the following properties hold for every $\L\in \mathrm{Loc}^{\mathrm{dR}}_{\Z_p}(X)_{\Q}$:
\begin{enumerate}
\item $\mathbf{N}_{\mathrm{dR}} (\L)[1/t]$ coincides with $D^{\dagger}_{\mathrm{rig}} (\L) [1/t]$ after inverting Fontaine's period $t$;

\item $\mathbf{N}_{\mathrm{dR}}(\L)$ is equipped with commuting differential operators $D_0, D_1, \ldots, D_d$, making $(\mathbf{N}_{\mathrm{dR}}(\L),\, D_i)$ a differential module over the differential ring $(\mathbf{C}^{\dagger}_{\Lambda},\, d_i)$ for all $i=0, \ldots, d$, where $d_i$ is the ``normalized'' differential operator on $\mathbf{C}^{\dagger}_{\Lambda}$ defined in Definition~\ref{defn: normalized differential operators}.
\end{enumerate}
\end{mainthm} 

\begin{mainrmk}
In the case where $X$ is a closed poly-annulus over $K$, the imperfect relative period ring $\mathbf{C}^{\dagger}_{\Lambda}$ consists of formal power series in $(d+1)$ variables $\pi_K$, $[T_1^{\flat}], \ldots, [T_d^{\flat}]$, subject to some convergence conditions (cf. \S \ref{subsec: closed-annuli}); it can be viewed as a ``relative Robba ring''. Here, $[T_i^\flat]$ are elements in $\mathbf{C}_{\Lambda}^{\dagger}$ that are Teichm\"{u}ller lifts of $T_i^\flat$ (determined by a compatible system of $p$-power roots of $T_i$), while $\pi_K$ is the variable in the classical Robba ring $\mathbf{C}^{\dagger}_K$ (see Construction~\ref{construction:pi_K}). In this case, the normalized differential operators on $\mathbf{C}^{\dagger}_{\Lambda}$ are $d_0=\frac{d}{d\pi_K}$ and $d_i=\frac{d}{d[T^{\flat}_i]}$ for $i=1, \ldots, d$.
\end{mainrmk}

\subsection{Related works}

To orient the reader, let us mention some relevant results in the literature. 
The theory of relative de Rham period rings, relative Tate--Sen theory, and (prototypes of)  relative $(\varphi, \Gamma)$-modules were studied by Andreatta and Brinon in a series of works \cite{Brinon, Andreatta-Brinon_Overconvergence, Andreatta-Brinon}. These works mostly center around the case of good reduction, that is, when $X$ admits a smooth formal model over $\spf \mO_K$. 
In a more general setting, Kedlaya and Liu \cite{Kedlaya-LiuI} studied $(\varphi, \Gamma)$-modules over \emph{perfect relative period rings} (cf. \S \ref{sec:relative_period_rings:perfect}) via the formalism of perfectoid spaces. In particular, they show that $p$-adic local systems give rise to \'etale $(\varphi, \Gamma)$-modules over these perfect  period rings (for example, see \S \ref{subsec: (phi, gamma)-modules-over-perfect-period-rings}). The arithmetic information one ultimately seeks, however, are encoded over the \emph{imperfect} period rings. In the classical setting, passage from the perfect period rings to the imperfect ones is the procedure known as \emph{decompletion}.

In a subsequent preprint \cite{Kedlaya-LiuII} (also see Kedlaya  \cite{Kedlaya-new-methods-phi-Gamma-modules}), 
a decompletion formalism is proposed over certain imperfect relative period rings, including a relative analogue of the overconvergence theorem of Cherbonnier--Colmez similar to Theorem~\ref{intro:theorem_equiv_of_categories}.
However, despite our best effort, some of the descent assertions claimed in \cite[in particular, \S 5]{Kedlaya-LiuII} necessary for this purpose seem to be incomplete and require significant further justification. In particular, it is unclear to us whether the proposed definition of ``$\mathbf{A}_{\psi, X}^r$'' in \cite[Definition~5.2.1]{Kedlaya-LiuII} by 
\[
\mathbf{A}_{\psi, X}^r\coloneqq \{ x \in \widetilde{\mathbf{A}}^r_{\psi,X} ~|~ \theta\circ \varphi^{-n}(x)\in A_{\psi, X, n}\textrm{ for all }n\ge -\log_p r\}
\]
gives suitable imperfect period rings in general, in order to decomplete relative $(\varphi, \Gamma)$-modules. Here $\psi$ stands for the perfectoid tower $X_{\infty}\rightarrow X$ as in \cite{Kedlaya-LiuII}, and $\widetilde{\mathbf{A}}^r_{\psi, X}$ is the perfect relative period ring for the tower $\psi$ constructed in \cite{Kedlaya-LiuI} (cf. \S \ref{sec:relative_period_rings:perfect}).

As it is rather unclear to us how to verify the claims in \cite{Kedlaya-LiuII}, we take a slightly different approach. One innovation of this article is to axiomatize the formal properties of imperfect period rings for the decompletion mechanism to work. These correspond to the terms \emph{decompleting data} and \emph{strongly decompleting data} in Theorem~\ref{intro:theorem_equiv_of_categories} and Theorem~\ref{intro:theorem_relative_Berger}, respectively. The axioms we impose on the imperfect period rings are inspired by, but are \textit{a priori} more restrictive than, the ones from \cite{Kedlaya-LiuII}. In the rest of the introduction, we will briefly describe these notions, as well as other ingredients in the proof of the main results.

\subsection{Axiomatization of imperfect period rings}

As mentioned earlier, in order to have a coherent theory to work with, our first goal is to initiate a systematic study of imperfect period rings. A central notion we introduce in this article is that of \emph{(weakly) decompleting data} (Definitions~\ref{defn: weakly decompleting datum} and \ref{definition:decompleting_datum}), which is engineered to axiomatize the idea of (weakly) decompleting towers introduced in \cite[Definition~5.2.3]{Kedlaya-LiuII}. 

A \emph{weakly decompleting datum} consists of a pair $\Lambda=(R_{\Lambda}, \mathbf{A}_{\Lambda})$, where $R_{\Lambda}$ is a subring of the tilt $\widehat A_{X, \infty}^{\flat}$  of the perfectoid algebra $\widehat A_{X, \infty}$ induced from perfectoid tower $X_{\infty}\rightarrow X$, while $\mathbf{A}_{\Lambda}$ is a subring of the Witt vectors $W(\widehat{A}_{X, \infty}^{\flat})$. They are required to satisfy a list of axioms; in particular, they acquire $(\varphi, \Gamma)$-actions. Intuitively, $R_{\Lambda}$ and $\mathbf{A}_{\Lambda}$ play the roles of $\mathbf{E}_K$ and $\mathbf{A}_K$ in the classical setting (for example, see \cite{Berger-differential}). Since most of the axioms are technical, we refer the reader to Definition~\ref{defn: weakly decompleting datum} for the precise definition. From such a pair $(R_{\Lambda}, \mathbf{A}_{\Lambda})$, one can further build imperfect relative period rings
\[
\mathbf{A}^r_{\Lambda}, \mathbf{A}^{\dagger}_{\Lambda}, \mathbf{B}^r_{\Lambda},
\mathbf{B}^{\dagger}_{\Lambda}, \mathbf{C}^{[s,r]}_{\Lambda}, \mathbf{C}^{r}_{\Lambda},  \mathbf{C}^{\dagger}_{\Lambda}
\]
for $r>0$ sufficiently small and for $0 < s \le r$, all equipped with suitable $(\varphi, \Gamma)$-actions. These form relative analogues of the imperfect period rings $\mathbf{A}^r_{K}, \mathbf{A}^{\dagger}_{K},  \mathbf{B}^r_{K}, 
\mathbf{B}^{\dagger}_{K}, \mathbf{C}^{[s,r]}_K, \mathbf{C}^r_K, \mathbf{C}^{\dagger}_K$  
in the classical setting.\footnote{We remark that $\mathbf{A}^r_{K}$, $\mathbf{B}^r_{K}$, $\mathbf{C}^r_K$ and $\mathbf{C}^{\dagger}_K$ are denoted by $\mathbf{A}^{\dagger, r'}_{K}$, $\mathbf{B}^{\dagger, r'}_{K}$, $\mathbf{B}^{\dagger, r'}_{\mathrm{rig}, K}$, and $\mathbf{B}^{\dagger}_{\mathrm{rig}, K}$ in \cite{Berger-differential} with $r'=(p-1)/pr$.} The axioms also ensure that there is a robust theory of finite projective (\'etale) $(\varphi, \Gamma)$-modules over these imperfect relative period rings.

For the purpose of decompletion and descent of $(\varphi, \Gamma)$-modules, we need to impose additional conditions. We axiomatize these conditions via the notion of \emph{decompleting data} (Definition~\ref{definition:decompleting_datum}). Specifically, we require
\begin{equation}\label{eq: vanishing of continuous group cohomology}
H^i_{\mathrm{cont}}(\Gamma, \varphi^{-1}(R_{\Lambda})/R_{\Lambda}) = 0
\end{equation}
for $i=0,1$, where $\varphi^{-1} (R_{\Lambda})$ stands for the pre-image of $R_{\Lambda}$ under the Frobenius map taken inside $\widehat{A}^{\flat}_{X,\infty}$. In fact, the precise definition of decompleting data asks for more: (i) the group cochain complex that computes the continuous group cohomology is required to be  \emph{strict exact} instead of simply being exact, which incorporates some subtle topological constraints; (ii) we ask that these conditions hold for all open normal subgroups of $\Gamma$, and for all finite \'etale covers of $X$.

These conditions allow us to use the continuous $\Gamma$-actions to decomplete $(\varphi, \Gamma)$-modules from perfect relative period rings to the imperfect ones. Indeed,  decompletion boils down to the vanishing of certain continuous group cohomology with coefficients in $\mathbf{A}^{\dagger}_{\Lambda} $ and $\mathbf{C}^{\dagger}_{\Lambda}$. A key observation is that these cohomological vanishings in characteristic 0 can be deduced from vanishings in characteristic $p$ (such as Condition \eqref{eq: vanishing of continuous group cohomology}) through systematic norm-compatible liftings of cochains ---  we call them ``good lifts'', see Corollary~\ref{cor: good lift} --- together with delicate $p$-adic approximation arguments. Combined with a further bootstrapping argument (which relies on a series of intricate d\'evissage steps), we establish the desired cohomological vanishings in Corollary~\ref{cor: complex-A-strict-exact} and Proposition~\ref{prop: C-complex-strict-exact}. Let us emphasize that these norm-compatible liftings and the bootstrapping procedure are by no means formal; they constitute the technical core of \S \ref{sec:decompleting_data}. These results eventually yield a collection of decompletion theorems in \S \ref{subsec: fully-faithfulness-indfree-(phi, Gamma)-modules} and \S \ref{subsec: descent-relative-etale-(phi, Gamma)-modules}, including Theorem~\ref{intro:theorem_equiv_of_categories}.

\subsection{Relative $p$-adic differential equations}

Before proceeding further, let us briefly recall the construction of the $p$-adic differential equation $\mathbf{N}_{\mathrm{dR}}(V)$ in the classical setting. Recall the classical Robba ring $\mathbf{C}^{\dagger}_K$ with coordinate in $\pi_K$. For any $p$-adic de Rham representation $V$ of $\Gal_K$, one can associate an \'etale $(\varphi, \Gamma_K)$-module $D^{\dagger}_{\mathrm{rig},K}(V)$ over $\mathbf{C}^{\dagger}_K$ via the Fontaine--Cherbonnier--Colmez equivalence. Since the differential operator $\nabla_0$ induced by the infinitesimal action of $\Gamma_K$ on $\mathbf{C}^{\dagger}_K$ differs from the standard differential operator $d/d\pi_K$ by a unit multiple of $t$, the ``correct differential operator'' on $D^{\dagger}_{\mathrm{rig},K} (V) $ should be the one induced by $\mathrm{Lie} (\Gamma_K)$ divided by $t$. In general, this introduces simple poles along vanishing locus of $t$. When $V$ is de Rham, the idea of Berger is to construct $\mathbf{N}_{\mathrm{dR}} (V)$ by systematically ``removing'' these singularities along the locus $t = 0$ via a modification procedure. 

To generalize to the relative setting, we need to axiomatize the process of ``removing singularities''. In fact, moving from the classical theory to the relative theory requires substantial effort. 
One difficulty is that the expected $\mathbf{N}_{\mathrm{dR}}(\L)$ is not \'etale in general, and unlike in the classical situation, we do not have a general decompletion theory for non-\'etale $(\varphi, \Gamma)$-modules over $\mathbf{C}^{\dagger}_{\Lambda}$. Therefore, although there is a standard modification procedure for $(\varphi, \Gamma)$-modules over perfect relative period rings (using either relative $B$-pairs or the relative Fargues--Fontaine curve), one cannot construct $\mathbf{N}_{\mathrm{dR}}(\L)$ by descent (cf. Remark~\ref{rem: relative-NdR-FF-curve}). 

To overcome this issue, we avoid such a decompletion theory by performing the modification procedure directly on the level of imperfect period rings. In other words, we directly modify $D^{\dagger}_{\mathrm{rig}}(\L)$ along the ``$t=0$'' locus. For this, we introduce a subclass of decompleting data called \emph{strongly decompleting data} (Definition~\ref{defn: strongly decompleting data}). Some of the conditions for strongly decompleting data (Condition (a) and (b) in Definition~\ref{defn: strongly decompleting data}) essentially require that the $t$-adic completion of $\mathbf{C}^{[s,r]}_{\Lambda}$ (for some suitable $0<s\le r$) can be identified with a subring of the form
\begin{equation}\label{eq:intro_inclusion}
A_{X,n} [\![t, u_1, \ldots, u_d]\!]^{\nabla=0} \subset \mathbb{B}_{\mathrm{dR},X}^+ (X_\infty).
\end{equation} 
Here, $A_{X,n}$ is the ring of analytic functions on the $n$-th cover in the perfectoid tower $X_{\infty}\rightarrow X$ and $\mathbb{B}_{\mathrm{dR},X}^+$ is the horizontal de Rham period sheaf on the pro-\'etale site $X_{\proet}$ defined in \cite{Scholze_p_adic_Hodge}. Given this identification, we can apply Beauville--Laszlo gluing (see Proposition~\ref{prop:BL gluing}) to glue $D^{\dagger}_{\mathrm{rig}}(\L)[1/t]$ with $\left(D_{\mathrm{dR}}(\L)\otimes A_{X,n} [\![t, u_1, \ldots, u_d]\!]\right)^{\nabla=0}$ to obtain a new finite projective module over $\mathbf{C}^{[s, r]}_{\Lambda}$. For this to work, we need a gluing isomorphism over $A_{X,n} [\![t, u_1, \ldots, u_d]\!]^{\nabla=0}[1/t]$.

In fact, such an isomorphism exists over $\mathbb{B}_{\mathrm{dR},X}^+ (X_\infty)[1/t]$, essentially by construction. The question then becomes how to descend $\Gamma$-modules along the inclusion \eqref{eq:intro_inclusion}. Note that this inclusion is a relative version of the $\Gamma_K$-equivariant inclusion $K_\infty [\![t]\!] \subset \mathbf{B}_{\mathrm{dR}}^+$ studied by Fontaine \cite{Fontaine_arithmetic}. Therefore, what we are looking for is a relative version of Fontaine's decompletion theory over $\mathbf{B}_{\mathrm{dR}}^+$. To this end, we prove the following descent result, generalizing previous works of Andreatta--Brinon \cite{Andreatta-Brinon_Overconvergence, Andreatta-Brinon}. 

\begin{mainthm} \label{intro-theorem-decompletion-BdR+}
Let $A_{X, \infty}$ denote the union $\bigcup_{n \ge 0} A_{X,n}$. Then the base change functor  
\[ 
\mathrm{Mod}_{/A_{X, \infty}[\![t, u_1, \ldots, u_d]\!]^{\nabla=0}}^{\Gamma} \lra  \mathrm{Mod}_{/ \mathbb{B}_{\mathrm{dR},X}^+ (X_\infty)}^{\Gamma} 
\] 
from finite projective $\Gamma$-modules over $A_{X, \infty}[\![t, u_1, \ldots, u_d]\!]^{\nabla=0}$ to those over $\mathbb{B}_{\mathrm{dR},X}^+ (X_\infty)$ is an equivalence of categories.
\end{mainthm}

We refer the reader to \S \ref{sec:relative_Tate_Sen} for the precise statement (Theorem~\ref{thm:equivalence_relative_B_dR^+}(2)). Combining this result and a further delicate descent argument (to descend from $A_{X,\infty}$ to $A_{X,n}$ using Condition (c) of Definition~\ref{defn: strongly decompleting data}), we complete the proof of Theorem \ref{intro:theorem_relative_Berger} in \S \ref{sec:relative_N_dR}.

\begin{mainrmk}
 \label{remark:intro_remark_AB_example}
When $X$ admits a smooth formal integral model over $\mathcal{O}_K$ as in the set-up of \cite{Andreatta-Brinon_Overconvergence, Andreatta-Brinon}, the imperfect period rings in \emph{loc. cit.} indeed correspond to strongly decompleting data (see \S \ref{subsec: relation-with-AB}). However, the framework of \emph{loc. cit.} does not cover the important case of a closed (poly)annulus (or a finite \'etale cover thereof). We construct explicit strongly decompleting data in such cases (see \S \ref{subsec: closed-annuli}).  
\end{mainrmk}

\begin{mainrmk}
Both decompleting and strongly decompleting data are ``stable under finite \'etale extensions'' in the sense that, if $\Lambda=(R_{\Lambda}, \mathbf{A}_{\Lambda})$ is a (strongly) decompleting data for $X$ and $X' \ra X$ is a finite \'etale cover, then one can canonically produce a (strongly) decompleting data $\Lambda'$ for $X'$ (see Theorem \ref{thm: weakly decompleting data finite etale} and \ref{thm: strongly decompleting under finite etale extensions}). This in turn implies that, Theorem \ref{intro:theorem_equiv_of_categories} and \ref{intro:theorem_relative_Berger} can be applied to study any $p$-adic local system \textit{locally} on a smooth rigid analytic space $X$ over $K$. More specifically, for any such $X$, there exists a covering of $X$ by affinoid open subspaces, each equipped with a finite \'etale map to a polyannulus (cf. \cite{Achinger}). Therefore, we have an explicit (strongly) decompleting data for each member of the covering by Remark \ref{remark:intro_remark_AB_example}.  We caution the reader that the notion of (strongly) decompleting data depends on the choice of the chart $h$ and we do not claim any global construction of $p$-adic differential equations over an arbitrary $X$. 
\end{mainrmk}

\begin{mainrmk}
In the one-dimensional situation, we explicitly construct strongly decompleting data for two examples in \S \ref{sec:example}: (i) closed annuli (and their finite \'etale covers); (ii) closed annuli with finitely many open disks removed (and their finite \'etale covers). In fact, these include the most important affinoid open subsets of a smooth curve. For instance, up to a base field extension, every type I\!I\!I point on a smooth curve admits a basis of open neighborhoods consisting of sets of the form (i), whereas every type V point admits a basis of open neighborhoods consisting of sets of the form (ii). Consequently, the constructions in this article allow us to prove the $p$-adic monodromy theorem locally around type I\!I\!I and V points in \cite{p_adic_monodromy_curve}.
\end{mainrmk}

\begin{mainrmk}
It is natural to wonder the relation between the techniques developed in the current manuscript and the theory of (pro-)locally analytic vectors, which provides another framework to decompletion. One central difficulty of the latter is that we do not know in general how to compute the (pro)-locally analytic vectors in perfect relative period rings. Furthermore, for the construction of relative $p$-adic differential equations, it seems unavoidable to study conditions of the same flavor as in the notion of strongly decompleting data. 
\end{mainrmk}

\subsection{Outline of contents}

In Section~\ref{sec: Banach-rings-modules}, we recall some definitions and preliminary facts on (semi)norms, Banach rings, and Banach modules that will be used in this article.

In Section~\ref{sec:classical_phi_Gamma}, we review the classical theory of $(\varphi, \Gamma_K)$-modules and $p$-adic differential equations associated with (de Rham) $p$-adic Galois representations, following \cite{Cherbonnier_Colmez, Berger-differential, Colmez08}. The perfect and imperfect period rings in the classical setting are recalled in \S \ref{ss:point perfect period ring} and \S \ref{ss:point imperfect period ring}, respectively. In \S \ref{ss:classical_pDE}, we recall the theory of $(\varphi, \Gamma_K)$-modules associated with $p$-adic representations of $\Gal_K$. Then in \S \ref{ss:NdR classical setting}, we provide an interpretation of Berger's construction of $p$-adic differential equations in terms of Beauville--Laszlo gluing.

In the remainder of the paper, we consider the relative situation. 

In Section~\ref{sec:relative_period_rings:perfect}, we review perfect relative period rings and the theory of $(\varphi, \Gamma)$-modules over them, mainly following \cite{Kedlaya-LiuI}. The notion of toric charts $X \rightarrow \D_K^d$ and toric towers $X_{\infty} \rightarrow X$ is introduced in \S \ref{ss:toric_tower}; this is the relative setting we work with throughout the article. In \S \ref{ss:perfect relative period rings} and \S \ref{ss:fet_over_tilde_A}, we review from \cite{Kedlaya-LiuI} the construction of perfect relative period rings and their behavior under finite \'etale extensions of $X$. Then in \S \ref{subsec: (phi, gamma)-modules-over-perfect-period-rings}, we recall from \cite{Kedlaya-LiuI} the theory of $(\varphi, \Gamma)$-modules over perfect relative period rings and their relation to \'etale $\Z_p$-local systems on $X$.

In Section~\ref{sec:weakly_decompleting_data}, we introduce and study the notions of weakly decompleting data and the associated imperfect relative period rings. We introduce these notions in \S \ref{subsection: weakly decompleting data}. In \S \ref{subsec: norm-compatible-expansions-good-lifts} - \S \ref{subsec: intersection-period-rings}, we study some important properties of weakly decompleting data and associated period rings, which will be used in later sections. In \S \ref{ss:weakly decompleting finite etale}, we explain that given a weakly decompleting datum for $X$, there is a natural way to construct a weakly decompleting datum for any finite \'etale extension of $X$.

In Section~\ref{sec:decompleting_data}, we introduce the notion of decompleting data and study their properties; this will allow us to decomplete and descend \'etale $(\varphi, \Gamma)$-modules over perfect relative period rings to the imperfect ones. We introduce the notion of decompleting data in \S \ref{subsec: definition-decompleting-data}. Then in \S \ref{subsec: strict-exactness-char-p} - \S \ref{subsec: strict-exactness-type-C}, we prove a sequence of results on the strict-exactness of certain group cochain complexes, valued in period rings of characteristic $p$ and of type $\fA$ and $\fC$.

In Section~\ref{sec: decompletion-relative-etale-phiGamma-modules}, we study decompletion of \'etale $(\varphi, \Gamma)$-modules under the assumption of decompleting data. We introduce relative $\varphi$-modules and $(\varphi, \Gamma)$-modules over imperfect period rings in \S \ref{subsec: definition-relative-(phi, Gamma)-modules}, and explain some preliminary results on them in \S \ref{subsec: preliminary-results-relative-(phi, Gamma)-modules}. Then in \S \ref{subsec: descent-relative-phi-modules}, we study some descent questions for relative $\varphi$-modules. The main results of this section are in \S \ref{subsec: fully-faithfulness-indfree-(phi, Gamma)-modules} and \S \ref{subsec: descent-relative-etale-(phi, Gamma)-modules}: in \S \ref{subsec: fully-faithfulness-indfree-(phi, Gamma)-modules}, we prove a full faithfulness result for relative $(\varphi, \Gamma)$-modules along the base change from the imperfect period rings to the perfect ones, while in \S \ref{subsec: descent-relative-etale-(phi, Gamma)-modules}, we establish the decompletion result of \'etale $(\varphi, \Gamma)$-modules. In \S \ref{ss:relative_D_rig}, we summarize some of the main results in this section and deduce the equivalence between $p$-adic local systems on $X$ and \'etale $(\varphi, \Gamma)$-modules over imperfect relative period rings. In particular, Theorem~\ref{intro:theorem_equiv_of_categories} is proved.

In Section~\ref{sec:strongly_decompleting_data}, we introduce and study the notion of strongly decompleting data, which will be used in later sections to construct a relative analogue of Berger's $\mathbf{N}_{\mathrm{dR}}$-functor. The definition of strongly decompleting data is given in \S \ref{subsec: definition-strongly-decompleting-data}, together with some preliminary facts. In \S \ref{subsec: strongly-decompleting-finite-etale}, we analyze the behavior of strongly decompleting data under finite \'etale extensions of $X$.

In Section~\ref{sec:relative_Tate_Sen}, we study a relative analogue of Fontaine's decompletion theory over $\mathbf{B}_{\mathrm{dR}}^+$, generalizing some results in \cite{Andreatta-Brinon}. We explain Tate's normalized trace in the relative setting in \S \ref{ss:relative_TS_setup}. Then in \S \ref{ss:relative_decompletion_degree_0} - \S \ref{subsec: decompletion-BdR+}, we prove the decompletion result claimed in Theorem~\ref{intro-theorem-decompletion-BdR+}.

In Section~\ref{sec:relative_N_dR}, given any strongly decompleting datum for $X$, we construct a relative analogue of Berger's $\mathbf{N}_{\mathrm{dR}}$-functor from the category of de Rham local systems on $X$ to that of $(\varphi, \Gamma)$-modules over $\fC_{\Lambda}^{\dagger}$; this establishes Theorem~\ref{intro:theorem_relative_Berger}. We begin by constructing the functor to $(\varphi, \Gamma)$-modules over the perfect period ring in \S \ref{subsec: tilde-NdR}. In \S \ref{subsec: NdR-preparation} - \S \ref{subsec: NdR-BL gluing}, we construct the relative $\mathbf{N}_{\mathrm{dR}}$ functor based on Beauville--Laszlo gluing and the descent result established in \S \ref{sec:relative_Tate_Sen}, and then study its properties as a $p$-adic differential equation in \S \ref{subsec: NdR-pdE}. In \S \ref{subsec: compatibility-NdR}, we show the compatibility of $\fN_{\mathrm{dR}}$ functor when we vary $X$.

In Section~\ref{sec:example}, our main goal is to construct a natural example of strongly decompleting datum when $X$ is either a closed annulus or a closed annulus with finitely many open disks removed (or a finite \'etale cover thereof). These examples are carefully studied in \S \ref{subsec: closed-annuli} and \S \ref{subsec: punctured-closed-annuli}. Moreover, we explicitly compute the corresponding normalized differential operators in \S \ref{subsec: normalized-differential-operators-closed-annuli}, and study the compatibility of $\fN_{\mathrm{dR}}$ with base field extensions in \S \ref{subsec: compatibility-NdR-base-field-extension}. Lastly in \S \ref{subsec: relation-with-AB}, we briefly discuss the relation of strongly decompleting data with the works in \cite{Andreatta-Brinon_Overconvergence, Andreatta-Brinon}; in particular, we explain that the good reduction cases studied in \emph{loc. cit.} are indeed strongly decompleting.

In Appendix~\ref{section: finitely generated modules}, we study finitely generated (resp. finite projective) modules over imperfect relative period rings associated with weakly decompleting data, and analyze their relations with coherent sheaves (resp. vector bundles) on the associated quasi-Stein spaces. This is used in the study of $\varphi$-modules over imperfect relative period rings in \S \ref{sec: decompletion-relative-etale-phiGamma-modules}.

\subsection{Notation and Convention}\label{ss: notation and convention}

We adopt the following notation and convention throughout this article.
\begin{itemize}
    \item Let $K$ be a $p$-adic field\footnote{Again, by a \emph{$p$-adic field} we mean a complete discrete valuation field of mixed characteristic $(0,p)$ with perfect residue field.} with ring of integers $\mO_K$ and residue field $k=k_K$. Let $\varpi$ be a uniformizer of $K$. Let $K_0=W(k)[1/p]$. 
    \item Let $\Gal_K$ denote the absolute Galois group of $K$. 
    \item Fix an algebraic closure $\overline{K}$ of $K$ and let $\C_K$ denote the completion of $\overline{K}$. This is a perfectoid field; let $\C_K^{\flat}$ denote its tilt. Let $\mO_{\C_K}$ (resp. $\mO_{\C_K}^{\flat}$) denote the ring of integers in $\C_K$ (resp. $\C_K^{\flat}$).\footnote{The ring $\mO^{\flat}_{\C_K}$ (resp. $ \C_K^{\flat}$) is denoted by $\widetilde{\mathbf{E}}^+$ (resp. $\widetilde{\mathbf{E}}$) in \cite{Berger-differential} and \cite{Colmez08}.} 
    \item  Fix a compatible system of $p$-power roots  $\{\zeta_{p^n} \}_{n \ge 0}$  of unity and a compatible system of $p$-power roots $\{p^{1/p^n}\}$ of $p$ in $\mO_{\C_K}$. This gives rise to elements $\varepsilon = (1, \zeta_p, \zeta_{p^2}, \ldots)$ and $p^{\flat}=(p,p^{1/p}, p^{1/p^2}\ldots)$ in $\mO_{\C_K}^\flat$. 
    \item Let $\cl \pi = \varepsilon - 1$. Write $\pi = [\varepsilon] -1 \in W(\mO_{\C_K}^{\flat})$, which is a lift of $\cl \pi$. Let $t=\log [\varepsilon]=\log(1+\pi)$ viewed as an element in Fontaine's de Rham period ring $\mathbf{B}_{\mathrm{dR}}$.
    \item For each $n\ge 1$, let $K_n=K(\zeta_{p^n})$. Let $K_{\infty}=\bigcup_{n\ge 0} K_n$ be the cyclotomic extension. Let $\widehat{K}_{\infty}$ be the $p$-adic completion of $K_{\infty}$ and let $\widehat{K}_{\infty}^\flat$ be its tilt. Write $\mO_{\widehat{K}_{\infty}}$ (resp.  $\mO^\flat_{\widehat{K}_{\infty}}$) for the ring of integers in $\widehat{K}_{\infty}$ (resp.  $\widehat{K}_{\infty}^\flat$).\footnote{The ring $\mO^\flat_{\widehat{K}_{\infty}}$ (resp. $ {\widehat{K}_{\infty}}^\flat$) is denoted by $\widetilde{\mathbf{E}}^+_K$ (resp. $\widetilde{\mathbf{E}}_K$) in \cite{Berger-differential} and \cite{Colmez08}.} 
    \item Let $H_K = \Gal_{K_{\infty}}$ be the Galois group of $K_\infty$ and let $\Gamma_K\coloneqq \mathrm{Gal}(K_{\infty}/K)\cong \Gal_K /H_K$. 
    \item Let $e_K = [K_{\infty}: (K_0)_{\infty}]$.
    \item Let $K'_0$ denote the maximal unramified extension of $K_0$ inside ${K}_{\infty}$. Denote by $\kappa=\kappa_K$ the residue field of $K'_0$. 
    \item The standard $p$-adic norm $|\cdot|$ on $\mathbb{C}_K$, normalized by $|p|=\frac{1}{p}$, defines a norm $|\cdot|$ on $\mathbb{C}_K^{\flat}$. In particular, $|\cdot|$ restricts to a norm on $\widehat{K}_{\infty}^{\flat}$. Note that $|\overline{\pi}| = |\varepsilon-1|=p^{-\frac{p}{p-1}}$. We then define a valuation $v\colon \mathbb{C}_K\rightarrow \mathbb{R}\cup \{+\infty\}$  by 
    \[
    v(y)\coloneqq - \log_p|y|.
    \] 
    This in turn determines a valuation $v\colon \mathbb{C}_K^\flat \rightarrow \mathbb{R}\cup \{+\infty\}$.\footnote{The valuation $v$ on $\mathbb{C}_K^\flat$ is denoted by $v_{\mathbf{E}}$ in \cite{Berger-differential} and \cite{Colmez08}.} In particular, we have $v(\overline{\pi}) = v(\varepsilon - 1)=\frac{p}{p-1}$.
    \item Let $k(\!(\cl \pi)\!)^{\mathrm{sep}}$ denote the separable closure of $k(\!(\cl \pi)\!)$ inside $\C_K^\flat$,\footnote{The ring $k(\!(\cl \pi)\!)^{\mathrm{sep}}$ is denoted by $\mathbf{E}$ in \cite{Berger-differential} and \cite{Colmez08}.} and write $\mathbf{E}_K\coloneqq (k(\!(\cl \pi)\!)^{\mathrm{sep}})^{H_K}$. Let  $\mathbf{E}_K^+$  be its ring of integers, which is a complete discretely valued ring with residue field $\kappa_K$. By the theory of fields of norms, we have an isomorphism 
\[
\varprojlim_{n} \mO_{K_n}/\mathfrak a \isom \mathbf{E}_K^+ 
\]
    where $\mathfrak a \subset \mO_{\C_K}$ is the ideal consisting of elements $x$ satisfying $|x| < p^{-\frac{1}{p-1}}$ and the inverse limit is taken over the Frobenius maps  $ \mO_{K_{n+1}}/\mathfrak a \ra  \mO_{K_n}/\mathfrak a$ sending $x \mapsto x^p$.   
    \item Fix a choice of uniformizer $\cl \pi_K$ of $\mathbf E_K^+$. Note that 
    $\mathbf E_K^+ \cong \kappa_K [\![\cl \pi_K]\!]$ and $v(\cl \pi_K)=\frac{p}{p-1}\cdot\frac{1}{e_K}$.  
    \item Throughout the article, for all the period rings $\widetilde{\mathbf{A}}^r_K, \widetilde{\mathbf{B}}^r_K, \widetilde{\mathbf{C}}^r_K, \widetilde{\mathbf{C}}^{[s,r]}_K$ and $\mathbf{A}^r_K, \mathbf{B}^r_K,\mathbf{C}^r_K, \mathbf{C}^{[s,r]}_K$, etc., as well as their relative counterparts, we always assume $r,s\in \Q$ in the superscripts.
\end{itemize}
Moreover, we follow the conventions on adic geometry in \cite{Huber_etale} and \cite{perfectoid}. For a smooth rigid analytic variety $X$ over $K$, we often view it as an adic space over $\spa(K,\mathcal{O}_K)$. Let $\mathrm{Loc}_{\Z_p}(X)$ denote the category of \'etale $\Z_p$-local systems on $X$ and let $\mathrm{Loc}_{\Z_p}(X)_{\Q}$ denote the corresponding isogeny category. Let $\mathrm{Loc}^{\mathrm{dR}}_{\Z_p}(X)$ denote the full subcategory of $\mathrm{Loc}_{\Z_p}(X)$ consisting of de Rham local systems in the sense of \cite{Scholze_p_adic_Hodge, Liu_Zhu}.

\subsection*{Acknowledgement}  

The authors are grateful to Kiran Kedlaya and Ruochuan Liu for their preprint \cite{Kedlaya-LiuII}. Evidently, some of the ideas behind the notion of (weakly) decompleting data stemmed from our initial attempts to reconstruct the proofs outlined in \textit{loc. cit.}  We would like to thank Shizhang Li and Lie Qian for discussions and comments on earlier drafts of the manuscript. We thank Heng Du for discussions around the contents of \cite{Kedlaya-LiuII} and on modifications of vector bundles on relative Fargues--Fontaine curves, which should lead to an alternative approach to the construction of $\widetilde{\mathbf{N}}_{\mathrm{dR}}(\L)$ in \S \ref{subsec: tilde-NdR}. 
During the project, H.D. was partially supported by the National Natural Science Foundation of China (Grant No. 12422101) and the National Key R\&D Program of China (Grant Nos. 2023YFA1009703 and 2021YFA1000704), Y.S.M. was partially supported by the Beijing Natural Science Foundation (Grant No. IS25030).

\vspace{0.1in}
\section{Banach rings and modules} \label{sec: Banach-rings-modules}

We begin by recalling some terminology and conventions concerning non-archimedean seminorms, Banach rings, and Banach modules, as well as the notion of rational localizations of Banach rings.

\subsection{Seminorms}
The (semi)norms in this article will always be the non-archimedean ones in the following sense.

\begin{definition}
Let $G$ be an abelian group. A (non-archimedean) \emph{seminorm} on $G$ is a function $|\cdot|\colon G \rightarrow \mathbb{R}_{\ge 0}$ such that $|0| = 0$ and $|g+h| \le \max\{|g|, |h|\}$ for all $g, h \in G$. If, in addition, $|g|\neq 0$ whenever $g\neq 0$, we say that $|\cdot|$ is a \emph{norm}.
\end{definition}

\begin{definition}\label{defn: equivalent seminorms} \hfill
\begin{enumerate}
\item For two (semi)norms $|\cdot|$ and $|\cdot|'$ on a given abelian group $G$, we say $|\cdot|$ and $|\cdot|'$ are \emph{equivalent} if there exist constants $c_1, c_2 > 0$ such that $|g| \le c_1 |g|' \textrm{ and } |g|' \le c_2 |g|$ for all $g \in G$.

\item Let $G$ and $H$ be abelian groups equipped with seminorms $|\cdot|$ and $|\cdot|'$ respectively, and let $f\colon G \rightarrow H$ be a group homomorphism. 
\begin{itemize}
\item We say that $f$ is \emph{bounded} if there is a constant $c > 0$ such that $|f(g)|' \le c|g|$ for all $g \in G$.

\item We say that $f$ is \emph{submetric} if $|f(g)|' \le |g|$ for all $g \in G$, and we say that $f$ is \emph{isometric} if $|f(g)|'=|g|$ for all $g \in G$.

\item Define the \emph{quotient seminorm} $|\cdot|_{\mathrm{quot}}$ on $\mathrm{image}(f)$ by
\[
|h|_{\mathrm{quot}}\coloneqq \inf\{|g|~|~ g \in G, ~f(g) = h\}.
\]
We say that $f$ is \emph{strict} if $|\cdot|_{\mathrm{quot}}$ and $|\cdot|'$ are equivalent on $\mathrm{image}(f)$.
\end{itemize}
\end{enumerate}
\end{definition}

\begin{definition}
Let $A$ be a ring and let $|\cdot|$ be a seminorm on $A$ (considered as an additive group).
\begin{enumerate}
\item We say that $|\cdot|$ is \emph{submultiplicative} if $|xy| \le |x|\cdot |y|$ for any $x, y \in A$.

\item We say that $|\cdot|$ is \emph{power-multiplicative} if it is submultiplicative and $|x^n| = |x|^n$ for all $x\in A$ and all integers $n \ge 1$.

\item We say that $|\cdot|$ is \emph{multiplicative} if $|1| = 1$ and $|xy| = |x|\cdot |y|$ for all $x, y \in A$.
\end{enumerate}
\end{definition}

\subsection{Banach rings and modules}

\begin{definition} \hfill
\begin{enumerate}
\item A \emph{Banach ring} is a commutative ring equipped with a submultiplicative norm under which it is complete.

\item Let $S$ be a Banach ring with a norm $|\cdot|$. A \emph{Banach} $S$-\emph{module} is an $S$-module $M$ equipped with a norm $|\cdot|_M$ under which $M$ is complete such that
\[
|s\cdot m|_M \le |s|\cdot |m|_M
\]
for all $s\in S$ and $m \in M$.

\item Let $S$ be a Banach ring. A \emph{Banach} $S$-\emph{algebra} is an $S$-algebra $S'$ equipped with a norm $|\cdot|'$ such that $(S', |\cdot|')$ is both a Banach ring and a Banach $S$-module.
\end{enumerate}
\end{definition}

\begin{remark}
Our conventions for Banach modules and Banach algebras agree with those in some standard references (for example, \cite{BGR}), but are more restrictive than those in some other literature (for example, \cite[Definition~2.2.7]{Kedlaya-LiuI}).
\end{remark}

We now recall some basic facts about finite projective modules over Banach rings (cf. \cite[\S 2]{Kedlaya-LiuI}).

\begin{lemma}[{\cite[Lemma ~2.2.6(a)]{Kedlaya-LiuI}}] \label{lem: finite module seminorm}
Let $S$ be a Banach ring and $M$ be a finite $S$-module. For any surjection $S^n \twoheadrightarrow M$ from a finite free $S$-module, we can equip $M$ with the quotient seminorm induced from the supremum norm on $S^n$ with respect to the canonical basis. Then, up to equivalence, the quotient seminorm does not depend on the choice of the surjection.
\end{lemma}

\begin{definition}\label{defn: finite-module seminorm}
Let $S$ be a Banach ring and $M$ be a finite $S$-module. A seminorm on $M$ as in Lemma~\ref{lem: finite module seminorm} is called a \emph{finite-module seminorm}. All finite-module seminorms on $M$ are equivalent, and so define the same topology on $M$.
\end{definition}

\begin{remark}[cf. {\cite[Remark~2.2.11]{Kedlaya-LiuI}}]\label{rmk: finite-module seminorm}
When $S$ contains a topologically nilpotent unit, a finite-module seminorm on $M$ is a complete norm. In general, a finite-module seminorm is not necessarily a norm, neither complete. On the other hand, finite-module (semi)norms behave well for finite projective $S$-modules as explained in the following.
\end{remark}

\begin{lemma}[{\cite[Lemma~2.2.12(a)(b)]{Kedlaya-LiuI}}] \label{lem: finite projective norm}
Let $S$ be a Banach ring and $P$ be a finite projective $S$-module. Choose a finite projective $S$-module $Q$ with an isomorphism $P\oplus Q\cong S^{\oplus n}$ of $S$-modules for some $n$. Equip $S^{\oplus n}$ with the sup-norm defined by the canonical basis. Then the subspace norm on $P$ inherited from $S^{\oplus n}$ is complete, and is equivalent to the quotient norm with respect to the projection from $S^{\oplus n}$. In particular, any finite-module seminorm on $P$ is a complete norm.
\end{lemma}

We also recall the following version of the open mapping theorem.

\begin{theorem}[see {\cite[Theorem~2.2.8]{Kedlaya-LiuI}} and the references therein] \label{thm: open mapping thm}
Let $S$ be a Banach ring containing a topologically nilpotent unit. Let $g\colon M\rightarrow N$ be a bounded surjective homomorphism of Banach $S$-modules. Then $g$ is open and strict.
\end{theorem}

Throughout this article, we restrict our attention to Banach rings satisfying the following additional condition.

\begin{assumption} \label{assump: banach-ring-unit-z}
Let $S$ be a Banach ring with a norm $|\cdot|$. We always assume that $S$ contains a unit $z$ such that $|z| < 1$ and $|z|\cdot |z^{-1}| = 1$.   
\end{assumption}

\begin{lemma} \label{lem: Banach-module-unit-z}
Let $S$ be a Banach ring and let $z$ be a unit in $S$ satisfying the conditions in Assumption~\ref{assump: banach-ring-unit-z}. Let $M$ be a Banach $S$-module, equipped with a norm $|\cdot|_M$. Then 
\[
|z\cdot m|_M = |z|\cdot |m|_M \textrm{ and } |z^{-1}\cdot m|_M = |z^{-1}|\cdot|m|_M
\]
for all $m \in M$.
\end{lemma}

\begin{proof}
For any $m \in M$, we have
\[
|m|_M = |zz^{-1}m|_M \le |z^{-1}|\cdot |z\cdot m|_M \le |z^{-1}|\cdot |z|\cdot |m|_M = |m|_M.
\]
and
\[
|m|_M = |zz^{-1}m|_M \le |z|\cdot |z^{-1}\cdot m|_M \le |z|\cdot |z^{-1}|\cdot |m|_M = |m|_M.
\]
Thus, the above inequalities are all equalities. 
\end{proof}

\begin{lemma}
Let $S$ be a Banach ring satisfying Assumption~\ref{assump: banach-ring-unit-z} and let $M$ be an $S$-module. Let $|\cdot|_1$ and $|\cdot|_2$ be two norms on $M$ such that $(M, |\cdot|_i)$ is a Banach $S$-module for each $i = 1, 2$. If $|\cdot|_1$ and $|\cdot|_2$ induce the same topology on $M$, then they are equivalent.
\end{lemma}

\begin{proof}
Since $|\cdot|_1$ and $|\cdot|_2$ induce the same topology on $M$, there is a constant $c$ such that for any $m \in M$ with $|m|_2 < c$, we have $|m|_1 < 1$. Write $|z| = b < 1$.

We will repeatedly apply Lemma~\ref{lem: Banach-module-unit-z} in what follows. For example, note that
$|z^n| = |z|^n = b^n$ 
for all integers $n$. Let $x \in M$ be a non-zero element. Choose $n$ such that
$b\cdot c|x|_2^{-1} \le |z^n|= b^n < c|x|_2^{-1}.$ 
Then we have $ |z^n\cdot x|_2 = |z^n|\cdot |x|_2 < c,$ 
so we have 
$bc|x|_2^{-1}\cdot |x|_1 \le  |z^n|\cdot |x|_1 = |z^n\cdot x|_1 < 1.$ 
Thus, we conclude that $|x|_1 \le (bc)^{-1}\cdot |x|_2.$ By symmetry, this implies that $|\cdot|_1$ and $|\cdot|_2$ are equivalent.
\end{proof}

\vspace{0.1in}
\section{$(\varphi, \Gamma_K)$-modules and Galois representations} \label{sec:classical_phi_Gamma}

In this section, we review the classical theory of $(\varphi, \Gamma_K)$-modules and $p$-adic  differential equations associated with (de Rham) $p$-adic Galois representations, following \cite{Cherbonnier_Colmez, Berger-differential, Colmez08}. The goals of this section are as follows. 
\begin{itemize}
    \item We first recall the classical theory of perfect period rings $\widetilde{\mathbf{A}}^\dagger_K, \widetilde{\mathbf{B}}^\dagger_K$, etc., and then introduce the imperfect period rings as their subrings.  
    \item We then review the theory of $(\varphi, \Gamma_K)$-modules and $p$-adic differential equations over these period rings. 
    \item In particular, we recall (and give another interpretation of) Berger's construction that associates a $p$-adic differential equation $\mathbf{N}_{\mathrm{dR}} (V)$ over the Robba ring to every de Rham representation.    
\end{itemize}
As explained in the introduction, we shall follow the notation of Kedlaya--Liu \cite{Kedlaya-LiuI} on period rings although the content in this section predates \textit{loc. cit..} and is mostly developed in \cite{Cherbonnier_Colmez} and \cite{Berger-differential}. 

We adopt the notation from the introduction \S \ref{ss: notation and convention}; in particular, recall the notation $K$, $\mathcal{O}_K$, $k$, $\varpi$, $\Gal_K$, $\C_K$, $\varepsilon$, $\overline{\pi}$, $\pi$, $K_n$, $K_{\infty}$, $\widehat{K}_{\infty}$, $\widehat{K}^{\flat}_{\infty}$, $H_K$, $\Gamma_K$, $e_K$, $K'_0$, $\kappa$, $\overline{\pi}_K$, as well as the standard norms (resp. valuations) $|\cdot|$ (resp. $v$) on $\C_K$ and $\C_K^{\flat}$. Moreover, we remind the reader that, for all the period rings in this article, the numbers $r,s$ in the superscripts are always assumed to be rational numbers.

\subsection{Perfect period rings} \label{ss:point perfect period ring}

In \S \ref{ss:point perfect period ring} and \S \ref{ss:point imperfect period ring}, we review the constructions of perfect and imperfect period rings in the classical setting.

\subsubsection{Perfect period rings of type $\mathbf{A}$ and $\mathbf{B}$}

Let $\widetilde{\mathbf{A}}=W(\C_K^{\flat})$ and $\widetilde{\mathbf{A}}^+ =W(\mO_{\C_K}^{\flat})$, equipped with natural $\Gal_K$- and Frobenius actions. Taking the invariant subrings under the action of $H_K$, we obtain
\[ 
\widetilde{\mathbf{A}}_K= (\widetilde{\mathbf{A}})^{H_K} = W(\widehat{K}_{\infty}^{\flat}) \quad \textup{and } \,\,\, \widetilde{\mathbf{A}}^+_K = (\widetilde{\mathbf{A}}^+)^{H_K}=W(\mO_{\widehat{K}_{\infty}^{\flat}}).
\] 
Note that every element $x\in \widetilde{\mathbf{A}}_K$ can be uniquely written as $x=\sum_{n=0}^{\infty}p^n[y_n]$ with $y_n\in \widehat{K}_{\infty}^{\flat}$. 
The ring $\widetilde{\mathbf{A}}^+$ is denoted by $\mathbf{A}_{\tu{inf}}$ in \cite{Fontaine_94}. We refer the reader to \textit{loc. cit..} for the construction of period rings such as $\mathbf{B}_{\mathrm{dR}}, \mathbf{B}_{\tu{crys}}, \mathbf{B}_{\tu{st}}$ and their variants.

\begin{definition}[{\cite[\S 2.1]{Berger-differential}, \cite[\S 5.2]{Colmez08}}]
\label{def:perfect_period_rings_type_A} \hfill
\begin{enumerate}
\item  For each $r>0$, define 
\begin{align*}
\widetilde{\mathbf{A}}^{r} & \coloneqq \left\{x=\sum_{n=0}^{\infty}p^n[y_n] \in \widetilde{\mathbf{A}}\,\Big|\,\lim_{n\rightarrow\infty}v(y_n)+\frac{n}{r}=+\infty\right\},  \\ 
\widetilde{\mathbf{A}}_K^{r} & \coloneqq \left\{x=\sum_{n=0}^{\infty}p^n[y_n] \in \widetilde{\mathbf{A}}_K\,\Big|\,\lim_{n\rightarrow\infty}v(y_n)+\frac{n}{r}=+\infty\right\}.
\end{align*}

\item On $\widetilde{\mathbf{A}}^{r}$, we define a valuation $v_{r}\colon \widetilde{\mathbf{A}}^{r}\rightarrow \mathbb{R}\cup \{+\infty\}$ by \footnote{The valuation $v_r$ is denoted as $v^{(0, r]}$ in \cite{Colmez08}.}
\[
v_{r}(x)\coloneqq \inf_{n\ge 0}\{v(y_n)+\frac{n}{r}\}.
\]
In fact, for each $0<s\le r$, we have a valuation $v_{s}$ on $\widetilde{\mathbf{A}}^{r}$ given by 
\[
v_{s}(x)\coloneqq \inf_{n\ge 0}\{v(y_n)+\frac{n}{s}\}.
\]
Similarly, we define valuations $v_s$ on $\widetilde{\mathbf{A}}_K^{r}$ for each $s \in (0, r]$.

\item  Define 
\[
\widetilde{\mathbf{A}}^{\dagger} \coloneqq \varinjlim_{r>0} \widetilde{\mathbf{A}}^{r}, \qquad \widetilde{\mathbf{A}}^{\dagger}_K \coloneqq \varinjlim_{r>0} \widetilde{\mathbf{A}}^{r}_K,
\]
where the transition maps are natural inclusions. Notice that $\widetilde{\mathbf{A}}$ (resp. $\widetilde{\mathbf{A}}_K$) is the $p$-adic completion of $\widetilde{\mathbf{A}}^{\dagger}$ (resp. $\widetilde{\mathbf{A}}^{\dagger}_K$).

\item The period rings $\widetilde{\mathbf{B}}, \widetilde{\mathbf{B}}^+, \widetilde{\mathbf{B}}^{r}, \widetilde{\mathbf{B}}^{\dagger},
    \widetilde{\mathbf{B}}_K$, $\widetilde{\mathbf{B}}^+_K$, $\widetilde{\mathbf{B}}^{r}_K$, and $\widetilde{\mathbf{B}}^{\dagger}_K$ are obtained from the corresponding period rings of type $\mathbf{A}$ by inverting $p$.  
\end{enumerate}
\end{definition}

\begin{remark}[The norm $\lambda_r$] \label{remark:norm_lambda_r} \hfill
\begin{enumerate}
\item
One can translate the definitions above into the language of norms. Consider norms $\lambda_{r}$ on $\widetilde{\mathbf{A}}^{r}$ and $\widetilde{\mathbf{A}}_K^{r}$ defined by
\begin{align*}
\lambda_{r}(x)\coloneqq \sup_{n\ge 0} \{p^{-n}\cdot |y_n|^r\}  = p^{-r\cdot v_{r}(x)} 
\end{align*}
Similarly, for each $s \in (0, r]$, we have norms $\lambda_{s}$ on $\widetilde{\mathbf{A}}^{r}$ and $\widetilde{\mathbf{A}}_K^{r}$ given by
\[
\lambda_{s}(x)\coloneqq \sup_{n\ge 0} \{p^{-n}\cdot |y_n|^s\} = p^{-s \cdot v_s (x)}.
\]
\item For example, consider $a^\flat = (a, a^{1/p}, a^{1/p^2}, \ldots) \in \mO_{\C_K}^\flat$ where $\{a^{1/p^m}\}_{m \ge 0}$ is a compatible choice of $p$-power roots of $a$. Then $[a^\flat] \in \widetilde{\mathbf{A}}^r$ for every $r > 0$, and we have $v_r ([a^\flat]) = v(a^\flat) = v (a)$ is independent of $r$; while $ \lambda_r ([a^\flat]) = p^{-r \cdot v(a)}$. On the other hand, note that $v_r (p) = \frac{1}{r}$; and 
$\lambda_r (p) = p^{-1} = |p|$ is independent of $r$.    
\item We have 
\[ 
\sq{\mathbf{A}}_{K}^{r}  = \left\{x=\sum_{n=0}^{\infty}p^n[y_n] \in \widetilde{\mathbf{A}}_K \,\Big|\,  p^{-n} |y_n|^r \ra 0 \tu{ as } n \ra +\infty \right\}. 
\]
Furthermore, $ \sq{\mathbf{A}}_{K}^{r} $ is complete with respect to the norm $\lambda_{r}$.  
\end{enumerate}
\end{remark}

\subsubsection{Perfect period rings of type $\mathbf{C}$}

\begin{definition}[{\cite[\S 2.3, Definition~2.16]{Berger-differential}, \cite[\S 5.4]{Colmez08}
}] \label{def:perfect_robba_ring_type_C_over_K} \noindent 
\begin{enumerate}
\item For each $0 < s \le r$,  define $\widetilde{\mathbf{C}}^{[s, r]}$  (resp. $\widetilde{\mathbf{C}}^{[s, r]}_K$ ) to be the completion of $\widetilde{\mathbf{B}}^r$  (resp. of $\widetilde{\mathbf{B}}^r_K$) with respect to \(\max\{\lambda_{s}, \lambda_{r}\}\).  

\item For each $r > 0$, define $\widetilde{\mathbf{C}}^{r}$ (resp.  $\widetilde{\mathbf{C}}^{r}_K$) to be the Fr\'echet completion of  $\widetilde{\mathbf{B}}^r$ (resp.  $\widetilde{\mathbf{B}}^r_K$) with respect to the family of norms $\lambda_{s}$ with $s \in  (0, r]$. 

\item Define $\sq{\mathbf C}^{\dagger}$ and $\sq{\mathbf C}_K^{\dagger}$ to be the unions  
$
\sq{\mathbf{C}}^{\dagger}\coloneqq \varinjlim_{r>0} \widetilde{\mathbf{C}}^{r}$ and $ \sq{\mathbf{C}}_K^{\dagger}\coloneqq  \varinjlim_{r>0} \widetilde{\mathbf{C}}^{r}_K.
$ 
\end{enumerate} 
 
\end{definition}

\begin{remark} \label{remark:subscript_K_as_invariant_perfect}
Note that for each $r$, we have  
$\widetilde{\mathbf{A}}_K^{r} =   (\widetilde{\mathbf{A}}^{r})^{H_K} = \widetilde{\mathbf{A}}^{r} \cap \widetilde{\mathbf{A}}_K$ as subrings of $\widetilde{\mathbf{A}}$. Similarly we have 
$\widetilde{\mathbf{B}}_K^{r} =   (\widetilde{\mathbf{B}}^{r})^{H_K}$ and $ \widetilde{\mathbf{C}}_K^{r} =   (\widetilde{\mathbf{C}}^{r})^{H_K}$.
  Consequently, 
we have 
$\widetilde{\mathbf{A}}_K^{\dagger} =   (\widetilde{\mathbf{A}}^{\dagger})^{H_K}$, $\widetilde{\mathbf{B}}_K^{\dagger} =   (\widetilde{\mathbf{B}}^{\dagger})^{H_K},$ and $  \widetilde{\mathbf{C}}_K^{\dagger} =   (\widetilde{\mathbf{C}}^{\dagger})^{H_K}.$  In particular, our definitions agree with the ones in 
 \cite{Berger-differential} and \cite{Colmez08}. 
\end{remark}

\begin{remark}
For each $r > 0$, the natural map $\widetilde{\mathbf{C}}^{r} \rightarrow \varprojlim_{0<s\le r} \widetilde{\mathbf{C}}^{[s,r]}$ is an isomorphism. A similar isomorphism holds for $\widetilde{\mathbf{C}}_K^{r}$ and $\widetilde{\mathbf{C}}_K^{[s,r]}$.
\end{remark}

\begin{remark}[Frobenius structures] \label{remark:Frobenius_classical}
The Frobenius maps on $\widetilde{\mathbf{A}}^+$ and $\widetilde{\mathbf{A}}$ uniquely extend to Frobenius maps on $\widetilde{\mathbf{B}}$ and $\widetilde{\mathbf{B}}^+$, which we denote by $\varphi$. Moreover, for $0 < s \le r$, we have Frobenius maps
\[
\varphi\colon \widetilde{\mathbf{A}}^{r} \xrightarrow[]{\sim} \widetilde{\mathbf{A}}^{r/p}, \qquad \varphi\colon \widetilde{\mathbf{B}}^{r} \xrightarrow[]{\sim} \widetilde{\mathbf{B}}^{r/p}, \qquad \varphi: \widetilde{\mathbf{C}}^{r} \xrightarrow[]{\sim} \widetilde{\mathbf{C}}^{r/p}, \qquad \varphi\colon \widetilde{\mathbf{C}}^{[s, r]}\xrightarrow[]{\sim} \widetilde{\mathbf{C}}^{[s/p,\,r/p]}
\]
extending uniquely the Frobenius on $\widetilde{\mathbf{A}}^+$. Passing to the colimits, we obtain Frobenius maps $\varphi$ on the period rings $\widetilde{\mathbf{A}}^\dagger, \widetilde{\mathbf{B}}^\dagger$, and $\widetilde{\mathbf{C}}^\dagger$. 
Similarly, there are natural Frobenius structures on $\widetilde{\mathbf{A}}_K^r, \widetilde{\mathbf{B}}_K^r$,  $\widetilde{\mathbf{C}}_K^r$, $\widetilde{\mathbf{C}}_K^{[s,r]}$, $\widetilde{\mathbf{A}}_K^\dagger, \widetilde{\mathbf{B}}_K^\dagger$, and $\widetilde{\mathbf{C}}_K^\dagger$.
\end{remark}

\begin{remark} \label{remark:notation_iota_n}
For later use, let us also recall the natural inclusion $\iota_0\colon \widetilde{\mathbf{A}}^r\hookrightarrow \mathbf{B}_{\mathrm{dR}}^+$ (cf. \cite[Proposition~2.11]{Berger-differential}) for any $r\ge 1$. In general, for each $r>0$ and non-negative integer $n\ge -\log_pr$, we have an injection
\[
\iota_n\coloneqq \iota_0\circ \varphi^{-n}\colon \widetilde{\mathbf{A}}^r\xrightarrow[]{\varphi^{-n}} \widetilde{\mathbf{A}}^{p^nr}\xrightarrow[]{\iota_0}\mathbf{B}_{\mathrm{dR}}^+.
\]
By continuity, for each $0<s\le r$ and non-negative integer $n$ such that $p^{-n}\in [s,r]$, the map $\iota_n$ above extends to injections
$\iota_n\colon \sq{\mathbf{C}}^{[s,r]} \ra \mathbf{B}_{\mathrm{dR}}^+, $ and $\iota_n\colon \sq{\mathbf{C}}^r \ra \mathbf{B}_{\mathrm{dR}}^+.$ 
\end{remark}

\subsection{Imperfect period rings} \label{ss:point imperfect period ring}

\subsubsection{Imperfect period rings of type $\mathbf{A}, \mathbf{B}, \mathbf{C}$} 
First,  recall that there are period rings $\mathbf{A} \subset \widetilde{\mathbf{A}}$ and $\mathbf{A}_K\subset \widetilde{\mathbf{A}}_K$, which are both $p$-adically complete and satisfy 
\[ 
\mathbf{A}/(p)=\mathbf{E}, \qquad \mathbf{A}_K/(p)=\mathbf{E}_K\cong \kappa_K(\!(\overline{\pi}_K)\!).
\] 
Furthermore, $\mathbf{A}$ is stable under the Frobenius $\varphi$ and $\Gal_K$-action on $\widetilde{\mathbf{A}}$. We have $(\mathbf{A})^{H_K} = \mathbf{A}_K,$  
which thus carries an action of $\Gamma_K$.  We omit the construction here and refer the reader to \cite[\S 6.1]{Colmez08} for details, noting only that for $K_0 = W(k)[\frac{1}{p}]$, $\mathbf{A}_{K_0}$ coincides with the $p$-adic completion of $W(k) [\![\pi]\!] [\frac{1}{\pi}]$ in $\widetilde{\mathbf{A}}$. Write $\mathbf A_{K_0}^+ \subset \mathbf A_{K_0}$ for the subring $W(k) [\![\pi]\!]$.  

\begin{construction}[{\cite[\S 6.2]{Colmez08}}] \label{construction:pi_K}
We construct a lift $\pi_K\in \mathbf{A}_K $ of $\overline{\pi}_K\in \mathbf{E}_K$ as follows. First consider the minimal polynomial $\overline{P}_K\in \mathbf{E}_{K'_0}^+[X]=\kappa_K [\![\cl \pi]\!][X]$ of $\overline{\pi}_K$ over $\mathbf{E}_{K'_0} = \kappa_K (\!(\cl \pi)\!)$. Choose a monic polynomial $P_K\in \mathbf{A}^+_{K'_0}[X]$ of the same degree that reduces to $\overline{P}_K$ modulo $p$. By Hensel's Lemma, $P_K$ has a unique root $\pi_K$ in $\mathbf{A}_K$ lifting $\overline{\pi}_K$. Moreover, $\mathbf{A}_K$ is a free $\mathbf{A}_{K'_0}$-module with basis $\{1, \pi_K, \ldots, \pi_K^{e_K-1}\}$.
\end{construction} 

We will need following equivalent description of $\mathbf{A}_K$. 
\begin{proposition} \label{prop:A_K_as_completion}
$\mathbf{A}_K$ is the $p$-adic completion of $\mathcal{O}_{K'_0}[\![\pi_K]\!][\frac{1}{\pi_K}]$. That is, we have  
\[\mathbf{A}_K=\left\{\sum_{i\in \mathbb{Z}}c_i\pi_K^i\,\Big|\,c_i\in \mathcal{O}_{K'_0}\; \textrm{ and }c_i\rightarrow 0 \textrm{ $p$-adically as }i\rightarrow -\infty\right\}.\]
\end{proposition}
 
\begin{proof}
This follows from \cite[Proposition~7.5]{Colmez08}.
\end{proof}

\begin{definition} \label{definition:mathbf_A^r_as_intersection} \hfill
\begin{enumerate}
\item For each $r>0$, define 
$\mathbf{A}^{r}\coloneqq \mathbf{A} \cap \widetilde{\mathbf{A}}^{r},$ and $ \mathbf{A}_K^{r}\coloneqq \mathbf{A}_K\cap \widetilde{\mathbf{A}}^{r}_K$   
where the intersection is taken inside $\widetilde{\mathbf{A}}$ and $\widetilde{\mathbf{A}}_K$ respectively.

\item Define 
${\mathbf{A}}^{\dagger} \coloneqq \varinjlim_{r>0}  {\mathbf{A}}^{r}$ and ${\mathbf{A}}^{\dagger}_K \coloneqq \varinjlim_{r>0}  {\mathbf{A}}^{r}_K$, 
where the transition maps are natural inclusions.

\item The period rings ${\mathbf{B}},  {\mathbf{B}}_K$, $\mathbf{B}^r$,  $ {\mathbf{B}}^{r}_K$, $\mathbf{B}^{\dagger}$, and $ {\mathbf{B}}^{\dagger}_K$ are obtained from the corresponding period rings of type $\mathbf{A}$ by inverting $p$.      

\item For each $0 < s \le r$, define ${\mathbf{C}}^{[s, r]}$ (resp.  ${\mathbf{C}}^{[s, r]}_K$) to be the completion of ${\mathbf{B}}^r$ (resp. ${\mathbf{B}}^r_K$) with respect to \(\max\{\lambda_{s}, \lambda_{r}\}\) inside $\widetilde{\mathbf{C}}^{[s, r]}$ (resp. $\widetilde{\mathbf{C}}^{[s, r]}_K$). 

\item For each $r > 0$, define ${\mathbf{C}}^{r}$ (resp. ${\mathbf{C}}^{r}_K$) to be the Fr\'echet completion of  ${\mathbf{B}}^r$ (resp. ${\mathbf{B}}^r_K$) inside $\widetilde{\mathbf{C}}^r$ (resp. $\widetilde{\mathbf{C}}^r_K$) with respect to the family of norms $\lambda_s$ with $s \in (0, r]$. 

\item Lastly, define ${\mathbf C}^{\dagger}$ and ${\mathbf C}_K^{\dagger}$ to be the union  
$
{\mathbf{C}}^{\dagger}\coloneqq \varinjlim_{r>0} {\mathbf{C}}^{r}$ and $
{\mathbf{C}}_K^{\dagger}\coloneqq \varinjlim_{r>0} {\mathbf{C}}^{r}_K.$ 
\end{enumerate}
\end{definition}

\begin{remark} \label{remark:completion_of_A_dagger} \hfill
\begin{enumerate}
\item  Since $\widetilde{\mathbf{A}}_K$ is the $p$-adic completion of $\widetilde{\mathbf{A}}_K^{\dagger}$, we deduce that $\mathbf{A}_K$ is the $p$-adic completion of $\mathbf{A}_K^{\dagger}$. 

\item  For each $r > 0$, the natural map $\mathbf{C}^{r} \rightarrow \varprojlim_{0<s\le r} \mathbf{C}^{[s,r]}$ is an isomorphism. A similar isomorphism holds for $\mathbf{C}^{r}_K$ and $\mathbf{C}_K^{[s,r]}$. 

\item  Remark~\ref{remark:subscript_K_as_invariant_perfect} continues to hold for imperfect period rings defined above. In particular, we have  $ {\mathbf{A}}_K^{\dagger} =   ( {\mathbf{A}}^{\dagger})^{H_K}$, $  {\mathbf{B}}_K^{\dagger} =   ( {\mathbf{B}}^{\dagger})^{H_K}$, and $ {\mathbf{C}}_K^{\dagger} =   ( {\mathbf{C}}^{\dagger})^{H_K}$. 
\end{enumerate}
\end{remark}

\begin{remark}[Topology]\label{remark: topology--classical}
We endow $\widetilde{\mathbf{A}}^r_K$ with the topology induced by the $\lambda_r$-norm and endow $\widetilde{\mathbf{A}}^{\dagger}_K$ with the inductive limit topology (i.e., the finest topology making $\widetilde{\mathbf{A}}^r_K\hookrightarrow \widetilde{\mathbf{A}}^{\dagger}_K$ continuous for all $r>0$). Inverting $p$, we obtain natural topologies on $\widetilde{\mathbf{B}}^r_K$ and $\widetilde{\mathbf{B}}^{\dagger}_K$. The ring $\widetilde{\mathbf{C}}^{[s,r]}_K$ is endowed with the topology induced by $\max\{\lambda_r, \lambda_s\}$, which then induces a Fr\'echet topology on $\widetilde{\mathbf{C}}^r_K$. Finally, we equip $\widetilde{\mathbf{C}}^{\dagger}_K$ with the inductive limit topology (i.e., the LF topology). The topologies on the imperfect period rings are defined similarly.
\end{remark}

\begin{remark}[Frobenius structures]\label{remark: Frobenius--classical}
The Frobenius $\varphi$ on $\widetilde{\mathbf{A}}_K$ (resp. on $\widetilde{\mathbf{A}}^{\dagger}_K$) preserves ${\mathbf{A}}_K$ (resp. ${\mathbf{A}}^{\dagger}_K$). For $0<s\le r$, we also have Frobenius maps
\[\varphi: \mathbf{A}^r_K\rightarrow \mathbf{A}^{r/p}_K, \qquad\varphi: \mathbf{B}^r_K\rightarrow \mathbf{B}^{r/p}_K, \qquad\varphi: \mathbf{C}^r_K\rightarrow \mathbf{C}^{r/p}_K, \qquad\varphi: \mathbf{C}^{[s,r]}_K\rightarrow \mathbf{C}^{[s/p, \,r/p]}_K\]
as well as Frobenius maps on $\mathbf{B}^{\dagger}_K$ and $\mathbf{C}^{\dagger}_K$. We caution the reader that $\varphi: \mathbf{A}^r_K\rightarrow \mathbf{A}^{r/p}_K$ is no longer bijective.
\end{remark}

\begin{remark}[{$\Gamma_K$-actions}]\label{remark: Gamma action--classical}
All period rings $\widetilde{\mathbf{A}}^r_K$, $\widetilde{\mathbf{A}}^{\dagger}_K$, $\widetilde{\mathbf{B}}^r_K$, $\widetilde{\mathbf{B}}^{\dagger}_K$, $\widetilde{\mathbf{C}}^{[s,r]}_K$, $\widetilde{\mathbf{C}}^r_K$, $\widetilde{\mathbf{C}}^{\dagger}_K$, as well as their imperfect counterparts, are equipped with continuous actions of $\Gamma_K$, where $\Gamma_K$ is endowed with the profinite topology. The $\Gamma_K$-actions commute with the Frobenius actions.
\end{remark}

\begin{notation}
We let $t = \log (1 + \pi)$. This is an element that lies in $\mathbf{C}^{\dagger}_K$ but not in $\mathbf{B}^{\dagger}_K$.
\end{notation}

\begin{lemma}\label{lemma:pi vs bar-pi}
Let
\[
r_K =\begin{cases}
    (2v(\mathfrak{d}_{\mathbf{E}_K/\mathbf{E}_{K_0}}))^{-1} & \text{if } \mathbf{E}_K/\mathbf{E}_{K_0} \text{ is ramified,}\\
    1 & \text{if } \mathbf{E}_K/\mathbf{E}_{K_0} \text{ is unramified}
    \end{cases}
\]
where $\mathfrak{d}_{\mathbf{E}_K/\mathbf{E}_{K_0}}$ denotes the discriminant of $\mathbf{E}_K/\mathbf{E}_{K_0}$. Then for each $0<r<r_K$, we have $\pi_K \in \mathbf{A}_K ^r$ and the element $u=\frac{\pi_K}{[\overline{\pi}_K]}$ is a unit in the ring of integers of $\widetilde{\mathbf{A}}_K^r$ (with respect to the norm $\lambda_{r}$). More precisely, if we write 
$u=\sum_{n=0}^{\infty}p^n[u_n]$  with $u_n\in \widehat{K}_{\infty}^{\flat}$, then $u_0=1$ and $rv(u_n)+n>0$ for all $n\ge 1$. 
\end{lemma}

\begin{proof}
This is \cite[Lemma~6.5(1)]{Colmez08} and its proof in \textit{loc. cit}. 
\end{proof}

\begin{construction} \label{construction:section_s_from_E_K}  For later use, consider a section $s\colon \mathbf{E}_K\rightarrow \mathbf{A}_K$ of the mod  $p$ reduction map given by
\[
s\colon \sum_{i\in \mathbb{Z}} a_i \overline{\pi}_K^i  \mapsto  \sum_{i\in \mathbb{Z}} [a_i] \pi_K^i.
\]
Note that for every $x\in \mathbf{A}_K$, we can construct a sequence $x_0, x_1, \ldots\in \mathbf{A}_K$ by setting $x_0=x$ and $x_{n+1}\coloneqq \frac{1}{p}(x_n-s(\overline{x}_n))$, where $\overline{x}_n\in \mathbf{E}_K$ denotes the mod $p$ reduction of $x_n$. Then we have $x=\sum_{n=0}^{\infty} p^n\cdot s(\overline{x}_n)$.

\begin{lemma} \label{lemma:section}
If $y\in \mathbf{E}_K$ and $0<r<r_K$ (cf. Lemma~\ref{lemma:pi vs bar-pi}), then $s(y)\in \mathbf{A}_K^{r}$ and $v_r(s(y))=v(y)$.
\end{lemma}

\begin{proof}
This is \cite[Lemma~7.2]{Colmez08}.
\end{proof}

\begin{remark} \label{remark:lambda_r_of_pi_K}
In particular, for $0 < r < r_K$, we have $v_r (\pi_K) = v(\overline{\pi}_K) = \frac{p}{p-1} \cdot \frac{1}{e_K}$ which is independent of $r$, and $\lambda_r (\pi_K) = p^{-r \cdot v (\overline \pi_K)} = p^{-\frac{p}{p-1} \cdot \frac{r}{e_K}}$.
\end{remark}
\end{construction}

\subsubsection{Explicit (geometric) descriptions}
We have the following explicit description of $\mathbf{A}_K^{r}$.

\begin{proposition} \label{prop:identificaiton_for_mathbf_A^r}
If $0<r<r_K$, there is an isomorphism
\[
\left\{\sum_{i\in \mathbb{Z}}c_iT^i\,\Big|\,c_i\in \mathcal{O}_{K'_0}, \,\,\frac{v(c_i)}{r}+ \frac{p}{p-1} \cdot \frac{i}{e_K}\rightarrow +\infty \textrm{ as }i\rightarrow -\infty\right\}\xrightarrow[]{\sim} \mathbf{A}_K^r
\]
sending $f\mapsto f(\pi_K)$.\footnote{The valuation is normalized so that $v(p)=1$} Moreover, the isomorphism preserves the valuations on both sides, where on the left hand side the valuation is $\inf_{i\in \mathbb{Z}} \left\{\frac{v(c_i)}{r}+ \frac{p}{p-1} \cdot \frac{i}{e_K} \right\},$ while on the right hand side the valuation is $v_r$.
\end{proposition}

\begin{proof}
This is \cite[Proposition~II.2.1]{Cherbonnier_Colmez} (also see  \cite[Proposition~1.4]{Berger-differential} and \cite[Proposition~7.5(i)]{Colmez08}).
\end{proof}

\begin{remark}
As an immediate corollary, for each $0<s\le r$, the valuation $\inf_{i\in \mathbb{Z}}\{\frac{v(c_i)}{s}+ \frac{p}{p-1} \cdot  \frac{i}{e_K}\}$ on the left hand side agrees with the valuation $v_{s}$ on the right hand side.
\end{remark}

\begin{remark} \label{remark:identificaiton_for_A^r}
For later use, it is helpful to rewrite the isomorphism in Proposition~\ref{prop:identificaiton_for_mathbf_A^r} as 
  \[
 \left\{\sum_{i\in \mathbb{Z}}c_iT^i\,\Big|\,c_i\in \mathcal{O}_{K'_0}, \,\,|c_i| \cdot \alpha^i \rightarrow 0 \textrm{ as }i\rightarrow -\infty\right\} \xrightarrow[]{\sim} \mathbf{A}_K^r, \]
where $\alpha = p^{- \frac{p}{p-1} \cdot \frac{r}{e_K}}$. Then the isomorphism preserves the norms on both sides, where the norm on the left hand side is $\sup_{i\in \mathbb{Z}} \left\{|c_i|\cdot \alpha^r\right\}$ and the one on the right hand side is $\lambda_r$. 

\end{remark}

\begin{remark} \label{remark:completion_of_A_dagger_2} Another immediate consequence is that, for $0 < r < r_K$, the isomorphism in Proposition \ref{prop:identificaiton_for_mathbf_A^r} induces an isomorphism
\[ 
\mathbf{A}^r_K/(p) \cong \mathbf{A}^\dagger_K/(p) \cong \mathbf{A}_K/(p) \cong \kappa_K (\!(\overline \pi_K)\!) \cong \mathbf{E}_K  
\]
In particular, $\mathbf{A}_K$ is the $p$-adic completion of $\mathbf{A}^r_K$ for any $r \in (0, r_K)$. This gives another way to recover $\mathbf{A}_K$ as the $p$-adic completion of $\mathbf{A}^\dagger_K$ (also see Remark~\ref{remark:completion_of_A_dagger}). 
\end{remark}

Next, we turn to explicit descriptions of $\mathbf{C}_K^r$ and $\mathbf{C}_K^{[s,r]}$.  It is convenient to introduce the following notation for functions on annuli. 

\begin{notation} \label{notn: annuli}
Let $L$ be a complete non-archimedean field equipped with norm $|\cdot|$. We use $\D_L$ to denote the rigid analytic closed unit disk over $L$. For an interval $I \subset \R_{> 0}$, we denote by $A_L (I)$ the rigid analytic annulus over $L$ centered at $0$ with radius $|T| \in I$. In particular, for $\alpha, \beta \in p^{\Q} \cap \R_{>0}$ with $\alpha \le \beta$, $A_L (\alpha, \beta)$ (resp.  $A_L {[\alpha, \beta]}$) is the open (resp. closed) annulus  $\alpha < |T| < \beta $ (resp. $\alpha \le |T| \le \beta $). We write $\Gamma_L (I) = \mO(A_L (I))$ for the ring of analytic functions on $A_L (I)$.Sometimes we emphasize that the coordinate for the annulus in consideration is $T$ by writing $A_L (I; T)$ and $\Gamma_L (I; T)$ instead of $A_L(I)$ and $\Gamma_L (I)$. 
\end{notation}

\begin{proposition}\label{prop:identification_C_with_annulus}
For $0 < s \le r$ sufficiently close to $0$, the isomorphism in Proposition~\ref{prop:identificaiton_for_mathbf_A^r} induces an identification \[\left\{\sum_{i\in \mathbb{Z}}c_iT^i\,\Big|\,c_i\in K'_0, \,\,v(c_i)+\frac{p}{p-1}\cdot\frac{it}{e_K}\rightarrow +\infty \textrm{ as }i\rightarrow \pm\infty, \textrm{ for all } t\in[s,r]\right\}\xrightarrow[]{\sim} \mathbf{C}_K^{[s,r]}\]
sending $f\mapsto f(\pi_K)$, where the left hand side can be re-written as
\[\left\{\sum_{i\in \mathbb{Z}}c_iT^i\,\Big|\,c_i\in K'_0, \,\,|c_i| \cdot p^{- \frac{p}{p-1} \cdot \frac{it}{e_K}}\rightarrow 0 \textrm{ as }i\rightarrow \pm\infty, \textrm{ for all } t\in [s,r] \right\}.\]
In other words, we can identify $\mathbf{C}_K^{[s,r]}$ with the ring of analytic functions on the closed annulus $\Gamma_{K'_0}([\alpha, \beta]; \pi_K)$ where $\alpha = p^{- \frac{p}{p-1} \cdot \frac{r}{e_K}}, \beta = p^{-\frac{p}{p-1} \cdot \frac{s}{e_K}}$. Taking limit as $s\rightarrow 0^+$, we obtain an identification
\[ \Gamma_{K'_0}([\alpha, 1); \pi_K) \isom \mathbf{C}_K^r\]  
for $r$ sufficiently small; namely, $\mathbf{C}_K^r$ can be identified as the ring of functions on an half-open annulus with open outer boundary at radius $1$ and closed inner boundary with radius $\alpha = p^{- \frac{p}{p-1} \cdot \frac{r}{e_K}}$. Consequently, taking a colimit as $r\rightarrow 0^+$, the ring $\mathbf{C}_K^{\dagger}$ can be identified with the classical Robba ring with coefficients in $K'_0$. 
\end{proposition}

\begin{proof}
Both statements follow from Proposition~\ref{prop:identificaiton_for_mathbf_A^r} by \cite[Proposition~7.5 \& 7.6]{Colmez08} and their proofs. 
\end{proof}

\begin{remark} Let $\alpha = p^{- \frac{p}{p-1} \cdot \frac{r}{e_K}}, \beta = p^{-\frac{p}{p-1} \cdot \frac{s}{e_K}}$ be as above. 
For any $t \in [s, r]$, the norm $\lambda_t$ on $\mathbf{C}_K^{[s, r]}$ coincides with the $\rho$-Gauss norm on $\Gamma_{K'_0}([\alpha, \beta]; \pi_K)$  for $\rho = p^{- \frac{p}{p-1} \cdot \frac{t}{e_K}}$ under the isomorphism in Proposition~\ref{prop:identification_C_with_annulus}. 
\end{remark}

\subsection{$(\varphi, \Gamma_K)$-modules and $p$-adic differential equations} \label{ss:classical_pDE}

We now recall the theory of $(\varphi,\Gamma_K)$-modules associated with $p$-adic Galois representations of $\Gal_K$, together with Berger's construction of the $p$-adic differential equations associated with \emph{de Rham} representations (cf. \cite{Berger-differential}).

\subsubsection{From $p$-adic Galois representations to $(\varphi, \Gamma_K)$-modules} \label{sss:overconvergent classical}

\begin{definition}\label{defn: phi Gamma modules--classical case}\hfill
\begin{enumerate}   
\item A \emph{$(\varphi, \Gamma_K)$-module}\footnote{Note that our notion of $(\varphi,\Gamma_K)$-modules over $\mathbf{A}^{\dagger}_K$ is referred to as \emph{\'etale $(\varphi, \Gamma_K)$-modules} in \cite{Cherbonnier_Colmez}.} over $\mathbf{A}_K^{\dagger}$ is a finite free $\mathbf{A}_K^{\dagger}$-module $M$ equipped with a Frobenius semilinear action $\varphi$ that induces an isomorphism  
\[ 
\varphi \otimes \mathrm{id}\colon M \otimes_{\mathbf{A}_K^{\dagger}, \varphi} \mathbf{A}_K^{\dagger} \isom M,
\] 
together with a continuous semilinear $\Gamma_K$-action that commutes with $\varphi$. The notions of \emph{$(\varphi, \Gamma_K)$-modules} over $\mathbf{B}_K^{\dagger}$ and $\mathbf{C}_K^{\dagger}$ are defined in the same way.
\item  A $(\varphi, \Gamma_K)$-module $M$ over either $\mathbf{B}_K^{\dagger}$ or  $\mathbf{C}_K^{\dagger}$ is called \emph{\'etale} if there exists a $(\varphi, \Gamma_K)$-module $M_0$ over $\mathbf{A}_K^{\dagger}$ such that $M$ is the base change of $M_0$ as a $(\varphi, \Gamma_K)$-module.  
\end{enumerate}
\end{definition}

\begin{remark}\label{rmk: descent to r}
Giving a $(\varphi, \Gamma_K)$-module $M$ over $\mathbf{A}^{\dagger}_K$ is equivalent to giving a compatible family of finite free $\mathbf{A}^r_K$-modules $M^r$ (for all sufficiently small $r>0$), each of which is equipped with a continuous semilinear $\Gamma_K$-action together with an isomorphism
$\varphi\colon M^r \otimes_{\mathbf{A}_K^r, \varphi} \mathbf{A}_K^{r/p}\isom M^r \otimes_{\mathbf{A}_K^r} \mathbf{A}_K^{r/p}$
commuting with the $\Gamma_K$-actions, where the second tensor is along the natural inclusion $\mathbf{A}_K^r\subset \mathbf{A}_K^{r/p}$. Here, compatibility means that, for $0<r<r'$, the $\mathbf{A}_K^r$-module $M^{r}$ (together with the $\varphi$- and $\Gamma_K$-structures) is the base change of $M^{r'}$ along $\mathbf{A}^{r'}_K\rightarrow \mathbf{A}^r_K$. There are similar statements for $(\varphi, \Gamma_K)$-modules over $\mathbf{B}^{\dagger}_K$ and $\mathbf{C}^{\dagger}_K$.
\end{remark}

Let $\mathrm{Rep}_{\mathbb{Q}_p}(\Gal_{K})$ (reps. $\mathrm{Rep}_{\mathbb{Z}_p}(\Gal_{K})$) denote the category of continuous $\Gal_K$-representations over finite dimensional $\mathbb{Q}_p$-vector spaces (resp. finite free $\mathbb{Z}_p$-modules). Such objects are often referred to as \emph{$p$-adic Galois representations}. 

\begin{definition}\hfill
\begin{enumerate}
\item For every $V\in \mathrm{Rep}_{\mathbb{Z}_p}(\Gal_{K})$, consider the $\mathbf{A}_K$-module
$D^{\mathrm{int}}_{K}(V)\coloneqq \big(V\otimes_{\mathbb{Z}_p}\mathbf{A}\big)^{H_K}$
and the $\mathbf{A}^{\dagger}_K$-module
$D^{\mathrm{int}, \dagger}_{K}(V)\coloneqq \big(V\otimes_{\mathbb{Z}_p}\mathbf{A}^{\dagger}\big)^{H_K},$
both equipped with natural $\varphi$- and $\Gamma_K$-actions.  
\item For every $V\in \mathrm{Rep}_{\mathbb{Q}_p}(\Gal_{K})$, consider the $\mathbf{B}^{\dagger}_K$-module
$D^{\dagger}_{K}(V)\coloneqq \big(V\otimes_{\mathbb{Q}_p}\mathbf{B}^{\dagger}\big)^{H_K}$
and the $\mathbf{C}^{\dagger}_K$-module
$
D^{\dagger}_{\mathrm{rig},K}(V)\coloneqq D^{\dagger}_{K}(V)\otimes_{\mathbf{B}^{\dagger}_{K}} \mathbf{C}^{\dagger}_{K},  
$
both equipped with natural $\varphi$- and $\Gamma_K$-actions.  
\end{enumerate}
\end{definition}

\begin{theorem}\label{thm: Galois representations and phi Gamma modules}
The functors $D^{\mathrm{int}}_{K}$ and $D^{\mathrm{int}, \dagger}_{K}$ induce equivalences of categories
\begin{align*}
\: \: \:  \mathrm{Rep}_{\mathbb{Z}_p}(\Gal_{K}) & \isom \big\{(\varphi, \Gamma_{K})\textrm{-modules over }\mathbf{A}_{K}\big\} 
 \\
& \isom \big\{(\varphi, \Gamma_{K})\textrm{-modules over }\mathbf{A}^{\dagger}_{K}\big\},
\end{align*} 
while the functors $D^{\dagger}_K$ and $D^{\dagger}_{\mathrm{rig}, K}$ induce equivalences of categories
\begin{align*}
\: \: \:  \mathrm{Rep}_{\mathbb{Q}_p}(\Gal_{K}) & \isom \big\{\textrm{\'etale }(\varphi, \Gamma_{K})\textrm{-modules over }\mathbf{B}^{\dagger}_{K}\big\} 
 \\
& \isom \big\{\textrm{\'etale }(\varphi, \Gamma_{K})\textrm{-modules over }\mathbf{C}^{\dagger}_{K}\big\}.
\end{align*} 
\end{theorem}

\begin{proof}
See \cite[Theorem~9.5.4 \& 9.5.6]{Kedlaya-LiuI}.
\end{proof}

\begin{remark}\label{remark: lattice in etale phi Gamma modules--classical case}
Let $V\in \mathrm{Rep}_{\mathbb{Q}_p}(\Gal_{K})$ and let $V^{\circ}$ be any $\Gal_K$-stable $\Z_p$-lattice inside $V$. If $M$ denotes either the \'etale $(\varphi, \Gamma_K)$-module $D^{\dagger}_K(V)$ or $D^{\dagger}_{\mathrm{rig}, K}(V)$, then we can take $M_0$ (as in Definition~\ref{defn: phi Gamma modules--classical case}(2)) to be $D^{\mathrm{int}, \dagger}_K(V^{\circ})$.
\end{remark}

\begin{remark}\label{rmk: phi Gamma modules over perfect period rings--classical}
There are also $(\varphi, \Gamma_K)$-module theories over perfect period rings. Indeed, by \cite[Theorem~9.5.4 \& 9.5.6]{Kedlaya-LiuI}, there are equivalences of categories
\begin{align*}
\: \: \:  \mathrm{Rep}_{\mathbb{Z}_p}(\Gal_{K}) & \isom \big\{(\varphi, \Gamma_{K})\textrm{-modules over }\widetilde{\mathbf{A}}^{\dagger}_{K}\big\} 
 \\
& \isom \big\{(\varphi, \Gamma_{K})\textrm{-modules over }\widetilde{\mathbf{A}}_{K}\big\}
\end{align*} 
and
\begin{align*}
\: \: \:  \mathrm{Rep}_{\mathbb{Q}_p}(\Gal_{K}) & \isom \big\{\textrm{\'etale }(\varphi, \Gamma_{K})\textrm{-modules over }\widetilde{\mathbf{B}}^{\dagger}_{K}\big\} 
 \\
& \isom \big\{\textrm{\'etale }(\varphi, \Gamma_{K})\textrm{-modules over }\widetilde{\mathbf{C}}^{\dagger}_{K}\big\}
\end{align*} 
For $V\in \mathrm{Rep}_{\mathbb{Z}_p}(\Gal_{K})$ or $\mathrm{Rep}_{\mathbb{Q}_p}(\Gal_{K})$, let $\widetilde{D}^{\mathrm{int}}_K(V)$, $\widetilde{D}^{\mathrm{int}, \dagger}_K(V)$, $\widetilde{D}^{\dagger}_K(V)$, $\widetilde{D}^{\dagger}_{\mathrm{rig}, K}(V)$ denote the associated $(\varphi, \Gamma_K)$-modules over $\widetilde{\mathbf{A}}_K$, $\widetilde{\mathbf{A}}^{\dagger}_K$, $\widetilde{\mathbf{B}}^{\dagger}_K$, $\widetilde{\mathbf{C}}^{\dagger}_K$, respectively. In particular, we have
\begin{enumerate}
\item for any $V\in \mathrm{Rep}_{\mathbb{Z}_p}(\Gal_{K})$, the $(\varphi, \Gamma_K)$-module $\widetilde{D}^{\mathrm{int}, \dagger}_K(V)$ is the base change of $D^{\mathrm{int}, \dagger}_K(V)$ along $\mathbf{A}^{\dagger}_K \rightarrow \widetilde{\mathbf{A}}^{\dagger}_K$;
\item for any $V\in \mathrm{Rep}_{\mathbb{Q}_p}(\Gal_{K})$, the \'etale $(\varphi, \Gamma_K)$-module $\widetilde{D}^{\dagger}_K(V)$ (resp. $\widetilde{D}^{\dagger}_{\mathrm{rig}, K}(V)$) is the base change of $\widetilde{D}^{\mathrm{int}, \dagger}_K(V^{\circ})$ along $\widetilde{\mathbf{A}}^{\dagger}_K\rightarrow \widetilde{\mathbf{B}}^{\dagger}_K$ (resp. $\widetilde{\mathbf{A}}^{\dagger}_K\rightarrow \widetilde{\mathbf{C}}^{\dagger}_K$)
where $V^{\circ}$ is any $\Gal_K$-stable $\Z_p$-lattice in $V$; 
\item $\widetilde{D}^{\dagger}_K(V)$ (resp. $\widetilde{D}^{\dagger}_{\mathrm{rig}, K}(V)$) is the base change of $D^{\dagger}_K(V)$ (resp. $D^{\dagger}_{\mathrm{rig}, K}(V)$) along $\mathbf{B}^{\dagger}_K\rightarrow \widetilde{\mathbf{B}}^{\dagger}_K$ (resp. $\mathbf{C}^{\dagger}_K\rightarrow \widetilde{\mathbf{C}}^{\dagger}_K$).
\end{enumerate}
We refer the reader to \cite{Cherbonnier_Colmez} and \cite[Theorem~9.5.4 \& 9.5.6]{Kedlaya-LiuI} for more details.
\end{remark}

\subsubsection{Differential modules} \label{sss:pDE_over_K}

Before proceeding further, we recall the notion of differential rings and differential modules. First, a \emph{differential ring} is a commutative ring $R$ equipped with a derivation $d\colon R \ra R$ (cf. \cite[\S 5]{kedlaya_book}), that is, an additive map satisfying the Leibniz rule $d (ab) = a d(b) + b d(a)$ for all $a, b \in R$. 
A \emph{differential module} over a differential ring $(R, d)$ is an $R$-module $M$ equipped with a differential operator $D\colon M \ra M$ compatible with $d$, that is, an additive map satisfying $D(ax) = a D(x) + d(a) x$ for all $a \in R, x \in M$. When $R$ is equipped with more than one differential operator $d_1, \ldots, d_m$, we abuse the terminology and say that $(M, D_1, \ldots, D_m)$ is a \emph{differential module} over $(R, d_1, \ldots, d_m)$ if $(M, D_i)$ is a differential module over $(R, d_i)$ for all $i=1,\ldots, m$.

We are primarily interested in the Robba ring $\mathbf{C}^{\dagger}_K$ on which the action of the Lie algebra of $\Gamma_K$ induces a differential operator $\nabla$ defined by $\nabla = \log (\gamma)/\log (\chi (\gamma))$ for any $\gamma$ sufficiently close to $1$, where $\chi\colon \Gamma_K\rightarrow \Z_p^{\times}$ is the usual cyclotomic character. By \cite[Lemme~4.2]{Berger-differential}, one has 
\begin{equation} \label{eq:nabla_vs_partial}
\nabla = t \cdot \partial = t \cdot (1+\pi) \cdot \frac{d}{d\pi}  
\end{equation} 
where $t=\log (1+\pi)$ and $\partial =  (1+\pi) \cdot \frac{d}{d\pi}$ is a well-defined differential operator on $\mathbf{C}^{\dagger}_K$. Then a differential module $(M, \partial_M)$ over $(\mathbf{C}^{\dagger}_K, \partial)$ can be viewed as a connection 
\[
\nabla_{M}\colon M \ra M \otimes_{\mathbf{C}^{\dagger}_K} \Omega^1_{\mathbf{C}^{\dagger}_K}
\]
where $\nabla_M (x) = \partial_M (x) \cdot d \log (1+ \pi) =  \frac{\partial_M (x)}{1 + \pi} \cdot d \pi$.\footnote{Note that this is a connection extending $d x = \frac{\partial x}{1+ \pi} d\pi$.} A finite free differential module over $(\mathbf{C}^{\dagger}_K, \partial)$ is referred to as a \emph{$p$-adic differential equation} in \cite{Berger-differential}.

\begin{remark}\label{remark:partial vs d/d pi_K}
In practice, one may replace $(\mathbf{C}^{\dagger}_K, \partial)$ with the differential ring $(\mathbf{C}^{\dagger}_K, \frac{d}{d\pi_K})$ without changing the underlying theory of $p$-adic differential equations, as $\partial=u_0\cdot \frac{d}{d\pi_K}$ for some nonzero element $u_0\in \mathbf{B}^{\dagger}_K$ (hence, a unit in $\mathbf{C}^{\dagger}_K$). 
\end{remark} 

\begin{definition} \label{def:p-adic_diff_eqn_with_Frob}
 Let $(M, \partial_M)$ be a finite free differential module over $(\mathbf{C}^{\dagger}_K, \partial)$.  A \emph{Frobenius structure} on $M$ is a $\varphi$-semilinear map $\varPhi\colon M \ra  M$ such that 
\begin{enumerate}
    \item[(1)] $\varPhi$ induces an isomorphism $\varphi^*(M) \isom M$, i.e.,  $\varPhi(M)$ generates $M$ as a $\mathbf{C}_{K}^{\dagger}$-module; 
    \item[(2)] $\varPhi$ is compatible with the differential structures; namely, the diagram
     \[ \begin{tikzcd}
         M \arrow[r, "\nabla_M"] \arrow[d, swap, "\varPhi"] & M \otimes_{\mathbf{C}^{\dagger}_K} \Omega^1_{\mathbf{C}^{\dagger}_K} \arrow[d, "\varPhi \otimes d \varphi"]  \\
         M \arrow[r, "\nabla_M"]  & M  \otimes_{\mathbf{C}^{\dagger}_K} \Omega^1_{\mathbf{C}^{\dagger}_K}
     \end{tikzcd} \]
    commutes. Equivalently, $\varPhi$ satisfies $\partial_M \circ \varPhi = p \,\varPhi \circ \partial_M$.    
\end{enumerate}
\end{definition}

\subsubsection{$p$-adic differential equations attached to de Rham representations} \label{sss:Berger_N_dR}

Let $D^{\dagger}_{\mathrm{rig}, K} (V)$ be the \'etale $(\varphi, \Gamma_K)$-module over  $\mathbf{C}^{\dagger}_K$ associated with some $V\in \mathrm{Rep}_{\mathbb{Q}_p}(\Gal_{K})$.
The action of the Lie algebra of $\Gamma_K$ induces a differential operator $\nabla_V =   \log (\gamma)/\log (\chi (\gamma))$ on $D^{\dagger}_{\mathrm{rig}, K} (V)$. To obtain a differential module over $(\mathbf{C}_{K}^{\dagger}, \partial)$, we consider 
\[ 
\partial_V = \frac{1}{t} \nabla_V = \frac{1}{\log (1+\pi)} \nabla_V.
\] 
In general, this operator does not preserve $D^{\dagger}_{\mathrm{rig}, K} (V)$, as it might introduce poles at the vanishing locus of $t$. The main technical result of \cite{Berger-differential} states that, when $V$ is de Rham, these poles can be removed systematically to produce an actual $p$-adic differential equation over $(\mathbf{C}^{\dagger}_K, \partial)$. 

\begin{theorem}[Berger] \label{thm:Berger_N_dR}
    Let $V$ be a de Rham representation of $\Gal_K$ of dimension $d$ with non-negative Hodge--Tate weights.\footnote{For general de Rham representations, one can first twist the representations by powers of cyclotomic character to make the Hodge--Tate weights non-negative.} Then there exists a unique free $\mathbf{C}^{\dagger}_K$-submodule $\mathbf{N}_{\mathrm{dR}} (V) \subset D^{\dagger}_{\mathrm{rig}, K} (V)$ of rank $d$ that is stable under $\partial_{V}$.  Moreover, $\mathbf{N}_{\mathrm{dR}}(V)$ is stable under $\varphi$ and $ \Gamma_K$-actions,\footnote{Note that, as a $(\varphi, \Gamma_K)$-module, $\mathbf{N}_{\mathrm{dR}}(V)$ is not necessarily \'etale.} and, consequently, it is a finite free differential module over $(\mathbf{C}^{\dagger}_K, \partial)$ with a Frobenius structure. Furthermore, the inclusion $\mathbf{N}_{\mathrm{dR}} (V) \subset D^{\dagger}_{\mathrm{rig}, K} (V)$ induces an isomorphism after inverting $t$.
\end{theorem}

\begin{proof}
This is {\cite[Th\'eor\`eme 5.10]{Berger-differential}}. 
\end{proof}

\subsection{$p$-adic differential equations via Beauville--Laszlo gluing} \label{ss:NdR classical setting} 

To prepare for later generalizations to the relative setting, here we provide an interpretation of Berger's $\mathbf{N}_{\textnormal{dR}}(V)$ that works for all de Rham representations $V$ of $\Gal_K$.  When $V$ has non-negative Hodge--Tate weights, this construction coincides with the one in \cite{Berger-differential}.

\subsubsection{Beauville--Laszlo gluing} \label{sss:BL_glue}

One key ingredient we need is \emph{Beauville--Laszlo gluing} for finite projective modules. This is well-known to experts, but we briefly recall the notion below for completeness.

\begin{definition}\label{definition: BL data}
A \emph{Beauville--Laszlo datum} is a triple $(A, A_f, \widehat{A})$ where $A$ is a commutative ring, $f\in A$ is a non-zero-divisor, $A_f\coloneq A[1/f]$, and $\widehat{A}\coloneqq \varprojlim_n A/f^n A$.
\end{definition}

\begin{proposition}[{Beauville--Laszlo gluing \cite{beauvillelaszlo}}]\label{prop:BL gluing}
Let $(A, A_f, \widehat{A})$ be a Beauville--Laszlo datum, and let \[\widehat{A}_f\coloneqq \widehat{A}[\frac{1}{f}]\cong \widehat{A}\otimes_AA_f.\] Then the functor that sends a finite projective $A$-module $M$ to the triple $(M\otimes_A A_f, M\otimes_A \widehat{A}, \psi_0)$, where $\psi_0$ is the canonical identification, induces an equivalence of categories between the category of finite projective $A$-modules and the category of triples $(M_1, M_2, \psi)$ where
\begin{itemize}
\item $M_1$ is a finite projective $A_f$-module;
\item $M_2$ is a finite projective $\widehat{A}$-module of the same rank;
\item $\psi$ is an isomorphism
$
M_1\otimes_{A_f}\widehat{A}_f\cong M_2\otimes_{\widehat{A}}\widehat{A}_f$ 
of $\widehat{A}_f$-modules.
\end{itemize}
We will refer to such a triple as a \emph{Beauville--Laszlo triple} (over a fixed Beauville--Laszlo datum). 
\end{proposition}

\subsubsection{Construction of $\mathbf{N}_{\textnormal{dR}}(V)$} \label{sss: NdR classical setting}

Let $\mathrm{Rep}^{\mathrm{dR}}_{\mathbb{Q}_p}(\Gal_{K})$ denote the full subcategory of de Rham representations in $\mathrm{Rep}_{\mathbb{Q}_p}(\Gal_{K})$.

\begin{construction} \label{construction:N_dR_for_all_HT_wts}     
Let $V\in \mathrm{Rep}^{\mathrm{dR}}_{\mathbb{Q}_p}(\Gal_{K})$ be a de Rham representation and let $M^{\dagger}\coloneqq D^{\dagger}_{\mathrm{rig}, K}(V)$ be the associated \'etale $(\varphi, \Gamma_K)$-module over $\mathbf{C}^{\dagger}_K$. Then $M^{\dagger}$ corresponds to a compatible family of finite projective modules $M^{[s,r]}$ over $\mathbf{C}^{[s,r]}_K$ equipped with compatible Frobenius isomorphisms 
\begin{equation} \label{eq:frob_isom}
\psi=\psi^{[s, r]}\colon M^{[s, r]} \otimes_{\mathbf{C}^{[s,r]}_K, \varphi} \mathbf{C}^{[s/p, \,r/p]}_K \isom M^{[s/p, \,r/p]} 
\end{equation}
for all $r$ sufficiently small and all $0<s\le r$. We will modify each $M^{[s,r]}$ along the ``$t=0$ locus'' to obtain another compatible system of finite projective modules $N^{[s,r]}$ which glue to a $(\varphi, \Gamma_K)$-module $N^{\dagger}$ over $\mathbf{C}^{\dagger}_K$. (When $[s,r]\cap p^{\mathbb{Z}}=\emptyset$, we simply put $N^{[s,r]}=M^{[s,r]}$ as $t$ is invertible on $\mathbf{C}^{[s,r]}_K$.) Let us split the construction into several steps. \\

\noindent \textbf{Step 1.} We first consider the situation where $\#\left([s,r]\cap p^{\mathbb{Z}}\right)=1$, say $p^{-n}\in [s,r]$. Recall the injective map 
$
\iota_n\colon \widetilde{\mathbf{C}}^{[s,r]}\hookrightarrow \mathbf{B}^+_{\mathrm{dR}}
$
from Remark \ref{remark:notation_iota_n}. It follows from \cite[Lemme~5.11]{Berger-differential} that, when $r$ is sufficiently small, the image of $\mathbf{C}^r_K$ under $\iota_n$ lies in $K_n[\![t]\!]$. The resulting map $\iota_n\colon \mathbf{C}^r_K\hookrightarrow K_n[\![t]\!]$ has dense image with respect to the $t$-adic topology. By continuity, this map extends to $\iota_n\colon \mathbf{C}^{[s,r]}_K\hookrightarrow K_n[\![t]\!]$, again with dense image. Consequently, the triple $(\mathbf{C}^{[s,r]}_K, \mathbf{C}^{[s,r]}_K[1/t], K_n[\![t]\!])$ is a Beauville--Laszlo datum (cf. Definition~\ref{definition: BL data}). We shall consider the finite projective $\mathbf{C}^{[s,r]}_K[1/t]$-module $M^{[s,r]}[1/t]$ and the finite projective $K_n[\![t]\!]$-module $D_{\mathrm{dR}}(V)\otimes_K K_n[\![t]\!]$, and then glue them into a finite projective $\mathbf{C}^{[s,r]}_K$-module $N^{[s,r]}$ via Beauville--Laszlo gluing.\\

\noindent \textbf{Step 2.} By Fontaine's construction of $D_{\mathrm{dif}}(V)$, there is an equality 
$
D_{\mathrm{dR}}(V) \otimes_K K_\infty (\!(t)\!) = D_{\mathrm{dif}} (V)
$ 
as $K_\infty (\!(t)\!)$-subspaces of $D_{\mathrm{dR}}(V) \otimes_K \mathbf{B}_{\mathrm{dR}}$. By \cite[Proposition~5.7 \& Corollaire~5.8]{Berger-differential}, we have an isomorphism 
\[ M^{[s,r]} \otimes_{\mathbf{C}^{[s,r]}_K, \iota_n }K_{\infty}(\!(t)\!) \isom D_{\mathrm{dif}} (V)\]
for $r$ sufficiently small. Putting these together, we obtain a $\Gamma_{K}$-equivariant isomorphism 
\begin{equation}
 \label{eq:psi_infty_iso}
\iota_n\colon M^{[s,r]} \otimes_{\mathbf{C}^{[s,r]}_K, \iota_n }K_{\infty}(\!(t)\!) \isom D_{\mathrm{dR}}(V)\otimes_K K_{\infty}(\!(t)\!).
\end{equation}
By construction, these isomorphisms fit into the following commutative diagram as $n$ varies
\begin{equation}
\label{eq:commutative_diagram_for_psi_infty}
\begin{tikzcd}
    M^{[s,r]} \otimes_{\mathbf{C}^{[s,r]}_K, \iota_n }K_{\infty}(\!(t)\!)  
   \arrow[r, "\iota_{n}"] \arrow[d, swap, "\psi^{[s, r]} \otimes \tu{id}"] & D_{\mathrm{dR}}(V)\otimes_K K_{\infty}(\!(t)\!)  \arrow[d, equal] \\
  M^{[s/p, \,r/p]} \otimes_{\mathbf{C}^{[s/p, \,r/p]}_K, \iota_{n+1} }K_{\infty}(\!(t)\!)   \arrow[r, "\iota_{n+1}"]
  & 
     D_{\mathrm{dR}}(V)\otimes_K K_{\infty}(\!(t)\!)
\end{tikzcd}
\end{equation}
where the vertical arrow on the left is induced by the Frobenius map $\psi^{[s, r]}$ from (\ref{eq:frob_isom}). 

We will need to prove the following proposition. For the setup, consider the composition (which we continue to denote by $\iota_n$ by a slight abuse of notation)
\[ 
\iota_n\colon M^{[s,r]} \hookrightarrow M^{[s,r]} \otimes_{\mathbf{C}^{[s,r]}_K, \iota_n }K_{\infty}(\!(t)\!) \isom  D_{\mathrm{dR}}(V)\otimes_K K_{\infty}(\!(t)\!)
\]
where the second map is the isomorphism (\ref{eq:psi_infty_iso}).

\begin{proposition}\label{lemma:image_of_iota_n_contained_in_level_n}
There exists a sufficiently large integer $n_0$ such that, whenever $[s,r]\cap p^{\mathbb{Z}}= \{ p^{-n} \}$ and $n \ge n_0$, we have 
\[\iota_n (M^{[s,r]}) \subset D_{\mathrm{dR}}(V)\otimes_K K_n(\!(t)\!).\]
\end{proposition} 

This assertion seems to be known to experts (for example, it is implicitly claimed before \cite[Proposition~5.15]{Berger-differential}). Let us provide a proof for completeness. As preparation, we provide an estimation of the $\Gamma$-actions on the period rings. Let $r_n=p^{-n}$.

\begin{lemma}\label{lemma: estimation Gamma on C}
There exists an integer $n_1\ge 0$ such that for any $n\ge n_1$, we have
\[\lambda_{r_n}\left((\gamma-1)(x)\right)\le p^{-\frac{p+1}{p}}\cdot \lambda_{r_n}(x)\]
for all $x\in \mathbf{C}^{[r_n, r_n]}_K$ and $\gamma\in \Gamma_{K_n}$.
\end{lemma}

\begin{proof}
It suffices to make the inequality hold for $x=\pi_K^m$ with $m\in \mathbb{Z}$. When $m>0$, we have (cf. \cite[Lemme~4.1]{Berger-differential})
\[\gamma(\pi_K^m)-\pi_K^m=\pi_K^m\left(\frac{\gamma(\pi_K)}{\pi_K}-1\right)\cdot\left(\frac{\gamma(\pi_K^{m-1})}{\pi_K^{m-1}}+\cdots+1\right)\]
and
\[\gamma(\pi_K^{-m})-\pi_K^{-m}=\pi_K^{-m}\left(\frac{\pi_K}{\gamma(\pi_K)}-1\right)\cdot\left(\frac{\pi_K^{m-1}}{\gamma(\pi_K^{m-1})}+\cdots+1\right).\]
Notice that the $\lambda_{r_n}$-norms of the last term in both equations are less than or equal to 1. Hence, it remains to show that, for $n$ sufficiently large (in particular, $r_n<r_K$), we have 
\[\lambda_{r_n}\left(\frac{\gamma(\pi_K)}{\pi_K}-1\right)\le p^{-\frac{p+1}{p}} \quad \tu{ and } \quad\lambda_{r_n}\left(\frac{\pi_K}{\gamma(\pi_K)}-1\right)\le p^{-\frac{p+1}{p}}\]
for all $\gamma\in \Gamma_{K_n}$. This is equivalent to showing that 
\[\lambda_{r_n}\left(\gamma(\pi_K)-\pi_K\right)\le p^{-\frac{p+1}{p}}\cdot \lambda_{r_n}(\pi_K)=p^{-\frac{p+1}{p}-\frac{1}{(p-1)p^{n-1}}\cdot \frac{1}{e_K}}\]
for all $\gamma\in \Gamma_{K_n}$.

We first compute $\lambda_{r_n}(\gamma(\pi)-\pi)$. Since $\gamma\in \Gamma_{K_n}$, we have $\chi(\gamma)\in 1+p^n\Z_p$ where $\chi$ is the cyclotomic character. Write $\chi(\gamma)=1+p^n\cdot u$ for some $u\in \Z_p$. Then we have
\[
\gamma(\pi)-\pi= (1+\pi)^{\chi(\gamma)}-1-\pi=(1+\pi)\cdot\left((1+\pi)^{p^nu}-1\right).
\]
Hence,
\[
\lambda_{r_n}(\gamma(\pi)-\pi)\le \lambda_{r_n}\left((1+\pi)^{p^nu}-1\right)=\lambda_{r_n}\left(\varphi^n((1+\pi)^u-1)\right)=\lambda_1((1+\pi)^u-1).
\]
Since $\frac{(1+\pi)^u-1}{\pi}\in \Z_p[\![\pi]\!]$, we have
$
\lambda_1((1+\pi)^u-1)\le \lambda_1(\pi)=p^{-\frac{p}{p-1}}$ 
where the last equality follows from the computation in \cite[Lemme~6.2]{Colmez08}.

Now we compute $\lambda_{r_n}(\gamma(\pi_K)-\pi_K)$ for $\gamma\in \Gamma_{K_n}$. If $e_K=1$, then $\pi_K=\pi$ and the above computation (for $n \ge 2$) already proves the lemma. For the rest of the proof, we assume $e_K>1$. Recall from Construction~\ref{construction:pi_K} the polynomial $P_K\in \mathbf{A}^+_{K'_0}[X]=W(\kappa)[\![\pi]\!][X]$ such that $P_K(\pi_K)=0$. Let $P_K^{\gamma}$ be the polynomial obtained from applying $\gamma$ to the coefficients of $P_K$. Then $P_K^{\gamma}-P_K=(\gamma(\pi)-\pi)\cdot R$ for some $R\in \mathbf{A}^+_{K'_0}[X]$. In particular, 
\[\lambda_{r_n}\left((P_K^{\gamma}-P_K)(\gamma(\pi_K))\right)\le \lambda_{r_n}(\gamma(\pi)-\pi).\]
On the other hand, using $P_K^{\gamma}(\gamma(\pi_K))=0$, we have
\[\lambda_{r_n}\left((P_K^{\gamma}-P_K)(\gamma(\pi_K))\right)=\lambda_{r_n}(P_K(\gamma(\pi_K)))=\lambda_{r_n}(P_K(\gamma(\pi_K))-P_K(\pi_K)).\]
Let \[\beta\coloneqq \frac{P_K(\gamma(\pi_K))-P_K(\pi_K)}{\gamma(\pi_K)-\pi_K},\]
which is a polynomial in $\pi_K$ and $\gamma(\pi_K)$ with coefficients in $\mathbf{A}^+_{K'_0}$. According to the computation in the proof of \cite[Lemme~9.4]{Colmez08}, $\beta\cdot [P'_K(\overline{\pi}_K)]^{-1}$ is a unit in the ring of integers in $\widetilde{\mathbf{A}}^{r_n}_K$ (with respect to the norm $\lambda_{r_n}$), and we have 
$|P'_K(\overline{\pi}_K)|=p^{-\delta_K}$ 
where $\delta_K=v(\frak{d}_{\mathbf{E}_K/\mathbf{E}_{K_0}})$. Then
\[\lambda_{r_n}(\beta^{-1})\le \lambda_{r_n}\left([P'_K(\overline{\pi}_K)]^{-1}\right)=p^{p^{-n}\cdot \delta_K}.\]
Therefore,
\[
\lambda_{r_n}(\gamma(\pi_K)-\pi_K)\le \lambda_{r_n}\left((\gamma(\pi_K)-\pi_K)\cdot \beta\right)\cdot \lambda_{r_n}(\beta^{-1})\le \lambda_{r_n}(\gamma(\pi)-\pi)\cdot \lambda_{r_n}(\beta^{-1})\le p^{-\frac{p}{p-1}+\frac{\delta_K}{p^n}}.
\]
Hence, it suffices to require $n$ sufficiently large such that $p^{n-1}\ge (p-1)\delta_K+\frac{p}{e_K}$, so that
\[p^{-\frac{p}{p-1}+\frac{\delta_K}{p^n}}\le p^{-\frac{p+1}{p}-\frac{1}{(p-1)p^{n-1}}\cdot\frac{1}{e_K}}.\]
\end{proof}

To prove Proposition~\ref{lemma:image_of_iota_n_contained_in_level_n}, we may assume $[s,r]=[r_n, r_n]$. Recall that each $M^{[r_n, r_n]}$ is equipped with a differential operator 
\[\nabla_V\coloneqq \frac{\log\gamma}{\log \chi(\gamma)}=-\left(\log \chi(\gamma)\right)^{-1}\cdot\sum_{k=1}^{\infty}\frac{(1-\gamma)^k}{k}\]
for any $\gamma$ sufficiently close to 1, making it a differential module over $(\mathbf{C}^{[r_n, r_n]}_K, \nabla)$. Consider the formal sum
\begin{equation}\label{eq: gamma as exp}
\mathrm{exp}\left(\log\chi(\gamma)\cdot \nabla_V\right)\coloneqq \sum_{m=0}^{\infty}\frac{(-1)^m}{m!}\left(\sum_{k=1}^{\infty}\frac{(1-\gamma)^k}{k}\right)^m.
\end{equation}

\begin{lemma}\label{lemma: gamma as exp}
There exists an integer $n_0>0$ such that for any $n\ge n_0$ and any $\gamma\in \Gamma_{K_n}$, the formal sum (\ref{eq: gamma as exp}) converges\footnote{The convergence means both $\sum_{k=1}^{\infty}\frac{(1-\gamma)^k}{k}$ and $\sum_{m=0}^{\infty}\frac{(-1)^m}{m!}\left(\sum_{k=1}^{\infty}\frac{(1-\gamma)^k}{k}\right)^m$ converge.} as an operator on $M^{[r_n, r_n]}$, and we have
$\gamma=\mathrm{exp}\left(\log\chi(\gamma)\cdot \nabla_V\right).$
\end{lemma}

\begin{proof}
For the convergence condition and the identity in the statement of the lemma, it suffices to construct a Banach norm $\mu$ on $M^{[r_n, r_n]}$ that is equivalent to the finite-module norm over $\mathbf{C}^{[r_n, r_n]}_K$, such that \footnote{In fact, we only need $\mu((\gamma-1)(x))\le c\cdot \mu(x)$ for some $c<p^{-\frac{1}{p-1}}$.}
\begin{equation}\label{eq: estimation Gamma action M}
\mu((\gamma-1)(x))\le p^{-\frac{p+1}{p}}\cdot \mu(x)
\end{equation}
for all $x\in M^{[r_n, r_n]}$ and $\gamma\in \Gamma_{K_n}$ (for all $n\ge n_0$).

Let $n_1$ be as in Lemma~\ref{lemma: estimation Gamma on C} and consider $M^{[r_{n_1}, r_{n_1}]}$. Pick a basis $e_1, \ldots, e_d$ of the finite free $\mathbf{C}^{[r_{n_1}, r_{n_1}]}_K$-module $M^{[r_{n_1}, r_{n_1}]}$ and equip $M^{[r_{n_1}, r_{n_1}]}$ with the corresponding finite-module norm $\mu$. By continuity, there exists $n_0\ge n_1$ such that
$
\mu((\gamma-1)(x))\le p^{-\frac{p+1}{p}}\cdot \mu(x)
$
for all $x\in M^{[r_{n_1}, r_{n_1}]}$ and $\gamma\in \Gamma_{K_{n_0}}$. We claim that this $n_0$ works. Indeed, let $n\ge n_0$, and consider the $\Gamma$-equivariant isomorphism
\[M^{[r_{n_1}, r_{n_1}]}\otimes_{\mathbf{C}^{[r_{n_1}, r_{n_1}]}_K, \varphi^{n-n_1}} \mathbf{C}^{[r_n, r_n]}_K\cong M^{[r_n, r_n]}\]
given by the Frobenius structure. Notice that the map 
$
\varphi^{n-n_1}\colon \mathbf{C}^{[r_{n_1}, r_{n_1}]}_K\rightarrow \mathbf{C}^{[r_n, r_n]}_K
$ 
is isometric where the source and target are equipped with the $\lambda_{r_{n_1}}$-norm and the $\lambda_{r_n}$-norm, respectively. We equip the $\mathbf{C}^{[r_n, r_n]}_K$-module $M^{[r_n, r_n]}$ with the finite module norm $\mu'$ using the basis $e_1\otimes 1,\ldots, e_d\otimes 1$. By construction, if we identify $M^{[r_{n_1}, r_{n_1}]}=M^{[r_{n_1}, r_{n_1}]}\otimes 1\subset M^{[r_n, r_n]}$, then $\mu'$ extends the norm $\mu$ on $M^{[r_{n_1}, r_{n_1}]}$. It remains to show that 
\[\mu'((\gamma-1)(e_i\otimes f))\le p^{-\frac{p+1}{p}}\cdot \mu'(e_i\otimes f)\]
for all $i=1, \ldots, d$ and $f\in \mathbf{C}^{[r_n, r_n]}_K$. 
Indeed, since
\[(\gamma-1)(e_i\otimes f)=(\gamma(e_i)-e_i)\otimes \gamma(f)+e_i\otimes (\gamma(f)-f),\]
we have
\[
\mu'((\gamma-1)(e_i\otimes f))\le \max\left\{\mu'((\gamma(e_i)-e_i)\otimes \gamma(f)), \,\,\mu'(e_i\otimes (\gamma(f)-f))\right\}\le p^{-\frac{p+1}{p}}\cdot \mu'(e_i\otimes f)
\]
where the last inequality follows from \eqref{eq: estimation Gamma action M} and Lemma~\ref{lemma: estimation Gamma on C}.
\end{proof}

We now return to the proof of Proposition~\ref{lemma:image_of_iota_n_contained_in_level_n}.

\begin{proof}[Proof of Proposition \ref{lemma:image_of_iota_n_contained_in_level_n}]
Let $n_0$ be as in Lemma~\ref{lemma: gamma as exp}. As mentioned earlier, we may assume $[s,r]=[r_n, r_n]$. Since $\iota_n$ is $\Gamma_K$-equivariant, the differential $\nabla_V$ on $M^{[r_n, r_n]}$ is compatible with the differential operator $\nabla_V$ on $D_{\mathrm{dR}}(V)\otimes_K K_{\infty}(\!(t)\!)$ (where the latter is a differential module over $(K_\infty (\!(t)\!), t \frac{d}{dt})$ by \eqref{eq:nabla_vs_partial}). In other words, we have $\iota_n \circ \nabla_V = \nabla_V \circ \iota_n$.  
Let $x_1, \ldots, x_d$ be a $K$-basis of $D_{\mathrm{dR}}(V)$. Let $x \in M^{[r_n, r_n]}$ be an arbitrary element and write 
\[ \iota_n (x) = \sum_{j=1}^d x_j \otimes c_j =  \sum_{j=1}^d\sum_{i \in \Z}  x_j \otimes a_{i, j} t^i\]
with $c_j\in K_{\infty}(\!(t)\!)$ and $a_{i, j} \in K_\infty$. We compute $\iota_n (\gamma \cdot x)$ for an arbitrary $\gamma \in \Gamma_{K_{n}}$ (with $n \ge n_0$). By Lemma~\ref{lemma: gamma as exp}, this is 
\begin{align*}
    \iota_n (\gamma \cdot x) \: & = 
\: \iota_n \left(\exp (\log \chi(\gamma) \cdot \nabla_V) (x)\right) \\ 
& = \:  
  \iota_n \left(\sum_{m=0}^{\infty} \frac{1}{m!}(\log \chi(\gamma))^m  ( \nabla_V)^m (x)\right) 
  \\ 
& = \: \sum_{m=0}^{\infty} \frac{1}{m!} (\log \chi (\gamma))^m (t \frac{d}{dt})^m  \Big( \sum_{j=1}^d\sum_{i \in \Z} x_j\otimes a_{i, j} t^i  \Big) \\ & = \: 
\sum_{j=1}^d  \sum_{i \in \Z } a_{i,j} \sum_{m=0}^{\infty} x_j\otimes \Big( \frac{i^m}{m!} (\log \chi (\gamma))^m    \Big) t^i    
\\ & = \: 
\sum_{j=1}^d  \sum_{i \in \Z } a_{i,j}  \cdot x_j\otimes\chi(\gamma)^i t^i.    
\end{align*}
On the other hand, we have 
\[
\gamma \cdot \iota_n (x)  \:  = \: \sum_{j=1}^d\sum_{i \in \Z} x_j\otimes\gamma(a_{i, j} t^i) = \: \sum_{j=1}^d\sum_{i \in \Z} \gamma(a_{i, j}) \cdot  x_j\otimes\chi(\gamma)^it^i. 
\]
This implies that $\gamma (a_{i, j}) = a_{i, j}$ for all $\gamma \in \Gamma_{K_n}$, and thus $a_{i, j} \in K_n$. 
\end{proof}

\begin{corollary} \label{cor:the_isom_iota_n_descends_to_finite_level_classical} 
Let $n\ge n_0$ as in Proposition \ref{lemma:image_of_iota_n_contained_in_level_n}. Then there is a canonical $\Gamma_K$-equivariant isomorphism
\[\iota_n: M^{[s,r]}[\frac{1}{t}]\otimes_{\mathbf{C}^{[s,r]}_K[\frac{1}{t}]}K_n(\!(t)\!)\cong D_{\mathrm{dR}}(V)\otimes_K K_n(\!(t)\!).\]
\end{corollary}

\begin{proof}
Indeed, $\iota_n$ from the lemma above supplies a $\Gamma_K$-equivariant map 
\[\iota_n: M^{[s,r]} \otimes_{\mathbf{C}^{[s,r]}_K}K_n(\!(t)\!) \ra D_{\mathrm{dR}}(V)\otimes_K K_n(\!(t)\!).\] 
This map must be an isomorphism by (\ref{eq:psi_infty_iso}) and faithfully flat descent along $K_n(\!(t)\!)\rightarrow K_{\infty}(\!(t)\!)$.
\end{proof}

\noindent \textbf{Step 3.} Applying Beauville--Laszlo gluing (cf. Proposition \ref{prop:BL gluing}) to the triple 
\[ \left(M^{[s,r]}[1/t], \: D_{\mathrm{dR}}(V)\otimes_K K_n[\![t]\!], \:  \iota_n\right),
\] 
we obtain a finite projective $\mathbf{C}^{[s,r]}_K$-module $N^{[s,r]}$ such that $N^{[s,r]}[1/t]=M^{[s,r]}[1/t]$.  Moreover, the $\Gamma_K$-action on $M^{[s,r]}[1/t]$ is compatible with the $\Gamma_K$-action on $ D_{\mathrm{dR}}(V)\otimes_K K_n[\![t]\!]$, since the isomorphism $\iota_{n}$ is $\Gamma_K$-equivariant. This in turn implies that $N^{[s,r]}$ is equipped with a natural $\Gamma_K$-action.\\

\noindent \textbf{Step 4.} Next, by gluing these finite projective modules $N^{[s, r]}$ for varying $0< s \le r$ such that $\# ([s, r] \cap p^\Z) = 1$, we obtain finite projective $\Gamma_K$-modules $N^{[s,r]}$ over $\mathbf{C}^{[s,r]}_K$ for arbitrary $[s,r]$ with $r$ sufficiently small. These $N^{[s,r]}$'s form a compatible family for varying $s, r$ and thus produce a finite projective $\Gamma_K$-module $N^r$ over $\mathbf{C}^r_K$ (resp. $N^{\dagger}$ over $\mathbf{C}^{\dagger}_K$). \\

\noindent \textbf{Step 5.} (Frobenius) We claim that $N^r$ (for $r$ sufficiently small) and $N^{\dagger}$ also possess natural Frobenius structures. It suffices to construct a Frobenius structure on $N^r$ for $r$ sufficiently small, which is equivalent to giving a compatible family of isomorphisms 
\[
\phi^{[s, r]}:  N^{[s, r]} \otimes_{\mathbf{C}^{[s,r]}_K, \varphi} \mathbf{C}^{[s/p, \,r/p]}_K \isom N^{[s/p, \,r/p]} 
\]
for all $r$ sufficiently small and $0 < s \le r$. For this, it again suffices to produce a compatible family of such isomorphisms in the case where $\left([s,r]\cap p^{\mathbb{Z}}\right)= \{ p^{-n} \}$. Recall from Step 3 that $N^{[s, r]}$ is glued from the Beauville--Laszlo triple $ \left(M^{[s,r]}[1/t], \: D_{\mathrm{dR}}(V)\otimes_K K_n[\![t]\!], \:  \iota_n\right),$ while $N^{[s/p, \,r/p]}$ is glued from the Beauville--Laszlo triple $ \left(M^{[s/p, \,r/p]}[1/t], \: D_{\mathrm{dR}}(V)\otimes_K K_{n+1}[\![t]\!], \:  \iota_{n+1}\right)$ 
(note that $[s/p, \,r/p] \cap p^\Z = \{p^{-(n+1)}\}$). By Proposition \ref{prop:BL gluing}, to construct $\phi^{[s, r]}$, it suffices to construct 
\begin{itemize}
    \item 
an isomorphism \[ \phi^{[s, r]} [\frac{1}{t}]:  M^{[s, r]} \otimes_{\mathbf{C}^{[s,r]}_K, \varphi} \mathbf{C}^{[s/p, \,r/p]}_K [\frac{1}{t}] \isom M^{[s/p, \,r/p]} [\frac{1}{t}]\] over $\mathbf{C}^{[s/p, \,r/p]}_K[\frac{1}{t}]$, 
\item  an isomorphism 
\[ (\phi^{[s, r]})^{\wedge}_{t}: \left( D_{\mathrm{dR}}(V)\otimes_K K_n[\![t]\!] \right) \otimes_{ K_n[\![t]\!]}  K_{n+1}[\![t]\!] \isom D_{\mathrm{dR}}(V)\otimes_K K_{n+1}[\![t]\!] 
\]
over $K_{n+1} [\![t]\!]$, 
\end{itemize}  
such that they agree over $K_{n+1} (\!(t)\!)$ after taking base changes along $\iota_{n+1}: \mathbf{C}^{[s/p, \,r/p]}_K[1/t]\hookrightarrow K_{n+1}(\!(t)\!)$ and the natural inclusion $K_{n+1} [\![t]\!] \hookrightarrow K_{n+1} (\!(t)\!)$, respectively (via the isomorphisms $\iota_n$ and $\iota_{n+1}$ in (\ref{eq:psi_infty_iso})). Indeed, we take the map $\phi^{[s, r]} [\frac{1}{t}]$ to be 
\[ 
\psi^{[s, r]} [\frac{1}{t}]: M^{[s, r]} \otimes_{\mathbf{C}^{[s,r]}_K, \varphi} \mathbf{C}^{[s/p,\,r/p]}_K [\frac{1}{t}] \isom M^{[s/p,\,r/p ]} [\frac{1}{t}] 
\] 
induced from (\ref{eq:frob_isom}) and take $(\phi^{[s, r]})^{\wedge}_{t}$ to be the identity map. It remains to show that these two maps agree after taking base changes to $K_{n+1} (\!(t)\!) $. Now, via the commutative diagram 
\[
\begin{tikzcd}
    \mathbf{C}^{[s, r]}_K \arrow[r, "\iota_n"] \arrow[d, "\varphi"] & K_n [\![t]\!] \arrow[d, hook] \\ 
    \mathbf{C}^{[s/p, \,r/p]}_K \arrow[r, "\iota_{n+1}"]   & K_{n+1} [\![t]\!],
\end{tikzcd}
\] the isomorphism $\psi^{[s, r]}$ from (\ref{eq:frob_isom}) induces an isomorphism 
\[ 
 \left( M^{[s, r]}[\frac{1}{t}] \otimes_{ \mathbf{C}^{[s, r]}_K, \iota_n} K_{n} (\!(t)\!) \right) 
\otimes_{K_{n} (\!(t)\!)} K_{n+1} (\!(t)\!) \isom M^{[s/p, \,r/p]}[\frac{1}{t}] \otimes_{ \mathbf{C}^{[s/p, \,r/p]}_K, \iota_{n+1}} K_{n+1} (\!(t)\!).
\] Unwinding definitions, this is precisely the base change of $\phi^{[s, r]} [\frac{1}{t}]$ along $\iota_{n+1}$. The fact that this agrees with $(\phi^{[s, r]})^{\wedge}_{t} = \tu{id}$ on $D_{\mathrm{dR}}(V)\otimes_K K_{n+1} (\!(t)\!)$ follows from the commutativity of 
the diagram (\ref{eq:commutative_diagram_for_psi_infty}). This finishes the construction of the Frobenius isomorphism $\phi^{[s, r]}$, and thus $\phi^r:  N^{r} \otimes_{\mathbf{C}^{r}_K, \varphi} \mathbf{C}^{ r/p}_K \isom N^{r/p}$, as desired. This in turn completes our construction of $N^r$ (resp. $N^\dagger$) as a $(\varphi, \Gamma_K)$-module over $\mathbf{C}^{r}_K$ (resp. over $\mathbf{C}^{\dagger}_K$). \\

\noindent \textbf{Step 6.} Finally, we verify that our construction agrees with the one of Berger.

\begin{proposition}
Suppose that $V$ has non-negative Hodge--Tate weights. Then the $(\varphi, \Gamma_K)$-module $N^\dagger$ over $\mathbf{C}^{\dagger}_K$ constructed above agrees with Berger's $\mathbf{N}_{\mathrm{dR}}(V)$ from \cite{Berger-differential} (cf. Theorem \ref{thm:Berger_N_dR}). 
\end{proposition} 

\begin{proof}
We first claim that, when $V$ has non-negative Hodge--Tate weights, the $(\varphi, \Gamma_K)$-module $N^{\dagger}$ is a submodule of $M^{\dagger}$. Indeed, when $r$ is sufficiently small and $[s,r]\cap p^{\Z}=\{p^{-n}\}$, by the proof of \cite[Corollaire 5.7]{Berger-differential}, we have an isomorphism 
$ 
M^{[s, r]} \otimes_{\mathbf{C}_K^{[s, r]}, \iota_n} K_\infty [\![t]\!] \cong D^{+}_{\mathrm{dif}} (V)
$ 
where the latter contains $D_\mathrm{dR}(V) \otimes_{K} K_\infty [\![t]\!]$ by the assumption on Hodge--Tate weights; in particular, we have
\[ D_{\mathrm{dR}} \otimes_{K} K_{\infty} [\![t]\!] \subset M^{[s, r]} \otimes_{\mathbf{C}_K^{[s, r]}, \iota_n} K_{\infty} [\![t]\!].\]
Taking intersection with 
$
D_{\mathrm{dR}} \otimes_{K} K_n (\!(t)\!)=M^{[s, r]} \otimes_{\mathbf{C}_K^{[s, r]}, \iota_n} K_n (\!(t)\!)
$ 
(cf. Corollary \ref{cor:the_isom_iota_n_descends_to_finite_level_classical}), we obtain
\[ D_{\mathrm{dR}} \otimes_{K} K_n [\![t]\!] \subset M^{[s, r]} \otimes_{\mathbf{C}_K^{[s, r]}, \iota_n} K_n [\![t]\!].\] 
It then follows from the construction (i.e., Beauville--Laszlo gluing) that $N^{[s, r]}\subset M^{[s, r]}$. Then the same statement holds for general $[s, r]$ (with $r$ sufficiently small) by construction. This in turn implies that $N^r \subset M^r$ for each  sufficiently small $r$ and thus $N^\dagger \subset M^\dagger$ as desired.  
Note that $N^\dagger$ is a finite free module over $\mathbf{C}_K^\dagger$ of full rank (cf. \cite[Proposition 4.12 \& Lemme 4.13]{Berger-differential}) and that $N^\dagger$ is stable under the differential operator $\partial=t^{-1}\cdot \nabla$. Therefore, $N^{\dagger}$ must coincide with $\mathbf{N}_{\mathrm{dR}}(V)$ in \cite{Berger-differential} by the uniqueness theorem in {\cite[Th\'eor\`eme 5.10]{Berger-differential}}. 
\end{proof}

To sum up, we obtain a functor
\begin{equation}
    \mathbf{N}_{\mathrm{dR}}:  
\mathrm{Rep}^{\mathrm{dR}}_{\mathbb{Q}_p}(\Gal_{K})  \lra \big\{(\varphi, \Gamma_{K})\textrm{-modules over }\mathbf{C}^{\dagger}_{K}\big\},
\end{equation}
sending $V\mapsto N^{\dagger}$, which recovers Berger's construction when $V$ has non-negative Hodge--Tate weights (justifying our slight abuse of notation here). 
\end{construction}

\subsubsection{Compatibility under extension of base field}

Let $L/K$ be a finite extension. Then there are natural $(\varphi, \Gamma_L)$-equivariant injections $\mathbf{C}^{[s, r]}_K\rightarrow \mathbf{C}^{[s, r]}_L$ for $r$ sufficiently small and $0 < s \le r$. This gives rise to a natural $(\varphi, \Gamma_L)$-equivariant injection $\mathbf{C}^{\dagger}_K\rightarrow \mathbf{C}^{\dagger}_L$.

\begin{proposition}  \label{prop:compatibility_with_L/K}
Let $L/K$ be a finite extension, and let $V\in \mathrm{Rep}^{\mathrm{dR}}_{\Q_p}(\Gal_K)$. Then there is a natural isomorphism
\[\mathbf{N}_{\mathrm{dR}}(V)\otimes_{\mathbf{C}^{\dagger}_K}\mathbf{C}^{\dagger}_L\xrightarrow[]{\sim} \mathbf{N}_{\mathrm{dR}}(V|_{\Gal_L})\]
of $(\varphi, \Gamma_L)$-modules.
\end{proposition}

\begin{proof} 
Let us denote by $N^{[s,r]}_K$ (resp. $N^{[s, r]}_L$) for the finite projective $\mathbf{C}^{[s,r]}_K$-module (resp. $\mathbf{C}^{[s,r]}_L$-module) constructed from $V$ (resp. from $V|_{\Gal_L}$) in Step 3 of Construction \ref{construction:N_dR_for_all_HT_wts}. Similarly, we have $M^{[s, r]}_K$ and $M^{[s, r]}_L$. It suffices to produce compatible isomorphisms 
\[
N^{[s, r]}_K \otimes_{\mathbf{C}^{[s,r]}_K} \mathbf{C}^{[s,r]}_L \isom N^{[s, r]}_L
\]
for $0< s \le r$ with $r$ being sufficiently small, that are equivariant with respect to the Frobenius and $\Gamma_L$-actions. By construction (i.e., Beauville--Laszlo gluing), it suffices to produce two isomorphisms after taking base changes along $\mathbf{C}^{[s,r]}_L\rightarrow \mathbf{C}^{[s,r]}_L[\frac{1}{t}]$ and $\iota_n:\mathbf{C}^{[s,r]}_L\rightarrow L_n [\![t]\!]$, respectively, such that the two isomorphisms agree over $L_n (\!(t)\!)$. Indeed, over $\mathbf{C}^{[s,r]}_L[\frac{1}{t}]$, we have 
\[
N^{[s, r]}_K \otimes_{\mathbf{C}^{[s,r]}_K}  \mathbf{C}^{[s,r]}_L [\frac{1}{t}]  \cong  M^{[s, r]}_K \otimes_{\mathbf{C}^{[s,r]}_K}  \mathbf{C}^{[s,r]}_L [\frac{1}{t}]  \cong M^{[s, r]}_L [\frac{1}{t}] \cong  N^{[s, r]}_L [\frac{1}{t}].
\]
for $r$ sufficiently small, where the second isomorphism follows from the fact that $D^{\dagger}_{\mathrm{rig}, K}(V)$ is compatible with extension of base field. Over $L_n [\![t]\!]$, we have 
\[ 
 \left( D_{\mathrm{dR}, K}(V)\otimes_K K_n[\![t]\!] \right) \otimes_{K_n[\![t]\!]}  L_n [\![t]\!] \isom  D_{\mathrm{dR}, L}(V|_{\Gal_L}) \otimes_L L_n [\![t]\!]
\]
by the base change property of $D_{\mathrm{dR}}$ functor. The claim that these two isomorphisms agree over $L_n (\!(t)\!)$ reduces to the claim that the isomorphism $\iota_n$ in Corollary~\ref{cor:the_isom_iota_n_descends_to_finite_level_classical} is compatible with base field extension $L/K$, which follows from the construction of the map $\iota_n$. For completeness, let us supply the argument for this compatibility. For this, 
we need to show that the following diagram commutes 
\[
\begin{tikzcd}
    M_K^{[s,r]} \otimes_{\mathbf{C}^{[s,r]}_K, \iota_n }K_{n}(\!(t)\!) \arrow[r, "\iota_{n, K}"] \arrow[d]
    & D_{\mathrm{dR}, K}(V)\otimes_K K_{n}(\!(t)\!) \arrow[d] \\
    M_L^{[s,r]} \otimes_{\mathbf{C}^{[s,r]}_L, \iota_n }L_{n}(\!(t)\!)   \arrow[r, "\iota_{n, L}"]
    & D_{\mathrm{dR}, L}(V)\otimes_L L_{n}(\!(t)\!),
\end{tikzcd}
\]
where we have used $\iota_{n, K}$ and $\iota_{n, L}$ to denote the isomorphism in Corollary~\ref{cor:the_isom_iota_n_descends_to_finite_level_classical} respectively over $K$ and $L$. Pick a $\Gal_K$-stable $\Z_p$-lattice $V^{\circ} \subset V$ and let $D_K^{\mathrm{int}, \dagger} (V^{\circ})$ be the associated $(\varphi, \Gamma_K)$-module over $\mathbf{A}_{K}^{\dagger}$ as in Theorem \ref{thm: Galois representations and phi Gamma modules}. We may assume that $r$ is sufficiently small so that $D_K^{\mathrm{int}, \dagger} (V^{\circ})$ is the base change of $D_K^{\mathrm{int}, r} (V^{\circ})=\left(V^{\circ}\otimes_{\Z_p}\mathbf{A}^r \right)^{H_K}$ along $\mathbf{A}_K^r\rightarrow \mathbf{A}^{\dagger}_K$. In particular, we have 
\[ M_K^{[s,r]} = D^{\mathrm{int}, r}_K (V^{\circ}) \otimes_{\mathbf{A}_K^r} \mathbf{C}_K^{[s, r]},\]
compatible with the Frobenius and $\Gamma_K$-structures. Now, the map $\iota_n$ is induced from the composition $\iota_0 \circ \varphi^{-n}$, where $\iota_0$ is the map given by \cite[Proposition 2.11]{Berger-differential}  (see \cite[\S 5.3]{Berger-differential}), thus it suffices to observe the commutativity of the diagram
\[
\begin{tikzcd}
   D^{\mathrm{int}, r}_K (V^{\circ}) \arrow[r, "\varphi^{-n}"]  \arrow[d] & 
    (V^{\circ} \otimes_{\Z_p} \widetilde{\mathbf{A}}^{p^nr})^{H_K}
    \arrow[r, "\iota_{0}"]   \arrow[d]
    &  (V \otimes_{\Q_p} \mathbf{B}_{\mathrm{dR}})^{H_K} \arrow[d] \\
   D^{\mathrm{int}, r}_L (V^{\circ}|_{\Gal_L}) \arrow[r, "\varphi^{-n}"]
    &    (V^{\circ} \otimes_{\Z_p} \widetilde{\mathbf{A}}^{p^nr})^{H_L}
    \arrow[r, "\iota_{0}"]   
    & (V \otimes_{\Q_p} \mathbf{B}_{\mathrm{dR}})^{H_L},
\end{tikzcd}
\]
where the maps $\varphi^{-n}$ are induced from the isomorphism $\varphi^{-n}: \widetilde{\mathbf{A}}^{r} \isom  \widetilde{\mathbf{A}}^{p^nr}$. This justifies the claim.

Finally, unwinding the construction of the $(\varphi, \Gamma_L)$-structure on $N^{[s, r]}_K$ from Step 4 \& 5  of Construction \ref{construction:N_dR_for_all_HT_wts}, we see that the base change isomorphism is indeed equivariant with respect to the $\varphi$- and $\Gamma_L$-actions. This completes the proof.
\end{proof}

\subsubsection{Interpretations in terms of $B$-pairs and the Fargues--Fontaine curve} \label{sss:remark:FF_curve_over_point}

For convenience of the reader, we also recall an interpretation of ${\mathbf{N}}_{\textnormal{dR}}(V)$ using the language of $B$-pairs introduced by Berger \cite{Berger_B_pair}. The category of $B$-pairs (over $K$) consists of pairs $(W_e, W_{\mathrm{dR}}^+)$ where $W_e$ is a finite free $\mathbf{B}_e$-module with a continuous $\Gal_K$-action and $W_{\mathrm{dR}}^+$ is a $\mathbf{B}_{\mathrm{dR}}^+$-lattice in $W_e\otimes_{\mathbf{B}_e} \mathbf{B}_{\mathrm{dR}}$ stable under $\Gal_K$ (cf. \cite[\S 2]{Berger_B_pair}). By \cite[Th\'eor\`eme 2.2.7]{Berger_B_pair}, there is an equivalence of categories
\begin{equation}\label{eq: B pairs}
\big\{(\varphi, \Gamma_{K})\textrm{-modules over }\mathbf{C}^{\dagger}_{K}\big\} \isom \big\{B\textrm{-pairs over }K\big\}.
\end{equation}
For $V\in \mathrm{Rep}_{\mathbb{Q}_p}(\Gal_K)$, the associated \'etale $(\varphi, \Gamma_K)$-module $D^{\dagger}_{\mathrm{rig}, K}(V)$ over $\mathbf{C}_{K}^{\dagger}$ corresponds to the $B$-pair $(V\otimes_{\Q_p}\mathbf{B}_e , V\otimes_{\Q_p}\mathbf{B}_{\mathrm{dR}}^+)$ under the equivalence (\ref{eq: B pairs})
(cf. \cite[\S 2.2]{Berger_B_pair}). Let us further assume that $V$ is de Rham and let $(D_{\mathrm{dR}}(V), \mathrm{Fil}^\bullet D_{\mathrm{dR}} (V))$ be the corresponding filtered $K$-vector space. Then there is an isomorphism  
\[(V\otimes_{\Q_p}\mathbf{B}_e , V\otimes_{\Q_p}\mathbf{B}_{\mathrm{dR}}^+) \cong 
( D_{\mathrm{dR}} (V) \otimes_{K} \mathbf{B}_e, \mathrm{Fil}^0 ( D_{\mathrm{dR}} (V)\otimes_{K}\mathbf{B}_{\mathrm{dR}} ) )
\]
of $B$-pairs over $K$. Now, after changing the $\mathbf{B}_{\mathrm{dR}}^+$-lattice, we may consider another $B$-pair
$
(D_{\mathrm{dR}} (V)\otimes_{K}\mathbf{B}_e,  D_{\mathrm{dR}} (V) \otimes_{K} \mathbf{B}_{\mathrm{dR}}^+). 
$ 
Then this new $B$-pair corresponds to the (not necessarily \'etale) $(\varphi, \Gamma_K)$-module $\mathbf{N}_{\mathrm{dR}}(V)$ under the equivalence (\ref{eq: B pairs}).

\begin{remark} \label{remark:modification_remark}
   The category of $B$-pairs over $K$, and thus that of $(\varphi, \Gamma_K)$-modules over $\mathbf{C}^{\dagger}_{K}$, is also equivalent to the category of $\Gal_K$-equivariant vector bundles over the Fargues--Fontaine curve $\mX_{\tu{FF}, \C_K}$ over $\C_K$ (cf. \cite[Theorem 9.5.8]{Kedlaya-LiuI}). From this perspective,  $D^{\dagger}_{\mathrm{rig}, K}(V)$ corresponds to a vector bundle $\mathcal E_{V}$ over $\mX_{\mathrm{FF}, \C_K}$, and $ \mathbf{N}_{\mathrm{dR}} (V)$ corresponds to the modification of $\mathcal E_{V}$ at $\infty \in \mX_{\mathrm{FF}, \C_K}$ using the lattice $D_{\mathrm{dR}} (V) \otimes_{K}\mathbf{B}_{\mathrm{dR}}^+$. 
\end{remark}

\begin{remark} \label{rem: relative-NdR-FF-curve}
In \S \ref{sec:relative_N_dR}, we shall construct a relative version of Berger's $\mathbf{N}_{\mathrm{dR}}$ functor. In the relative setting, however, the methods of $B$-pairs and Fargues--Fontaine curves are not helpful.

In fact, by \cite[Theorem~9.3.12]{Kedlaya-LiuI}, there is an equivalence of categories between the category of relative $B$-pairs (or a suitable category of vector bundles on the relative Fargues--Fontaine curve) and the category of $(\varphi, \Gamma)$‑modules over perfect relative period rings. However, we lack general descent results for $(\varphi,\Gamma)$‑modules from perfect relative period rings to imperfect ones; descent is known only for the \'etale objects (see Corollary~\ref{cor: decompletion-etale-(phi, Gamma)-modules}). For this reason, we cannot construct the relative $\mathbf{N}_{\mathrm{dR}}$ functor from relative $B$-pairs nor from the relative Fargues--Fontaine curve.
\end{remark}

\subsubsection{Construction of $\widetilde{\mathbf{N}}_{\textnormal{dR}}(V)$}\label{sss:tilde NdR classical setting}

For any de Rham representation $V$, we have constructed a $(\varphi, \Gamma_K)$-module $\mathbf{N}_{\mathrm{dR}}(V)$ over the imperfect period ring $\mathbf{C}^{\dagger}_K$ by modifying $D^{\dagger}_{\mathrm{rig}}(V)$ along the ``$t=0$ locus''. We have seen that $\mathbf{N}_{\mathrm{dR}}(V)$ naturally carries the structure of a differential module. Let us point out that there is also a theory of ``modified $(\varphi, \Gamma_K)$-modules'' over the perfect period ring $\widetilde{\mathbf{C}}^{\dagger}_K$.

\begin{construction} \label{construction:tilde NdR classical setting}
This time, we consider the Beauville--Laszlo datum 
\[\left(\widetilde{\mathbf{C}}^{[s,r]}_K,\,\, \widetilde{\mathbf{C}}^{[s,r]}_K[\frac{1}{t}],\,\, \mathbf{B}^+_{\mathrm{dR}, K}\right)\]when $\#\left([s,r]\cap p^{\mathbb{Z}}\right)=1$, where $\mathbf{B}^+_{\mathrm{dR}, K}:=(\mathbf{B}^+_{\mathrm{dR}})^{H_K}$. For $V\in \mathrm{Rep}^{\mathrm{dR}}_{\Q_p}(\Gal_K)$, consider the associated \'etale $(\varphi, \Gamma_K)$-module $\widetilde{M}^{\dagger}:=\widetilde{D}^{\dagger}_{\mathrm{rig},K}(V)$ (cf. Remark \ref{rmk: phi Gamma modules over perfect period rings--classical}), which corresponds to a compatible system of finite projective modules $\widetilde{M}^{[s,r]}$ over $\widetilde{\mathbf{C}}^{[s,r]}_K$. Then, by Beauville--Laszlo gluing, the finite projective $\widetilde{\mathbf{C}}^{[s,r]}_K[\frac{1}{t}]$-module $\widetilde{M}^{[s,r]}[\frac{1}{t}]$ and the finite projective $\mathbf{B}^+_{\mathrm{dR},K}$-module $D_{\mathrm{dR}}(V)\otimes_K \mathbf{B}^+_{\mathrm{dR},K}$ glue to a finite projective $\widetilde{\mathbf{C}}^{[s,r]}_K$-module $\widetilde{N}^{[s,r]}$. By a similar argument as in Construction \ref{construction:N_dR_for_all_HT_wts}, these $\widetilde{N}^{[s,r]}$ glue to a $(\varphi, \Gamma_K)$-module $\widetilde{N}^r$ (resp. $\widetilde{N}^{\dagger}$) over $\widetilde{\mathbf{C}}^r_K$ (resp. $\widetilde{\mathbf{C}}^{\dagger}_K)$. This yields a functor 
\[\widetilde{\mathbf{N}}_{\mathrm{dR}}: \mathrm{Rep}^{\mathrm{dR}}_{\Q_p}(\Gal_K)\rightarrow \big\{(\varphi, \Gamma_{K})\textrm{-modules over }\widetilde{\mathbf{C}}^{\dagger}_{K}\big\}\]
sending $V\mapsto \widetilde{N}^{\dagger}$. By construction, it is straightforward to verify that 
$\widetilde{\mathbf{N}}_{\mathrm{dR}}(V) = \mathbf{N}_{\mathrm{dR}}(V) \otimes_{\mathbf{C}^{\dagger}_K } \widetilde{\mathbf{C}}^{\dagger}_K$ as $(\varphi, \Gamma_K)$-modules. Consequently, $\mathbf{N}_{\mathrm{dR}}(V)$ is the unique descent of $\widetilde{\mathbf{N}}_{\mathrm{dR}}(V)$ along $\mathbf{C}^{\dagger}_K\rightarrow \widetilde{\mathbf{C}}^{\dagger}_K$, via the equivalence of categories between $(\varphi, \Gamma_K)$-modules over $\mathbf{C}^{\dagger}_K$ and those over $\widetilde{\mathbf{C}}^{\dagger}_K$ (cf. \cite[Theorem 9.5.8]{Kedlaya-LiuI}).
\end{construction}

\vspace{0.1in}
\section{Perfect relative period rings} \label{sec:relative_period_rings:perfect}

In this section and the next, we introduce the ``relative period rings'' (perfect and imperfect ones) arising from toric towers. These period rings are the main objects of study in this article, over which we shall systematically develop the theories of relative $(\varphi,\Gamma)$-modules and $p$-adic differential equations. Our treatment is inspired by \cite{Kedlaya-LiuII}, but differs from \textit{loc. cit.} in a fundamental way, especially in how we treat \textit{imperfect period rings}. We shall discuss the main differences in more detail at the beginning of \S \ref{sec:weakly_decompleting_data}.

The goals of this section can be summarized as follows. 
\begin{itemize}
    \item Introduce the notion of toric charts and toric towers.
    \item Introduce perfect relative period rings arising from toric towers. 
    \item Recall the construction of (overconvergent) relative $(\varphi, \Gamma)$-modules over these perfect period rings attached to $p$-adic local systems. 
\end{itemize}
We shall mainly follow \cite{Kedlaya-LiuI} for the notion of perfect relative period rings and $(\varphi, \Gamma)$-modules over them. The treatment here parallels our treatment in the classical case (cf. \S \ref{ss:point perfect period ring}, \S \ref{ss:point imperfect period ring}, and \S \ref{ss:classical_pDE}). We retain the notation from \S \ref{ss: notation and convention}.

\subsection{Toric charts and toric towers}  \label{ss:toric_tower}

For any integer $d\ge 1$, let
$
\D_K^d\coloneqq \spa (K \gr{T_1, \ldots, T_d}, \mO_K \gr{T_1, \ldots, T_d})
$ 
denote the $d$-dimensional closed unit polydisk over $K$, equipped with the log structure associated with the normal crossing divisor $\{T_1\cdots T_d=0\}$. For each $n \ge 0$, let
\[
\D^d_{K, n}\coloneqq \spa \left(K_n \gr{T_1^{\frac{1}{p^n}}, \ldots, T_d^{\frac{1}{p^n}}}, \mO_{K_n} \gr{T_1^{\frac{1}{p^n}}, \ldots, T_d^{\frac{1}{p^n}}}\right). 
\]
This is the $d$-dimensional closed unit polydisk over $K_n$, with coordinates $T_1^{\frac{1}{p^n}}, \ldots, T_d^{\frac{1}{p^n}}$; we equip it with the log structure associated with the normal crossing divisor 
$
\{T_1^{\frac{1}{p^n}}\cdots T_d^{\frac{1}{p^n}}=0\}.
$  
When the base field $K$ and dimension $d$ are understood, we often write $\D_{n} = \D^d_{K, n}$ to simplify the notation. 
The natural maps $\psi_n\colon \D_n \ra \D^d_K$ (induced by the natural inclusion of rings) form a tower of Kummer \'etale maps
\[
\psi\colon \cdots \ra \D_{n+1}\ra \D_n\ra \cdots \ra\D_1\ra \D^d_K.
\]
Let $\D_{\infty}=\varprojlim_{n\ge 0} \D_n$, viewed as an element in the \emph{pro-Kummer \'etale site} $\left(\D^d_K\right)_{\mathrm{prok\acute{e}t}}$ studied in \cite{DLLZ1}. The tower is \textit{log perfectoid} with associated perfectoid space 
\[ 
\widehat{\D}_{\infty} = \spa (\widehat{K}_\infty \gr{T_1^{\frac{1}{p^\infty}}, \ldots, T_d^{\frac{1}{p^\infty}}}, \mO_{\widehat{K}_\infty} \gr{T_1^{\frac{1}{p^\infty}}, \ldots, T_d^{\frac{1}{p^\infty}}}).
\]Moreover, the tower is Galois in the sense of \cite[Definition ~6.1.2]{DLLZ1}, with Galois group $\Gamma \cong \Gamma_{\mathrm{geom}} \rtimes \Gamma_K$ where $\Gamma_{\mathrm{geom}}$ is the Galois group of 
\[\D_{\infty}\rightarrow \varprojlim_{n \ge 0}\D^d_K\times_{\spa(K, \mO_K)} \spa(K_n, \mO_{K_n}).\]
In particular, $\Gamma_{\mathrm{geom}}$ acts on $K_{\infty}$ trivially and permutes the $p$-power roots of $T_i$'s. On the other hand, $\Gamma_K = \Gal(K_\infty/K)$ acts on the $p$-power roots of unity via the usual cyclotomic character $\chi$ and acts trivially on the prescribed $p$-power roots $T_i^{1/p^n}$'s of $T_i$.\footnote{In fact, such a choice of compatible system of $p$-power roots of $T_i$'s determines a splitting $\Gamma\cong \Gamma_{\mathrm{geom}}\rtimes \Gamma_K$.} Note that $\Gamma_{\mathrm{geom}} \cong \Z_p^{\oplus d}$ and $\Gamma_K \subset \Z_p^\times$ is an open subgroup, and 
for any $\sigma \in \Gamma_K, \tau \in \Gamma_{\mathrm{geom}}$, we have  $\sigma \tau \sigma^{-1} = \tau^{\chi(\sigma)}.$ 
Finally, for any integer $n\ge 0$, the Galois group of $\D_{\infty}\rightarrow \D_n$ is an open subgroup $\Gamma_n\subset\Gamma$. Then $\Gamma_n\cong \Gamma_{\mathrm{geom},n}\rtimes \Gamma_{K_n}$ where $\Gamma_{\mathrm{geom},n}=p^n\Gamma_{\mathrm{geom}}\subset \Gamma_{\mathrm{geom}}$.
 
We now specify the setting of interest. Throughout the remainder of \S \ref{sec:relative_period_rings:perfect}, let $A_X$ be a reduced affinoid $K$-algebra (in the sense of Tate's rigid analytic geometry) and write $X=\spa (A_X, A_X^+)$ for the associated adic space over $\spa(K,\mO_K)$. In particular, $A_X^+ = A_X^{\circ}$. 
 
\begin{definition} \label{defn: toric chart}
A \emph{toric chart} of $X$ is an \'etale morphism $h\colon X \ra \D_K^d$ of adic spaces (for some $d\ge 1$) that is a composition of rational localizations and finite \'etale maps such that $h$ is away from the locus $T_1\cdots T_d = 0$.
\end{definition}
 
\begin{definition} \label{def:toric_chart}
Let $h\colon X\ra \D^d_K$ be a toric chart. For each $n\ge 0$, let $X_n\coloneqq X \times_{\D^d_K} \D_n$. The tower $\psi=\{\psi_n\colon \D_n\ra\D^d_K\}$ then pulls back to a tower $\psi = \{\psi_n\colon X_n \ra X\}$ over $X$, which we refer to as the \emph{toric tower} over $X$ relative to the toric chart $h$. 
\end{definition} 

Since the image of $h$ does not intersect with the locus $\{T_1\cdots T_d=0\}$, each $\psi_n$ is finite \'etale, and the tower $\psi\colon X_{\infty}=\varprojlim X_n \ra X$ is an affinoid perfectoid object in the pro-\'etale site $X_{\proet} $ of $X$. Write $\widehat{X}_\infty = X \times_{\D^d_K} \widehat{\D}_\infty$ for the associated perfectoid space.  Note that $\psi$ is Galois with the same Galois group $\Gamma \cong \Gamma_{\mathrm{geom}} \rtimes \Gamma_K$. Furthermore, each connected component of $\psi\colon X_{\infty} \rightarrow X$ is Galois with Galois group isomorphic to a finite index open subgroup of $\Gamma = \Gamma_{\mathrm{geom}} \rtimes \Gamma_K$. 

\begin{definition} \label{construction:A_psi_X_infty} \hfill
\begin{enumerate}
\item For every $n\ge 0$, write $X_n = \spa (A_{X, n}, A_{X, n}^+)$. In particular, we have $ A_{X,0} =A_X$ and $A_{X,0}^+ =A_X^+$. Each $A_{X, n}$ is equipped with the spectral norm $|\cdot|$, which extends to a norm $|\cdot|$ on  
$
A_{X, \infty} \coloneqq \bigcup_{n\ge 0} A_{X, n}.
$  

\item Let $\widehat A_{X, \infty}$ denote the completion of $A_{X, \infty}$ with respect to $|\cdot|$, and define 
\[
\widehat A_{X, \infty}^+\coloneqq \left\{x \in \widehat A_{X, \infty} \: \big| \: |x| \le 1 \right\}.
\]
Then $\widehat A_{X, \infty}$ is a perfectoid affinoid $\Q_p$-Banach algebra, $(\widehat A_{X, \infty}, \widehat A_{X, \infty}^+)$ is a perfectoid Huber pair, and we have 
$\widehat{X}_\infty = \spa (\widehat A_{X, \infty}, \widehat A_{X, \infty}^+)$. 

\item Denote the tilt of $(\widehat A_{X, \infty}, \widehat A_{X, \infty}^+)$ by 
$
(\sq R_X, \sq R_X^+)\coloneqq (\widehat A_{X, \infty}^\flat, \widehat A_{X, \infty}^{\flat, +}).$  
We write $|\cdot|$ also for the norm on $\sq R_X$ induced from the norm $|\cdot|$ on $\widehat A_{X, \infty}$.

\item Define a valuation $v\colon \sq R_X \ra \R \cup \{+\infty\}$ by $v(y)\coloneqq -\log_p |y|$, which extends the valuation $v$ on $\widehat K_\infty^\flat \subset \sq R_X$ defined in \S \ref{ss: notation and convention}. 

\item Note that $\varpi$ is a pseudouniformizer in $\widehat{A}_{X, \infty}^+$. 
We fix a tilt $\varpi^\flat \in  \sq R_X^+$ of $\varpi$. 
\end{enumerate} 
\end{definition}

\subsection{Perfect period rings} \label{ss:perfect relative period rings}   

We recall the constructions of perfect period rings in the relative setting following \cite{Kedlaya-LiuI}. These are natural extensions of the classical perfect period rings reviewed in \S \ref{ss:point perfect period ring}. 

\begin{definition} \label{definition:valuation_v_r_relative_perfect}
Let $h\colon X \ra \D^d_K$ be a toric chart.
\begin{enumerate}
\item  Define $\widetilde{\mathbf{A}}_X\coloneqq W(\sq R_X)$ and $\widetilde{\mathbf{A}}_X^+\coloneqq W(\sq R_X^+)$.

\item  For each $r > 0$, we define a valuation $v_{r}\colon \widetilde{\mathbf{A}}^+_X \rightarrow \mathbb{R}\cup \{+\infty\}$ by 
\[
v_{r}\left(\sum_{n=0}^{\infty} p^n [x_n] \right)\coloneqq \inf_{n\ge 0}\left\{v(x_n)+\frac{n}{r}\right\},
\] 
which extends the valuation $v_r$ on $\widetilde{\mathbf{A}}_{K}^+ \subset \widetilde{\mathbf{A}}_X^+$. This corresponds to the norm $\lambda_r$ given by 
\[
\lambda_r \big(\sum_{n=0}^{\infty} p^n [x_n]  \big) = \sup_n \{p^{-n} |x_n|^r\}.
\] 
Note that $\lambda_r(x) = p^{- r \cdot v_r(x)}$. 

\item For each $r > 0$, write $\widetilde{\mathbf{A}}_X^r$ for the completion of $\widetilde{\mathbf{A}}_X^+[[\varpi^\flat]^{-1}]$ under the norm $\lambda_r$. Equivalently,
\[
\widetilde{\mathbf{A}}_X^r\coloneqq \left \{ \sum_{n \ge 0} p^n [x_n] \in \widetilde{\mathbf{A}}_X \, \Bigg|   \,  p^{-n} |x_n|^r \ra 0 \tu{ as } n \ra \infty \right\}. 
\]

\item Let $\widetilde{\mathbf{A}}_X^{\dagger}$ be the union $\widetilde{\mathbf{A}}_X^{\dagger} \coloneqq \varinjlim_{r>0} \widetilde{\mathbf{A}}_X^r$.

\item The period rings $\widetilde{\mathbf{B}}_X$, $\widetilde{\mathbf{B}}_X^r$, and $\widetilde{\mathbf{B}}_X^\dagger$ are obtained from the corresponding period rings of type $\mathbf{A}$ by inverting $p$. 

\item  For each $0 < s \le r$, define $\widetilde{\mathbf{C}}_X^{[s, r]}$ to be the completion of the period ring $\widetilde{\mathbf{B}}_X^r$ with respect to the norm \(\max\{\lambda_{s}, \lambda_{r}\}\) (equivalently, completion with respect to the family of norms $\{\lambda_t\}_{t\in [s,r]}$).

\item For each $r > 0$, let $\widetilde{\mathbf{C}}_X^{r}$ be the Fr\'echet completion of $\widetilde{\mathbf{B}}_X^r$ with respect to the family of norms $\{\lambda_{s}\}_{s \in (0, r]}$.

\item Finally, define $\widetilde{\mathbf{C}}_X^\dagger$ to be the union  
$\widetilde{\mathbf{C}}_X^{\dagger} \coloneqq \varinjlim_{r>0} \widetilde{\mathbf{C}}_X^r.$
\end{enumerate}
\end{definition}     

\begin{remark}
Strictly speaking, the perfect period rings $\widetilde{\mathbf{A}}_X$, $\widetilde{\mathbf{A}}_X^r$, etc. depend on the choice of the toric chart $h\colon X \ra \D^d_K$. In most cases, the relevant toric chart will be evident from the context. For situations where this choice requires emphasis, we adopt the notations $\widetilde{\mathbf{A}}_{X,h}$, $\widetilde{\mathbf{A}}_{X,h}^r$, etc.
\end{remark}

\begin{remark}[Topology]
We endow $\widetilde{\mathbf{A}}^r_X$ with the topology induced by the $\lambda_r$-norm, and endow $\widetilde{\mathbf{A}}^{\dagger}_X$ with the inductive limit topology (i.e. the finest topology making $\widetilde{\mathbf{A}}^r_X\hookrightarrow \widetilde{\mathbf{A}}^{\dagger}_X$ continuous for all $r>0$). Inverting $p$, we obtain induced topologies on $\widetilde{\mathbf{B}}^r_X$ and $\widetilde{\mathbf{B}}^{\dagger}_X$. The ring $\widetilde{\mathbf{C}}^{[s,r]}_X$ is equipped with the topology given by $\max\{\lambda_r, \lambda_s\}$, which then induces a Fr\'echet topology on $\widetilde{\mathbf{C}}^r_X$. Finally, we equip $\widetilde{\mathbf{C}}^{\dagger}_X$ with the inductive limit topology (i.e. the LF topology).
\end{remark}

\begin{remark} \label{remark:Frobenius_relative_and_action_of_Gamma} \hfill
\begin{enumerate}
\item For more general adic spaces over $K$, the perfect period rings above globalize to period sheaves on $X_{\text{pro\'et}}$, which are denoted by $\widetilde{\mathbb{A}}_X, \widetilde{\mathbb{A}}^r_X, \widetilde{\mathbb{A}}^\dagger_X, \widetilde{\mathbb{B}}_X, \widetilde{\mathbb{B}}^r_X, \widetilde{\mathbb{B}}^\dagger_X,   \widetilde{\mathbb{C}}^r_X, \widetilde{\mathbb{C}}^\dagger_X$, etc. For more details, we refer the reader to \cite[\S 9.3]{Kedlaya-LiuI}. We will not make essential use of the sheaf version of these period rings in this article. 

\item Just as in the classical setting, there are natural surjective maps $\theta\colon \widetilde{\mathbf{A}}^+_X  \to \widehat A_{X, \infty}^+$ and $\theta\colon \widetilde{\mathbf{A}}^1_X  \to \widehat A_{X,\infty}$. 

\item Remark~\ref{remark:Frobenius_classical} also applies to the relative setting. In other words, the Frobenius maps on $\widetilde{\mathbf{A}}_X^+$ and $\widetilde{\mathbf{A}}_X$ extend to Frobenius structures on $\widetilde{\mathbf{A}}_X^r$, $\widetilde{\mathbf{A}}_X^\dagger$, $\widetilde{\mathbf{B}}_X$, $\widetilde{\mathbf{B}}_X^r$, $\widetilde{\mathbf{B}}_X^\dagger$, $\widetilde{\mathbf{C}}_X^{[s,r]}$, $\widetilde{\mathbf{C}}_X^r$, and $\widetilde{\mathbf{C}}_X^\dagger$, which we also denote by $\varphi$. 

\item The natural $\Gamma$-action on $\widetilde{R}_X$ induces an action on $\widetilde{\mathbf{A}}_X$. Note that the $\lambda_r$-norm is preserved under the $\Gamma$-action. So for each $r > 0$, the subspace $\widetilde{\mathbf{A}}_X^r\subset \widetilde{\mathbf{A}}_X$ is stable under the $\Gamma$-action, and the action is isometric with respect to $\lambda_r$. This in turn induces a continuous $\Gamma$-action on $\widetilde{\mathbf{A}}_X^{\dagger}$. Similarly, there are natural isometric (resp. continuous) $\Gamma$-actions on $\widetilde{\mathbf{B}}_X^r$ and $\widetilde{\mathbf{C}}_X^{[s, r]}$ (resp. on $\widetilde{\mathbf{B}}_X^{\dagger}$, $\widetilde{\mathbf{C}}_X^r$, and $\widetilde{\mathbf{C}}_X^{\dagger}$). 
\end{enumerate}
\end{remark}

The following definition will be used later. 

\begin{definition} \label{def:w_k_for_relative_A_tilde} 
For each $k\in \mathbb{Z}_{\ge 0}$, we define a map $
w_k\colon \widetilde{\mathbf{A}}_{X} \rightarrow \mathbb{R}\cup \{+\infty\}$ 
as follows. For $x=\sum_{n=0}^{\infty}p^n[x_n]\in \widetilde{\mathbf{A}}_{X}$ with $x_n\in \widetilde{R}_{X}$, set 
$w_k(x)\coloneqq \inf_{n\le k} \{v(x_n)\}.$  
In particular, if $x\in \widetilde{\mathbf{A}}^r_{X}$, we have $v_r(x)=\inf_{k\ge 0}\{w_k(x)+\frac{k}{r}\}$.
\end{definition}

\subsection{Perfect period rings under finite \'etale extensions} \label{ss:fet_over_tilde_A}

We now recall how the perfect period rings in \S \ref{ss:perfect relative period rings} behave under finite \'etale extensions of $X$. The following results are either explicitly stated in \cite[Propositions~5.5.3 and 5.5.4]{Kedlaya-LiuI} or follow implicitly from the proofs therein.

\begin{proposition}\label{prop: perfect theory} 

Let $(*)_{\mathrm{f\acute{e}t}}$ denote the category of finite \'etale $*$-algebras. \hfill

\begin{enumerate}
\item There are equivalences of categories
\[
\left(\widetilde{\mathbf{A}}^{\dagger}_X\right)_{\mathrm{f\acute{e}t}} \cong \left(\widetilde{\mathbf{A}}_X\right)_{\mathrm{f\acute{e}t}} \cong \left(\widetilde{R}_X\right)_{\mathrm{f\acute{e}t}} \cong \left(\widehat{A}_{X, \infty}\right)_{\mathrm{f\acute{e}t}}
\]
where the first two equivalences are induced from the natural inclusion $\widetilde{\mathbf{A}}^{\dagger}_X\hookrightarrow \widetilde{\mathbf{A}}_X$ and the modulo $p$ map $\widetilde{\mathbf{A}}_X\rightarrow \widetilde{R}_X$, while the last equivalence is the tilting equivalence.

\item Let $f\colon X'\rightarrow X$ be a finite \'etale morphism, equipped with the toric chart $h'=h\circ f\colon X'\rightarrow \mathbb{D}_K.$  Then $f$ induces finite \'etale morphisms 
\[
\widetilde{\mathbf{A}}^{\dagger}_X\rightarrow \widetilde{\mathbf{A}}^{\dagger}_{X'}, \quad \widetilde{\mathbf{A}}_X\rightarrow \widetilde{\mathbf{A}}_{X'}, \quad \widetilde{R}_{X}\rightarrow \widetilde{R}_{X'}, \quad \widehat{A}_{X, \infty}\rightarrow \widehat{A}_{X', \infty}
\] 
which correspond to each other through the equivalences in part (1). 

\item The map $f$ in part (2) induces a compatible system of finite \'etale morphism $\widetilde{\mathbf{A}}^r_X\rightarrow \widetilde{\mathbf{A}}^r_{X'}$ for all $r > 0$ sufficiently small; namely,
\begin{enumerate}
\item[(i)] The map $\widetilde{\mathbf{A}}^{\dagger}_X \rightarrow \widetilde{\mathbf{A}}^{\dagger}_{X'}$ is the base change of $\widetilde{\mathbf{A}}^r_X\rightarrow \widetilde{\mathbf{A}}^r_{X'}$ along the inclusion $\widetilde{\mathbf{A}}^r_X\hookrightarrow \widetilde{\mathbf{A}}^{\dagger}_X$;

\item[(ii)] for any $0<s\le r$, $\widetilde{\mathbf{A}}^s_{X}\rightarrow \widetilde{\mathbf{A}}^s_{X'}$ is the base change of $\widetilde{\mathbf{A}}^r_X\rightarrow \widetilde{\mathbf{A}}^r_{X'}$ along $\widetilde{\mathbf{A}}^r_X\hookrightarrow \widetilde{\mathbf{A}}^s_{X}$.
\end{enumerate}
In addition, for $r$ sufficiently small and for any $0<s\le r$, we have
\[
\widetilde{\mathbf{C}}^{[s,r]}_{X'}\cong \widetilde{\mathbf{C}}^{[s,r]}_X\otimes_{\widetilde{\mathbf{A}}^r_X}\widetilde{\mathbf{A}}^r_{X'},\qquad \widetilde{\mathbf{C}}^r_{X'}\cong \widetilde{\mathbf{C}}^r_X\otimes_{\widetilde{\mathbf{A}}^r_X}\widetilde{\mathbf{A}}^r_{X'}.
\]
In particular, $\widetilde{\mathbf{C}}^{[s,r]}_X\rightarrow \widetilde{\mathbf{C}}^{[s,r]}_{X'}$ and $\widetilde{\mathbf{C}}^r_X\rightarrow \widetilde{\mathbf{C}}^r_{X'}$ are finite \'etale.

\end{enumerate}
\end{proposition}

\begin{proof}
For the reader’s convenience, we either recall some of the arguments from \cite[\S 5.5]{Kedlaya-LiuI} or indicate more precisely where the corresponding statements are proved in \textit{loc. cit.} We remark that Theorem~\ref{thm: open mapping thm} is used in \textit{loc. cit.}; note that $\pi_K$ is a topologically nilpotent unit in $\widetilde{R}_X$.

\begin{enumerate}
\item Since $\widetilde{\mathbf{A}}_X$ is $p$-adically complete, $(\widetilde{\mathbf{A}}_X, (p))$ is henselian. So the mod $p$ map induces an equivalence
\[ 
\left(\widetilde{\mathbf{A}}_X\right)_{\mathrm{f\acute{e}t}} \cong \left(\widetilde{R}_X\right)_{\mathrm{f\acute{e}t}}.
\]
Since the $p$-adic completion of $\widetilde{\mathbf{A}}^{\dagger}_X$ is $\widetilde{\mathbf{A}}_X$, it suffices to show that $(\widetilde{\mathbf{A}}^{\dagger}_X, (p))$ is henselian. For each $r > 0$, consider the norm $\mu_r$ on $\widetilde{\mathbf{A}}_X^r$ given by the maximum of $\lambda_r$ and the $p$-adic norm (from that on $\widetilde{\mathbf{A}}_X$). Note that for each $0 < s \le r$, the map
$
(\widetilde{\mathbf{A}}_X^r, \mu_r) \rightarrow (\widetilde{\mathbf{A}}_X^s, \mu_s)
$ 
is submetric. For $x \in \widetilde{\mathbf{A}}^{\dagger}_X$, let $\mu(x) = \inf_{r > 0}\{ \mu_r(x) \}$ where the infimum runs over all sufficiently small $r > 0$ such that $x \in \widetilde{\mathbf{A}}^r_X$. Then $\mu(x) \le 1$ for any $x \in \widetilde{\mathbf{A}}^{\dagger}_X$, and $\mu(x) < 1$ if and only if $x \in p\widetilde{\mathbf{A}}^{\dagger}_X$. Thus, $(\widetilde{\mathbf{A}}^{\dagger}_X, (p))$ is henselian by \cite[Lemma ~2.2.3(b)]{Kedlaya-LiuI}.

\item The finite \'etale morphism $f\colon X' \rightarrow X$ naturally induces a finite \'etale map $\widehat{A}_{X, \infty} \rightarrow \widehat{A}_{X', \infty}$. This further induces the other finite \'etale morphisms in the statement via part (1).

\item Let $x_1, \ldots, x_n$ be any choice of elements in $\widetilde{\mathbf{A}}^{\dagger}_{X'}$ whose mod $p$ reductions generate $\widetilde{R}_{X'}$ as an $\widetilde{R}_X$-module. Then by \cite[Lemma ~5.5.2]{Kedlaya-LiuI}, $x_1, \ldots, x_n$ generate $\widetilde{\mathbf{A}}^r_{X'}$ as an $\widetilde{\mathbf{A}}^r_X$-module for any sufficiently small $r > 0$. Furthermore, it is shown in the proof of \cite[Proposition ~5.5.3]{Kedlaya-LiuI} that the map $\widetilde{\mathbf{A}}^r_X \rightarrow \widetilde{\mathbf{A}}^r_{X'}$ is finite \'etale for each sufficiently small $r > 0$. This justifies (i) and (ii).
Lastly, it is shown in the proof of \cite[Proposition ~5.5.4]{Kedlaya-LiuI} that for any sufficiently small $r > 0$, the map 
\[
\widetilde{\mathbf{C}}^r_X\otimes_{\widetilde{\mathbf{B}}^r_X}\widetilde{\mathbf{B}}^r_{X'} \rightarrow \widetilde{\mathbf{C}}^r_{X'}
\]
is an isometric isomorphism with respect to $\lambda_r$. This implies the remaining statements.
\end{enumerate}
\end{proof}

\begin{remark}
We caution the reader that the perfectoid tower $X_{\infty}\rightarrow X$ need not be connected; consequently, the associated period rings may admit nontrivial idempotents. We demonstrate this phenomenon with the following example. Assume that $K$ contains $\zeta_{p^n}$ for some $n\ge 1$ and let $h\colon X\rightarrow \D_K^d$ be a toric chart. Consider the $p^n$-th Kummer cover $\D_n\rightarrow \D_K^d$. Taking base change along $h$, we obtain a finite \'etale map $g\colon X'\rightarrow X$. Consider $X'$ equipped with the toric chart $h'\colon X'\rightarrow \D_K^d$ given by $h'=h\circ g$. Note that 
\[(X')_{\infty}=X'\times_{h', \D_K^d} \D_{\infty}\cong X_{\infty}\times (\Z/p^n\Z)^d\]
is isomorphic to $p^{nd}$ copies of $X_{\infty}$. By Proposition~\ref{prop: perfect theory}, the resulting perfect period ring $\widetilde{\mathbf{A}}_{X'}$ (resp. $\widetilde{\mathbf{A}}_{X'}^{\dagger}$, $\widetilde{\mathbf{B}}_{X'}^{\dagger}$, $\widetilde{\mathbf{C}}_{X'}^{\dagger}$) is isomorphic to the direct product of $p^{nd}$ copies of $\widetilde{\mathbf{A}}_{X}$ (resp. $\widetilde{\mathbf{A}}_{X}^{\dagger}$, $\widetilde{\mathbf{B}}_{X}^{\dagger}$, $\widetilde{\mathbf{C}}_{X}^{\dagger}$). Moreover, $\Gamma_n=\Gal(X_{\infty}/X_n)$ is an open normal subgroup of $\Gamma$ of index $p^{nd}$. The $\Gamma$-action on $(X')_{\infty}$ is given by the $\Gamma_n$-action on $X_{\infty}$ and the finite group $\Gamma/\Gamma_n\cong (\Z/p^{n}\Z)^d$ permutes the $p^{nd}$ copies of $X_{\infty}$. The $\Gamma$-actions on the period rings follow a similar pattern.
\end{remark}

\subsection{$(\varphi, \Gamma)$-modules over perfect period rings} \label{subsec: (phi, gamma)-modules-over-perfect-period-rings}

Next, let us recall the notions of $\varphi$-modules and $(\varphi,\Gamma)$-modules over perfect period rings, following \cite{Kedlaya-LiuI}. The notions of $\varphi$-modules and $(\varphi,\Gamma)$-modules introduced in \S\ref{ss:classical_pDE} naturally extend to the setting of perfect relative period rings. A key difference in the relative setting is that we must work with finite projective modules rather than merely finite free modules.

\begin{definition} \label{def:varphi_mod_relative_perfect}
Let $h\colon X \ra \D^d_K$ be a toric chart as before.
\begin{enumerate}
\item A \emph{$\varphi$-module} over $\widetilde{\mathbf{A}}^{\dagger}_{X}$ is a finite projective $\widetilde{\mathbf{A}}^{\dagger}_{X}$-module $M$ together with an isomorphism
\[
\varphi_{M}\colon M\otimes_{\widetilde{\mathbf{A}}^{\dagger}_{X},\varphi} \widetilde{\mathbf{A}}^{\dagger}_{X} \isom M 
\]
of $\widetilde{\mathbf{A}}^{\dagger}_{X}$-modules. There are similar notions over the period rings $\widetilde{\mathbf{A}}_{X}$, $\widetilde{\mathbf{B}}_{X}$, $\widetilde{\mathbf{B}}^{\dagger}_{X}$, and $\widetilde{\mathbf{C}}^{\dagger}_{X}$. 
\item  Let $\mathbf{D}$ be one of the period rings considered above. 
A \emph{$(\varphi, \Gamma)$-module} over $\mathbf{D}$ is a $\varphi$-module $M$ over $\mathbf{D}$, equipped with a continuous semilinear action of $\Gamma$ that commutes with the $\varphi$-action. 
\item A $\varphi$-module (resp. $(\varphi, \Gamma)$-module) over  $\widetilde{\mathbf{B}}_{X}^{\dagger}$ or $ \widetilde{\mathbf{C}}_{X}^{\dagger}$ is called \emph{\'etale} if it comes from a $\varphi$-module (resp. $(\varphi, \Gamma)$-module) over  $\widetilde{\mathbf{A}}^{\dagger}_{X}$ via base change. There is a similar notion over $\widetilde{\mathbf{B}}_{X}$ using the lattice $\widetilde{\mathbf{A}}_X$. 
\end{enumerate}
\end{definition}
 
It is worth pointing out that our $\varphi$-modules over perfect period rings of type $\mathbf A$ are referred to as \'etale $\varphi$-modules in \cite{Cherbonnier_Colmez}. 
We have the following descent results of $\varphi$-modules and $(\varphi, \Gamma)$-modules.

\begin{proposition}[Kedlaya--Liu] \label{prop:descending_phi_mod_along_4_and_10} \hfill
\begin{enumerate}
\item The map $\widetilde{\mathbf{A}}_{X}^{\dagger} \ra \widetilde{\mathbf{A}}_{X}$ induces an equivalence of categories of $\varphi$-modules 
\[ 
\mathrm{Mod}_{/\widetilde{\mathbf{A}}_{X}^{\dagger}}^{\varphi} \isom \mathrm{Mod}_{/\widetilde{\mathbf{A}}_{X} }^{\varphi}.
\]
Consequently, the map $\widetilde{\mathbf{B}}_{X}^{\dagger} \ra \widetilde{\mathbf{B}}_{X}$ induces an equivalence of categories of \'etale $\varphi$-modules
\[ 
\mathrm{Mod}_{/\widetilde{\mathbf{B}}_{X}^{\dagger}}^{\varphi, \ett} \isom \mathrm{Mod}_{/\widetilde{\mathbf{B}}_{X} }^{\varphi, \ett}.
\]
    
\item The map $\widetilde{\mathbf{B}}_{X}^{\dagger}  \ra \widetilde{\mathbf{C}}_{X}^{\dagger}$ induces an equivalence of categories of \'etale $\varphi$-modules 
\[  
\mathrm{Mod}_{/\widetilde{\mathbf{B}}_{X}^{\dagger}}^{\varphi, \ett}  \isom  \mathrm{Mod}_{/\widetilde{\mathbf{C}}_{X}^{\dagger}}^{\varphi, \ett}. 
\]
\end{enumerate} 
\end{proposition}

\begin{proof}
Part (1) is the equivalence between (d) and (e) in \cite[Theorem ~8.5.3]{Kedlaya-LiuI}. Part (2) is the equivalence between (e) and (f) in \cite[Theorem ~8.5.6]{Kedlaya-LiuI}. 
\end{proof}

\begin{corollary} \label{cor: equiv-etale-(phi, Gamma)-mods-over-perfect-period-rings}
The map $\widetilde{\mathbf{A}}_{X}^{\dagger} \ra \widetilde{\mathbf{A}}_{X}$ (resp. $\widetilde{\mathbf{B}}_{X}^{\dagger} \ra \widetilde{\mathbf{B}}_{X}$, $\widetilde{\mathbf{B}}_{X}^{\dagger}  \ra \widetilde{\mathbf{C}}_{X}^{\dagger}$) induces an equivalence
\[
\mathrm{Mod}_{/\widetilde{\mathbf{A}}_{X}^{\dagger}}^{(\varphi, \Gamma)} \isom \mathrm{Mod}_{/\widetilde{\mathbf{A}}_{X} }^{(\varphi, \Gamma)} \quad (\text{resp.} \quad \mathrm{Mod}_{/\widetilde{\mathbf{B}}_{X}^{\dagger}}^{(\varphi, \Gamma), \ett} \isom \mathrm{Mod}_{/\widetilde{\mathbf{B}}_{X} }^{(\varphi, \Gamma), \ett}, \quad \mathrm{Mod}_{/\widetilde{\mathbf{B}}_{X}^{\dagger}}^{(\varphi, \Gamma), \ett}  \isom  \mathrm{Mod}_{/\widetilde{\mathbf{C}}_{X}^{\dagger}}^{(\varphi, \Gamma), \ett}).
\]
\end{corollary}

\begin{proof}
This follows from applying Lemma~\ref{lem: general-fact-equivalence-phi-Gamma-modules} below to the equivalences in Proposition~\ref{prop:descending_phi_mod_along_4_and_10}.    
\end{proof}

\begin{remark}
By definition, the base change functors in Proposition~\ref{prop:descending_phi_mod_along_4_and_10} and Corollary~\ref{cor: equiv-etale-(phi, Gamma)-mods-over-perfect-period-rings} are compatible with tensor structures.  
\end{remark}

Finally, we recall how to associate $(\varphi, \Gamma)$-modules to $p$-adic local systems. Let $\mathrm{Loc}_{\Z_p}(X)$ denote the category of \'etale $\Z_p$-local systems on $X$, and write $\mathrm{Loc}_{\Z_p}(X)_{\Q}$ for its isogeny category. 

\begin{theorem}[Kedlaya--Liu] \label{thm: equivalence-local-system-(phi, Gamma)-modules-perfect-period-rings}
There is a natural equivalence of $\otimes$-categories
\[
\widetilde{D}^{\mathrm{int}}\colon \mathrm{Loc}_{\Z_p} (X)\isom\mathrm{Mod}_{/\widetilde{\mathbf{A}}_X }^{(\varphi,\Gamma)}.
\]
\end{theorem}

\begin{proof}
This follows from the equivalence of (a) and (e) in \cite[Theorem ~8.5.5]{Kedlaya-LiuI}. Let us note that this equivalence reduces to the Riemann--Hilbert correspondence of Katz after reducing modulo $p$ (cf. \cite[Proposition ~4.1.1]{Katz}, also see \cite[Proposition ~3.4]{BS_crystal}). By the construction, this equivalence is compatible with tensor structures.
\end{proof}

\begin{corollary} \label{cor: equiv-local-system-overconvergent-etale-(phi, Gamma)-mods}
There is a natural equivalence of $\otimes$-categories
\[
\widetilde{D}^{\mathrm{int}, \dagger}\colon \mathrm{Loc}_{\Z_p}(X) \isom 
     \mathrm{Mod}^{(\varphi, \Gamma)}_{/\widetilde{\mathbf{A}}_X^{\dagger}}.
\]
Furthermore, this induces natural equivalences of $\otimes$-categories
\[
\widetilde{D}^{\dagger} \colon \mathrm{Loc}_{\Z_p}(X)_{\Q} \isom \mathrm{Mod}^{(\varphi, \Gamma), \ett}_{/\widetilde{\mathbf{B}}_X^{\dagger}}
\]
and
\[
\widetilde{D}_{\mathrm{rig}}^{\dagger} \colon \mathrm{Loc}_{\Z_p}(X)_{\Q} \isom \mathrm{Mod}^{(\varphi, \Gamma), \ett}_{/\widetilde{\mathbf{C}}_X^{\dagger}}.
\]

\end{corollary}

\begin{proof}
This follows from Theorem~\ref{thm: equivalence-local-system-(phi, Gamma)-modules-perfect-period-rings}, Proposition~\ref{prop:descending_phi_mod_along_4_and_10}, and Corollary~\ref{cor: equiv-etale-(phi, Gamma)-mods-over-perfect-period-rings}.    
\end{proof}

\vspace{0.1in}
\section{Weakly decompleting data and imperfect period rings} \label{sec:weakly_decompleting_data} 

In this section, we introduce and study the notion of \emph{weakly decompleting data}. This is a central notion in this article, designed to axiomatize the conditions needed for the construction and study of imperfect relative period rings. Some of the axioms are borrowed from, or inspired by, the notion of \emph{weakly decompleting towers} introduced in \cite[Definition~5.2.3]{Kedlaya-LiuII}. As mentioned in the introduction, several key arguments in \emph{loc. cit.} appear to be incomplete, and we do not know how to fill the gaps in those arguments. Instead, we introduce an \emph{a priori} more restrictive notion and show that these conditions are indeed satisfied in some important examples. We treat these examples in \S\ref{sec:example}.

\begin{remark}
The term \emph{decompleting} is coined to indicate that these conditions allow for decompletion (or \textit{deperfection}) of relative $(\varphi,\Gamma)$-modules. In particular, the axioms are designed to ensure that the resulting imperfect period rings generalize those in the classical setting (cf. \S\ref{ss:point imperfect period ring}) and enjoy the desired properties (one key property is that the ``completed perfections'' of these imperfect period rings should coincide with the corresponding perfect period rings).  
\end{remark}

\begin{remark} 
Later in this article, we refine the notion of weakly decompleting data by introducing the notions of \emph{decompleting data} and \emph{strongly decompleting data} (cf. \S\ref{sec:decompleting_data} and \S\ref{sec:strongly_decompleting_data}). These stronger conditions are designed to yield better descent results and to provide a robust framework for the study of $p$-adic differential equations over the associated imperfect period rings.
\end{remark}

\subsection{Weakly decompleting data and imperfect period rings} \label{subsection: weakly decompleting data}

\noindent Let $h\colon X \ra \D_K^d$ be a toric chart and $\psi = \{\psi_n\colon X_n \ra X\}$ be the corresponding toric tower, as in Definition~\ref{defn: toric chart}. We begin with a pair $\Lambda\coloneqq (R_{\Lambda}, \fA_{\Lambda})$ where
\begin{itemize}
\item $R_{\Lambda}\subset \widetilde{R}_{X}$ is a $\kappa(\!(\overline{\pi}_K)\!)$-subalgebra such that $R_{\Lambda}$ is complete with respect to the norm $|\cdot|$ inherited from $\widetilde{R}_{X}$ (cf. Definition~\ref{construction:A_psi_X_infty}(3));

\item $\fA_{\Lambda}\subset \widetilde{\mathbf{A}}_{X}$ is a $\varphi$-stable $\mathbf{A}_K$-subalgebra which is complete with respect to the $p$-adic topology such that the image of the composite
$
\fA_{\Lambda}/p \rightarrow \widetilde{\mathbf{A}}_{X} /p \cong \widetilde{R}_{X}
$ 
lies in $R_{\Lambda}$. 
\end{itemize}

\begin{definition}[Imperfect period rings] \label{def: imperfect-period-rings-B-C}
Given a pair $\Lambda = (R_{\Lambda}, \mathbf{A}_{\Lambda})$ as above, we consider the following rings.
\begin{enumerate}
\item For any $r>0$, let $\mathbf{A}^r_{\Lambda}\coloneqq \mathbf{A}_{\Lambda}\cap \widetilde{\mathbf{A}}_X^r$ equipped with the $\lambda_r$-norm inherited from $\widetilde{\mathbf{A}}_X^r$. Write $\fB_{\Lambda}\coloneqq \fA_{\Lambda}[p^{-1}]$ and $\mathbf{B}^r_{\Lambda}\coloneqq \mathbf{A}^r_{\Lambda}[p^{-1}]$.

\item Let $\mathbf{A}^{\dagger}_{\Lambda}$ be the union $\mathbf{A}^{\dagger}_{\Lambda}\coloneqq \varinjlim_{r>0} \mathbf{A}^r_{\Lambda}$, and $\mathbf{B}^{\dagger}_{\Lambda}\coloneqq \mathbf{A}^{\dagger}_{\Lambda}[p^{-1}]$. It follows that $\mathbf{A}^{\dagger}_{\Lambda}=\mathbf{A}_{\Lambda}\cap \widetilde{\mathbf{A}}^{\dagger}_X$.

\item For $0<s\le r$, let $\mathbf{C}^{[s,r]}_{\Lambda}$ be the completion of $\mathbf{B}^r_{\Lambda}$ with respect to the norm $\max\{\lambda_r, \lambda_s\}$ (equivalently, completion with respect to the family of norms $\{\lambda_t\}_{t\in [s,r]}$). Let $\mathbf{C}^r_{\Lambda}$ be the Fr\'echet completion of $\mathbf{B}^r_{\Lambda}$ with respect to the family of norms $\{\lambda_s\}_{s\in (0,r]}$, and denote $\mathbf{C}^{\dagger}_{\Lambda}\coloneqq  \varinjlim_{r>0}\mathbf{C}^r_{\Lambda}$.
\end{enumerate}
\end{definition}

By definition, $\mathbf{A}^r_{\Lambda}$ (resp. $\mathbf{A}^{\dagger}_{\Lambda}$, $\mathbf{B}^r_{\Lambda}$, $\mathbf{B}^{\dagger}_{\Lambda}$, $\mathbf{C}^{[s,r]}_{\Lambda}$, $\mathbf{C}^r_{\Lambda}$, $\mathbf{C}^{\dagger}_{\Lambda}$) can naturally be identified as a subring of $\widetilde{\mathbf{A}}^r_X$ (resp. $\widetilde{\mathbf{A}}^{\dagger}_X$, $\widetilde{\mathbf{B}}^r_X$, $\widetilde{\mathbf{B}}^{\dagger}_X$, $\widetilde{\mathbf{C}}^{[s,r]}_X$, $\widetilde{\mathbf{C}}^r_X$, $\widetilde{\mathbf{C}}^{\dagger}_X$). We equip them with the corresponding subspace topologies. In light of this, we introduce the following notation.

\begin{definition} \label{def: Frobenius-pullback-imperfect-period-rings} \hfill
\begin{enumerate}
\item For each $r>0$ and integer $m\ge 0$, write
$ 
R_{\Lambda, (m)}\coloneqq \varphi^{-m}(R_{\Lambda}) \subset \widetilde{R}_X,$ and $ \mathbf{A}^r_{\Lambda, (m)}\coloneqq \varphi^{-m}(\fA^{r/p^m}_{\Lambda}) \subset \widetilde{\mathbf{A}}^r_X
$, 
where $\varphi^m\colon \widetilde{\mathbf{A}}^r_X\xrightarrow[]{\sim} \widetilde{\mathbf{A}}^{r/p^m}_X$ is the $m$-th power of the Frobenius map as in Remark~\ref{remark:Frobenius_relative_and_action_of_Gamma}(3). Let
$ 
\breve{\mathbf{A}}^r_{\Lambda}\coloneqq \bigcup_{m=0}^{\infty}\mathbf{A}^r_{\Lambda, (m)}\subset \widetilde{\mathbf{A}}^r_X$ and $  \breve{\mathbf{A}}^{\dagger}_{\Lambda}\coloneqq \bigcup_{r>0}\breve{\mathbf{A}}^r_{\Lambda}\subset \widetilde{\mathbf{A}}^{\dagger}_X.$ 

\item Denote $\mathbf{B}^r_{\Lambda, (m)}\coloneqq \mathbf{A}^r_{\Lambda, (m)}[p^{-1}] \subset \widetilde{\mathbf{B}}^r_X$, $\breve{\mathbf{B}}^r_{\Lambda}\coloneqq \breve{\mathbf{A}}^r_{\Lambda}[p^{-1}] \subset \widetilde{\mathbf{B}}^r_X$, and $\breve{\mathbf{B}}^{\dagger}_{\Lambda}\coloneqq \bigcup_{r>0}\breve{\mathbf{B}}^r_{\Lambda}\subset \widetilde{\mathbf{B}}^{\dagger}_X$.

\item For each $0 < s \le r$ and integer $m \ge 0$, let
$
\mathbf{C}^{[s,r]}_{\Lambda, (m)}\coloneqq \varphi^{-m}(\mathbf{C}^{[s/p^m, r/p^m]}_{\Lambda}) \subset \widetilde{\mathbf{C}}^{[s,r]}_X
$ 
and 
$
\breve{\mathbf{C}}^{[s,r]}_{\Lambda}\coloneqq \bigcup_{m=0}^{\infty}\mathbf{C}^{[s,r]}_{\Lambda, (m)}\subset \widetilde{\mathbf{C}}^{[s,r]}_X.$ 
Write $\mathbf{C}^r_{\Lambda, (m)}\coloneqq \varprojlim_{0<s\le r} \mathbf{C}^{[s,r]}_{\Lambda, (m)} = \varphi^{-m}(\mathbf{C}^{r/p^m}_{\Lambda}) \subset \widetilde{\mathbf{C}}^r_X$. Finally, let 
$\breve{\mathbf{C}}^r_{\Lambda}\coloneqq \bigcup_{m=0}^{\infty} \mathbf{C}^r_{\Lambda, (m)}\subset \widetilde{\mathbf{C}}^r_X $ and $ \breve{\mathbf{C}}^{\dagger}_{\Lambda}\coloneqq \bigcup_{r>0}\breve{\mathbf{C}}^r_{\Lambda}\subset \widetilde{\mathbf{C}}^{\dagger}_X.$ 
\end{enumerate}
\end{definition}

\begin{remark}[{Frobenius}] \label{remark: Frobenius_wdc_datum}
The Frobenius map $\varphi\colon \widetilde{\mathbf{A}}^r_X\xrightarrow[]{\sim}\widetilde{\mathbf{A}}^{r/p}_X$ induces a Frobenius map $\varphi\colon \mathbf{A}^r_{\Lambda}\rightarrow \mathbf{A}^{r/p}_{\Lambda}$ for all $r$. Similarly, we have Frobenius maps 
$
\varphi\colon \mathbf{B}^r_{\Lambda} \ra \mathbf{B}^{r/p}_{\Lambda},  \varphi\colon \mathbf{C}^r_{\Lambda} \ra \mathbf{C}^{r/p}_{\Lambda}, $ and $ \varphi\colon \mathbf{C}^{[s, r]}_{\Lambda} \ra \mathbf{C}^{[s/p,\,r/p]}_{\Lambda},$ 
as well as on their breve variants $\breve{\mathbf{B}}^r_{\Lambda}, \breve{\mathbf{C}}^r_{\Lambda}, \breve{\mathbf{C}}^{[s,r]}_{\Lambda}$. Note that $\mathbf{B}^{\dagger}_{\Lambda}, \mathbf{C}^{\dagger}_{\Lambda}, \breve{\mathbf{A}}^{\dagger}_{\Lambda}, \breve{\mathbf{B}}^{\dagger}_{\Lambda}, \breve{\mathbf{C}}^{\dagger}_{\Lambda}$ are stable under the Frobenius.
\end{remark}

We now introduce the notion of weakly decompleting datum. 

\begin{definition}[Weakly decompleting datum] \label{defn: weakly decompleting datum}
A pair $\Lambda=(R_{\Lambda}, \mathbf{A}_{\Lambda})$ as above is called a \emph{weakly decompleting datum} for $(X,h)$ if there exists $r_* > 0$ such that the following conditions hold.
\begin{enumerate}
\item The ring $R_{\Lambda}$ admits a finite $p$-basis, i.e., $\varphi^{-1}(R_{\Lambda})$ (considered as a subring of $\widetilde{R}_X$) is a finite free $R_{\Lambda}$-module.

\item The colimit perfection $R^{\mathrm{perf}}_{\Lambda}\coloneqq \varinjlim \; (R_{\Lambda}\xrightarrow[]{\varphi}R_{\Lambda}\xrightarrow[]{\varphi}\cdots )$ is dense in $\widetilde{R}_X$.

\item The subring $\mathbf{A}_{\Lambda}\subset \widetilde{\mathbf{A}}_X$ is stable under the $\Gamma$-action, and coincides with the $p$-adic completion of $\mathbf{A}^{\dagger}_{\Lambda}$.

\item For each $r\in (0, r_*]\cap\Q$, $\mathbf{A}^r_{\Lambda}$ is complete with respect to the $\lambda_r$-norm.

\item For each $r\in (0, r_*]\cap\Q$, the composite
\[\mathbf{A}^r_{\Lambda}/p\rightarrow \mathbf{A}_{\Lambda}/p\rightarrow R_{\Lambda}\]
is an isomorphism (so $\mathbf{A}_{\Lambda}/p\rightarrow R_{\Lambda}$ is also an isomorphism).

\item For each $r\in (0, r_*]\cap\Q$, the mod $p$ reduction map $\mathbf{A}^r_{\Lambda}\rightarrow R_{\Lambda}$ is strict; namely, the quotient (semi)norm on $R_{\Lambda}$ induced from the $\lambda_r$-norm on $\mathbf{A}^r_{\Lambda}$ is equivalent to the norm $|\cdot|^r$ on $R_{\Lambda}$ induced from $\widetilde{R}_X$.

\item The ring $R_{\Lambda}$ is an affinoid $\kappa(\!(\overline{\pi}_K)\!)$-algebra (in the sense of Tate's rigid analytic geometry).

\item For each $0<s \le r \le r_*$ with $s, r \in \mathbb{Q}$, the Banach $K'_0$-algebra $\mathbf{C}^{[s,r]}_{\Lambda}$ is a reduced affinoid $K'_0$-algebra (in the sense of Tate's rigid analytic geometry).
\end{enumerate}
\end{definition}

\begin{lemma} \label{lem: breve-A-dense-weakly-decompleting-datum}
Let $\Lambda = (R_{\Lambda}, \fA_{\Lambda})$ be a weakly decompleting datum. For each $r \in (0, r_*] \cap \Q$, $\breve{\mathbf{A}}^r_{\Lambda}$ is dense in $\widetilde{\mathbf{A}}^r_X$ with respect to the $\lambda_r$-norm.    
\end{lemma}

\begin{proof}
Note that the mod $p$ reduction map $\fA_{\Lambda}^r \rightarrow R_{\Lambda}$ is surjective (see Definition~\ref{defn: weakly decompleting datum}(5)). Thus, the induced map
$
\fA_{\Lambda, (m)}^r \rightarrow R_{\Lambda, (m)}
$ 
is surjective for each $m \ge 0$, and so does the map $\breve{\mathbf{A}}^r_{\Lambda} \rightarrow R_{\Lambda}^{\mathrm{perf}}$. Since $R_{\Lambda}^{\mathrm{perf}}$ is dense in $\widetilde{R}_X$ (Definition~\ref{defn: weakly decompleting datum}(2)), this implies the statement by the definition of $\widetilde{\fA}_X^r$.
\end{proof}

\begin{notation} \label{notn: rs-rational}
Throughout the article, all $r$ and $s$ in the superscripts of the period rings will be assumed to be in $\Q_{>0}$.
\end{notation}

\begin{remark}[{$\Gamma$-actions}] \label{rem: stability-under-Gamma-actions}
The condition (3) above implies $R_{\Lambda} \subset \widetilde{R}_X$ is stable under the $\Gamma$-action since $R_{\Lambda} = \fA_{\Lambda}/p$. Furthermore, by Remark~\ref{remark:Frobenius_relative_and_action_of_Gamma}(4), it also implies that all of the imperfect period rings
\[
\mathbf{A}^r_{\Lambda}, \mathbf{A}^{\dagger}_{\Lambda}, \mathbf{B}_{\Lambda}, \mathbf{B}^r_{\Lambda}, \mathbf{B}^{\dagger}_{\Lambda}, \mathbf{C}^{[s,r]}_{\Lambda}, \mathbf{C}^r_{\Lambda}, \mathbf{C}^{\dagger}_{\Lambda}
\]
and their breve variants
\[
\breve{\mathbf{A}}^r_{\Lambda}, \breve{\mathbf{A}}^{\dagger}_{\Lambda}, \breve{\mathbf{B}}^r_{\Lambda}, \breve{\mathbf{B}}^{\dagger}_{\Lambda}, \breve{\mathbf{C}}^{[s,r]}_{\Lambda}, \breve{\mathbf{C}}^r_{\Lambda}, \breve{\mathbf{C}}^{\dagger}_{\Lambda}
\]
are stable under the $\Gamma$-actions.
\end{remark}

\begin{remark} \label{rem: imperfect-period-rings-affinoid}
Suppose $\Lambda=(R_{\Lambda}, \mathbf{A}_{\Lambda})$ is weakly decompleting.
\begin{enumerate}
\item Since $R_{\Lambda}$ is a subring of the perfectoid ring $\widetilde{R}_{X}$, it is reduced. Note that $|\cdot|$ is a complete power-multiplicative $\kappa(\!(\overline{\pi}_K)\!)$-algebra norm on $R_{\Lambda}$. By \cite[\S 6.2.4, Theorem ~1]{BGR}, the norm $|\cdot|$ coincides with the spectral norm on $R_{\Lambda}$ considered as a reduced affinoid algebra.

\item Similarly, $\mathrm{max}\{\lambda_r, \lambda_s\}$ is a complete power-multiplicative $K'_0$-algebra norm on $\mathbf{C}^{[s,r]}_{\Lambda}$. By the same reason as in (1), $\mathrm{max}\{\lambda_r, \lambda_s\}$ agrees with the spectral norm on $\mathbf{C}^{[s,r]}_{\Lambda}$ viewed as a reduced affinoid $K'_0$-algebra.

\item We remind the reader that $R_{\Lambda}$, $\mathbf{A}^r_{\Lambda}$, and $\mathbf{C}^{[s,r]}_{\Lambda}$ are Banach rings, equipped with norms $|\cdot|$, $\lambda_r$, and $\mathrm{max}\{\lambda_r, \lambda_s\}$, respectively. They contain topologically nilpotent elements $\overline{\pi}_K$, $\pi_K$, and $\pi_K$ respectively, which play the role of $z$ in Assumption \ref{assump: banach-ring-unit-z}.
\end{enumerate}
\end{remark}

\begin{lemma} \label{lem: C breve dense in C tilde}
Suppose $\Lambda = (R_{\Lambda}, \mathbf{A}_{\Lambda})$ is weakly decompleting. 
\begin{enumerate}
\item For each $r\in (0,r_*]$, the map $\mathbf{C}^r_{\Lambda} \rightarrow  \varprojlim_{0<s\le r} \mathbf{C}^{[s,r]}_{\Lambda}$ is an isomorphism.

\item For each $r\in (0,r_*]$ and $s\in (0,r]$, the subring $\breve{\mathbf{C}}^{[s,r]}_{\Lambda}$ is dense inside $\widetilde{\mathbf{C}}^{[s,r]}_X$ with respect to the $\max\{\lambda_r, \lambda_s\}$-norm.
\end{enumerate}    
\end{lemma}

\begin{proof} \hfill
\begin{enumerate}
\item The map is injective by definition. For surjectivity, note that $\fC_{\Lambda}^r$ has a dense image in $\fC_{\Lambda}^{[s, r]}$  for each $0 < s \le r$. Let $x \in \varprojlim_{0<s\le r} \mathbf{C}^{[s,r]}_{\Lambda}$. For each integer $n \ge 1$, we have $x_n \in \fC_{\Lambda}^r$ such that $\lambda_t(x-x_n) \le \frac{1}{n}$ for all $t \in [\frac{r}{n}, r]$. The sequence $\{x_n\}_{n=1}^{\infty}$ is Cauchy for each $\lambda_t$ with $t \in (0, r]$. Thus, its limit $x$ lies in $\fC_{\Lambda}^r$.

\item It suffices to show that $\breve{\fB}_{\Lambda}^r$ is dense in $\widetilde{\mathbf{C}}^{[s,r]}_{X}$ with respect to the $\max\{\lambda_s,\lambda_r\}$-norm. So it suffices to show that $\breve{\fB}_{\Lambda}^r$ is dense in $\widetilde{\mathbf{B}}^r_{X}$ with respect to $\max\{\lambda_s,\lambda_r\}$. Let $x\in \widetilde{\mathbf{B}}^r_{X}$ and take any $0<\varepsilon<1$. We claim that there exists $x'\in \breve{\fB}_{\Lambda}^r$ with $\max\{\lambda_s(x-x'),\lambda_r(x-x')\} < \varepsilon$. By scaling, we may assume $x \in \widetilde{\mathbf{A}}^r_{X}$. Since $\breve{\fA}_{\Lambda}^r$ is dense in $\widetilde{\mathbf{A}}^r_{X}$ with respect to $\lambda_r$ by Lemma~\ref{lem: breve-A-dense-weakly-decompleting-datum}, there exists $x'\in \breve{\fA}_{\Lambda}^r$ with $\lambda_r(x-x')<\varepsilon^{\frac{r}{s}}$. This implies 
$
\lambda_s(x-x')\le (\lambda_r(x-x'))^{\frac{s}{r}} < \varepsilon
$ 
as claimed.
\end{enumerate}
\end{proof}

\subsection{Norm compatible expansions and good lifts} \label{subsec: norm-compatible-expansions-good-lifts}

In the remainder of \S \ref{sec:weakly_decompleting_data}, we study properties of weakly decompleting data. Let $\Lambda = (R_{\Lambda}, \fA_{\Lambda})$ be a weakly decompleting datum for $(X, h)$.

We will often need to express an element in $\widetilde{\mathbf{A}}^r_{X}$ (resp. $\mathbf{A}_{\Lambda}^r$) as a convergent sum of certain ``nice elements'' in $\widetilde{\mathbf{A}}^{r_0}_X$ (resp. $\fA_{\Lambda}^{r_0}$) for some $r_0 > r$. For this purpose, we introduce the notion of \emph{norm compatible expansions}.

\begin{definition}[Norm compatible expansions] \label{defn: norm compatible expansion}
Let $0 < r\le r_0$ and $x\in \widetilde{\mathbf{B}}^r_X$. A \emph{norm compatible expansion} of $x$ \emph{relative to} $\widetilde{\mathbf{A}}^{r_0}_X$ is a $\lambda_r$-convergent sum 
\[ 
x=\sum_{n=n_0}^{\infty} p^n x_n
\] 
for some $n_0\in \mathbb{Z}$ and $x_n\in \widetilde{\mathbf{A}}^{r_0}_X$ satisfying $\lambda_{r_0}(x_n) = |\overline{x}_n|^{r_0}$ for all $n$. Here, $\overline{x}_n \in \widetilde{R}_X$ denotes the reduction of $x_n$ modulo $p$.
\end{definition}

\begin{remark} \label{rmk: equivalent condition}
The condition $\lambda_{r_0}(x_n) = |\overline{x}_n|^{r_0}$ in Definition~\ref{defn: norm compatible expansion} is equivalent to requiring $\lambda_{r'}(x_n)=|\overline{x}_n|^{r'}$ for all $0 < r' \le r_0$. Indeed, for any $y \in \widetilde{\mathbf{A}}^{r_0}_X$, we can uniquely write $y = \sum_{n=0}^{\infty}p^n[y_n]$ with $y_n\in \widetilde{R}_{X}$; in particular, $y_0 = \overline{y}$. Suppose $\lambda_{r_0}(y) = |\overline{y}|^{r_0}$. Then $p^{-n}|y_n|^{r_0}\le |\overline{y}|^{r_0}$ for each $n\ge 0$. So for any $0<r'\le r_0$, we have
$
p^{-n}|y_n|^{r'} \le p^{-\frac{nr'}{r_0}}|y_n|^{r'} \le |\overline{y}|^{r'}.
$ 
This implies $\lambda_{r'}(y)=\sup_{n\ge 0}\{p^{-n}|y_n|^{r'}\} = |\overline{y}|^{r'}$.
\end{remark}

\begin{example} \label{example: norm compatible expansion}
For any $x\in \widetilde{\mathbf{B}}^r_X$, the Witt expansion $x=\sum_{n=n_0}^{\infty}p^n[z_n]$ with $z_n\in \widetilde{R}_{X}$ is a norm compatible expansion of $x$ relative to  $\widetilde{\mathbf{A}}^{r_0}_X$ for any $r_0 \ge r$.
\end{example}

The following fact will be used later.

\begin{lemma} \label{lem: norm-product-of-reductions}
Let $w$ and $w'$ be two non-zero elements in $\widetilde{\mathbf{A}}^{r_0}_X$ such that $\lambda_{r_0}(w) = |\overline{w}|^{r_0}$, $\lambda_{r_0}(w') = |\overline{w'}|^{r_0}$ and $|\overline{w}\cdot \overline{w'}| = |\overline{w}|\cdot |\overline{w'}|$. Then, for any $0< r <r_0$, we have $\lambda_r(w\cdot w') = |\overline{w}\cdot \overline{w'}|^r$.    
\end{lemma}

\begin{proof}
Write $w=\sum_{n=0}^{\infty}p^n[w_n]$ and $w=\sum_{n=0}^{\infty} p^n[w'_n]$ with $w_n, w'_n\in \widetilde{R}_{X}$. In particular, $w_0 = \overline{w}$ and $w'_0 = \overline{w'}$. By the computation in Remark~\ref{rmk: equivalent condition}, for any $0 < r < r_0$, we have $p^{-n}|w_n|^r < |\overline{w}|^r =\lambda_r(w)$ for all $n\ge 1$. Similarly for $w'_n$ and $w'$. Thus, in the expansion 
\[
w\cdot w'=[w_0\cdot w'_0]+\textrm{other terms},
\]
the term
$
\lambda_r([w_0\cdot w'_0]) = |w_0\cdot w'_0|^r = |\overline{w}\cdot \overline{w'}|^r = |\overline{w}|^r \cdot |\overline{w'}|^r
$ 
strictly dominates all the other terms, and so $\lambda_r(w\cdot w')=|\overline{w}\cdot \overline{w'}|^r$ as desired.     
\end{proof}

We have the following nice property for norm compatible expansions.

\begin{proposition} \label{prop: tech lemma norm compatible expansion}
Let $0 < r \le r_0$ and $x\in \widetilde{\mathbf{B}}^r_X$. Suppose $x=\sum_{n=n_0}^{\infty}p^n x_n$ is a norm compatible expansion of $x$ relative to $\widetilde{\mathbf{A}}^{r_0}_X$. Then
\[
\lambda_s(x) = \max_{n\ge n_0}\{p^{-n}\cdot\lambda_s(x_n)\}
\] 
for each $0 <s \le r$.
\end{proposition}

To prove Proposition~\ref{prop: tech lemma norm compatible expansion}, we need the following fact (see Definition~\ref{def:w_k_for_relative_A_tilde} for $w_k$).

\begin{lemma} \label{lemma: w_k}
Let $x\in \widetilde{\mathbf{A}}_{X}$, and write $\overline{x} \in \widetilde{R}_{X}$ for its reduction modulo $p$. Let $r_0 > 0$ and $\widehat{\overline{x}}\in \widetilde{\mathbf{A}}^{r_0}_X$ be a lift of $\overline{x}$ satisfying $\lambda_{r_0}(\widehat{\overline{x}})=|\overline{x}|^{r_0}$. Let $k\in \mathbb{Z}_{\ge 0}$. Then we have 
\[
w_k\left(\frac{x-\widehat{\overline{x}}}{p}\right)\ge \inf\left\{w_{k+1}(x), w_0(x)-\frac{k+1}{r_0}\right\}.
\]
\end{lemma}
 
\begin{proof}
This follows from a similar argument as in the proof of \cite[Lemma ~7.3]{Colmez08}. We recall the proof for the sake of completeness. We have
\[
w_k\left(\frac{x-\widehat{\overline{x}}}{p}\right)=w_{k+1}(x-\widehat{\overline{x}})\ge \inf\left\{w_{k+1}(x), w_{k+1}(\widehat{\overline{x}})\right\}.
\]
Since $w_0(x)=v(\overline{x})$, it remains to show $w_{k+1}(\widehat{\overline{x}})\ge v(\overline{x})-\frac{k+1}{r_0}$. Indeed, we have 
$
v(\overline{x})=v_{r_0}(\widehat{\overline{x}})\le w_{k+1}(\widehat{\overline{x}})+\frac{k+1}{r_0}.$ 
\end{proof}

\begin{proof}[Proof of Proposition~\ref{prop: tech lemma norm compatible expansion}]
We follow a similar argument as in \cite[\S 7]{Colmez08}. By multiplying a $p$-power, we may assume that $x\in \widetilde{\mathbf{A}}^r_{X}$ and $n_0=0$. Consider a sequence $y_0, y_1, \ldots\in \widetilde{\mathbf{A}}^r_{X}$ given by $y_0=x$ and 
\[
y_{n+1}\coloneqq \frac{1}{p^{n+1}}(x-x_0-p\cdot x_1-p^2\cdot x_2-\cdots-p^n\cdot x_n).
\]
Then $y_{n+1}=\frac{y_n-x_n}{p}$ and $\overline{y}_n = \overline{x}_n$. By Lemma~\ref{lemma: w_k}, we have
\[
w_k(y_{n+1})\ge \inf\left\{w_{k+1}(y_n), w_0(y_n)-\frac{k+1}{r_0}\right\}.
\]
By induction on $n$, we deduce that
\[
w_k(y_n)\ge \inf\left\{w_{k+n}(x), \inf_{0\le i\le n-1}\left\{w_i(x)-\frac{k+n-i}{r_0}\right\}\right\}.
\]
In particular, we obtain by taking $k=0$ that 
\[
v(\overline{y}_n)=w_0(y_n)\ge \inf_{0\le i\le n}\left\{w_i(x)-\frac{n-i}{r_0}\right\}.
\] 
Now, we prove that $\lambda_s(x)=\max_{n\ge 0}\{p^{-n}\lambda_s(x_n)\}$ for each $0 <s \le r$; or equivalently, 
$
v_s(x) = \min_{n\ge 0}\{v_s(p^n x_n)\}.
$ 
On the one hand, since $x=\sum_{n=0}^{\infty}p^nx_n$, we have $v_s(x)\ge \min_{n\ge 0}\{v_s(p^nx_n)\}$. On the other hand, for each $n\ge 0$, we have
\begin{align*}
&v_s(p^nx_n)=\frac{n}{s}+v(\overline{x}_n) = \frac{n}{s}+v(\overline{y}_n)\ge \frac{n}{s}+\inf_{0\le i\le n}\left\{w_i(x)-\frac{n-i}{r_0}\right\}\\
&=\inf_{0\le i\le n}\left\{w_i(x)+\frac{i}{s}+(n-i)(\frac{1}{s}-\frac{1}{r_0})\right\}\ge \inf_{0\le i\le n}\left\{w_i(x)+\frac{i}{s}\right\}\ge v_s(x).
\end{align*}
\end{proof}

\begin{proposition}[cf. {\cite[Lemma ~5.2.7, Corollary ~5.2.8]{Kedlaya-LiuII}}] \label{prop: good lift weakly decompleting}
Let $\Lambda = (R_{\Lambda}, \fA_{\Lambda})$ be a weakly decompleting datum for $(X,h)$.
\begin{enumerate}
\item For a sufficiently small $r_0 > 0$, we have that every $y \in R_{\Lambda}$ lifts to some $\widetilde{y} \in \fA_{\Lambda}^{r_0}$ satisfying
\[
\lambda_r(\widetilde{y}-[y])\le p^{-\frac{1}{2}}\cdot|y|^r
\]
for each $0<r\le r_0$. In particular, $\lambda_r(\widetilde{y}) = \lambda_r([y]) = |y|^r$ for each $0<r\le r_0$. We call such $\widetilde{y}$ a \emph{good lift} of $y$.

\item Let $r_0$ be as in part (1). Then for any $0 < r < r_0$ and $x\in \mathbf{A}_{\Lambda}^r$, we have a norm compatible expansion
$x=\sum_{n=0}^{\infty} p^n x_n$ with $x_n\in \fA_{\Lambda}^{r_0}$, such that $x_n$ is a good lift of $\overline{x}_n$ for each $n$. In particular, we have \[\lambda_s(x)=\max_{n\ge 0}\{p^{-n}\lambda_s(x_n)\}\] for each $0<s\le r$.

\item Let $r_0$ be as in part (1). Then for any $0 < s < r\le r_0$, the ring $\fA_{\Lambda}^s$ is the completion of $\mathbf{A}_{\Lambda}^r$ with respect to $\lambda_s$.
\end{enumerate}
\end{proposition}

\begin{proof} 
Part (1) is essentially \cite[Lemma ~5.2.7]{Kedlaya-LiuII} and we recall the proof below for completeness. Part (2) is claimed in \cite[Corollary ~5.2.8]{Kedlaya-LiuII}, but the proof is omitted. We shall explain the proof below in more detail.

\begin{enumerate}
\item Choose any $0 < r_1 < r_*$. By Definition~\ref{defn: weakly decompleting datum}(6), there exists a constant $c \ge 1$ such that every $y\in R$ lifts to some $\widetilde{y}\in \fA_{\Lambda}^{r_1}$ with $\lambda_{r_1}(\widetilde{y})\le c|y|^{r_1}$. In particular, $\lambda_{r_1}(\widetilde{y}-[y])\le c|y|^{r_1}$. By writing $\widetilde{y}-[y]=\sum_{n=1}^{\infty}p^n[z_n]$ with $z_n \in \widetilde{R}_X$, we observe that
\[
p^{-n}\cdot|z_n|^r\le p^{-1+\frac{r}{r_1}}\cdot\left(p^{-n}\cdot |z_n|^{r_1}\right)^{\frac{r}{r_1}}
\]
for any $0 < r \le r_1$ and $n\ge 1$. Hence
\[
\lambda_r(\widetilde{y}-[y])\le p^{-1+\frac{r}{r_1}}\cdot \lambda_{r_1}(\widetilde{y}-[y])^{\frac{r}{r_1}}\le p^{-1+\frac{r}{r_1}}\cdot c^{\frac{r}{r_1}}\cdot |y|^r.
\]
Take any $r_0\in (0,r_1)$ such that $(pc)^{\frac{r_0}{r_1}} \le p^{\frac{1}{2}}$, and the statement holds.

\item We construct two sequences $\{x_n\}$, $\{y_n\}$ in $\fA_{\Lambda}^{r_0}$ inductively. Let $y_0 = x$ and write $\overline{y}_0\in R_{\Lambda}$ for its reduction modulo $p$. Let $x_0\in \fA_{\Lambda}^{r_0}$ be a good lift of $\overline{y}_0$ as in part (1). Suppose we already have $x_0, \ldots, x_n$ and $y_0, \ldots, y_n$. Set $y_{n+1}\coloneqq \frac{1}{p}(y_n-x_n)\in \fA_{\Lambda}$, and let $\overline{y}_{n+1}\in R_{\Lambda}$ be its reduction modulo $p$. Take a good lift $x_{n+1}\in \fA_{\Lambda}^{r_0}$ of $\overline{y}_{n+1}$ as in part (1). 

The construction above provides an expansion $x=\sum_{n=0}^{\infty}p^n x_n$ in $\fA_{\Lambda}$. We need to show that this converges with respect to $\lambda_r$. Note that $\overline{y}_n = \overline{x}_n$ for all $n\ge 0$. As in the proof of Proposition~\ref{prop: tech lemma norm compatible expansion}, we have 
\[
v(\overline{x}_n) = v(\overline{y}_n)\ge \inf_{0\le i\le n}\left\{w_i(x)-\frac{n-i}{r_0}\right\}.
\]
Thus,
\[
v_r(p^nx_n)=\frac{n}{r}+v(\overline{x}_n)\ge \frac{n}{r}+\inf_{0\le i\le n}\left\{w_i(x)-\frac{n-i}{r_0}\right\} =\inf_{0\le i\le n}\left\{w_i(x)+\frac{i}{r}+(n-i)(\frac{1}{r}-\frac{1}{r_0})\right\}.
\]
We claim that the last term goes to $+\infty$ as $n\rightarrow +\infty$. Indeed, since $v(x_i)+\frac{i}{r}\rightarrow +\infty$ as $i\rightarrow +\infty$, we deduce that $w_i(x)+\frac{i}{r}\rightarrow +\infty$ as $i\rightarrow +\infty$. The claim then follows from $\frac{1}{r}-\frac{1}{r_0}>0$. The second statement holds by Proposition~\ref{prop: tech lemma norm compatible expansion}. 

\item This follows directly from part (2).
\end{enumerate} 
\end{proof}

\begin{remark}
In Proposition~\ref{prop: good lift weakly decompleting}(1), we can replace $p^{-\frac{1}{2}}$ by any constant in $(\frac{1}{p}, 1)$ (via a similar proof). But we will not need this generality.
\end{remark}

Recall the notation $R_{\Lambda, (m)}$ and $\mathbf{A}^r_{\Lambda, (m)}$ from Definition~\ref{def: Frobenius-pullback-imperfect-period-rings}.  

\begin{corollary} \label{cor: good lift}
Let $r_0$ be as in Proposition~\ref{prop: good lift weakly decompleting}, and $m \ge 0$ be an integer. 
\begin{enumerate}
\item Every $y \in R_{\Lambda, (m)}$ lifts to some $\widetilde{y} \in \fA_{\Lambda, (m)}^{r_0}$ such that
\[
\lambda_r(\widetilde{y} - [y]) \le p^{-\frac{1}{2}}\cdot |y|^r 
\]
for any $0 < r \le r_0$. In particular, $\lambda_r(\widetilde{y}) = \lambda_r([y]) = |y|^r$ for each $0<r\le r_0$. We still refer to such $\widetilde{y}$ as a \emph{good lift} of $y$.

\item For any $0 < r < r_0$ and $x \in \fA_{\Lambda, (m)}^r$, there is a norm compatible expansion $x = \sum_{n=0}^{\infty} p^n x_n$ with $x_n \in \fA_{\Lambda, (m)}^{r_0}$, such that $x_n$ is a good lift of $\overline{x}_n$ for each $n$. In particular, 
\[
\lambda_s(x) = \max_{n\ge 0}\{p^{-n}\lambda_s(x_n)\}
\]
for each $0 < s \le r$.

\item The statements analogous to (1) and (2) remain true if we replace $R_{\Lambda, (m)}$ and $\mathbf{A}^r_{\Lambda, (m)}$ by $R_{\Lambda, (m)} / R_{\Lambda}$ and $\fA_{\Lambda, (m)}^r / \fA_{\Lambda}^r$, respectively. Here, $R_{\Lambda, (m)} / R_{\Lambda}$ and $\fA_{\Lambda, (m)}^r / \fA_{\Lambda}^r$ are equipped with the quotient norms (still denoted by $|\cdot|$ and $\lambda_r$ respectively, by abuse of notation).
\end{enumerate}
\end{corollary}

\begin{proof} \hfill
\begin{enumerate}
\item Let $y \in R_{\Lambda, (m)}$. By Proposition~\ref{prop: good lift weakly decompleting}(1), $z\coloneqq \varphi^m(y) \in R_{\Lambda}$ lifts to some $\widetilde{z} \in \fA_{\Lambda}^{\frac{r_0}{p^m}}$ with $\lambda_{\frac{r}{p^m}}(\widetilde{z}-[z]) \le p^{-\frac{1}{2}}\cdot |z|^{\frac{r}{p^m}}$ for any $0 < r \le r_0$. Thus, letting $\widetilde{y} = \varphi^{-m}(\widetilde{z}) \in \fA_{\Lambda, (m)}^{r_0}$, we have
\[
\lambda_r(\widetilde{y} - [y]) = \lambda_{\frac{r}{p^m}}(\widetilde{z} - [z]) \le  p^{-\frac{1}{2}}\cdot |z|^{\frac{r}{p^m}} = p^{-\frac{1}{2}}\cdot |y|^r.
\]

\item Given part (1), this follows from a similar argument as in the proof of Proposition~\ref{prop: good lift weakly decompleting}(2).

\item Let $y \in R_{\Lambda, (m)}/R_{\Lambda}$. Choose a lift $y_1 \in R_{\Lambda, (m)}$ of $y$ such that $|y_1|$ is equal to the quotient norm of $y$. Choose a good lift $\widetilde{y}_1 \in \fA_{\Lambda, (m)}^{r_0}$ of $y_1$ given by part (1), so that $\lambda_r(\widetilde{y}_1-[y_1]) \le p^{-\frac{1}{2}}\cdot |y_1|^r$ for any $0 < r \le r_0$. It follows directly from part (2) that $\lambda_r(\widetilde{y}_1)$ is equal to the quotient norm of the image of $\widetilde{y}_1$ in $\fA_{\Lambda, (m)}^r / \fA_{\Lambda}^r$ for any $0 < r < r_0$. This implies the statement analogous to part (1), and that analogous to part (2) can be deduced similarly.  
\end{enumerate}
\end{proof}

\begin{remark} \label{rem: variation-of-good-lifts}
In the situation of Corollary~\ref{cor: good lift}, let $y \in R_{\Lambda, (m)}$ and a good lift $\widetilde{y} \in \mathbf{A}_{\Lambda, (m)}^{r_0}$ of $y$ be given. For any $y' \in R_{\Lambda, (m)}$ with $|y'-y| < |y|$, choose a good lift $\widetilde{z} \in \mathbf{A}_{\Lambda, (m)}^{r_0}$ of $y'-y$. Then $\widetilde{y}+\widetilde{z} \in \mathbf{A}_{\Lambda, (m)}^{r_0}$ is a good lift of $y'$ with $\lambda_r((\widetilde{y}+\widetilde{z})-\widetilde{y}) = \lambda_r(\widetilde{z}) = |y'-y|^r$ for any $0 < r \le r_0$.     
\end{remark}

\subsection{Strict decompositions of imperfect period rings} \label{subsec: strict-decompositions-imperfect-period-rings}

To prepare for later sections, we study certain strict decompositions of Frobenius pullbacks of imperfect period rings. Some of the results are stated in \cite[\S 5.8]{Kedlaya-LiuII} in a slightly more general setting, but with certain proofs omitted. We provide here details specialized to our setting. Keep the setup that $\Lambda = (R_{\Lambda}, \fA_{\Lambda})$ is a weakly decompleting datum for $(X, h)$. Fix a choice of $p$-basis $\{e_1 = 1, e_2, \ldots, e_{\ell}\}$ of $R_{\Lambda, (1)} = \varphi^{-1}(R_{\Lambda})$ as an $R_{\Lambda}$-module, so that
\[
R_{\Lambda, (1)} \cong \bigoplus_{i=1}^{\ell} R_{\Lambda}\cdot e_i.
\]
This induces an isomorphism of $R_{\Lambda}$-modules
\begin{equation} \label{eq: p-basis-R}
R_{\Lambda, (1)} \cong R_{\Lambda}\oplus (R_{\Lambda, (1)}/R_{\Lambda}).
\end{equation}
Since $\overline{\pi}_K \in R_{\Lambda}$ is a topologically nilpotent unit, the above isomorphism is strict by the open mapping theorem (cf. Theorem \ref{thm: open mapping thm}). Here, the left hand side is equipped with the norm $|\cdot|$ inherited from $\widetilde{R}_X$. For the right hand side, $R_{\Lambda}$ is equipped with the given norm $|\cdot|$, and $R_{\Lambda, (1)}/R_{\Lambda}$ is equipped with the quotient norm. For each integer $j \ge 0$, we have the induced isomorphism
\[
R_{\Lambda, (j+1)}\cong \bigoplus_{i=1}^{\ell} R_{\Lambda, (j)}\cdot \varphi^{-j}(e_i) \cong  R_{\Lambda, (j)}\oplus (R_{\Lambda, (j+1)})/R_{\Lambda, (j)}))
\]
of $R_{\Lambda, (j)}$-modules. By induction, these further induce an isomorphism of $R_{\Lambda}$-modules 
\[
R_{\Lambda, (m+1)} \cong R_{\Lambda}\oplus \bigoplus_{i=0}^m R_{\Lambda, (i+1)} / R_{\Lambda, (i)}
\]
for each $m \ge 0$.

\begin{lemma} \label{lem: Frob decomp R uniform strict}
The above isomorphism of $R_{\Lambda}$-modules
\begin{equation} \label{eq: Frobenius-decomposition-R}
R_{\Lambda, (m+1)} \cong R_{\Lambda}\oplus \bigoplus_{i=0}^m R_{\Lambda, (i+1)} / R_{\Lambda, (i)}    
\end{equation}
is uniformly strict with respect to $m$; namely, the constants $c_1, c_2$ as in Definition \ref{defn: equivalent seminorms}(1) can be taken to be independent of $m$. Here, the left-hand side is equipped with the norm $|\cdot|$ inherited from $\widetilde{R}_X$. The right-hand side is equipped with the sup norm given by the direct sum presentation, where each $R_{\Lambda, (i+1)} / R_{\Lambda, (i)}$ is equipped with the quotient norm.
\end{lemma}

\begin{proof}
Let $C\ge 1$ be a constant such that for the isomorphism~\eqref{eq: p-basis-R}, the norm on the right is bounded above (resp. below) by that on the left by $C$ (resp. $C^{-1}$). Then for the induced isomorphism
\[
R_{\Lambda, (i+1)} \cong R_{\Lambda, (i)} \oplus (R_{\Lambda, (i+1)}/R_{\Lambda, (i)})
\]
of $R_{\Lambda, (i)}$-modules, the norm on the right is bounded above (resp. below) by that on the left by $C^{\frac{1}{p^i}}$ (resp. $C^{-\frac{1}{p^i}}$). Since $\sum_{i=0}^{\infty}\frac{1}{p^i} = \frac{p}{p-1}$, for any $m \ge 0$ and the isomorphism
\[
R_{\Lambda, (m+1)} \cong R_{\Lambda}\oplus \bigoplus_{i=0}^m R_{\Lambda, (i+1)} / R_{\Lambda, (i)}, 
\]
the norm on the right is bounded above (resp. below) by that on the left by $C^{\frac{p}{p-1}}$ (resp. $C^{-\frac{p}{p-1}}$).
\end{proof}

Let $r_0$ be as in Proposition~\ref{prop: good lift weakly decompleting}. For each $1 \le i \le \ell$, let $\widetilde{e}_i \in \mathbf{A}^{r_0}_{\Lambda, (1)}$ be a good lift of $e_i$ as in Corollary~\ref{cor: good lift}(1).

\begin{lemma} \label{lem: Frob decomp A}
Let $r_0$ be as in Proposition~\ref{prop: good lift weakly decompleting}. There exists $0 < r_1 \le r_0$ such that for any $0 < r \le r_1$, the natural map
\[
\bigoplus_{i=1}^{\ell} \mathbf{A}_{\Lambda}^r\cdot \widetilde{e}_i \rightarrow \mathbf{A}^r_{\Lambda, (1)}
\]
is an isomorphism of $\mathbf{A}_{\Lambda}^r$-modules.
\end{lemma}

\begin{proof}
First recall that we have an isomorphism of $R_{\Lambda}$-modules
\begin{equation} \label{eq: R-frobenius-decomposition-finite-module-norm}
R_{\Lambda, (1)} \cong \bigoplus_{i=1}^{\ell} R_{\Lambda}\cdot e_i.    
\end{equation}
Equip the right hand side with the finite-module norm given by the presentation (cf. Definition~\ref{defn: finite-module seminorm} and Remark~\ref{rmk: finite-module seminorm}). Since $\overline{\pi}_K \in R_{\Lambda}$ is a topologically nilpotent unit, the above isomorphism is strict by Theorem~\ref{thm: open mapping thm}. Let $C' \ge 1$ be a constant such that the norm on the right in \eqref{eq: R-frobenius-decomposition-finite-module-norm} is bounded above by that on the left by $C'$. Write
$
D\coloneqq C'\cdot \max\{1, |e_1|, \ldots, |e_{\ell}|\}.
$ 
We now consider the map
\[
\bigoplus_{i=1}^{\ell} \mathbf{A}_{\Lambda}^r\cdot \widetilde{e}_i \rightarrow \mathbf{A}^r_{\Lambda, (1)}.
\]
Since $\widetilde{\mathbf{A}}_X$ is $p$-torsion free and $p$-adically separated, the injectivity of the map follows from the isomorphism~\eqref{eq: R-frobenius-decomposition-finite-module-norm}. 

For surjectivity, we use a similar argument as in the proof of \cite[Lemma~5.5.2]{Kedlaya-LiuI}. Let $0 < r \le r_0$ and $z \in \mathbf{A}^r_{\Lambda, (1)}$. For each integer $n \ge 0$, we construct $z_n \in \mathbf{A}^r_{\Lambda, (1)}$ and $\widetilde{x}_{n, 1}, \ldots, \widetilde{x}_{n, \ell} \in \mathbf{A}_{\Lambda}^r$ inductively as follows. Let $z_0 = z$. Given $z_n$, write
$
\overline{z}_n = \sum_{i=1}^{\ell} x_{n, i}\cdot e_i
$ 
via the isomorphism~\eqref{eq: R-frobenius-decomposition-finite-module-norm}. For each $1 \le i \le \ell$, choose a good lift $\widetilde{x}_{n, i} \in \mathbf{A}_{\Lambda}^{r_0}$ as in Proposition~\ref{prop: good lift weakly decompleting}(1). Then $z_n - (\sum_{i=1}^d \widetilde{x}_{n, i}\cdot \widetilde{e}_i) \in p\mathbf{A}^r_{\Lambda, (1)}$. Set
\[
z_{n+1} = \frac{1}{p}(z_n - (\sum_{i=1}^{\ell} \widetilde{x}_{n, i}\cdot \widetilde{e}_i)).
\]
Note that
\[
\lambda_{r_0}(\widetilde{x}_{n, i}\cdot \widetilde{e}_i) \le \lambda_{r_0}(\widetilde{x}_{n, i})\cdot \lambda_{r_0}(\widetilde{e}_i) = |x_{n, i}|^{r_0} |e_i|^{r_0} \le D^{r_0} |\overline{z}_n|^{r_0}.
\]
So we have $\lambda_{r_0}(z_{n+1}) \le pD^{r_0}\lambda_{r_0}(z_n)$, which implies $\lambda_{r_0}(z_n) \le (pD^{r_0})^n \lambda_{r_0}(z_0)$ for each $n \ge 0$. We deduce
\[
|x_{n, i}| \le D|e_i|^{-1} (\lambda_{r_0}(z_0))^{\frac{1}{r_0}}\cdot (Dp^{\frac{1}{r_0}})^n.
\]
Thus, by choosing $0 < r_1 \le r_0$ sufficiently small so that $D^{r_1}p^{\frac{r_1}{r_0}-1} < 1$, we see that $\sum_{n = 0}^{\infty} p^n \widetilde{x}_{n, i}$ converges under $\lambda_r$-norm to an element in $\mathbf{A}_{\Lambda}^r$ for any $0 < r \le r_1$. By construction, we have
$
z = \sum_{i=1}^{\ell} (\sum_{n=0}^{\infty} p^n \widetilde{x}_{n, i})\cdot \widetilde{e}_i
$ as desired.  
\end{proof}

By Lemma~\ref{lem: Frob decomp A}, we have an isomorphism of $\fA_{\Lambda}^r$-modules
\[
\mathbf{A}^r_{\Lambda, (1)} \cong \mathbf{A}^r_{\Lambda} \oplus (\mathbf{A}^r_{\Lambda, (1)} / \mathbf{A}^r_{\Lambda})
\]
for each $0 < r \le r_1$. For each $m \ge 0$, this induces an isomorphism of $\fA_{\Lambda, (m)}^r$-modules
\[
\mathbf{A}^r_{\Lambda, (m+1)} \cong \mathbf{A}^r_{\Lambda, (m)} \oplus (\mathbf{A}^r_{\Lambda, (m+1)} / \mathbf{A}^r_{\Lambda, (m)}).
\]

\begin{corollary} \label{cor: Frob decomp A uniform strict}
Let $r_1$ be as in Lemma~\ref{lem: Frob decomp A}. There exists $0 < r_2 \le r_1$ such that for any $0 < r \le r_2$, the induced isomorphism of $\mathbf{A}_{\Lambda}^r$-modules
\begin{equation} \label{eq: Frobenius-decomposition-A}
\mathbf{A}^r_{\Lambda, (m+1)} \cong \fA_{\Lambda}^r\oplus \bigoplus_{i=0}^m \fA_{\Lambda, (i+1)}^r/\fA_{\Lambda, (i)}^r 
\end{equation}
is uniformly strict with respect to $m$. Here, the norm on the left-hand side is the $\lambda_r$-norm, while the right-hand side is equipped with the sup-norm given by the direct sum presentation, where each $\fA_{\Lambda, (i+1)}^r/\fA_{\Lambda, (i)}^r$ is equipped with the quotient norm.
In particular, for any $0 < s \le r \le r_2$, the induced isomorphism  
\[
\fC^{[s, r]}_{\Lambda, (m+1)} \cong \fC_{\Lambda}^{[s, r]}\oplus \bigoplus_{i=0}^m \fC_{\Lambda, (i+1)}^{[s, r]}/\fC_{\Lambda, (i)}^{[s, r]}
\]
of $\mathbf{C}_{\Lambda}^{[s, r]}$-modules 
is uniformly strict with respect to $m$. 
\end{corollary}

\begin{proof}
First note that in the isomorphism~\eqref{eq: Frobenius-decomposition-A}, the $\lambda_r$-norm on the left hand side is less than or equal to the norm on the right hand side. So it suffices to find a uniform constant $C_1 \ge 1$ such that the norm on the right is bounded above by $C_1$ times the $\lambda_r$-norm on the left.

The mod $p$ reduction of \eqref{eq: Frobenius-decomposition-A} gives the isomorphism (\ref{eq: Frobenius-decomposition-R}) in Lemma~\ref{lem: Frob decomp R uniform strict}. Let $C$ be a uniform bound as in the proof of Lemma~\ref{lem: Frob decomp R uniform strict}.
Let $z \in \fA_{\Lambda, (m+1)}^r$. By Corollary~\ref{cor: good lift}(2), we may assume that $z$ is a good lift of $\overline{z}$ without loss of generality. We first write down the decomposition of $z$ in \eqref{eq: Frobenius-decomposition-A}, by a similar construction as in the proof of Lemma~\ref{lem: Frob decomp A}. For each $n \ge 0$, we set $z_n \in \fA_{\Lambda, (m+1)}^r$, $\widetilde{x}_n \in \fA_{\Lambda}^r$, and $\widetilde{x}_{n, i} \in \fA_{\Lambda, (i+1)}^r/\fA_{\Lambda, (i)}^r$ for $0 \le i \le m$ inductively as follows. Let $z_0 = z$. Given $z_n$, write
$
\overline{z}_n = x_n+\sum_{i=0}^m x_{n, i}
$ 
via the isomorphism~\eqref{eq: Frobenius-decomposition-R}. Choose a good lift $\widetilde{x}_n$ of $x_n$ and a good lift $\widetilde{x}_{n, i}$ of $x_{n, i}$ as in Corollary~\ref{cor: good lift}. Then $z_n-(\widetilde{x}_n+\sum_{i=0}^m \widetilde{x}_{n, i}) \in p\fA_{\Lambda, (m+1)}^r$. Let
\[
z_{n+1} = \frac{1}{p}(z_n-(\widetilde{x}_n+\sum_{i=0}^m \widetilde{x}_{n, i})).
\]

Now, by a similar computation as in the proof of Lemma~\ref{lem: Frob decomp A}, we have
\[
|x_n|, ~~|x_{n, i}| \le C (\lambda_{r_1}(z))^{\frac{1}{r_1}} \cdot (Cp^{\frac{1}{r_1}})^n
\]
for each $n$ and $i$. Choose $0 < r_2 \le \min\{r_1, 1\}$ so that $C^{r_2}p^{\frac{r_2}{r_1}-1} < 1$. Then for any $0 < r \le r_2$, we have
\[
\lambda_r(p^n\widetilde{x}_n) = p^{-n}|x_n|^r \le C^r (\lambda_{r_1}(z))^{\frac{r}{r_1}} = C^r \lambda_r(z) \le C \lambda_r(z).
\]
Here, we used the assumption that $z$ is a good lift of $\overline{z}$, so $\lambda_s(z) = |z|^s$ for any $0 < s \le r_0$. Similarly, we obtain
$
\lambda_r(p^n\widetilde{x}_{n, i}) \le C\lambda_r(z)
$ 
for each $0 \le i \le m$. Since the decomposition of $z$ is given by
\[
z = \sum_{n=0}^{\infty} (p^n\widetilde{x}_n+\sum_{i=0}^m p^n\widetilde{x}_{n, i}),
\]
we deduce that the norm on the right of \eqref{eq: Frobenius-decomposition-A} is bounded above by $C$ times the $\lambda_r$-norm on the left.
\end{proof}

\subsection{Intersection of period rings} \label{subsec: intersection-period-rings}

In this subsection, we verify that, under the weakly decompleting assumption, the intersections of period rings behave as expected. These results are stated in \cite[Lemma ~5.2.10, 5.2.12]{Kedlaya-LiuII}, but the proofs outlined in \emph{loc. cit.} seem to be incomplete. We provide here necessary details.

\begin{proposition}[cf. {\cite[Lemma ~5.2.10]{Kedlaya-LiuII}}]  \label{prop: intersection A}
Let $r_0$ be as in Proposition~\ref{prop: good lift weakly decompleting}(1). For any $0<s\le r<r_0$, we have 
$
\widetilde{\mathbf{A}}^s_X\cap \mathbf{C}_{\Lambda}^{[s,r]}=\mathbf{A}_{\Lambda}^r,
$  
where the intersection is taken inside $\widetilde{\mathbf{C}}^{[s,s]}_X$. Consequently, we have 
\[
\widetilde{\mathbf{A}}^{\dagger}_X \cap \mathbf{C}^{\dagger}_{\Lambda} = \mathbf{A}^{\dagger}_{\Lambda}
\]
where the intersection takes place in $\widetilde{\mathbf{C}}^{\dagger}_X$.
\end{proposition}

\begin{proof} 
It is clear that $\mathbf{A}_{\Lambda}^r\subset \widetilde{\mathbf{A}}^s_{X} \cap \mathbf{C}_{\Lambda}^{[s,r]}$. For the converse inclusion, let $x\in \widetilde{\mathbf{A}}^s_{X} \cap \mathbf{C}_{\Lambda}^{[s,r]}$. By definition, there exists a sequence $x_0, x_1,\ldots$ in $\mathbf{B}^r_{\Lambda}$ such that $\lambda_t(x_i-x)\le p^{-i}$ for all $i\ge 0$ and all $t \in [s,r]$. Applying Proposition~\ref{prop: good lift weakly decompleting}(2) to each $x_i$, we can write $x_i$ as a $\lambda_r$-convergent sum 
\[
x_i = \sum_{n\ge n_i}p^n x_{i, n}
\] for some integer $n_i \le -1$ and $x_{i, n}\in \fA_{\Lambda}^{r_0}$ satisfying $\lambda_{r_0}(x_{i, n})=|\overline{x}_{i,n}|^{r_0}$ for all $n$ (so by Remark~\ref{rmk: equivalent condition}, $\lambda_{r'}(x_{i,n})=|\bar{x}_{i,n}|^{r'}$ for any $n$ and $0< r' \le r_0$).

For each $i$, let $x'_i\coloneqq \sum_{n=0}^{\infty} p^n x_{i,n}$ and $x''_i\coloneqq \sum_{n=n_i}^{-1}p^n x_{i, n}$. Then $x_i-x = x''_i+(x'_i-x)$. Since $x\in \widetilde{\mathbf{A}}^s_{X}$ and $x'_i\in \mathbf{A}_{\Lambda}^r \subset \widetilde{\mathbf{A}}^r_{X}\subset \widetilde{\mathbf{A}}^s_{X}$, we have $x'_i-x\in \widetilde{\mathbf{A}}^s_{X}$. By Example~\ref{example: norm compatible expansion}, the Witt expansion $x'_i-x=\sum_{n=0}^{\infty}p^n[z_n]$ is a norm compatible expansion of $x'_i-x$. Thus,
\[
x_i-x=\sum_{n=n_i}^{-1}p^nx_{i,n}+\sum_{n=0}^{\infty}p^n[z_n]
\]
is a norm compatible expansion of $x_i-x$ relative to $\widetilde{\mathbf{A}}^{r_0}_X$. By Proposition~\ref{prop: tech lemma norm compatible expansion}, we have
\[
p^{-n}\cdot\lambda_s(x_{i,n})\le \lambda_s(x-x_i)\le p^{-i}
\]
for each $i\ge 0$ and $n_i\le n\le -1$. In particular, $\lambda_s(x_{i,n})\le p^{n-i}$, and so $|\overline{x}_{i,n}|\le p^{\frac{n-i}{s}}$. Note that we have 
$
p^{-n}\cdot p^{\frac{(n-i)r}{s}}\le p^{1-\frac{(i+1)r}{s}}
$ 
when $n\le -1$. We thus obtain
\[
\lambda_r(x''_i)\le \max_{n_i\le n\le -1}\left\{p^{-n}\cdot \lambda_r(x_{i,n})\right\}=\max_{n_i\le n\le -1}\left\{p^{-n}\cdot |\overline{x}_{i,n}|^r\right\}\le p^{1-\frac{(i+1)r}{s}}
\]
for all $i\ge 0$. Hence, the sequence $x'_0, x'_1, \ldots \in \mathbf{A}_{\Lambda}^r$ converges to $x$ under the $\lambda_r$-norm, and $x\in \mathbf{A}_{\Lambda}^r$.
\end{proof}

\begin{proposition}[cf. {\cite[Lemma ~5.2.12]{Kedlaya-LiuII}}] \label{prop: intersection C}
Let $r_0$ be as in Proposition~\ref{prop: good lift weakly decompleting}(1). For any $0<s\le s'\le r\le r'<r_0$, we have 
$
\mathbf{C}_{\Lambda}^{[s, r]} \cap \fC_{\Lambda}^{[s', r']} = \fC_{\Lambda}^{[s, r']}
$ 
where the intersection is taken inside $\fC_{\Lambda}^{[s',r]}$.
\end{proposition}

We first show the following lemma, which is a slightly enhanced version of \cite[Lemma ~5.2.11]{Kedlaya-LiuII}.

\begin{lemma}[cf. {\cite[Lemma ~5.2.11]{Kedlaya-LiuII}}] \label{lemma: C=A+C}
Let $r_0$ be as in Proposition~\ref{prop: good lift weakly decompleting}(1). For any $0<s\le r\le r'<r_0$, each $x\in \mathbf{C}_{\Lambda}^{[s,r]}$ can be written as $x=x'+x''$ with $x'\in \mathbf{A}_{\Lambda}^r$ and $x''\in \fC_{\Lambda}^{[s,r']}$.
\end{lemma}

\begin{proof} 
The proof is similar to the one for Proposition~\ref{prop: intersection A}. First, note that $\mathbf{C}_{\Lambda}^{[s,r]}$ is the completion of $\mathbf{B}^r_{\Lambda}$ with respect to $\{\lambda_t\}_{t\in [s,r]}$, and $\mathbf{A}_{\Lambda}^r$ is the completion of $\fA_{\Lambda}^{r'}$ with respect to $\lambda_r$ by Proposition~\ref{prop: good lift weakly decompleting}(3). So each $x\in \mathbf{C}_{\Lambda}^{[s,r]}$ can be written as a sum $x=x_0+x_1+\cdots$ with $x_0, x_1, \ldots$ in $\fB_{\Lambda}^{r'}$ such that the sum is convergent with respect to $\{\lambda_t\}_{t\in [s,r]}$, and
$
\lambda_t\left(x-(x_0+x_1+\cdots+x_i)\right)\le p^{-i-1}
$
for all $i \ge 0$ and $t \in [s, r]$. In particular, $\lambda_t(x_i)\le p^{-i}$ for all $i\ge 1$ and $t \in [s, r]$.
By Proposition~\ref{prop: good lift weakly decompleting}(2), each $x_i$ admits a norm compatible expansion  
$
x_i=\sum_{n=n_i}^{\infty}p^n x_{i,n}
$ 
for some $n_i\le -1$ and $x_{i,n}\in \fA_{\Lambda}^{r_0}$, which is convergent under $\lambda_{r'}$. 
By Remark~\ref{rmk: equivalent condition}, we have $\lambda_{r'}(x_{i,n}) = |\overline{x}_{i, n}|^{r'}$ for each $0 < r' \le r_0$, $i\ge 0$, and $n\ge n_i$. By Proposition~\ref{prop: tech lemma norm compatible expansion}, we have 
\[
\lambda_t(x_i)=\max_{n\ge n_i} \{p^{-n}\cdot\lambda_t(x_{i,n})\}
\] 
for all $t\in [s,r']$. Let $x'_i\coloneqq \sum_{n=0}^{\infty} p^n x_{i,n} \in \fA_{\Lambda}^{r'}$ and $x''_i\coloneqq \sum_{n=n_i}^{-1} p^n x_{i,n} \in \fB_{\Lambda}^{r_0}$.
We claim that the sum $\sum_{j=0}^{\infty} x''_j$ converges to an element $x'' \in \fC_{\Lambda}^{[s, r']}$. For any $t \in [s,r]$, we have
\[
p^{-n}\cdot |\overline{x}_{i,n}|^t=p^{-n}\cdot \lambda_t(x_{i,n})\le \lambda_t(x_i)\le p^{-i}
\]
for each $i \ge 1$ and $n_i\le n\le -1$. In particular, $|\overline{x}_{i,n}|\le p^{\frac{n-i}{s}}$ for all $i\ge 1$ and $n_i\le n\le -1$. Thus for any $t \in [s, r']$, 
\[
\lambda_t(x''_i)\le \max_{n_i\le n\le -1}\{p^{-n}\lambda_t(x_{i,n})\} = \max_{n_i\le n\le -1}\{p^{-n}|\overline{x}_{i,n}|^t\}\le \max_{n_i\le n\le -1}\{p^{-n}\cdot p^{\frac{(n-i)t}{s}}\}\le  p^{1-\frac{(i+1)t}{s}}.
\]
This yields the claim.

On the other hand, the sum $\sum_{j=0}^{\infty} x'_j$ converges to an element $x'\in \mathbf{A}_{\Lambda}^r$, since 
\[
\lambda_r(x'_i) \le \sup_{n\ge 0}\{p^{-n}\lambda_r(x_{i,n})\} \le \lambda_r(x_i)
\]
which goes to $0$ as $i\rightarrow +\infty$.
\end{proof}

\begin{proof}[Proof of Proposition~\ref{prop: intersection C}]
It is clear that $\fC_{\Lambda}^{[s,r']} \subset \mathbf{C}_{\Lambda}^{[s,r]}\cap \fC_{\Lambda}^{[s', r']}$. For the converse inclusion, let $x\in \mathbf{C}_{\Lambda}^{[s,r]}\cap \fC_{\Lambda}^{[s', r']}$. By Lemma~\ref{lemma: C=A+C}, we can write $x=x'+x''$ with $x'\in \mathbf{A}_{\Lambda}^r$ and $x''\in \fC_{\Lambda}^{[s, r']}$. By Proposition~\ref{prop: intersection A}, we obtain
\[
x'\in \mathbf{A}_{\Lambda}^r\cap \fC_{\Lambda}^{[s', r']}\subset \widetilde{\mathbf{A}}^r_{X}\cap \fC_{\Lambda}^{[r, r']} = \fA_{\Lambda}^{r'}\subset \fC_{\Lambda}^{[s,r']}.
\]
Hence, $x = x'+x''\in \fC_{\Lambda}^{[s, r']}$.
\end{proof}

\subsection{Weakly decompleting data under finite \'etale extensions} \label{ss:weakly decompleting finite etale}

Let $h\colon X\rightarrow \mathbb{D}_K^d$ be a toric chart as before. Let $g\colon X'\rightarrow X$ be a finite \'etale morphism, and let $h'\colon X'\rightarrow \mathbb{D}_K$ be the composition $h\circ g$. The goal of this section is to show that, given any weakly decompleting datum for $(X, h)$, there is a natural way to obtain a weakly decompleting datum for $(X', h')$. We first observe the following equivalences of categories between finite \'etale algebras.

\begin{proposition} \label{prop:equivalence of categories finite etale}
Let $\Lambda = (R_{\Lambda}, \fA_{\Lambda})$ be a weakly decompleting datum for $(X,h)$. Then the equivalences in Proposition~\ref{prop: perfect theory}(1) extend to a natural diagram of equivalences of categories:
\begin{equation} \label{eq: equivalence of categories finite etale}
\begin{tikzcd}
\left(\mathbf{A}^{\dagger}_{\Lambda}\right)_{\mathrm{f\acute{e}t}} \arrow[rr, "\cong"] \arrow[d, "\cong"'] && \left(\widetilde{\mathbf{A}}^{\dagger}_X \right)_{\mathrm{f\acute{e}t}} \arrow[d, "\cong"] &\\
\left(\fA_{\Lambda}\right)_{\mathrm{f\acute{e}t}} \arrow[d, "\cong"'] \arrow[rr, "\cong"]&& \left(\widetilde{\mathbf{A}}_{X}\right)_{\mathrm{f\acute{e}t}} \arrow[d, "\cong"] &\\
\left(R_{\Lambda}\right)_{\mathrm{f\acute{e}t}} \arrow[r, "\cong"] & \left(R_{\Lambda}^{\mathrm{perf}}\right)_{\mathrm{f\acute{e}t}} \arrow[r, "\cong"]& \left(\widetilde{R}_{X}\right)_{\mathrm{f\acute{e}t}} \arrow[r, "\cong"] & \left(\widehat{A}_{X, \infty}\right)_{\mathrm{f\acute{e}t}}
\end{tikzcd}    
\end{equation}
\end{proposition}

\begin{proof}
Since $R_{\Lambda}$ is a uniform Banach algebra with respect to the norm $|\cdot|$ and $R_{\Lambda}^{\mathrm{perf}}$ is dense in $\widetilde{R}_{X}$, the equivalences on the bottom row follow from \cite[Theorem ~3.1.15]{Kedlaya-LiuI}. Note that  $(\fA_{\Lambda}, (p))$ is a henselian pair since $\fA_{\Lambda}$ is $p$-adically complete. So $\left(\fA_{\Lambda}\right)_{\mathrm{f\acute{e}t}}\cong \left(R_{\Lambda}\right)_{\mathrm{f\acute{e}t}}$.

Since the $p$-adic completion of $\fA_{\Lambda}^{\dagger}$ is equal to $\fA_{\Lambda}$, it suffices to show $(\mathbf{A}^{\dagger}_{\Lambda}, (p))$ is henselian. This follows from the same argument as in the proof of Proposition~\ref{prop: perfect theory}(1), by considering the norm $\mu_r$ on $\fA_{\Lambda}^r$ for any $r \in (0, r_*]$ given by the maximum of $\lambda_r$ and the $p$-adic norm (from that on $\widetilde{\fA}_X$).
\end{proof}

\begin{construction} \label{construction:wdc_finite_etale}
Let $\Lambda = (R_{\Lambda}, \fA_{\Lambda})$ be a weakly decompleting datum for $(X,h)$. The finite \'etale morphism $g\colon X'\rightarrow X$ induces a finite \'etale ring map $\widehat{A}_{X, \infty} \rightarrow \widehat{A}_{ X',\infty}$. Consider the diagram~\eqref{eq: equivalence of categories finite etale} and view $\widehat{A}_{X',\infty}$ as an object in $\left(\widehat{A}_{X, \infty}\right)_{\mathrm{f\acute{e}t}}$. Let ${\fA'}^{\dagger}$, $\fA'$, and $R'$ be the corresponding objects in $\left(\mathbf{A}_{\Lambda}^{\dagger}\right)_{\mathrm{f\acute{e}t}}$, $\left(\fA_{\Lambda}\right)_{\mathrm{f\acute{e}t}}$, and $\left(R_{\Lambda}\right)_{\mathrm{f\acute{e}t}}$ respectively. Note that the corresponding objects in the finite \'etale sites $\left(\widetilde{\mathbf{A}}^{\dagger}_X \right)_{\mathrm{f\acute{e}t}}$, $\left(\widetilde{\mathbf{A}}_X \right)_{\mathrm{f\acute{e}t}}$, and $\left(\widetilde{R}_X \right)_{\mathrm{f\acute{e}t}}$ are precisely $\widetilde{\mathbf{A}}^{\dagger}_{X'}$, $\widetilde{\mathbf{A}}_{X'}$ and $\widetilde{R}_{X'}$ respectively. Furthermore, ${\fA'}^{\dagger}$, $\fA'$, and $R'$ can respectively be identified as subrings of $
\widetilde{\mathbf{A}}^{\dagger}_{X'}$, $\widetilde{\mathbf{A}}_{X'},$ and $ \widetilde{R}_{X'}$. Since we have  
\[
\left(\mathbf{A}_{\Lambda}^{\dagger}\right)_{\mathrm{f\acute{e}t}}\cong \varinjlim_{r>0} \left(\mathbf{A}_{\Lambda}^r\right)_{\mathrm{f\acute{e}t}} \text{ and } \left(\widetilde{\mathbf{A}}^{\dagger}_{X}\right)_{\mathrm{f\acute{e}t}}\cong \varinjlim_{r>0}\left(\widetilde{\mathbf{A}}^r_{X}\right)_{\mathrm{f\acute{e}t}},
\] 
there are ${\fA'}^r\in \left(\mathbf{A}_{\Lambda}^r\right)_{\mathrm{f\acute{e}t}}$ for all $r > 0$ sufficiently small such that
\begin{itemize}
\item ${\fA'}^r$ is a subring of $\widetilde{\mathbf{A}}^r_{X'}$;
\item ${\fA'}^r\subset {\fA'}^s$ for $0 < s <r$, and ${\fA'}^s \cong {\fA'}^r\otimes_{\mathbf{A}_{\Lambda}^r} \fA_{\Lambda}^s$;
\item ${\fA'}^{\dagger}\cong {\fA'}^r\otimes_{\mathbf{A}_{\Lambda}^r}\mathbf{A}^{\dagger}_{\Lambda}$;
\item $\widetilde{\mathbf{A}}^r_{X'}\cong {\fA'}^r\otimes_{\mathbf{A}_{\Lambda}^r}\widetilde{\mathbf{A}}^r_{X}$;
\item ${\fA'}^{\dagger} = \bigcup_{r>0} {\fA'}^r$;
\item $\fA' = {\fA'}^r\otimes_{\mathbf{A}_{\Lambda}^r} \fA_{\Lambda}$.
\end{itemize}
We deduce that ${\fA'}^r = \widetilde{\mathbf{A}}^r_{X'}\cap \fA'$; indeed, since $\mathbf{A}_{\Lambda}^r\rightarrow {\fA'}^r$ is flat, we have
\[
\widetilde{\mathbf{A}}^r_{X'}\cap \fA'={\fA'}^r\otimes_{\mathbf{A}_{\Lambda}^r}\left(\widetilde{\mathbf{A}}^r_{X}\cap \fA_{\Lambda}\right)={\fA'}^r\otimes_{\mathbf{A}_{\Lambda}^r} \mathbf{A}_{\Lambda}^r = {\fA'}^r.
\]
For each $0 < s \le r$, write ${\fC'}^{[s, r]}$ for the completion of ${\fB'}^r\coloneqq {\fA'}^r[p^{-1}]$ with respect to the norm $\max\{\lambda_s, \lambda_r\}$, and let ${\fC'}^r$ be the Fr\'echet completion of ${\fB'}^r$ with respect to the family of norms $\{\lambda_s\}_{s\in (0,r]}$.
\end{construction}

The main result of this subsection is the following.

\begin{theorem}\label{thm: weakly decompleting data finite etale}
The pair $(R', \fA')$ is a weakly decompleting datum for $(X',h')$.
\end{theorem}

\begin{proof} 
We need to verify
\begin{enumerate}
\item[(a)] $R'$ is complete with respect to the norm $|\cdot|$ on $\widetilde{R}_{X'}$;

\item[(b)] $\fA'$ is $\varphi$-stable and $p$-adically complete, and the image of the composite
\[
\fA'/p \rightarrow \widetilde{\fA}_{X'}/p \cong \widetilde{R}_{X'}
\]
lies in $R'$,
\end{enumerate}
and then verify the axioms
\begin{enumerate}
\item $R'$ admits a finite $p$-basis;

\item $\left(R'\right)^{\mathrm{perf}}$ is dense in $\widetilde{R}_{X'}$;

\item $\fA'$ is stable under the $\Gamma$-action, and is equal to the $p$-adic completion of ${\fA'}^{\dagger} \coloneqq \bigcup_{r>0}{\fA'}^r$;

\item for each $r>0$ sufficiently small, ${\fA'}^r$ is complete with respect to the $\lambda_r$-norm on $\widetilde{\mathbf{A}}_{X'}$;

\item the composite ${\fA'}^r/p \rightarrow \fA'/p \rightarrow R'$ is an isomorphism for any sufficiently small $r > 0$;

\item the mod $p$ reduction ${\fA'}^r/p \rightarrow R'$ is strict for any sufficiently small $r > 0$;

\item $R'$ is an affinoid $\kappa(\!(\overline{\pi}_K)\!)$-algebra;

\item for each sufficiently small $r > 0$, ${\fC'}^{[s, r]}$ is a reduced affinoid $K_0'$-algebra for any $0 < s \le r$.
\end{enumerate}

We first check the conditions (a) and (b).

\begin{enumerate}
\item[(a)] Recall that $R_{\Lambda}$ is a Banach ring equipped with the norm $|\cdot|$ induced from $\widetilde{R}_{X}$. Consider $R'$ as a finite projective $R_{\Lambda}$-module, and choose a presentation $R'\oplus Q\cong R_{\Lambda}^{\oplus n}$; this defines a complete norm by Lemma~\ref{lem: finite projective norm}. It induces a presentation 
\[
\widetilde{R}_{X'}\oplus \widetilde{Q}\cong \left(\widetilde{R}_{X}\right)^{\oplus n}
\]
where $\widetilde{Q} = Q\otimes_{R_{\Lambda}} \widetilde{R}_{X}$, and so induces a complete norm on $\widetilde{R}_{X'}$ extending the one on $R'$. By the construction (based on Theorem~\ref{thm: open mapping thm}), this norm is equivalent to $|\cdot|$. Consequently, $R'$ is complete with respect to $|\cdot|$.

\item [(b)] By the equivalences of categories \eqref{eq: equivalence of categories finite etale}, the homomorphism $\varphi\colon R'\rightarrow R'$ uniquely lifts to $\varphi\colon \fA'\rightarrow \fA'$ which coincides with the Frobenius $\varphi$ on $\widetilde{\mathbf{A}}_{X'}$. Furthermore, applying Lemma~\ref{lem: finite projective norm} to the $p$-adic norm on $\fA_{\Lambda}$, we deduce that it extends to $\fA'$ (which is finite projective over $\fA_{\Lambda}$) and is complete. Again by the equivalences \eqref{eq: equivalence of categories finite etale}, the composite $\fA'/p \rightarrow \widetilde{\fA}_{X'}/p \cong \widetilde{R}_{X'}$ induces an isomorphism $\fA'/p \cong R'$. 
\end{enumerate}

Now we check the conditions (1)-(8).

\begin{enumerate}
\item By \cite[\S 1.1.2(ii)]{Berthelot-Messing}, any $p$-basis for $R_{\Lambda}$ is also a $p$-basis for $R'$.

\item From the above, 
\[
\begin{tikzcd}
R' \arrow[r, "\varphi"] & R'\\
R_{\Lambda} \arrow[u] \arrow[r, "\varphi"] & R_{\Lambda} \arrow[u]
\end{tikzcd}
\]
is a pushout diagram, and so
\[
\left(R'\right)^{\mathrm{perf}}\cong R_{\Lambda}^{\mathrm{perf}}\otimes_{R_{\Lambda}} R'.
\]
Thus, $\left(R'\right)^{\mathrm{perf}}$ is dense in $\widetilde{R}_{X'}$ as $R^{\mathrm{perf}}$ is dense in $\widetilde{R}_{X}$.

\item This follows from the equivalences \eqref{eq: equivalence of categories finite etale}; in fact, the stability under the $\Gamma$-action can be argued in the same way as for $\varphi$. Since the $p$-adic completion of $\fA_{\Lambda}^{\dagger}$ is equal to $\fA_{\Lambda}$, we have (for $r > 0$ sufficiently small)
\[
\varprojlim_{n}\left({\fA'}^{\dagger}/p^n\right) = \varprojlim_{n}\left(\mathbf{A}^{\dagger}_{\Lambda}/p^n\otimes_{\mathbf{A}_{\Lambda}^r}{\fA'}^r\right)\cong \left(\varprojlim_{n}\left(\mathbf{A}^{\dagger}_{\Lambda}/p^n\right)\right)\otimes_{\mathbf{A}_{\Lambda}^r}{\fA'}^r = \fA_{\Lambda}\otimes_{\mathbf{A}_{\Lambda}^r}{\fA'}^r = \fA'.
\]

\item Equip $\mathbf{A}_{\Lambda}^r$ with the $\lambda_r$-norm, and consider ${\fA'}^r$ as a finite projective $\mathbf{A}^r$-module. Choose a presentation ${\fA'}^r\oplus Q\cong (\mathbf{A}_{\Lambda}^r)^{\oplus n}$. Applying Lemma~\ref{lem: finite module seminorm}, we obtain a complete norm on ${\fA'}^r$ equivalent to the finite-module norm. We have an induced presentation 
\[
\widetilde{\mathbf{A}}^r_{X'}\oplus \widetilde{Q}\cong \left(\widetilde{\mathbf{A}}_{X}\right)^{\oplus n},
\] 
which yields a complete norm on $\widetilde{\mathbf{A}}^r_{X'}$ (equivalent to the finite-module norm). We need to show that the finite-module norm on $\widetilde{\mathbf{A}}^r_{X'}$ is equivalent to the $\lambda_r$-norm, for any sufficiently small $r>0$.

Let $x_1, \ldots, x_m$ be elements in $\widetilde{\mathbf{A}}_{X'}$ whose mod $p$ reductions $\overline{x}_1, \ldots, \overline{x}_m\in \widetilde{R}_{X'}$ generate $\widetilde{R}_{X'}$ as an $\widetilde{R}_{X}$-module. By \cite[Lemma ~5.5.2]{Kedlaya-LiuI} and its proof, for each sufficiently small $r > 0$, we have that 
\begin{itemize}
\item $x_1, \ldots, x_m$ generate $\widetilde{\mathbf{A}}_{X'}^r$ as a $\widetilde{\mathbf{A}}_X^r$-module;
\item the surjective map
\[
\left(\widetilde{\mathbf{A}}_X^r\right)^{\oplus m} \twoheadrightarrow \widetilde{\mathbf{A}}_{X'}^r, ~~(a_i)_{1 \le i \le m} \mapsto \sum_{i=1}^m a_i x_i,
\]
is strict with respect to $\max_{1\le i\le m}\{\lambda_r(a_i)\}$ on the source and $\lambda_r$-norm on the target. 
\end{itemize}
Thus, the two norms are equivalent.

\item This follows again from the equivalences \eqref{eq: equivalence of categories finite etale}.

\item Choose a presentation ${\fA'}^r\oplus Q\cong (\mathbf{A}_{\Lambda}^r)^{\oplus n}$, which induces a complete norm on ${\fA'}^r$ via Lemma~\ref{lem: finite projective norm}. This is equivalent to the $\lambda_r$-norm by above. Reduction modulo $p$ gives a presentation $R'\oplus \overline{Q}\cong R^{\oplus n}$, which induces a complete norm on $R'$ via Lemma~\ref{lem: finite projective norm} (equivalent to the norm $|\cdot|$). Thus, the desired strictness follows from the strictness of $(\mathbf{A}_{\Lambda}^r)^{\oplus n}\rightarrow R_{\Lambda}^{\oplus n}$.

\item Since $R'$ is finite \'etale over $R_{\Lambda}$, it is an affinoid $\kappa(\!(\overline{\pi}_K)\!)$-algebra.

\item It suffices to show that the map
\[
f\colon \fC_{\Lambda}^r\otimes_{\fB_{\Lambda}^r} {\fB'}^r \rightarrow {\fC'}^r
\]
is an isometric isomorphism with respect to $\lambda_r$-norm, for any sufficiently small $r > 0$. We follow a similar argument as in the proof of \cite[Proposition ~5.5.4]{Kedlaya-LiuI}. First, the map $f$ is a strict surjection for any $r > 0$ sufficiently small by a similar computation as in the proof of \cite[Lemma ~5.5.2]{Kedlaya-LiuI}, noticing that ${\fB'}^r$ is finite projective over $\fB_{\Lambda}^r$. In particular, there is a constant $c \ge 1$ (which may depend on $r$) such that each $z \in {\fC'}^r$ can be lifted to some $\sum_i x_i\otimes y_i \in \fC_{\Lambda}^r\otimes_{\fB_{\Lambda}^r} {\fB'}^r$ with $\max_i\{\lambda_r(x_i y_i)\} \le c\lambda_r(z)$. Note that the map $f$ is submetric. Given $z \in {\fC'}^r$, let $z_0 \in {\fB'}^r$ with $\lambda_r(z-z_0)\le \frac{1}{2}c^{-1} \lambda_r(z)$. In particular, $\lambda_r(z_0) = \lambda_r(z)$. Choose a lift $\sum_i x_i\otimes y_i \in \fC_{\Lambda}^r\otimes_{\fB_{\Lambda}^r} {\fB'}^r$ of $z-z_0$ such that $\max_i\{\lambda_r(x_i y_i)\} \le c\lambda_r(z-z_0)$. Then $\max_i\{\lambda_r(x_i y_i)\} \le \frac{1}{2}\lambda_r(z)$. So $1\otimes z_0+\sum_i x_i\otimes y_i$ is a lift of $z$ satisfying
\[
\max\{\lambda_r(z_0), \max_i\{\lambda_r(x_i y_i)\}\} = \lambda_r(z_0) = \lambda_r(z).
\]
Thus, the map $f$ is isometric. In particular, $f$ is an isomorphism.
\end{enumerate}
\end{proof}

\begin{corollary}
Let $\Lambda=(R_{\Lambda}, \mathbf{A}_{\Lambda})$ and $(R', \mathbf{A}')$ be as above. Then there are canonical $(\varphi, \Gamma)$-equivariant isomorphisms
\[{\mathbf{C}'}^r\cong \mathbf{C}^r_{\Lambda}\otimes_{\mathbf{A}^r_{\Lambda}} {\mathbf{A}'}^r\cong \mathbf{C}^r_{\Lambda}\otimes_{\mathbf{B}^r_{\Lambda}} {\mathbf{B}'}^r\]
(for $r$ sufficiently small) and 
\[{\mathbf{C}'}^{\dagger}\cong \mathbf{C}^{\dagger}_{\Lambda}\otimes_{\mathbf{A}^{\dagger}_{\Lambda}} {\mathbf{A}'}^{\dagger}\cong \mathbf{C}^{\dagger}_{\Lambda}\otimes_{\mathbf{B}^{\dagger}_{\Lambda}} {\mathbf{B}'}^{\dagger}\]
\end{corollary}

\begin{proof}
This is discussed in the proof of Theorem~\ref{thm: weakly decompleting data finite etale} when we check condition (8).
\end{proof}

\vspace{0.1in}
\section{Decompleting data}  \label{sec:decompleting_data}

In this section, we introduce an upgraded version of weakly decompleting data called \textit{decompleting data}. As the terminology suggests, this notion will be used in \S \ref{sec: decompletion-relative-etale-phiGamma-modules} when we consider the decompletion of certain relative $(\varphi, \Gamma)$-modules. The goals of this section are the following: 
\begin{itemize}
\item Define the notion of decompleting data (in terms of strict-exactness of certain group cochain complexes in characteristic $p$).

\item Study fundamental properties of decompleting data. In particular, we will derive the strict-exactness of certain group cochain complexes valued in period rings of type $\mathbf{A}$ and $\mathbf{C}$. 
\end{itemize}
The main results are Corollary~\ref{cor: complex-A-strict-exact} and Proposition~\ref{prop: C-complex-strict-exact}.

\subsection{Definition of decompleting data} \label{subsec: definition-decompleting-data}

Let $h\colon X \ra \D_K^d$ be a toric chart and $\psi = \{\psi_n\colon X_n \ra X\}$ be the toric tower relative to $h$ as in Definition~\ref{defn: toric chart}. Let $\Lambda = (R_{\Lambda}, \mathbf{A}_{\Lambda})$ be a weakly decompleting data for $(X, h)$ as in Definition~\ref{defn: weakly decompleting datum}. Recall the stability of imperfect period rings under the $\Gamma$-action in Remark~\ref{rem: stability-under-Gamma-actions}.

\begin{notation}
For any topological ring $S$ whose topology is induced by a (semi)norm (and more generally, a topological $S$-module $M$), which is equipped with a continuous $\Gamma$-action, we let 
\[
\mathcal{C}^{(\bullet)}(\Gamma, S), \qquad (\textup{resp. }  \mathcal{C}^{(\bullet)}(\Gamma, M)) 
\]
denote the complex of inhomogeneous continuous $\Gamma$-cochains with values in $S$ (resp. $M$), equipped with the supremum norm induced from the norm on $S$ (resp. on $M$). Namely, for $f \in \mathcal{C}^{(n)}(\Gamma, S)$ (resp. $f \in \mathcal{C}^{(n)}(\Gamma, M)$), define
\[
|f|\coloneqq \sup_{(\gamma_1, \ldots, \gamma_n) \in \Gamma^n} |f(\gamma_1, \ldots, \gamma_n)|.
\]
Note that since $\Gamma$ is compact, the supremum norm is well-defined. For any open subgroup $\Gamma_0 \subset \Gamma$, we similarly define $\mathcal{C}^{(\bullet)}(\Gamma_0, S)$ and $\mathcal{C}^{(\bullet)}(\Gamma_0, M)$ (note that since $\Gamma$ is compact, $\Gamma_0$ is necessarily of finite index, so it is closed in $\Gamma$). 
\end{notation}

\begin{notation}
Let $s, r \in \Q$ be rational numbers such that $0 < s \le r \le r_*$ and let $m'> m\ge 0$ be integers. Consider the Banach rings $\widetilde{R}_X$, $\widetilde{\mathbf{A}}^r_X$, and $\widetilde{\mathbf{C}}^{[s,r]}_X$, equipped with norms $|\cdot|$, $\lambda_r$, and $\max\{\lambda_r, \lambda_s\}$, respectively, as well as their closed subrings $R_{\Lambda, (m)}$, $\mathbf{A}^r_{\Lambda, (m)}$, and $\mathbf{C}^{[s,r]}_{\Lambda, (m)}$. Then we consider the quotients $\widetilde{R}_X/R_{\Lambda, (m)}$, $\widetilde{\mathbf{A}}^{r}_{X}/ \mathbf{A}^r_{\Lambda, (m)}$, and $\widetilde{\mathbf{C}}^{[s, r]}_{X}/ \fC^{[s, r]}_{\Lambda, (m)}$, equipped with the quotient (semi)norms. Similarly, we consider $R_{\Lambda, (m')}/R_{\Lambda, (m)}$, $\mathbf{A}^r_{\Lambda, (m')}/\mathbf{A}^r_{\Lambda, (m)}$, and $\fC^{[s, r]}_{\Lambda, (m')}/\fC^{[s, r]}_{\Lambda, (m)}$. Since all quotients above are taken with respect to closed subspaces, the quotient seminorms are actually complete norms. Also note that the $\Gamma$-actions are isometric on these quotient spaces. 

For any open subgroup $\Gamma_0\subset \Gamma$, we consider the cochain complexes
\[  
\mathcal{C}^{(\bullet)}(\Gamma_0, \widetilde{R}_X/R_{\Lambda, (m)}) , 
\qquad \mathcal{C}^{(\bullet)}(\Gamma_0, \widetilde{\mathbf{A}}^{r}_{X}/ \mathbf{A}^r_{\Lambda, (m)}), \qquad
\mathcal{C}^{(\bullet)}(\Gamma_0, \widetilde{\mathbf{C}}^{[s, r]}_{X}/ \fC^{[s, r]}_{\Lambda, (m)})
\]  
and
\[  
\mathcal{C}^{(\bullet)}(\Gamma_0, R_{\Lambda, (m')}/R_{\Lambda, (m)}) , 
\qquad \mathcal{C}^{(\bullet)}(\Gamma_0, \mathbf{A}^r_{\Lambda, (m')}/ \mathbf{A}^r_{\Lambda, (m)}), \qquad
\mathcal{C}^{(\bullet)}(\Gamma_0, \fC^{[s, r]}_{\Lambda, (m')}/ \fC^{[s, r]}_{\Lambda, (m)}).
\] 
It is clear that differentials of these cochain complexes are submetric. We say that such a complex is \emph{strict exact} at $d^{(i)}$ for the differential map $d^{(i)}\colon \mathcal{C}^{(i)}  \ra \mathcal{C}^{(i+1)}$, if $\ker(d^{(i)}) = \im(d^{(i-1)})$ and $d^{(i)}$ is strict in the sense of Definition \ref{defn: equivalent seminorms}, namely, the quotient (semi)norm on $\mathcal{C}^{i} /\ker(d^{(i)}) \cong \im (d^{(i)})$ is equivalent to the subspace norm on $\im(d^{(i)}) \subset \mathcal{C}^{(i+1)}$. For example, if we denote the first complex above as $\mathcal{C}^{(\bullet)}$, then strict exactness at $d^{(i)}$ means that there exist constants $\mathfrak{c}, \mathfrak{c}'>0$ such that for any $g \in \mathcal{C}^{(i)}$, we have $|d^{(i)} g| \le \mathfrak{c} \cdot |g|$, and there exists $f' \in \mathcal C^{(i-1)}$ satisfying
\[ 
|g + d^{(i-1)} f'| \le \mathfrak{c}' \cdot |d^{(i)} g|.
\]
\end{notation}

\begin{definition} \label{definition:decompleting_datum}
A weakly decompleting datum  $\Lambda=(R_{\Lambda}, \mathbf{A}_{\Lambda})$ for $(X,h)$ is called \emph{decompleting} if the following conditions hold.
\begin{enumerate}
\item For any open subgroup $\Gamma_0 \subset \Gamma$, the cochain complex
\[
\mathcal{C}^{(\bullet)}(\Gamma_0, \varphi^{-1}(R_{\Lambda}) / R_{\Lambda})
\]
is strict exact at $d^{(0)}$ and $d^{(1)}$.

\item For any finite \'etale map   $g\colon X'\rightarrow X$, Condition (1) holds with $R_{\Lambda}$ replaced by $R'$, where $R'\in \left(R_{\Lambda}\right)_{\mathrm{f\acute{e}t}}$ is given by Construction~\ref{construction:wdc_finite_etale}. 
\end{enumerate}
\end{definition}

\begin{remark} \label{remark:decompleting_finite_etale} In fact, Condition (2) is not needed for the results in \S \ref{sec:decompleting_data}. It is included to ensure that the notion of decompleting data is stable under finite \'etale extensions (cf. Construction~\ref{construction:wdc_finite_etale} and Theorem~\ref{thm: weakly decompleting data finite etale}). This is used in the proof of Theorem~\ref{thm: decompletion-etale-(phi, Gamma)-modules} (Step 2).
\end{remark}

\subsection{Strict-exactness in characteristic $p$} \label{subsec: strict-exactness-char-p}

In the remainder of this section, we assume $\Lambda = (R_{\Lambda}, \mathbf{A}_{\Lambda})$ is a decompleting data for $(X, h)$. 

\begin{lemma} \label{lem: complex-R-strict-exact}  The cochain complex
\[
\mathcal{C}^{(\bullet)}(\Gamma_0, \widetilde{R}_X  / R_{\Lambda, (m_0)})
\]
is strict exact at $d^{(0)}$ and $d^{(1)}$ for every open subgroup $\Gamma_0 \subset \Gamma$ and every $m_0 \ge 0$.
\end{lemma}

\begin{proof}
For simplicity, we consider the case $m_0 = 0$; the statement for general $m_0$ follows from a similar argument. Fix an open subgroup $\Gamma_0 \subset \Gamma$. \\

\noindent \textbf{Step 1.} We first claim that for each $j \ge 0$, the cochain complex 
\[
\mathcal{C}^{(\bullet)}(\Gamma_0, R_{\Lambda, (j+1)} / R_{\Lambda, (j)})
\] 
is strict exact at $d^{(0)}$ and $d^{(1)}$.
For this, note that for each $j \ge 0$, $\varphi^j$ on $\widetilde{R}_X$ induces a $\Gamma_0$-equivariant isomorphism
$
\varphi^j\colon R_{\Lambda, (j+1)} / R_{\Lambda, (j)} \stackrel{\sim}{\longrightarrow} R_{\Lambda, (1)} / R_{\Lambda}. 
$ 
Since $\Lambda$ is decompleting and $\varphi^j$ is compatible with the differential maps, the claim holds.\\

\noindent \textbf{Step 2.} Next, we show that the cochain complex in Step 1 is uniformly strict (with respect to $j$) at $d^{(0)}$ and $d^{(1)}$.
We only prove this for $d^{(1)}$, as that for $d^{(0)}$ follows from a similar argument. Note first that the differential maps are submetric. Let $C \ge 1$ be a constant such that 
$
|x|_{\mathrm{quot}} \le C|d^{(1)}(x)|
$ 
for any $x \in \mathcal{C}^{(1)}(\Gamma_0, R_{\Lambda, (1)} / R_{\Lambda})$, where $|\cdot|_{\mathrm{quot}}$ denotes the quotient norm on 
\[
\mathcal{C}^{(1)}(\Gamma_0, R_{\Lambda, (1)} / R_{\Lambda}) / \ker(d^{(1)}).
\]
Let $y \in \mathcal{C}^{(1)}(\Gamma_0, R_{\Lambda, (j+1)} / R_{\Lambda, (j)})$. Consider $x = \varphi^j (y) \in \mathcal{C}^{(1)}(\Gamma_0, R_{\Lambda, (1)} / R_{\Lambda})$. Then $|x+z| \le C|d^{(1)}(x)|$ for some $z \in \ker(d^{(1)}) \in \mathcal{C}^{(1)}(\Gamma_0, R_{\Lambda, (1)} / R_{\Lambda})$. Since the norms are power-multiplicative and differential maps are $p$-power multiplicative, we have
\[
|y+\varphi^{-j}(z)| \le C^{\frac{1}{p^j}}|d^{(1)}(y)|
\]
with $\varphi^{-j}(z) \in \ker(d^{(1)}) \subset \mathcal{C}^{(1)}(\Gamma_0, R_{\Lambda, (j+1)} / R_{\Lambda, (j)})$. This shows the uniform strictness at $d^{(1)}$.\\ 

\noindent \textbf{Step 3.} We claim that the complex
\[
\mathcal{C}^{(\bullet)}(\Gamma_0, R_{\Lambda, (m+1)} / R_{\Lambda})
\]
is exact and uniformly strict (with respect to $m$) at $d^{(0)}$ and $d^{(1)}$. We prove this for $d^{(1)}$; that for $d^{(0)}$ follows from a similar (and simpler) argument.

For each $j \ge 0$, consider the $\Gamma$-equivariant exact sequence
\begin{equation} \label{eq: filtration-Gamma-equiv}
0 \rightarrow R_{\Lambda, (j)}/R_{\Lambda} \rightarrow R_{\Lambda, (j+1)}/R_{\Lambda} \rightarrow R_{\Lambda, (j+1)}/R_{\Lambda, (j)} \rightarrow 0.    
\end{equation}
By the proof of Lemma~\ref{lem: Frob decomp R uniform strict}, there exists a constant $D \ge 1$ such that for each $j \ge 0$ and $r \in R_{\Lambda, (j+1)}/R_{\Lambda, (j)}$, the section $r' \in R_{\Lambda, (j+1)}/R_{\Lambda}$ given by the decomposition in Lemma~\ref{lem: Frob decomp R uniform strict} satisfies $|r'| \le D^{p^{-j}}|r|$.

For each $j \ge 0$, we will inductively construct a constant $B_j \ge 1$ such that for any $x \in \mathcal{C}^{(1)}(\Gamma_0, R_{\Lambda, (j+1)}/R_{\Lambda})$, we have $|x|_{\mathrm{quot}} \le B_j|d^{(1)}(x)|$. For $j = 0$, we take $B_0 = C$ as in Step 2. Suppose $B_{j-1}$ is given, and let $x \in \mathcal{C}^{(1)}(\Gamma_0, R_{\Lambda, (j+1)}/R_{\Lambda})$. Write $\overline{x} \in \mathcal{C}^{(1)}(\Gamma_0, R_{\Lambda, (j+1)}/R_{\Lambda, (j)})$ given by the pointwise projection of $x$ under the projection in \eqref{eq: filtration-Gamma-equiv}. By Step 2, there exists $\overline{z} \in \mathcal{C}^{(0)}(\Gamma_0, R_{\Lambda, (j+1)}/R_{\Lambda, (j)})$ such that $|\overline{x}+d^{(0)}(\overline{z})| \le C^{p^{-j}}|d^{(1)}(x)|$. Write $z \in \mathcal{C}^{(0)}(\Gamma_0, R_{\Lambda, (j+1)}/R_{\Lambda})$ for the $0$-cochain given by the pointwise section of $\overline{z}$ (via Lemma~\ref{lem: Frob decomp R uniform strict}). 

Let $x_1 = x+d^{(0)}(z) \in \mathcal{C}^{(1)}(\Gamma_0, R_{\Lambda, (j+1)}/R_{\Lambda})$, and consider its projection $\overline{x}_1 \in \mathcal{C}^{(1)}(\Gamma_0, R_{\Lambda, (j+1)}/R_{\Lambda, (j)})$. Note that $d^{(1)}(x) = d^{(1)}(x_1)$. Write $w \in \mathcal{C}^{(1)}(\Gamma_0, R_{\Lambda, (j+1)}/R_{\Lambda})$ for the $1$-cochain given by the pointwise section of $\overline{x}_1$. Then we have 
\[
|w| \le (CD)^{p^{-j}}|d^{(1)}(x)| \quad \text{and} \quad x_1-w \in \mathcal{C}^{(1)}(\Gamma_0, R_{\Lambda, (j)}/R_{\Lambda}).
\]
Thus, $|d^{(1)}(x_1-w)| \le (CD)^{p^{-j}}|d^{(1)}(x)|$. By the induction hypothesis, there is $t \in \mathcal{C}^{(0)}(\Gamma_0, R_{\Lambda, (j)}/R_{\Lambda})$ such that 
\[
|x_1-w+d^{(0)}(t)| \le B_{j-1}(CD)^{p^{-j}}|d^{(1)}(x)|.
\]
Hence, we deduce
\[
|x+d^{(0)}(z+t)| \le B_{j-1}(CD)^{p^{-j}}|d^{(1)}(x)|.
\]
So we may take $B_j = B_{j-1}(CD)^{p^{-j}}$.

Since $B_j \le B_0(CD)^{\frac{p}{p-1}}$ for all $j \ge 0$, this proves the claim.\\

Write $C_1 \ge 1$ for a uniform bound for $\mathcal{C}^{(\bullet)}(\Gamma_0, R_{\Lambda, (m+1)} / R_{\Lambda})$ at $d^{(j)}$ for $j = 0, 1$; namely, $|y|_{\mathrm{quot}} \le C_1|d^{(j)}(y)|$ for all $m\ge 0$, $j=0,1$, and $y \in \mathcal{C}^{(j)}(\Gamma_0, R_{\Lambda, (m+1)} / R_{\Lambda})$.
Now, we consider the cochain complex
\begin{equation} \label{eq: complex-tilde R-over-R}
\mathcal{C}^{(\bullet)}(\Gamma_0, \widetilde{R}_X / R_{\Lambda}).    
\end{equation}

\noindent \textbf{Step 4.} We show that the complex \eqref{eq: complex-tilde R-over-R} is exact at $d^{(0)}$ and $d^{(1)}$.\\ 

We first explain the proof for $d^{(1)}$. It suffices to show that for any $y \in \ker(d^{(1)})$ and $\epsilon > 0$, there exists $z \in \mathcal{C}^{(0)}(\Gamma_0, R_{\Lambda, (m+1)} / R_{\Lambda})$ for some $m$ (which may depend on $y$ and $\epsilon$) such that $|y-d^{(0)}(z)| < \epsilon$. In fact, since $\mathcal{C}^{(\bullet)}(\Gamma_0, R_{\Lambda, (m+1)} / R_{\Lambda})$ is uniformly strict at $d^{(0)}$, this would allow us to inductively construct a sequence of elements $z_j \in \mathcal{C}^{(0)}(\Gamma_0, \widetilde{R}_X / R_{\Lambda})$ for $j \ge 0$ with $\lim_{j \rightarrow \infty} |z_j| = 0$ such that $y = \sum_{j=0}^{\infty} d^{(0)}(z_j)$. 

Since ${R}^{\mathrm{perf}}_{\Lambda}$ is dense inside $\widetilde{R}_X$, we can find $y_1 \in \mathcal{C}^{(1)}(\Gamma_0, R_{\Lambda, (m+1)} / R_{\Lambda})$ for some sufficiently large $m$ such that $|y-y_1| < C_1^{-1} \epsilon$. Since differential maps are submetric, we have
$
|d^{(1)}(y-y_1)| = |d^{(1)}(y_1)| < C_1^{-1} \epsilon.
$ 
By Step 3, there exists $z \in \mathcal{C}^{(0)}(\Gamma_0, R_{\Lambda, (m+1)} / R_{\Lambda})$ such that
\[
|y_1-d^{(0)}(z)| \le C_1|d^{(1)}(y_1)| < C_1\cdot C_1^{-1} \epsilon = \epsilon.
\]
Thus, we have
\[
|y-d^{(0)}(z)| \le \max\{|y-y_1|, |y_1-d^{(0)}(z)|\} < \epsilon.
\]
This shows the exactness at $d^{(1)}$.

The exactness at $d^{(0)}$ (namely, injectivity of $d^{(0)}$) follows from a similar argument. Indeed, one can take $z=0$ in this case; hence, for any $y\in \ker (d^{(0)})$, we deduce that $|y|<\epsilon$ for all $\epsilon>0$, namely, $y=0$.\\

\noindent \textbf{Step 5.} It remains to show that the complex \eqref{eq: complex-tilde R-over-R} is strict at $d^{(0)}$ and $d^{(1)}$. We only prove this for $d^{(1)}$, as that for $d^{(0)}$ follows from a similar argument. Let $y \in \mathcal{C}^{(1)}(\Gamma_0, \widetilde{R}_X / R_{\Lambda})$. We may assume that $|d^{(1)}(y)| > 0$. Choose $y_1 \in \mathcal{C}^{(1)}(\Gamma_0, R_{\Lambda, (m+1)} / R_{\Lambda})$ for some $m$ such that $|y-y_1| \le \frac{1}{2} |d^{(1)}(y)|$. Since $|d^{(1)}(y-y_1)| \le |y-y_1|$, we have $|d^{(1)}(y_1)| = |d^{(1)}(y)|$. By Step 3, there exists $z \in \mathcal{C}^{(0)}(\Gamma_0, R_{,\Lambda, (m+1)} / R_{\Lambda})$ such that 
\[
|y_1-d^{(0)}(z)| \le C_1|d^{(1)}(y_1)| = C_1|d^{(1)}(y)|.
\]
Therefore,
\[
|y-d^{(0)}(z)| \le \max\{|y-y_1|, |y_1-d^{(0)}(z)| \} \le C_1|d^{(1)}(y)|,
\]
as desired.
\end{proof}

\subsection{Strict-exactness for type $\mathbf{A}$} \label{subsec: strict-exactness-type-A}

To prove the strict-exactness of cochain complexes of type $\mathbf{A}$ (see Proposition~\ref{prop: strict exactness for A(1)/A} and Corollary~\ref{cor: complex-A-strict-exact} below), we work with the following set-up.

\begin{setup} \label{setup: r0-Gamma0}
Let $r_0 > 0$ be sufficiently small such that Corollaries~\ref{cor: good lift} and \ref{cor: Frob decomp A uniform strict} are applicable for any $0 < r \le r_0$. Let $\Gamma_0 \subset \Gamma$ be an open normal subgroup of the form $\Gamma_0 \cong (p^a\Z_p)^d \rtimes (1+p^b\Z_p) \subset \Z_p^d \rtimes \Z_p^{\times}$ for some integers $a \ge 0$ and $b \ge 1$ (note that normality holds when $b \ge a$). By Lemma~\ref{lem: complex-R-strict-exact}, after further shrinking $r_0$ if necessary, we may further assume that
\[
|z|_{\mathrm{quot}}^r \le \frac{3}{2}|d^{(\ell)}(z)|^r  
\]
for all $z \in \mathcal{C}^{(\ell)}(\Gamma_0, \widetilde{R}_X / R_{\Lambda, (m_0)})$, $\ell = 0, 1$, and $0 < r \le r_0$.
\end{setup}

\begin{lemma} \label{lem: norm-compatible-expansion-for-continuous-functions}
Assume Set-up~\ref{setup: r0-Gamma0}. Let $0 < r < r_0$. For any element $x \in \mathcal{C}^{(1)}(\Gamma_0, \fA^{r}_{\Lambda, (1)}/ \fA^{r}_{\Lambda})$, we can express $x$ as a $\lambda_r$-convergent norm compatible expansion
\begin{equation} \label{eq: norm compatible expansion}
x = \sum_{j=0}^{\infty} p^j x_j 
\end{equation}
with each $x_j \in \mathcal{C}^{(1)}(\Gamma_0, \fA^{r_0}_{\Lambda, (1)}/ \fA^{r_0}_{\Lambda})$ such that $\lambda_s(x_j) = |\overline{x}_j|^s$ and 
$
\lambda_{s}(x) = \max_j\{p^{-j}|\overline{x}_j|^s\}
$ 
for any $0 < s \le r$.    
\end{lemma}

\begin{proof}

Let us start with the following claim.\\

\noindent \textbf{Claim.} For any element $z \in \mathcal{C}^{(1)}(\Gamma_0, R_{\Lambda, (1)}/R_{\Lambda})$, there exists a ``good lift'' $\widetilde{z} \in \mathcal{C}^{(1)}(\Gamma_0, \fA^{r_0}_{\Lambda, (1)}/ \fA^{r_0}_{\Lambda})$ in the sense that for each $\gamma \in \Gamma_0$, $\widetilde{z}(\gamma) \in \fA^{r_0}_{\Lambda, (1)}/ \fA^{r_0}_{\Lambda}$ is a good lift of $z(\gamma) \in R_{\Lambda, (1)}/R_{\Lambda}$ in the sense of Proposition \ref{prop: good lift weakly decompleting} and Corollary \ref{cor: good lift}.\\

\noindent \textit{Proof of Claim.} Note that $z$ is a continuous function on $\Gamma_0$ valued in $R_{\Lambda, (1)}/R_{\Lambda}$. Since $\Gamma_0$ is compact, $z$ is uniformly continuous. For notational simplicity, we assume $d = 1$, as the statement for general $d$ would follow from the same argument.

We represent elements in $\Gamma_0 = p^a\Z_p\rtimes (1+p^b\Z_p)$ by pairs $(\alpha, \beta)$ with $\alpha, \beta \in \Z_p$ so that $(\alpha, \beta)$ corresponds to $(p^a \alpha, 1+p^b \beta) \in p^a\Z_p\rtimes (1+p^b\Z_p)$. For such a pair $(\alpha, \beta)$, we further write $\alpha = \sum_{k=0}^{\infty} [\alpha_k] p^k$ and $\beta = \sum_{k=0}^{\infty} [\beta_k] p^k$ with $\alpha_k, \beta_k \in \F_p$. Pick a sequence of integers $1 \le m_1 < m_2 < m_3 < \cdots$ inductively so that for any $(\alpha, \beta), (\delta, \sigma) \in \Z_p^2$ with $(\alpha, \beta) - (\delta, \sigma) \in (p^{m_j}\Z_p)^2$, we have
\[
|z(\alpha, \beta) - z(\delta, \sigma)| \le \left(\frac{1}{2}\right)^j. 
\]

We shall construct $\widetilde{z}(\alpha, \beta)$ for all $(\alpha, \beta)\in \Z_p^2$. Firstly, pick a good lift (in the sense of Corollary~\ref{cor: good lift}(3)) $\widetilde{z}(\sum_{k = 0}^{m_1-1} [\alpha_k] p^k, \sum_{k = 0}^{m_1-1} [\beta_k] p^k)$ of $z(\sum_{k = 0}^{m_1-1} [\alpha_k] p^k, \sum_{k = 0}^{m_1-1} [\beta_k] p^k)$ for every choice of $\alpha_k, \beta_k \in \F_p$ ($0 \le k \le m_1-1$). Next, we construct good lifts
$\widetilde{z}(\sum_{k = 0}^{m_2-1} [\alpha_k] p^k, \sum_{k = 0}^{m_2-1} [\beta_k] p^k)$ of 
$z(\sum_{k = 0}^{m_2-1} [\alpha_k] p^k, \sum_{k = 0}^{m_2-1} [\beta_k] p^k)$ for every choice of $\alpha_k, \beta_k \in \F_p$ ($0 \le k \le m_2-1$) satisfying
\begin{equation}\label{eq: inductive step}
\lambda_r\left(\widetilde{z}(\sum_{k = 0}^{m_2-1} [\alpha_k] p^k, \sum_{k = 0}^{m_2-1} [\beta_k] p^k) - \widetilde{z}(\sum_{k = 0}^{m_1-1} [\alpha_k] p^k, \sum_{k = 0}^{m_1-1} [\beta_k] p^k)\right) \le \left(\frac{1}{2}\right)^r.
\end{equation}
Indeed, if 
\[
\left|z(\sum_{k = 0}^{m_2-1} [\alpha_k] p^k, \sum_{k = 0}^{m_2-1} [\beta_k] p^k) - z(\sum_{k = 0}^{m_1-1} [\alpha_k] p^k, \sum_{k = 0}^{m_1-1} [\beta_k] p^k)\right|<\left|z(\sum_{k = 0}^{m_1-1} [\alpha_k] p^k, \sum_{k = 0}^{m_1-1} [\beta_k] p^k)\right|,
\]
we can construct the desired good lift using Remark~\ref{rem: variation-of-good-lifts}. Otherwise, we simply pick any good lift and the inequality \eqref{eq: inductive step} automatically holds.
Repeating this process, for all $j\ge 1$, we obtain good lifts $\widetilde{z}(\sum_{k = 0}^{m_j-1} [\alpha_k] p^k, \sum_{k = 0}^{m_j-1} [\beta_k] p^k)$ of $z(\sum_{k = 0}^{m_j-1} [\alpha_k] p^k, \sum_{k = 0}^{m_j-1} [\beta_k] p^k)$ for all choices of $\alpha_k, \beta_k \in \F_p$ ($0 \le k \le m_j-1$) satisfying
\[
\lambda_r\left(\widetilde{z}(\sum_{k = 0}^{m_{j+1}-1} [\alpha_k] p^k, \sum_{k = 0}^{m_{j+1}-1} [\beta_k] p^k) - \widetilde{z}(\sum_{k = 0}^{m_j-1} [\alpha_k] p^k, \sum_{k = 0}^{m_j-1} [\beta_k] p^k)\right) \le \left(\frac{1}{2}\right)^{jr}.
\]
For a general $(\alpha, \beta)\in \Z_p^2$ with $\alpha=\sum_{k = 0}^{\infty} [\alpha_k] p^k$
and $\beta=\sum_{k = 0}^{\infty} [\beta_k] p^k$, we put
\[\widetilde{z}(\alpha, \beta)\coloneqq \lim_{j\rightarrow\infty}\widetilde{z}(\sum_{k = 0}^{m_j-1} [\alpha_k] p^k, \sum_{k = 0}^{m_j-1} [\beta_k] p^k).\]
By construction, $\widetilde{z}$ is continuous. This finishes the proof of the claim.\\

Returning to the proof of the lemma, we employ the same construction as in the proof of Proposition~\ref{prop: good lift weakly decompleting}(2). Let $y_0 = x$, and let $x_0 \in \mathcal{C}^{(1)}(\Gamma_0, \fA^{r_0}_{\Lambda, (1)}/ \fA^{r_0}_{\Lambda})$ be a good lift of $\overline{y}_0$ given by the above claim. Given $x_0, \ldots, x_j$ in $\mathcal{C}^{(1)}(\Gamma_0, \fA^{r_0}_{\Lambda, (1)}/ \fA^{r_0}_{\Lambda})$ and $y_0, \ldots, y_j$ in $\mathcal{C}^{(1)}(\Gamma_0, \fA^{r}_{\Lambda, (1)}/ \fA^{r}_{\Lambda})$, set
\[
y_{j+1}\coloneqq \frac{1}{p}(y_j-x_j) \in \mathcal{C}^{(1)}(\Gamma_0, \fA^{r}_{\Lambda, (1)}/ \fA^{r}_{\Lambda}),
\]
and pick a good lift $x_{j+1} \in \mathcal{C}^{(1)}(\Gamma_0, \fA^{r_0}_{\Lambda, (1)}/ \fA^{r_0}_{\Lambda})$ of $\overline{y}_{j+1}$ given by the above claim. Then it follows from the construction that $x = \sum_{j=0}^{\infty} p^j x_j$ is a desired norm compatible expansion.
\end{proof}

\begin{proposition} \label{prop: strict exactness for A(1)/A}
Assume Set-up~\ref{setup: r0-Gamma0}. There exists $0 < r_1 < r_0$ such that for any $0 < r \le r_1$ (here $r_1$ depends on the choice of $\Gamma_0$), the cochain complex
\[
\mathcal{C}^{(\bullet)}(\Gamma_0, \fA^{r}_{\Lambda, (m_0+1)} / \fA^{r}_{\Lambda, (m_0)})
\]
is strict exact at $d^{(0)}$ and $d^{(1)}$ for each $m_0 \ge 0$. Furthermore, we have
\[
\lambda_{r, \mathrm{quot}}(x) \le \frac{3}{2}\lambda_r(d^{(\ell)}(x))
\]
for any $\ell = 0, 1$, $0 < r \le r_1$, and any $x \in \mathcal{C}^{(\ell)}(\Gamma_0, \fA^{r}_{\Lambda, (m_0+1)} / \fA^{r}_{\Lambda, (m_0)})$. Here, $\lambda_{r, \mathrm{quot}}$ stands for the quotient norm on $\mathcal{C}^{(\ell)}(\Gamma_0, \fA^{r}_{\Lambda, (m_0+1)} / \fA^{r}_{\Lambda, (m_0)})/\ker(d^{(\ell)})$.
\end{proposition}

\begin{proof}
To simplify notation, we assume $m_0 = 0$; the case for general $m_0$ follows from a similar argument. Choose $r_0'$ and $r_1$ with $0 < r_1 < r_0' < r_0$ such that $(\frac{3}{2}p)^{\frac{r_0'}{r_0}} \le p^{\frac{1}{2}}$ and $(\frac{3}{2}p)^{\frac{r_1}{r_0'}} \le p^{\frac{1}{2}}$. Let $0 < r \le r_1$.\\

\noindent \textbf{Step 1.} (exactness at $d^{(0)}$). \\
Let $x \in \mathcal{C}^{(0)}(\Gamma_0, \fA^{r}_{\Lambda, (1)} / \fA^{r}_{\Lambda})=\fA^{r}_{\Lambda, (1)} / \fA^{r}_{\Lambda}$ such that $d^{(0)}(x) = 0$. Choose a norm compatible expansion $x=\sum_{j=0}^{\infty}p^jx_j$ as in Corollary~\ref{cor: good lift}(3). Since $d^{(0)}(\overline{x}) = d^{(0)}(\overline{x}_0) = 0$, we have $\overline{x}_0 = 0$ by Lemma~\ref{lem: complex-R-strict-exact}. Then $\lambda_r(x_0) = |\overline{x}_0|^r = 0$, so $x_0 = 0$ and
\[
d^{(0)}(px_1+p^2x_2+\cdots) = pd^{(0)}(x_1+px_2+\cdots) = 0.
\]
Since $\mathcal{C}^{(1)}(\Gamma_0, \fA^{r}_{\Lambda, (1)} / \fA^{r}_{\Lambda})$ is $p$-torsion free, we deduce that $d^{(0)}(\overline{x}_1) = 0$, which then implies $x_1 = 0$ by a similar argument as above. Proceeding inductively, we obtain $x_j = 0$ for all $j$, and hence $x = 0$. This proves the exactness at $d^{(0)}$.\\

\noindent \textbf{Step 2.} (strictness at $d^{(0)}$)\\
Let $0 < r \le r_1$ and $x \in \mathcal{C}^{(0)}(\Gamma_0, \fA^{r}_{\Lambda, (1)}/\mathbf{A}^r_{\Lambda})$. Once again, we choose a norm compatible expansion $x=\sum_{j=0}^{\infty}p^jx_j$ as in Corollary \ref{cor: good lift}(3). For each $j$, we have $|\overline{x}_j|^{r_0} \le \frac{3}{2}|d^{(0)}(\overline{x}_j)|^{r_0}$ (see Set-up~\ref{setup: r0-Gamma0}). Hence, since $\overline{d^{(0)}(x_j)} = d^{(0)}(\overline{x}_j)$, we have
\begin{equation} \label{eq: strictness at d0 for A(1)/A}
\lambda_s (x_j) = |\overline{x}_j|^s \le \frac{3}{2}|d^{(0)}(\overline{x}_j)|^s \le \frac{3}{2}\lambda_s (d^{(0)}(x_j))    
\end{equation}
for any $0 < s \le r_0$. Note that for each $j \ge 0$,
\[
\lambda_{r_0}(d^{(0)}(x_j)) \le \lambda_{r_0}(x_j) = |\overline{x}_j|^{r_0} \le \frac{3}{2}|d^{(0)}(\overline{x}_j)|^{r_0}.
\]
Since $\overline{d^{(0)}(x_j)} = d^{(0)}(\overline{x}_j)$ and $(\frac{3}{2} p)^{\frac{r_1}{r_0}} \le p^{\frac{1}{2}}$, we deduce from the same argument as in the proof of Proposition~\ref{prop: good lift weakly decompleting}(1) that $\lambda_s(d^{(0)}(x_j)) = |d^{(0)}(\overline{x}_j)|^s$ for all $0 < s \le r_1$. In particular, the presentation
$
d^{(0)}(x) = \sum_{j=0}^{\infty} p^j d^{(0)}(x_j)
$ 
is a norm compatible expansion. By Proposition~\ref{prop: tech lemma norm compatible expansion}, we have
\[
\lambda_r(d^{(0)}(x)) = \max_{j \ge 0}\{p^{-j} \lambda_r(d^{(0)}(x_j))\}.
\]
Consequently, we deduce from the inequality~\eqref{eq: strictness at d0 for A(1)/A} that $\lambda_r(x) \le \frac{3}{2}\lambda_r(d^{(0)}(x))$.\\

\noindent \textbf{Step 3.} (exactness at $d^{(1)}$)\\
Let $x \in \mathcal{C}^{(1)}(\Gamma_0, \fA^{r}_{\Lambda, (1)} / \fA^{r}_{\Lambda})$ such that $d^{(1)}(x) = 0$. Write $x=\sum_{j=0}^{\infty}p^jx_j$ as in \eqref{eq: norm compatible expansion}. In particular, $\lambda_r(x)=p^{-j_0}|x_{j_0}|^r$ for some $j_0$. Moreover, 
there exists $j_1>j_0$ such that 
$
p^{-j}|\overline{x}_j|^r \le p^{-\frac{1}{2}}\lambda_r(x)
$ 
for all $j \ge j_1$. We claim that 
$
p^{-j}|\overline{x}_j|^s \le p^{-\frac{1}{2}}\lambda_s(x) 
$
for all $0 < s \le r$ and all $j\ge j_1$. Indeed, since 
\[p^{-j}|\overline{x}_j|^r\le p^{-\frac{1}{2}}\lambda_r(x)=p^{-\frac{1}{2}}\cdot p^{-j_0}|x_{j_0}|^r,\]
we have
\[
\left(\frac{|\overline{x}_{j}|}{|\overline{x}_{j_0}|}\right)^r \le p^{j-j_0-\frac{1}{2}}.
\]
Hence, for any $0 < s \le r$, we have
\[
\left(\frac{|\overline{x}_{j}|}{|\overline{x}_{j_0}|}\right)^s \le (p^{j-j_0-\frac{1}{2}})^{\frac{s}{r}} \le p^{j-j_0-\frac{1}{2}}.
\]
Therefore, 
\[p^{-j}|\overline{x}_{j}|^s \le p^{-\frac{1}{2}}p^{-{j_0}}|\overline{x}_{j_0}|^s\le p^{-\frac{1}{2}}\lambda_s(x),
\]
as desired.
To prove the exactness at $d^{(1)}$, we need to find $y'\in \mathcal{C}^{(0)}(\Gamma_0, \fA^{r_0}_{\Lambda, (1)} / \fA^{r_0}_{\Lambda})$ such that $d^{(0)}(y')=x$. Such an element can be constructed by applying the following claim repeatedly.\\

\noindent \textbf{Claim 1.} There exists $y \in \mathcal{C}^{(0)}(\Gamma_0, \fA^{r_0}_{\Lambda, (1)} / \fA^{r_0}_{\Lambda})$ such that $\lambda_r(x-d^{(0)}(y)) \le p^{-\frac{1}{2}}\lambda_r(x)$.\\

\noindent \textit{Proof of Claim 1.} Since $d^{(1)}(\overline{x}) = d^{(1)}(\overline{x}_0) = 0$, there exists $\overline{y}_0 \in \mathcal{C}^{(0)}(\Gamma_0, R_{\Lambda, (1)}/R_{\Lambda})$ with $d^{(0)}(\overline{y}_0) = \overline{x}_0$ by Lemma~\ref{lem: complex-R-strict-exact}. By Corollary~\ref{cor: good lift}(3), there exists a lift $y_0 \in \mathcal{C}^{(0)}(\Gamma_0, \fA^{r_0}_{\Lambda, (1)} / \fA^{r_0}_{\Lambda})$ of $\overline{y}_0$ such that
$
\lambda_{r_0}(y_0) = |\overline{y}_0|^{r_0}.
$ 
Since $d^{(0)}$ is submetric, we have
\[
\lambda_{r_0}(d^{(0)}(y_0)) \le \lambda_{r_0}(y_0) = |\overline{y}_0|^{r_0} \le \frac{3}{2}|d^{(0)}(\overline{y}_0)|^{r_0} = \frac{3}{2}|\overline{x}_0|^{r_0}.
\]
Since $\overline{x}_0 = \overline{d^{(0)}(y_0)}$, $(\frac{3}{2} p)^{\frac{r_1}{r_0}} \le p^{\frac{1}{2}}$, and $r \le r_1$, we deduce from the same argument as in the proof of Proposition~\ref{prop: good lift weakly decompleting}(1) that
\[
\lambda_r(x_0-d^{(0)}(y_0)) \le p^{-\frac{1}{2}} |\overline{x}_0|^r = p^{-\frac{1}{2}}\lambda_r(x_0).
\]

If $j_1 = 1$, then 
$
\lambda_r(x-d^{(0)}(y_0)) \le p^{-\frac{1}{2}}\lambda_r(x).
$ 
In this case, we can simply take $y = y_0$.
Suppose $j_1 \ge 2$. Note that $x-d^{(0)}(y_0) = px^{(1)}$ for some $x^{(1)} \in \mathcal{C}^{(1)}(\Gamma_0, \fA^{r}_{\Lambda, (1)}/\mathbf{A}^r_{\Lambda})$. We claim that we can express $px^{(1)}$ as a $\lambda_r$-convergent norm compatible expansion
$
px^{(1)} = \sum_{j =1}^{\infty} p^j x_j^{(1)},
$ 
with each $x_j^{(1)} \in \mathcal{C}^{(1)}(\Gamma_0, \fA^{r_0}_{\Lambda, (1)} / \fA^{r_0}_{\Lambda})$, such that 
$
p^{-j}|\overline{x}^{(1)}_j|^r \le p^{-\frac{1}{2}}\lambda_r(x)
$ 
for all $j \ge j_1$. Let $S$ be the set of indices $j \ge 1$ such that $p^{-j}|x_j|^r > p^{-\frac{1}{2}}\lambda_r(x)$. Then $S$ is finite, and $\max_{j\in S} \{{j\}} \le j_1-1$. Consider
\[
x'\coloneqq (x_0-d^{(0)}(y_0))+\sum_{j\ge 1, j \notin S} p^j x_j.
\]
Note that $\lambda_r(x') \le p^{-\frac{1}{2}}\lambda_r(x)$, since $\lambda_r(x_0-d^{(0)}(y_0)) \le p^{-\frac{1}{2}}\lambda_r(x_0) \le p^{-\frac{1}{2}}\lambda_r(x)$. By Lemma~\ref{lem: norm-compatible-expansion-for-continuous-functions}, $x'$ admits a $\lambda_r$-convergent norm compatible expansion $x' = \sum_{j=1}^{\infty} p^j x_j'$. Then
\[
px^{(1)} = x'+\sum_{j \in S} p^j x_j = \sum_{j \notin S} p^j x_j'+ \sum_{j \in S} p^j(x_j+x_j')
\]
is a desired norm compatible expansion.

Note that $pd^{(1)}(x^{(1)})=d^{(1)}(x-d^{(0)}(y_0)) = 0$ and $\mathcal{C}^{(2)}(\Gamma_0, \mathbf{A}^r_{\Lambda, (1)}/\mathbf{A}^r_{\Lambda})$ is $p$-torsion free. So $d^{(1)}(x^{(1)}) = 0$. Applying the same procedure above to $x^{(1)}$ (the same way as in obtaining $y_0$ from $x$), we obtain an element $y_1 \in \mathcal{C}^{(0)}(\Gamma_0, \fA^{r_0}_{\Lambda, (1)} / \fA^{r_0}_{\Lambda})$ such that
\[
d^{(0)}(\overline{y}_1) = \overline{x}^{(1)}_1, \quad \lambda_{r_0}(y_1) = |\overline{y}_1|^{r_0},
\]
and
\[
\lambda_r(x^{(1)}_1-d^{(0)}(y_1)) \le p^{-\frac{1}{2}}\lambda_r(x^{(1)}_1).
\]
We now have $x-d^{(0)}(y_0+py_1) = p^2 x^{(2)}$ for some $x^{(2)} \in \mathcal{C}^{(1)}(\Gamma_0, \fA^{r}_{\Lambda, (1)}/\mathbf{A}^r_{\Lambda})$. Iterating this process yields $x^{(j)} \in \mathcal{C}^{(1)}(\Gamma_0, \fA^{r}_{\Lambda, (1)}/\mathbf{A}^r_{\Lambda})$ and $y_j \in \mathcal{C}^{(0)}(\Gamma_0, \fA^{r}_{\Lambda, (1)}/\mathbf{A}^r_{\Lambda})$ for $0 \le j \le j_1-1$, satisfying analogous conditions as above. By construction and the definition of $j_1$, we deduce that
\[
\lambda_r(x-d^{(0)}(y)) \le p^{-\frac{1}{2}}\lambda_r(x).
\]
where $y \coloneqq y_0+py_1+\cdots +p^{j_1-1}y_{j_1-1}$. This finishes the proof of Claim 1.\\

We now return to the construction of $y'$ such that $d^{(0)}(y')=x$. Let $y$ be the element from Claim 1. Let $y^{(0)} = 0$ and let $y^{(1)} = y$. Next, applying Claim 1 to $x-d^{(0)}(y^{(0)}+y^{(1)})$ instead of $x$, we obtain an element $y^{(2)} \in \mathcal{C}^{(0)}(\Gamma_0, \fA^{r}_{\Lambda, (1)}/\mathbf{A}^r_{\Lambda})$ such that
\[
\lambda_r \left(x-d^{(0)}(y^{0}+y^{(1)}+y^{(2)})\right) \le p^{-\frac{1}{2}}\lambda_r \left(x-d^{(0)}(y^{0}+y^{(1)})\right).
\]
Proceeding inductively, we obtain $y^{(j)} \in \mathcal{C}^{(0)}(\Gamma_0, A^{r}_{\Lambda, (1)}/\mathbf{A}^r_{\Lambda})$ for each $j \ge 0$ satisfying
\[
\lambda_r (x-d^{(0)}(y^{(0)}+\cdots y^{(j)})) \le p^{-\frac{j}{2}}\lambda_r (x). 
\]
Furthermore, it follows from the construction above that $\lambda_r (y^{(j)})$ goes to $0$ as $j \rightarrow \infty$. Since $\mathcal{C}^{(0)}(\Gamma_0, \fA^{r}_{\Lambda, (1)}/\mathbf{A}^r_{\Lambda})$ is complete with respect to $\lambda_r$, we obtain a well-defined element
\[
y'\coloneqq \sum_{j=0}^{\infty} y^{(j)} \in \mathcal{C}^{(0)}(\Gamma_0, \fA^{r}_{\Lambda, (1)}/\mathbf{A}^r_{\Lambda})
\]
satisfying $x = d^{(0)}(y')$. This proves the exactness at $d^{(1)}$.\\

\noindent \textbf{Step 4.} (strictness at $d^{(1)}$)\\
In fact, we will show that for any $0 < r \le r_1$ and any $x \in \mathcal{C}^{(1)}(\Gamma_0, \fA^{r}_{\Lambda, (1)}/\mathbf{A}^r_{\Lambda})$, there exists $z\in \mathcal{C}^{(0)}(\Gamma_0, \fA^{r}_{\Lambda, (1)}/\mathbf{A}^r_{\Lambda})$ satisfying $\lambda_r(x+d^{(0)}(z))\le \frac{3}{2}\lambda_r(d^{(1)}(x))$. Without loss of generality, we may assume that $\lambda_r(x) \le \lambda_r(x+d^{(0)}(z))$ for all $z \in \mathcal{C}^{(0)}(\Gamma_0, \fA^{r}_{\Lambda, (1)}/\mathbf{A}^r_{\Lambda})$. Write $x=\sum_{j=0}^{\infty}p^jx_j$ as in the norm compatible expansion~\eqref{eq: norm compatible expansion}. By Set-up~\ref{setup: r0-Gamma0}, for each $j \ge 0$, we have
\[
|\overline{x}_j|_{\mathrm{quot}}^{r_0} \le \frac{3}{2}|d^{(1)}(\overline{x}_j)|^{r_0}.
\]
If $|\overline{x}_j| = |\overline{x}_j|_{\mathrm{quot}}$, let $z_j = 0$. If $|\overline{x}_j| > |\overline{x}_j|_{\mathrm{quot}}$, we choose $\overline{z}_j \in \mathcal{C}^{0}(\Gamma_0, R_{\Lambda, (1), i}/R_{\Lambda, i})$ satisfying $|\overline{x}_j+d^{(0)}(\overline{z}_j)| = |\overline{x}_j|_{\mathrm{quot}}$ (in particular, $|\overline{x}_j| = |d^{(0)}(\overline{z}_j)|$). By Corollary~\ref{cor: good lift}(3), there exists a lift $z_j \in \mathcal{C}^{(0)}(\Gamma_0, \fA^{r_0}_{\Lambda, (1)} / \fA^{r_0}_{\Lambda})$ of $\overline{z}_j$ such that
\[
\lambda_s(z_j-[\overline{z}_j]) \le p^{-\frac{1}{2}}|\overline{z}_j|^s
\]
for all $0 < s \le r_0$.\\

\noindent \textbf{Claim 2.} For any $0 < s \le r_0'$, we have
\[
\lambda_s(x_j+d^{(0)}(z_j)) \le \max\{p^{-\frac{1}{2}}|\overline{x}_j|^s, |\overline{x}_j+d^{(0)}(\overline{z}_j)|^s \} \le \max\left\{p^{-\frac{1}{2}}|\overline{x}_j|^s, \left(\frac{3}{2}\right)^{\frac{s}{r_0}}|d^{(1)}(\overline{x}_j)|^s \right\}
\]
for every $j\ge 0$.\\

\noindent \textit{Proof of Claim 2.} We have 
\[
\lambda_{r_0}(d^{(0)}(z_j)) \le \lambda_{r_0}(z_j) = |\overline{z}_j|^{r_0} \le \frac{3}{2}|d^{(0)}(\overline{z}_j)|^{r_0},
\]
where the first inequality follows from that $d^{(0)}$ is submetric, and the last inequality is due to Set-up~\ref{setup: r0-Gamma0}. Since $\overline{d^{(0)}(z_j)} = d^{(0)}(\overline{z}_j)$, $(\frac{3}{2}p)^{\frac{r_0'}{r_0}} \le p^{\frac{1}{2}}$, and $s \le r_0'$, we deduce from the same argument as in the proof of Proposition~\ref{prop: good lift weakly decompleting}(1) that
\[
\lambda_s(d^{(0)}(z_j)-[d^{(0)}(\overline{z}_j)]) \le p^{-\frac{1}{2}}|d^{(0)}(\overline{z}_j)|^s \le p^{-\frac{1}{2}}|\overline{x}_j|^s.
\]
Note that 
\[
(x_j+d^{(0)}(z_j))- [\overline{x}_j+d^{(0)}(\overline{z}_j)] = (x_j-[\overline{x}_j])+(d^{(0)}(z_j)-[d^{(0)}(\overline{z}_j)])+([\overline{x}_j]+[d^{(0)}(\overline{z}_j)]-[\overline{x}_j+d^{(0)}(\overline{z}_j)]).
\]
Since $x_j$ is a good lift of $\overline{x}_j$ by construction, we have $\lambda_s(x_j-[\overline{x}_j]) \le p^{-\frac{1}{2}}|\overline{x}_j|^s$. Furthermore, using Witt vector expansion, we deduce that
\[
\lambda_s([\overline{x}_j]+[d^{(0)}(\overline{z}_j)]-[\overline{x}_j+d^{(0)}(\overline{z}_j)]) \le p^{-1}\max\{|\overline{x}_j|^s, |d^{(0)}(\overline{z}_j)|^s\} = p^{-1}|\overline{x}_j|^s.
\]
Combining these inequalities, we complete the proof of Claim 2.\\ 

Next, we consider the following cases separately.\\

\noindent \textbf{Case 4.1:} $\frac{3}{2}|d^{(1)}(\overline{x}_j)|^{r_0} < p^{-\frac{1}{2}}|\overline{x}_j|^{r_0}$.\\

In this case, we have
\[
\left(\frac{3}{2}\right)^{\frac{s}{r_0}}|d^{(1)}(\overline{x}_j)|^s < (p^{-\frac{1}{2}})^{\frac{s}{r_0}}|\overline{x}_j|^s
\]
for each $0 < s \le r_0'$. So by Claim 2,
\[
\lambda_s(d^{(1)}(x_j)) \le \lambda_s(x_j+d^{(0)}(z_j)) \le (p^{-\frac{1}{2}})^{\frac{s}{r_0}}|\overline{x}_j|^s
\]
for each $0 < s \le r_0'$.\\ 

\noindent \textbf{Case 4.2:} $\frac{3}{2}|d^{(1)}(\overline{x}_j)|^{r_0} \ge p^{-\frac{1}{2}}|\overline{x}_j|^{r_0}$.\\

For this case, we have $\frac{3}{2}|d^{(1)}(\overline{x}_j)|^{r_0'} \ge p^{-\frac{1}{2}}|\overline{x}_j|^{r_0'}$ because $r_0' < r_0$. Hence,
\[
\lambda_{r_0'}(d^{(1)}(x_j)) \le \lambda_{r_0'}(x_j+d^{(0)}(z_j)) \le \frac{3}{2}|d(\overline{x}_j)|^{r_0'},
\]
where the first inequality holds as $d^{(1)}$ is submetric, and the second inequality follows from Claim 2. Since $\overline{d^{(1)}(x_j)} = d^{(1)}(\overline{x}_j)$ and $(\frac{3}{2} p)^{\frac{r_1}{r_0'}} \le p^{\frac{1}{2}}$, we deduce using the same argument as in the proof of Proposition~\ref{prop: good lift weakly decompleting}(1) that $\lambda_s(d^{(1)}(x_j)) = |d^{(1)}(\overline{x}_j)|^s$ for all $0 < s \le r_1$.\\

To proceed, consider
\[
z\coloneqq \sum_{j \ge 0} p^j z_j \in \mathcal{C}^{(0)}(\Gamma_0, \fA^{r}_{\Lambda, (1)}/\mathbf{A}^r_{\Lambda}),
\]
and consider the expansion
\[
x+d^{(0)}(z) = \sum_{j \ge 0} p^j(x_j+d^{(0)}(z_j)),
\]
which induces
\[
d^{(1)}(x+d^{(0)}(z)) = d^{(1)}(x) = \sum_{j = 0}^{\infty} p^j d^{(1)}(x_j).
\]
Let $S_1$ be the set of $j \ge 0$ such that Case 1 holds, and $S_2$ be the set of $j \ge 0$ such that Case 2 holds. The expansion of $d^{(1)}(x)$ decomposes into
\[
d^{(1)}(x) = \sum_{j \in S_1} p^j d^{(1)}(x_j) + \sum_{j \in S_2} p^j d^{(1)}(x_j).
\]

By the above observations for Case 4.1 and Case 4.2, we have
\[
\lambda_r(d^{(1)}(x_j)) \le \lambda_r(x_j+d^{(0)}(z_j)) \le (p^{-\frac{1}{2}})^{\frac{r}{r_0}}|\overline{x}_j|^r
\]
for each $j \in S_1$, and $\lambda_r(d^{(1)}(x_j)) = |d^{(1)}(\overline{x}_j)|^r$ for each $j \in S_2$ so that the sum $\sum_{j \in S_2} p^j d^{(1)}(x_j)$ is a norm compatible expansion.\\

\noindent \textbf{Claim 3.} There exists $j_1 \in S_2$ such that $p^{-j_1}|\overline{x}_{j_1}|^r = \lambda_r(x)$ and $\lambda_r(x_{j_1}) \le \lambda_r(x_{j_1}+d^{(0)}(z_{j_1}))$.\\

\noindent \textit{Proof of Claim 3.} Recall that
$
\lambda_r(x) = \max_j \{p^{-j}|\overline{x}_j|^r\}.
$ 
So the set of indices $j$ satisfying $p^{-j}|\overline{x}_j|^r = \lambda_r(x)$ is non-empty. Note that among such $j$'s, there exists an index $j_1$ which further satisfies $\lambda_r(x_{j_1}) \le \lambda_r(x_{j_1}+d^{(0)}(z_{j_1}))$; otherwise, we have a contradiction to the assumption that $\lambda_r(x) \le \lambda_r(x+d^{(0)}(z'))$ for all $z' \in \mathcal{C}^{(0)}(\Gamma_0, \fA^{r}_{\Lambda, (1)}/\mathbf{A}^r_{\Lambda})$. By Claim 2, we have
\[
|\overline{x}_{j_1}|^r = \lambda_r(x_{j_1}) \le \left(\frac{3}{2}\right)^{\frac{r}{r_0}}|d^{(1)}(\overline{x}_{j_1})|^r,
\]
so $\frac{3}{2}|d^{(1)}(\overline{x}_{j_1})|^{r_0} \ge |\overline{x}_{j_1}|^{r_0} \ge p^{-\frac{1}{2}}|\overline{x}_{j_1}|^{r_0}$. Thus, $j_1 \in S_2$ and Claim 3 holds.\\

In summary, we have
\[
\lambda_r\left(\sum_{j \in S_1} p^j d^{(1)}(x_j)\right) \le p^{-\frac{r}{2 r_0}} \lambda_r(x) < \lambda_r(x),
\]
and
\[
\lambda_r\left(\sum_{j \in S_2} p^j d^{(1)}(x_j)\right) = \max_{j \in S_2}\{p^{-j}\lambda_r(d^{(1)}(x_j))\}
\]
by Proposition~\ref{prop: tech lemma norm compatible expansion} and using that $\sum_{j \in S_2} p^j d^{(1)}(x_j)$ is a norm compatible expansion. Consequently, by Claim 3, we conclude that 
\[
\lambda_r(x) \le \frac{3}{2}\lambda_r(d^{(1)}(x)).
\]
\end{proof}

The following statement will be used in Theorem~\ref{thm: decompletion-etale-(phi, Gamma)-modules}.

\begin{corollary} \label{cor: complex-A-strict-exact}
Assume Set-up~\ref{setup: r0-Gamma0}, and let $0 < r_1 < r_0$ be as in Proposition~\ref{prop: strict exactness for A(1)/A}. For any $0 < r \le r_1$, the cochain complex
\[
\mathcal{C}^{(\bullet)}(\Gamma_0, \widetilde{\mathbf{A}}^r_X / \fA^{r}_{\Lambda, (m_0)})
\]
is strict exact at $d^{(0)}$ and $d^{(1)}$ for each $m_0 \ge 0$. Furthermore, we have
\[
\lambda_{r, \mathrm{quot}}(x) \le \frac{3}{2}\lambda_r(d^{(\ell)}(x))
\]
for any $\ell = 0, 1$, $0 < r \le r_1$, and any $x \in \mathcal{C}^{(\ell)}(\Gamma_0, \widetilde{\mathbf{A}}^r_X / \fA^{r}_{\Lambda, (m_0)}))$. 
\end{corollary}

\begin{proof}
Let $0 < r \le r_1$. For any $m > m_0$, consider the cochain complex
\[
\mathcal{C}^{(\bullet)}(\Gamma_0, \mathbf{A}^r_{\Lambda, (m)} / \mathbf{A}^r_{\Lambda, (m_0)}).
\]
It follows from the same argument as in the proof of Proposition~\ref{prop: strict exactness for A(1)/A} that this complex is exact at $d^{(0)}$ and $d^{(1)}$, and that 
\[
\lambda_{r, \mathrm{quot}}(x) \le \frac{3}{2}\lambda_r(d^{(\ell)}(x))
\]
for any $x \in \mathcal{C}^{(\ell)}(\Gamma_0, \mathbf{A}^r_{\Lambda, (m)} / \mathbf{A}^r_{\Lambda, (m_0)})$ with $\ell = 0, 1$. Thus, we can proceed similarly as in the proof of Lemma~\ref{lem: complex-R-strict-exact} (the part after Step 3 therein) to deduce the statement.
\end{proof}

\subsection{Strict-exactness for type $\mathbf{C}$} \label{subsec: strict-exactness-type-C}

We now analyze the cochain complexes of type $\fC$. The following result will be used in Theorem~\ref{thm: fully-faithfulness-for-free-modules}.

\begin{proposition} \label{prop: C-complex-strict-exact}
Assume Set-up~\ref{setup: r0-Gamma0}, and let $0 < r_1 < r_0$ be as in Proposition~\ref{prop: strict exactness for A(1)/A}. For any $0 < s \le r \le r_1$, the cochain complex
\[
\mathcal{C}^{(\bullet)}(\Gamma_0, \widetilde{\mathbf{C}}^{[s, r]}_X / \fC^{[s, r]}_{\Lambda, (m_0)})
\]
is strict exact at $d^{(0)}$ for each $m_0 \ge 0$. Furthermore, we have
\[
\max\{\lambda_s, \lambda_r\}(x) \le \frac{3}{2}\max\{\lambda_s, \lambda_r\}(d^{(0)}(x))
\]
for any $0 < s \le r \le r_1$ and any $x \in \mathcal{C}^{(0)}(\Gamma_0, \widetilde{\mathbf{C}}^{[s, r]}_X / \fC^{[s, r]}_{\Lambda, (m_0)})$.
\end{proposition}

\begin{proof}
To simplify notation, we set $m_0 = 0$. The case for general $m_0$ follows from a similar argument. Write $\lambda_{s, r}(\cdot)\coloneqq \max\{\lambda_s(\cdot), \lambda_r(\cdot)\}$. Recall from the proof of Corollary~\ref{cor: complex-A-strict-exact} that for any $0 < t \le r_1$, $m \ge 0$, and $x \in \mathcal{C}^{(0)}(\Gamma_0, \fA^t_{\Lambda, (m)} / \fA^t_{\Lambda})[p^{-1}]$, we have
\begin{equation} \label{eq: inequality-for-C-complex}
\lambda_t(x) \le \frac{3}{2}\lambda_t(d^{(0)}(x)). 
\end{equation}

We first consider the exactness at $d^{(0)}$. Let $y \in \mathcal{C}^{(0)}(\Gamma_0, \widetilde{\mathbf{C}}^{[s, r]}_X / \fC^{[s, r]}_{\Lambda})$ such that $d^{(0)}(y) = 0$. It suffices to show that for any $\epsilon > 0$, we have $\lambda_{s, r}(y) < \epsilon$. Let $y_1 \in \mathcal{C}^{(0)}(\Gamma_0, \mathbf{A}^r_{\Lambda, (m)} / \mathbf{A}^r_{\Lambda})[p^{-1}]$ for some $m$ such that $\lambda_{s, r}(y-y_1) \le \frac{1}{3}\epsilon$. Then
\[
\lambda_{s, r}(d^{(0)}(y-y_1)) = \lambda_{s, r}(d^{(0)}(y_1)) \le \frac{1}{3}\epsilon.
\] 
Choose $j \ge 0$ sufficiently large so that $p^j y_1 \in \mathcal{C}^{(0)}(\Gamma_0, \mathbf{A}^r_{\Lambda, (m)} / \mathbf{A}^r_{\Lambda})$. By the inequality~\eqref{eq: inequality-for-C-complex}, we have
\[
\lambda_{s, r}(p^j y_1) \le \frac{3}{2}\lambda_{s, r}(d^{(0)}(p^j y_1)) = \frac{3}{2}\cdot p^{-j} \lambda_{s, r}(d^{(0)}(y_1)) \le \frac{1}{2}p^{-j}\cdot \epsilon. 
\]
Thus, we have
\[
\lambda_{s, r}(y) \le \max\{\lambda_{s, r}(y-y_1), \lambda_{s, r}(y_1) \} \le \frac{1}{2}\epsilon.
\]
This shows the exactness at $d^{(0)}$. 

Finally, we can proceed similarly as in the proof of Lemma~\ref{lem: complex-R-strict-exact} to deduce from \eqref{eq: inequality-for-C-complex} that for any $x \in \mathcal{C}^{(0)}(\Gamma_0, \widetilde{\mathbf{C}}^{[s, r]}_X / \fC^{[s, r]}_{\Lambda})$, we have
$
\lambda_{s, r}(x) \le \frac{3}{2}\lambda_{s, r}(d^{(0)}(x)).
$ 
\end{proof}

\vspace{0.1in}
\section{Decompletion of relative \'etale $(\varphi, \Gamma)$-modules} \label{sec: decompletion-relative-etale-phiGamma-modules}

Let $h\colon X \rightarrow \D_K^d$ be a toric chart, and $\Lambda = (R_{\Lambda}, \mathbf{A}_{\Lambda})$ be a decompleting datum throughout the section. Recall the (perfect and imperfect) relative period rings associated with $\Lambda$ introduced in Definitions~\ref{definition:valuation_v_r_relative_perfect}, \ref{def: imperfect-period-rings-B-C}, and \ref{def: Frobenius-pullback-imperfect-period-rings}.
In this section, we study the \textit{decompletion} and \textit{deperfection} of relative $\varphi$-modules and $(\varphi, \Gamma)$-modules from perfect relative period rings to those over imperfect relative period rings. More precisely, we will consider \textit{descent} of $\varphi$-modules and $(\varphi, \Gamma)$-modules along maps of period rings as in the following diagram
\begin{equation} \label{diagram:descent_along_all_the_rings}
\begin{tikzcd}[column sep=2em, row sep=1em]
    \mathbf{A}_{\Lambda}^{\dagger} \arrow[rr, "(1)"]   \arrow[ddd] & &   \widetilde{\mathbf{A}}_{X}^{\dagger} \arrow[rrrdd, "(4)"] \arrow[ddd]   
    \\   \\ 
    &  &  &
    \mathbf{A}_{\Lambda} \arrow[from=llluu, crossing over, swap,  "(3)"] \arrow[rr, swap, "(2)"] & &\widetilde{\mathbf{A}}_{X}   \arrow[ddd] \\ 
    \mathbf{B}_{\Lambda}^{\dagger} \arrow[rr, "(5)"] \arrow[ddd, "(9)"]  & &   \widetilde{\mathbf{B}}_{X}^{\dagger} \arrow[ddd,  near end, "(10)"] \arrow[rrrdd, "(8)"] \\  \\ 
       &  & &  \mathbf{B}_{\Lambda} \arrow[from=uuu, crossing over] \arrow[from=llluu, swap,  "(7)", crossing over] \arrow[rr, swap, "(6)"] & & \widetilde{\mathbf{B}}_{X}   \\
      \mathbf{C}_{\Lambda}^{\dagger} \arrow[rr, "(11)"]& & \widetilde{\mathbf{C}}_{X}^{\dagger}  
\end{tikzcd}
\end{equation}
constructed in \S \ref{sec:relative_period_rings:perfect} and \S \ref{sec:weakly_decompleting_data}.
The main results are Theorems~\ref{thm: fully-faithfulness-for-free-modules} and \ref{thm: equivalence-isogeny-local-system-overconvergent-etale-(phi, Gamma)-modules-C}.

\subsection{Relative $(\varphi, \Gamma)$-modules} \label{subsec: definition-relative-(phi, Gamma)-modules}

We first introduce the notions of $\varphi$-modules and $(\varphi, \Gamma)$-modules, which are analogues of Definition \ref{def:varphi_mod_relative_perfect}. 

\begin{definition}[$\varphi$-modules and \'etale $\varphi$-modules]  \label{def:phi_modules_and_etale_phi_modules} \hfill
\begin{enumerate}
\item  Let $r>0$ and let $\mathbf{D}^r$ be one of the following rings
\begin{equation} \label{eq:some_rings_with_convergence_condition} 
 \{
  \mathbf{A}_{\Lambda}^r, \mathbf{B}_{\Lambda}^r, \mathbf{C}_{\Lambda}^r, \breve{\mathbf{A}}_{\Lambda}^r, \breve{\mathbf{B}}_{\Lambda}^r, \breve{\mathbf{C}}_{\Lambda}^r,
 \widetilde{\mathbf{A}}^r_{X}, \widetilde{\mathbf{B}}^r_{X}, \widetilde{\mathbf{C}}^r_{X}
 \}.
 \end{equation} 
A \emph{$\varphi$-module} over $\mathbf{D}^r $ is a finite projective $\mathbf{D}^r$-module $M$ together with an isomorphism
\[
\varphi_{M}\colon M\otimes_{{\mathbf{D}}^r,\varphi} {\mathbf{D}}^{r/p} \isom M\otimes_{{\mathbf{D}}^r,\iota}{\mathbf{D}}^{r/p}
\]
of ${\mathbf{D}}^{r/p}$-modules, where $\varphi\colon \mathbf{D}^{r} \ra \mathbf{D}^{r/p}$ is the Frobenius map (see Remarks~\ref{remark:Frobenius_classical} and \ref{remark:Frobenius_relative_and_action_of_Gamma}(3)), and $\iota$ is the inclusion map. \item Suppose $0 < s \le \frac{r}{p}$ and let $\mathbf{D}^{[s, r]}$ be one of the following rings
\begin{equation} \label{eq:some_rings_with_convergence_condition_1}
    \{\mathbf{C}_{\Lambda}^{[s, r]}, \breve{\mathbf{C}}_{\Lambda}^{[s, r]},  \widetilde{\mathbf{C}}_{X}^{[s, r]} \}
\end{equation}   
A \emph{$\varphi$-module} over $\mathbf{D}^{[s, r]}$ is a finite projective $\mathbf{D}^{[s, r]}$-module $M$ with an isomorphism
\[
\varphi_M\colon (M\otimes_{\mathbf{D}^{[s, r]}, \varphi} \mathbf{D}^{[s/p, \,r/p]})\otimes_{\mathbf{D}^{[s/p, \,r/p]}, \iota} \mathbf{D}^{[s, \frac{r}{p}]} \isom M\otimes_{\mathbf{D}^{[s, r]}, \iota} \mathbf{D}^{[s, \frac{r}{p}]}
\]
of $\mathbf{D}^{[s, \frac{r}{p}]}$-modules.
\item Let $\mathbf{D}$ be one of the rings from the set 
\begin{equation} \label{eq:some_rings_with_convergence_condition_2}
\{ \mathbf{A}_{\Lambda}^{\dagger}, \mathbf{B}_{\Lambda}^{\dagger}, \mathbf{C}_{\Lambda}^{\dagger}, \breve{\mathbf{A}}_{\Lambda}^{\dagger}, \breve{\mathbf{B}}_{\Lambda}^{\dagger}, \breve{\mathbf{C}}_{\Lambda}^{\dagger}
\} \tu{ or } \{ \mathbf{A}_{\Lambda}, \mathbf{B}_{\Lambda}, \breve{\mathbf{A}}_{\Lambda}, \breve{\mathbf{B}}_{\Lambda}
 \}.
\end{equation}
A \emph{$\varphi$-module} over $\mathbf{D}$ is a finite projective $\mathbf{D}$-module $M$ together with an isomorphism
\[
\varphi_{M}\colon M\otimes_{{\mathbf{D}},\varphi} \mathbf{D} \isom M 
\]
of $\mathbf{D}$-modules. 
 
\item A $\varphi$-module over $\mathbf{B}_{\Lambda}$ (resp. $\breve{\mathbf{B}}_{\Lambda}$
) is called \emph{\'etale} if it comes from a $\varphi$-module over $\mathbf{A}_{\Lambda}$ (resp. $\breve{\mathbf{A}}_{\Lambda}$) via base change. Similarly, a $\varphi$-module over $\mathbf{B}_{\Lambda}^{\dagger}$ or $ \mathbf{C}_{\Lambda}^{\dagger}$ (resp. $\breve{\mathbf{B}}_{\Lambda}^{\dagger}$ or $\breve{\mathbf{C}}_{\Lambda}^{\dagger}$) is called \emph{\'etale} if it comes from a $\varphi$-module over $\mathbf{A}_{\Lambda}^{\dagger} $ (resp. $\breve{\mathbf{A}}_{\Lambda}^{\dagger}$) via base change. 
\end{enumerate}
\end{definition}

\begin{definition}[$(\varphi, \Gamma)$-modules] \hfill
\label{def:relative_phi_gamma_modules}   
\begin{enumerate}
\item Let $\mathbf{D}$ be one of the rings in (\ref{eq:some_rings_with_convergence_condition}), (\ref{eq:some_rings_with_convergence_condition_1}), or  (\ref{eq:some_rings_with_convergence_condition_2}) from  Definition~\ref{def:phi_modules_and_etale_phi_modules}. A \emph{$\Gamma$-module} over $\mathbf{D}$ is a finite projective $\mathbf{D}$-module equipped with a continuous semilinear action of $\Gamma$.\footnote{Semilinearity here means that $\gamma (\alpha \cdot x) = \gamma (\alpha) \cdot \gamma (x)$ for all $\alpha \in \mathbf{D}$ and $x \in M$. We refer the reader to \cite[Definition~1.2.2]{Kedlaya-LiuII} for a discussion of topology on $M$.}  A \emph{$(\varphi, \Gamma)$-module} over $\mathbf{D}$ is a $\varphi$-module $M$ over $\mathbf{D}$, equipped with a continuous semilinear action of $\Gamma$ that commutes with the action of $\varphi$. 

\item  A $(\varphi, \Gamma)$-module over one of the following rings 
\begin{equation} \label{eq:some_rings_with_convergence_condition_3} 
    \mathbf{B}_{\Lambda}, \breve{\mathbf{B}}_{\Lambda}, 
    \mathbf{B}_{\Lambda}^{\dagger}, \breve{\mathbf{B}}_{\Lambda}^\dagger,    
    \mathbf{C}_{\Lambda}^{\dagger}, \breve{\mathbf{C}}_{\Lambda}^{\dagger}, 
     \widetilde{\mathbf{B}}_{X}, 
   \widetilde{\mathbf{B}}_{X}^{\dagger}, \widetilde{\mathbf{C}}_{X}^{\dagger}, 
  \end{equation} 
is called \emph{\'etale} if it is the base change of a $(\varphi, \Gamma)$-module over the rings 
\[    
 \mathbf{A}_{\Lambda}, \breve{\mathbf{A}}_{\Lambda}, 
    \mathbf{A}_{\Lambda}^{\dagger}, \breve{\mathbf{A}}_{\Lambda}^\dagger,     
    \mathbf{A}_{\Lambda}^{\dagger}, \breve{\mathbf{A}}_{\Lambda}^{\dagger}, 
     \widetilde{\mathbf{A}}_{X}, 
   \widetilde{\mathbf{A}}_{X}^{\dagger}, \widetilde{\mathbf{A}}_{X}^{\dagger} 
\]
respectively. 
\end{enumerate}
\end{definition}

\begin{remark} \hfill
\begin{enumerate}
    \item We point out a slight deviation from the standard terminology in literature: in the definitions above, all $\varphi$-modules and $(\varphi, \Gamma)$-modules are assumed to be finite projective. 
    
    \item Our notion of \'etale $(\varphi, \Gamma)$-modules requires the corresponding lattices over period rings of type $\mathbf{A}$ to be stable under both $\varphi$ and $\Gamma$. This is similar to \cite[Definition~I.4.2]{Cherbonnier_Colmez}. On the other hand, \cite{Kedlaya-LiuII} only requires the lattice to be stable under $\varphi$, so the definition we use is more restrictive. 
\end{enumerate} 
\end{remark}

\begin{notation} \hfill
\begin{enumerate}
    \item Let $\mathbf{D}$ be one of the rings in (\ref{eq:some_rings_with_convergence_condition}), (\ref{eq:some_rings_with_convergence_condition_1}), or (\ref{eq:some_rings_with_convergence_condition_2}). We write 
    $
    \mathrm{Mod}_{/\mathbf{D}}^{\varphi}$, (resp. $ \mathrm{Mod}_{/\mathbf{D}}^{\Gamma}, $ resp. $\mathrm{Mod}_{/\mathbf{D}}^{(\varphi, \Gamma)}$)   for the category of $\varphi$-modules (resp. the category of $\Gamma$-modules, resp. the category of $(\varphi, \Gamma)$-modules over $\mathbf{D}$.  
    
    \item Let $\mathbf{D}$ be one of the rings in (\ref{eq:some_rings_with_convergence_condition_3}). We write $ \mathrm{Mod}_{/\mathbf{D}}^{\varphi, \ett} 
     $ (resp. $  \mathrm{Mod}_{/\mathbf{D}}^{(\varphi, \Gamma), \ett}$)
    for the category of \'etale $\varphi$-modules (resp. the category of \'etale $(\varphi, \Gamma)$-modules) over $\mathbf{D}$.
\end{enumerate} 
\end{notation}

\begin{remark}
Let $\mathbf{D}$ be one of the rings in  \eqref{diagram:descent_along_all_the_rings}. Let $M$ and $N$ be $\varphi$-modules over $\mathbf{D}$. Let 
\[ L\coloneqq \Hom_{\mathbf{D}}(M, N)\cong M^{\vee}\otimes_{\mathbf{D}} N, \] which is a finite projective module over $\mathbf{D}$. The maps  $\varphi_M\colon M\otimes_{\mathbf{D}, \varphi} \mathbf{D} \stackrel{\cong}{\rightarrow} M$ and $\varphi_N\colon N\otimes_{\mathbf{D}, \varphi} \mathbf{D} \stackrel{\cong}{\rightarrow} N$ induce an isomorphism
\[
\varphi_L\colon L\otimes_{\mathbf{D}, \varphi} \mathbf{D} \cong \Hom_{\mathbf{D}}(M\otimes_{\mathbf{D}, \varphi} \mathbf{D},  N\otimes_{\mathbf{D}, \varphi} \mathbf{D}) \stackrel{\cong}{\longrightarrow} L.
\]
Also note that $f \in L = \Hom_{\mathbf{D}}(M, N)$ is compatible with $\varphi_M$ and $\varphi_N$ if and only if $\varphi(f) = f$, where $\varphi(f)\coloneqq \varphi_L(f\otimes 1)$ by a slight abuse of notation.  
\end{remark}

\begin{remark} \label{rem: base-change-functor-compatible-with-tensor-dual}
Let $D_1 \rightarrow D_2$ be one of the maps in the diagram~\eqref{diagram:descent_along_all_the_rings}. By definition, the base change functor from the category of $\varphi$-modules (resp. $\Gamma$-modules) over $D_1$ to that over $D_2$ is compatible with tensor and dual structures.  
\end{remark}

\subsection{Preliminary results} \label{subsec: preliminary-results-relative-(phi, Gamma)-modules}

We begin with some preliminary facts about $(\varphi, \Gamma)$-modules. The following lemma is already used in the proof of Corollary~\ref{cor: equiv-etale-(phi, Gamma)-mods-over-perfect-period-rings}. 

\begin{lemma} \label{lem: general-fact-equivalence-phi-Gamma-modules}
Let $D_1 \rightarrow D_2$ be one of the maps in the diagram~\eqref{diagram:descent_along_all_the_rings}. 
\begin{enumerate}
    \item Suppose this map induces an equivalence of categories of either $\varphi$-modules or of $\Gamma$-modules over $D_1$ and $D_2$. Then it also induces an equivalence of categories of $(\varphi, \Gamma)$-modules over $D_1$ and $D_2$. 
    \item 
Similarly, if this map induces a fully faithful functor from the category of $\varphi$-modules (resp. $\Gamma$-modules) over $D_1$ to the category of $\varphi$-modules (resp. $\Gamma$-modules) over $D_2$, then it induces a fully faithful functor of $(\varphi, \Gamma)$-modules between them.     
    \item 
    To show that the base change functor from the category of $\Gamma$-modules over $D_1$ to that over $D_2$ is an equivalence, it suffices to replace $\Gamma$ by an(y) open normal subgroup $\Gamma_0 \subset \Gamma$. 
\end{enumerate}

\end{lemma}

\begin{proof}
First, suppose that the map $D_1 \rightarrow D_2$ induces a fully faithful functor from the category of $\Gamma$-modules over $D_1$ to that over $D_2$. Let $M$ and $N$ be $(\varphi, \Gamma)$-modules over $D_1$. Write 
$L\coloneqq \Hom_{D_1}(M, N) \cong M^{\vee}\otimes_{D_1} N.$  
We have $(L\otimes_{D_1} D_2)^{\Gamma = 1} \subset L$ by the assumption, so $(L\otimes_{D_1} D_2)^{\Gamma = 1, \varphi = 1} \subset L$. Thus, the map also induces a fully faithful functor of $(\varphi, \Gamma)$-modules between them. 

Next, suppose that the map induces an equivalence of $\Gamma$-modules between them. Let $M_2$ be a $(\varphi, \Gamma)$-module over $D_2$. Then we have a $\Gamma$-module $M_1$ over $D_1$ such that $M_1\otimes_{D_1} D_2 \cong M_2$ is compatible with $\Gamma$-actions. Consider the isomorphism  
$
\varphi_{M_2}\colon M_2\otimes_{D_2, \varphi} D_2 \stackrel{\cong}{\longrightarrow} M_2
$ of $\Gamma$-modules over $D_2$. 
Since the isomorphism $M_2\otimes_{D_2, \varphi} D_2 \cong (M_1\otimes_{D_1, \varphi} D_1)\otimes_{D_1} D_2$ is compatible with $\Gamma$-actions, $\varphi_{M_2}$ descends to an isomorphism  
$
\varphi_{M_1}\colon M_1\otimes_{D_1, \varphi} D_1 \stackrel{\cong}{\longrightarrow} M_1 
$ over $D_1$. 
Thus, $M_1$ is a $(\varphi, \Gamma)$-module over $D_1$ such that $M_1\otimes_{D_1} D_2 \cong M_2$ is compatible with the $(\varphi, \Gamma)$-action. This implies that the map also induces an equivalence of $(\varphi, \Gamma)$-modules between them. 

The other statements in parts (1) and (2) follow from similar arguments. 

Finally, let us prove part (3). The full faithfulness is clear. It remains to show essential surjectivity. For this, suppose that $M_1$ is a finite projective $D_1$-submodule of $M_2$ stable under $\Gamma_0$ such that $M_2 \cong M_1 \otimes_{D_1} D_2$ as $\Gamma_0$-modules. We have to show that for any $g \in \Gamma$, $M_1$ is also stable under $g$. To this end, let 
\[
M_1 ' \coloneqq g \cdot M_1 = \{g \cdot x ~|~ x \in M_1\} \subset M_2.
\] 
Note that, since $\Gamma$ acts on $D_1$,  $ M_1' \subset M_2$ is a $D_1$-submodule, which is finite projective. Since $\Gamma_0$ is a normal subgroup of $\Gamma$, $M_1'$ is stable under the action of $\Gamma_0$. Furthermore, the map 
$
M_1' \otimes_{D_1} D_2 \lra M_2 \cong M_1 \otimes_{D_1} D_2
$  
given by $(g \cdot x) \otimes \alpha \mapsto \alpha (g \cdot x)$ is an isomorphism between finite projective $\Gamma_0$-modules over $D_2$. Since the base change functor from finite projective $\Gamma_0$-modules from $D_1$ to $D_2$ is fully faithful by assumption,  we conclude that the natural inclusion $M_1' \subset  M_2$ induces an identification $M_1' = M_1$ of  submodules as desired. 
\end{proof}

\begin{lemma} \label{lem: equivalence-phi-mod-base-change-to-breve}
The base change functor from the category of $\varphi$-modules over $\mathbf{A}_{\Lambda}^{\dagger}$ (resp. $\mathbf{B}_{\Lambda}^{\dagger}$, $\mathbf{C}_{\Lambda}^{\dagger}$) to the category of $\varphi$-modules over $\breve{\mathbf{A}}_{\Lambda}^{\dagger}$ (resp. $\breve{\mathbf{B}}_{\Lambda}^{\dagger}$, $\breve{\mathbf{C}}_{\Lambda}^{\dagger}$) is an equivalence.    
\end{lemma}

\begin{proof}
We only consider the case for $\mathbf{A}_{\Lambda}^{\dagger}$; the other cases follow from the same argument. We will construct a quasi-inverse to the base change functor. Let $M$ be a $\varphi$-module over $\breve{\mathbf{A}}_{\Lambda}^{\dagger}$. Then there exist $r > 0$ and $n \ge 0$ such that $M$ comes from a $\varphi$-module $N$ over $\varphi^{-n}(\mathbf{A}_{\Lambda}^{\frac{r}{p^n}})$, i.e. $N$ is a finite projective $\varphi^{-n}(\mathbf{A}_{\Lambda}^{\frac{r}{p^n}})$-module equipped with an isomorphism
\[
\varphi_N\colon N\otimes_{\varphi^{-n}(\mathbf{A}_{\Lambda}^{\frac{r}{p^n}}), \varphi} \varphi^{-n}(\mathbf{A}_{\Lambda}^{\frac{r}{p^{n+1}}}) \stackrel{\cong}{\longrightarrow} N\otimes_{\varphi^{-n}(\mathbf{A}_{\Lambda}^{\frac{r}{p^n}})} \varphi^{-n}(\mathbf{A}_{\Lambda}^{\frac{r}{p^{n+1}}}),
\]
and $M = N\otimes_{\varphi^{-n}(\mathbf{A}_{\Lambda}^{\frac{r}{p^n}})} \breve{\mathbf{A}}_{\Lambda}^{\dagger}$, compatible with the Frobenius structures. 
Let $N_0\coloneqq N\otimes_{\varphi^{-n}(\mathbf{A}_{\Lambda}^{\frac{r}{p^n}}), \varphi^n} \mathbf{A}_{\Lambda}^{\frac{r}{p^n}}$. Then $N_0$ is a finite projective module over $\mathbf{A}_{\Lambda}^{\frac{r}{p^n}}$, and $\varphi_N$ induces an isomorphism
\[
\varphi_{N_0}\colon N_0\otimes_{\mathbf{A}_{\Lambda}^{\frac{r}{p^n}}, \varphi} \mathbf{A}_{\Lambda}^{\frac{r}{p^{n+1}}} \stackrel{\cong}{\longrightarrow} N_0\otimes_{\mathbf{A}_{\Lambda}^{\frac{r}{p^n}}} \mathbf{A}_{\Lambda}^{\frac{r}{p^{n+1}}}.
\]
It is straightforward to check that a quasi-inverse is given by  $N_0\otimes_{\mathbf{A}_{\Lambda}^{\frac{r}{p^n}}} \mathbf{A}_{\Lambda}^{\dagger}$ with Frobenius induced from $\varphi_{N_0}$. 
\end{proof}

The following result is proved in \cite{Kedlaya-LiuI}.

\begin{proposition}
  [cf. {\cite[Lemma~6.1.5, Theorem~6.3.12]{Kedlaya-LiuI}}]
Let $r > 0$ and $0 < s \le \frac{r}{p}$. The category of $\varphi$-modules over $\widetilde{\mathbf{C}}_X^{\dagger}$ is naturally equivalent to the category of $\varphi$-modules over $\widetilde{\mathbf{C}}_{X}^{[s, r]}$.        
\end{proposition}

We have an analogous result for $\varphi$-modules over imperfect period rings of $\fC$-type. 

\begin{lemma} \label{lem: equiv-phi-module-base-change-r-to-[s, r]}
Let $r_0$ be as in Proposition~\ref{prop: good lift weakly decompleting}. Let $0 < s < r < r_0$ such that $s, r \in \mathbb{Q}$ and $s \le \frac{r}{p}$. Then the base change functor from the category of $\varphi$-modules over $\mathbf{C}_{\Lambda}^r$ to the category of $\varphi$-modules over $\mathbf{C}_{\Lambda}^{[s, r]}$ is an equivalence.     
\end{lemma}

\begin{proof}
We will construct a quasi-inverse to the base change functor. Let $(M_0, \varphi_{M_0})$ be a $\varphi$-module over $\mathbf{C}_{\Lambda}^{[s, r]}$. If we write $M_1 \coloneqq M_0\otimes_{\mathbf{C}_{\Lambda}^{[s, r]}, \varphi} \mathbf{C}_{\Lambda}^{[s/p, \,r/p]}$, then $\varphi_{M_0}$ is an isomorphism
\[
\varphi_{M_0}\colon M_1\otimes_{\mathbf{C}_{\Lambda}^{[s/p, \,r/p]}} \mathbf{C}_{\Lambda}^{[s, \frac{r}{p}]} \stackrel{\cong}{\longrightarrow} M_0\otimes_{\mathbf{C}_{\Lambda}^{[s, r]}} \mathbf{C}_{\Lambda}^{[s, \frac{r}{p}]}. 
\]
Let $M_2\coloneqq  M_0\otimes_{\mathbf{C}_{\Lambda}^{[s, r]}, \varphi^2} \mathbf{C}_{\Lambda}^{[\frac{s}{p^2}, \frac{r}{p^2}]}$. Then $\varphi_{M_0}$ above induces an isomorphism
\[
\varphi_{M_1}\colon M_2\otimes_{\mathbf{C}_{\Lambda}^{[\frac{s}{p^2}, \frac{r}{p^2}]}} \mathbf{C}_{\Lambda}^{[\frac{s}{p}, \frac{r}{p^2}]} \stackrel{\cong}{\longrightarrow} M_1\otimes_{\mathbf{C}_{\Lambda}^{[s/p, \,r/p]}} \mathbf{C}_{\Lambda}^{[\frac{s}{p}, \frac{r}{p^2}]}
\]
via pullback along the map $\varphi\colon \mathbf{C}_{\Lambda}^{[s, \frac{r}{p}]} \rightarrow \mathbf{C}_{\Lambda}^{[\frac{s}{p}, \frac{r}{p^2}]}$. Iterating this process, we obtain
\[
M_n\coloneqq  M_0\otimes_{\mathbf{C}_{\Lambda}^{[s, r]}, \varphi^n} \mathbf{C}_{\Lambda}^{[\frac{s}{p^n}, \frac{r}{p^n}]}
\]
for each $n \ge 0$ with an isomorphism
\[
\varphi_{M_n}\colon M_{n+1}\otimes_{\mathbf{C}_{\Lambda}^{[\frac{s}{p^{n+1}}, \frac{r}{p^{n+1}}]}} \mathbf{C}_{\Lambda}^{[\frac{s}{p^n}, \frac{r}{p^{n+1}}]} \stackrel{\cong}{\longrightarrow} M_n\otimes_{\mathbf{C}_{\Lambda}^{[\frac{s}{p^n}, \frac{r}{p^n}]}} \mathbf{C}_{\Lambda}^{[\frac{s}{p^n}, \frac{r}{p^{n+1}}]}.
\]
By the construction, the isomorphisms $\varphi_{M_n}$'s satisfy the cocycle condition. Thus, we have a finite projective $\mathbf{C}_{\Lambda}^r$-module $M$ given by Corollary~\ref{cor: fin-proj-vec-bun-Cr} such that $M\otimes_{\mathbf{C}_{\Lambda}^r} \mathbf{C}_{\Lambda}^{[\frac{s}{p^n}, \frac{r}{p^n}]} \cong M_n$ for each $n \ge 0$. Furthermore, $\varphi_{M_n}$'s induce an isomorphism
\[
\varphi_M\colon M\otimes_{\mathbf{C}_{\Lambda}^r, \varphi} \mathbf{C}_{\Lambda}^{\frac{r}{p}} \stackrel{\cong}{\longrightarrow} M\otimes_{\mathbf{C}_{\Lambda}^r, \iota} \mathbf{C}_{\Lambda}^{\frac{r}{p}}.
\]
It is straightforward to check that this gives a quasi-inverse.
\end{proof}

\begin{corollary} \label{cor: free-(phi, Gamma)-modules}
Let $r > 0$, and let $M$ be a $(\varphi, \Gamma)$-module over $\fC_{\Lambda}^r$ such that $M\otimes_{\fC_{\Lambda}^r} \widetilde{\fC}_X^{\dagger}$ is free over $\widetilde{\fC}_X^{\dagger}$. Then there exists $0 < r' \le r$ such that for any $0 < s \le r'$, the induced $(\varphi, \Gamma)$-module $M\otimes_{\fC_{\Lambda}^r} \fC_{\Lambda}^{[s, r']}$ is free over $\fC_{\Lambda}^{[s, r']}$.
\end{corollary}

\begin{proof}
Choose $0 < r_1 \le r$ sufficiently small so that $M\otimes_{\fC_{\Lambda}^r} \widetilde{\fC}_X^{r_1}$ is free over $\widetilde{\fC}_X^{r_1}$. Let $s_1 = \frac{r_1}{p}$, and consider the induced $(\varphi, \Gamma)$-module
\[
\widetilde{M}^{[s_1, r_1]} \coloneqq (M\otimes_{\fC_{\Lambda}^r} \widetilde{\fC}_X^{r_1})\otimes_{\widetilde{\fC}_X^{r_1}} \widetilde{\fC}_X^{[s_1, r_1]}.
\]

Write $\lambda_{s_1, r_1}(\cdot) \coloneqq \max\{\lambda_{s_1}(\cdot), \lambda_{r_1}(\cdot)\}$. Let $x_1, \ldots, x_{\ell} \in M$ be a set of generators of $M$ as a finite $\fC_{\Lambda}^r$-module. Choose a $\widetilde{\fC}_X^{[s_1, r_1]}$-basis $\{e_1, \ldots, e_n\}$ of $\widetilde{M}^{[s_1, r_1]}$. Equip $\widetilde{M}^{[s_1, r_1]}$ with the finite-module norm given by the presentation
\[
\widetilde{M}^{[s_1, r_1]} \cong \bigoplus_{i=1}^n \widetilde{\fC}_X^{[s_1, r_1]}\cdot e_i.
\]
For each $1 \le i \le n$, we can write $e_i = \sum_{j=1}^{\ell} a_{ij} x_j$ for some $a_{ij} \in \widetilde{\fC}_X^{[s_1, r_1]}$. Recall from Lemma~\ref{lem: C breve dense in C tilde}(2) that $\breve{\fC}_{\Lambda}^{[s_1, r_1]}$ is dense in $\widetilde{\fC}_X^{[s_1, r_1]}$. So we can choose an integer $m \ge 0$ sufficiently large such that for each $1 \le i \le n$ and $1 \le j \le \ell$, there exists $b_{ij} \in \fC_{\Lambda, (m)}^{[s_1, r_1]}$ satisfying $\lambda_{s_1, r_1}(e_i-f_i) \le \frac{1}{2}$ with $f_i\coloneqq \sum_{j=1}^{\ell} b_{ij} x_j$. Since $\widetilde{\fC}_X^{[s_1, r_1]}$ is complete under the $\lambda_{s_1, r_1}$-norm, $\{f_1, \ldots, f_n\}$ gives another basis of $\widetilde{M}^{[s_1, r_1]}$. We have $f_1, \ldots, f_n \in M^{[s_1, r_1]}_{(m)}\coloneqq M\otimes_{\fC_{\Lambda}^r} \fC_{\Lambda, (m)}^{[s_1, r_1]}$.

Note that $\fC_{\Lambda, (k)}^{[s_1, r_1]}$ is complete under the $\lambda_{s_1, r_1}$-norm for each $k \ge 0$. Using a similar approximation argument as above, we deduce that 
$
x_j \in \bigoplus_{i=1}^n \fC_{\Lambda, (m)}^{[s_1, r_1]}\cdot f_i$ for each $1 \le j \le \ell$  (after enlarging $m$ if necessary). Thus, $M^{[s_1, r_1]}_{(m)}$
is free over $\fC_{\Lambda, (m)}^{[s_1, r_1]}$. 

Lastly, let $r' = \frac{r_1}{p^m}$ and $s' = \frac{s_1}{p^m}$. Then the induced $(\varphi, \Gamma)$-module
\[
M^{[s', r']}\coloneqq M^{[s_1, r_1]}_{(m)}\otimes_{\fC_{\Lambda, (m)}^{[s_1, r_1]}, \varphi^m} \fC_{\Lambda}^{[s', r']}
\]
is free over $\fC_{\Lambda}^{[s', r']}$. Since $s' = \frac{r'}{p}$, this implies the desired statement by Lemma~\ref{lem: equiv-phi-module-base-change-r-to-[s, r]}.
\end{proof}

\subsection{Descending relative $\varphi$-modules} \label{subsec: descent-relative-phi-modules}

We first study some descent results for $\varphi$-modules. We will refer to the arrows in the diagram~\eqref{diagram:descent_along_all_the_rings}. Recall that Proposition~\ref{prop:descending_phi_mod_along_4_and_10} induces descent of \'etale $\varphi$-modules along the arrows (4), (8), and (10).

\subsubsection{Descent of $\varphi$-modules along the arrows (2) and (6)} 

Recall from Definition~\ref{defn: weakly decompleting datum}(5) that $\fA_{\Lambda}/p \cong R_{\Lambda}$.  The essence of the next two results is the Riemann--Hilbert correspondence between $\F_p$-local systems in characteristic $p$ and vector bundles with Frobenius structures, due to Katz \cite[Proposition~4.1.1]{Katz} (see also \cite[Proposition~3.2.7]{Kedlaya-LiuI} and \cite[Proposition~3.4]{BS_crystal}).

\begin{lemma} \label{lemma:KL2_5.4.5}
Let $M$ be a finite projective $\varphi$-module over $\mathbf{A}_{\Lambda}$. Let $m \ge 1$ be a positive integer.
\begin{enumerate}
    \item There exists a faithfully finite \'etale $\mathbf{A}_{\Lambda}/p^m$-algebra $S_m$ such that $M \otimes_{\mathbf{A}_{\Lambda}} S_m$ is finite free over $S_m$ and admits a $\varphi$-fixed basis.
    
    \item If $m_0 < m$ is another integer such that $M \otimes_{\mathbf{A}_{\Lambda}} S_{m_0}$ admits a $\varphi$-fixed basis for some faithfully finite \'etale $\mathbf{A}_{\Lambda}/p^{m_0}$-algebra $S_{m_0}$, then $S_m$ can be chosen in a way such that $S_{m}/p^{m_0}$ is faithfully finite \'etale over $S_{m_0}$, and that this basis over $S_{m_0}$ lifts to a $\varphi$-fixed basis of  $M \otimes_{\mathbf{A}_{\Lambda}} S_m$. 
    \end{enumerate}
\end{lemma}

\begin{proof} The same proof of \cite[Lemma~3.2.6]{Kedlaya-LiuI} carries over.    
\end{proof}

\begin{corollary} \label{cor:Katz_RH_mod_p^m}
Let $m \ge 1$ be a positive integer. There is a natural tensor equivalence between the category of $\Z/p^m$-local systems on $\spec R_{\Lambda}$ and finite projective $\varphi$-modules over $\mathbf{A}_{\Lambda}/p^m$. 
\end{corollary}

\begin{proof}
The same proof of \cite[Proposition~3.2.7]{Kedlaya-LiuI} carries over. 
\end{proof}

The following is the main result of this subsection.

\begin{proposition}[cf. {\cite[Corollary~5.4.6]{Kedlaya-LiuII}}] \label{prop:descending_phi_mod_along_2} 
The base change functor along $\mathbf{A}_{\Lambda} \ra \widetilde{\mathbf{A}}_{X}$ induces an equivalence of categories of $\varphi$-modules 
\[ 
\mathrm{Mod}_{/\mathbf{A}_{\Lambda}}^{\varphi} \isom \mathrm{Mod}_{/\widetilde{\mathbf{A}}_X}^{\varphi}.
\]
Consequently, the base change functor along $\mathbf{B}_{\Lambda} \rightarrow \widetilde{\mathbf{B}}_{X}$ induces an equivalence of categories of \'etale $\varphi$-modules
\[
\mathrm{Mod}_{/\mathbf{B}_{\Lambda}}^{\varphi, \ett} \isom \mathrm{Mod}_{/\widetilde{\mathbf{B}}_X}^{\varphi, \ett}.
\]
\end{proposition}

\begin{proof}  
This is essentially  \cite[Corollary~5.4.6]{Kedlaya-LiuII}. We recall the argument for completeness. 
We first show full faithfulness. Since $\Hom_{\varphi} (M_1, M_2) \cong (M_1^\vee \otimes M_2)^{\varphi = 1}$ for (finite projective) $\varphi$-modules $M_1$, $M_2$ over $\mathbf{A}_{\Lambda}$, it suffices to show that for any finite projective $\varphi$-module $M \in \mathrm{Mod}_{/\mathbf{A}_{\Lambda}}^{\varphi}$, we have 
$
(M \otimes_{\mathbf{A}_{\Lambda}} \widetilde{\mathbf{A}}_X)^{\varphi = 1} \subset M.
$  
Note that since $\fA_{\Lambda}$ and $\widetilde{\fA}_X$ are $p$-adically complete and $\widetilde{\fA}_X$ is $p$-torsion free, if any element $x \in \widetilde{\fA}_X$ satisfies the condition that its image in $\widetilde{\fA}_X/p^m$ lies in $\fA_{\Lambda}/p^m$ for every $m \ge 1$, then $x \in \fA_{\Lambda}$. Since $M$ is finite projective over $\mathbf{A}_{\Lambda}$, we further reduce to showing that for each integer $m \ge 1$, we have 
\[ 
(M \otimes_{\mathbf{A}_{\Lambda}} \mathbf{A}_{\Lambda}/p^m)^{\varphi = 1} = (M \otimes_{A}  \widetilde{\mathbf{A}}_{X}/p^m )^{\varphi = 1}.
\]

Let $v \in (M \otimes_{\mathbf{A}_{\Lambda}} \widetilde{\mathbf{A}}_{X}/p^m )^{\varphi = 1}$. Let $S_m$ be a finite faithfully \'etale algebra over $\mathbf{A}_{\Lambda}/p^m$ over which $M$ admits a $\varphi$-fixed basis as in Lemma~\ref{lemma:KL2_5.4.5}. Write 
\[ 
\sq S_m \coloneqq S_m \otimes_{\mathbf{A}_{\Lambda}} \widetilde{\mathbf{A}}_{X} = S_m \otimes_{\mathbf{A}_{\Lambda}/p^m}  \widetilde{\mathbf{A}}_{X}/p^m.
\]  
View $v$ as an element in $M \otimes_{\mathbf{A}_{\Lambda}} \widetilde{S}_m$, which has a $\varphi$-fixed basis given by a $\varphi$-fixed basis in $M \otimes_{\mathbf{A}_{\Lambda}} S_m$. By writing $v$ in terms of this basis and noting that $(\sq S_m)^{\varphi = 1} = (S_m)^{\varphi = 1}$, we deduce $v \in (M \otimes_{A} S_m)^{\varphi =1}$. Thus,
\[
v \in (M \otimes_{\mathbf{A}_{\Lambda}} S_m) \cap  (M \otimes_{\mathbf{A}_{\Lambda}}  \widetilde{\mathbf{A}}_{X}/p^m ).
\] 
We claim that 
\[  
(M \otimes_{\mathbf{A}_{\Lambda}} S_m) \cap  (M \otimes_{\mathbf{A}_{\Lambda}}  \widetilde{\mathbf{A}}_{X}/p^m ) = M \otimes_{\mathbf{A}_{\Lambda}} (S_m  \cap   \widetilde{\mathbf{A}}_{X}/p^m ) = M/p^m.
\]
The first equality above follows from that $M$ is flat over $\mathbf{A}_{\Lambda}$. Since $S_m$ is faithfully flat over $\mathbf{A}_{\Lambda}/p^m$, we have $S_m \cap \widetilde{\mathbf{A}}_{X}/p^m = \mathbf{A}_{\Lambda}/p^m$ as subrings of $\sq S_m$, which yields the second equality. This proves the desired full faithfulness. 

Now we show essential surjectivity. Let $\sq M$ be a finite projective $\varphi$-module over $\widetilde{\mathbf{A}}_{X}$ of rank $d$. By the modulo-$p^m$ full faithfulness proved above, it suffices to descend (for each $m \ge 1$) the finite projective $\varphi$-module $\sq M_m \coloneqq \sq M /p^m$ over $\widetilde{\mathbf{A}}_{X}/p^m$ to a finite projective $\varphi$-module $M_m$ over $\mathbf{A}_{\Lambda}/p^m$. Note that the $\varphi$-module $\sq M$ corresponds to a $\Z_p$-local system $\widetilde{\T}$ on $\spec \widetilde{R}_{X}$ of rank $d$ by \cite[Theorem~8.5.3]{Kedlaya-LiuI}. Since the colimit perfection ${R}_{\Lambda}^{\mathrm{perf}}$ is dense in $\widetilde{R}_{X}$ (see Definition~\ref{defn: weakly decompleting datum}(2)), $\widetilde{\T}$ comes from a $\Z_p$-local system $\T$ on $\spec R_{\Lambda}$ again by \cite[Theorem~8.5.3]{Kedlaya-LiuI}. By Corollary~\ref{cor:Katz_RH_mod_p^m}, the $\Z/p^m$-local system $\T/p^m$ corresponds to a finite projective $\varphi$-module $M_m$ over $\fA_{\Lambda}/p^m$, and it follows from the constructions that $M_m\otimes_{\mathbf{A}_{\Lambda}} \widetilde{\mathbf{A}}_{X} \cong \sq M_m$ is compatible with $\varphi$. 
\end{proof}

\subsubsection{full faithfulness along the arrows (1), (3), (5), (7), (9), and (11)}

\begin{proposition}[cf. {\cite[Proposition~5.4.8]{Kedlaya-LiuII}}]  \label{prop:descending_phi_mod_along_3_and_9} \hfill
\begin{enumerate}
\item The base change functor of $\varphi$-modules 
\[ 
\mathrm{Mod}_{/\mathbf{A}_{\Lambda}^{\dagger} }^{\varphi} \lra \mathrm{Mod}_{/ \mathbf{A}_{\Lambda}}^{\varphi}.
\] 
along the map $\mathbf{A}_{\Lambda}^{\dagger} \ra \mathbf{A}_{\Lambda}$ is fully faithful. Consequently, the base change functor of \'etale $\varphi$-modules
\[
\mathrm{Mod}_{/\mathbf{B}_{\Lambda}^{\dagger}}^{\varphi, \ett} \lra \mathrm{Mod}_{/ \mathbf{B}_{\Lambda}}^{\varphi, \ett}
\]
along $\mathbf{B}_{\Lambda}^{\dagger} \rightarrow \mathbf{B}_{\Lambda}$ is fully faithful.
    
\item The base change functor of \'etale $\varphi$-modules 
\[ 
\mathrm{Mod}_{/\mathbf{B}_{\Lambda}^{\dagger} }^{\varphi, \ett} \lra \mathrm{Mod}_{/ \mathbf{C}_{\Lambda}^{\dagger} }^{\varphi, \ett}.
\] 
along the map $\mathbf{B}_{\Lambda}^{\dagger} \ra \mathbf{C}_{\Lambda}^\dagger$ is an equivalence of categories. 

\item The base change functor of $\varphi$-modules
\[
\mathrm{Mod}_{/\mathbf{A}_{\Lambda}^{\dagger}}^{\varphi} \lra \mathrm{Mod}_{/ \widetilde{\mathbf{A}}_X^{\dagger}}^{\varphi}
\]
along the map $\mathbf{A}_{\Lambda}^{\dagger} \rightarrow \widetilde{\mathbf{A}}_X^{\dagger}$ is fully faithful. Consequently, the base change functor of \'etale $\varphi$-modules 
\[
\mathrm{Mod}_{/\mathbf{B}_{\Lambda}^{\dagger}}^{\varphi, \ett} \lra \mathrm{Mod}_{/ \widetilde{\mathbf{B}}_X^{\dagger}}^{\varphi, \ett}
\]
along $\mathbf{B}_{\Lambda}^{\dagger} \rightarrow \widetilde{\mathbf{B}}_X^{\dagger}$ is fully faithful. Furthermore, the base change functor of \'etale $\varphi$-modules \[
\mathrm{Mod}_{/\mathbf{C}_{\Lambda}^{\dagger}}^{\varphi, \ett} \lra \mathrm{Mod}_{/ \widetilde{\mathbf{C}}_X^{\dagger}}^{\varphi, \ett}
\]
along $\mathbf{C}_{\Lambda}^{\dagger} \rightarrow \widetilde{\mathbf{C}}_X^{\dagger}$ is also fully faithful.

\end{enumerate}
\end{proposition}

\begin{proof}
The proof is similar to that of full faithfulness in Proposition~\ref{prop:descending_phi_mod_along_2}. 
\begin{enumerate}
\item It suffices to show that for any finite projective $\varphi$-module $M^{\dagger}$ over $\mathbf{A}^{\dagger}_{\Lambda}$, every element $v \in (M^\dagger \otimes_{\mathbf{A}_{\Lambda}^{\dagger}} \mathbf{A}_{\Lambda})^{\varphi = 1}$ lies in $M^{\dagger}$. By Propositions~\ref{prop:descending_phi_mod_along_2} and \ref{prop:descending_phi_mod_along_4_and_10}, we have 
\[
(M^\dagger \otimes_{\mathbf{A}_{\Lambda}^{\dagger}} \mathbf{A}_{\Lambda})^{\varphi = 1} = (M^\dagger \otimes_{\mathbf{A}_{\Lambda}^{\dagger}}  \widetilde{\mathbf{A}}_{X} )^{\varphi = 1} = (M^\dagger \otimes_{\mathbf{A}^{\dagger}_{\Lambda}} \widetilde{\mathbf{A}}_{X}^{\dagger})^{\varphi = 1}.
\]
Hence, $v$ lies in the intersection 
\[ 
(M^\dagger \otimes_{\mathbf{A}_{\Lambda}^{\dagger}} \mathbf{A}_{\Lambda}) \cap (M^\dagger \otimes_{\mathbf{A}_{\Lambda}^{\dagger}} \widetilde{\mathbf{A}}_{X}^{\dagger}) = M^\dagger \otimes_{\mathbf{A}_{\Lambda}^{\dagger}} ( \mathbf{A}_{\Lambda} \cap \widetilde{\mathbf{A}}_{X}^{\dagger})=M^{\dagger},
\]
as desired, where the last equality follows from $\mathbf{A}_{\Lambda}^{\dagger} = \mathbf{A}_{\Lambda} \cap \widetilde{\mathbf{A}}_{X}^{\dagger}$ (see Definition~\ref{def: imperfect-period-rings-B-C}(2)). 

\item It suffices to prove full faithfulness. Again, this amounts to showing that for any finite projective \'etale $\varphi$-module $M^{\dagger}$ over $\mathbf{B}^{\dagger}_{\Lambda}$, every element $v \in (M^\dagger \otimes_{\mathbf{B}_{\Lambda}^{\dagger}} \mathbf{C}_{\Lambda}^{\dagger})^{\varphi = 1}$ lies in $M^{\dagger}$. Via the inclusion $\mathbf{C}_{\Lambda}^{\dagger} \ra \widetilde{\mathbf{C}}_{X}^{\dagger}$, we view $v$ as an element in $M^\dagger \otimes_{\mathbf{B}_{\Lambda}^{\dagger}} \widetilde{\mathbf{C}}_{X}^{\dagger}$, which is an \'etale $\varphi$-module over $\widetilde{\mathbf{C}}_{X}^{\dagger}$. By Proposition~\ref{prop:descending_phi_mod_along_4_and_10}(2), we have
\[
v \in (M^\dagger \otimes_{\mathbf{B}_{\Lambda}^{\dagger}} \widetilde{\mathbf{C}}_{X}^{\dagger})^{\varphi = 1} =  (M^\dagger \otimes_{\mathbf{B}_{\Lambda}^{\dagger}} \widetilde{\mathbf{B}}_{X}^{\dagger})^{\varphi = 1}.
\]
Thus,
\[ 
v \in (M^\dagger \otimes_{\mathbf{B}_{\Lambda}^{\dagger}} \widetilde{\mathbf{B}}_{X}^{\dagger})  \cap (M^\dagger \otimes_{\mathbf{B}_{\Lambda}^{\dagger}}\mathbf{C}^{\dagger}_{\Lambda}) = 
 M^\dagger \otimes_{\mathbf{B}_{\Lambda}^{\dagger}}  (\widetilde{\mathbf{B}}_{X}^{\dagger} \cap \fC_{\Lambda}^\dagger)  = M^\dagger
\]
where the last equality follows from Proposition~\ref{prop: intersection A}.

\item The statement for $\fA_{\Lambda}^{\dagger} \rightarrow \widetilde{\fA}_X^{\dagger}$ (and therefore that for $\fB_{\Lambda}^{\dagger} \rightarrow \widetilde{\fB}_X^{\dagger}$) is a direct consequence of part (1), Propositions~\ref{prop:descending_phi_mod_along_4_and_10}(1) and \ref{prop:descending_phi_mod_along_2}. As a consequence, the statement for $\fC_{\Lambda}^{\dagger} \rightarrow \widetilde{\fC}_X^{\dagger}$ follows from Proposition~\ref{prop:descending_phi_mod_along_4_and_10}(2).
\end{enumerate}
\end{proof}

\begin{remark}\label{remark: descent of phi modules implies descent of phi Gamma modules}
Thanks to Lemma~\ref{lem: general-fact-equivalence-phi-Gamma-modules}, the statements in Propositions~\ref{prop:descending_phi_mod_along_2} and \ref{prop:descending_phi_mod_along_3_and_9} remain true if we replace $\varphi$-modules by $(\varphi, \Gamma)$-modules.
\end{remark}

\subsection{Descending relative $(\varphi, \Gamma)$-modules, I: full faithfulness along the arrow (11)} \label{subsec: fully-faithfulness-indfree-(phi, Gamma)-modules}

Notice that we have not analyzed the base change functor along the map $\mathbf{C}_{\Lambda}^{\dagger} \rightarrow \widetilde{\mathbf{C}}_X^{\dagger}$ (arrow (11)) in the non-\'etale situation. To state our result on full faithfulness along the arrow (11), we consider the full subcategory
\[\mathrm{Mod}^{(\varphi, \Gamma), \mathrm{ind-free}}_{/\mathbf{C}_{\Lambda}^{\dagger}} \subset \mathrm{Mod}_{/\fC_{\Lambda}^{\dagger}}^{(\varphi, \Gamma)}\]
consisting of $(\varphi, \Gamma)$-modules $M$ over $\fC_{\Lambda}^{\dagger}$ such that $M\otimes_{\fC_{\Lambda}^{\dagger}} \widetilde{\fC}_X^{\dagger}$ is free over $\widetilde{\fC}_X^{\dagger}$.
 
\begin{theorem} \label{thm: fully-faithfulness-for-free-modules}
Let $\Lambda = (R_{\Lambda}, \mathbf{A}_{\Lambda})$ be a decompleting datum for $(X, h)$. The functor
\[
\mathrm{Mod}^{(\varphi, \Gamma), \mathrm{ind-free}}_{/\mathbf{C}_{\Lambda}^{\dagger}} \rightarrow \mathrm{Mod}^{(\varphi, \Gamma)}_{/\widetilde{\mathbf{C}}_X^{\dagger}}
\]
given by the base change along $\fC_{\Lambda}^{\dagger} \rightarrow \widetilde{\mathbf{C}}_X^{\dagger}$ is fully faithful.
\end{theorem}

\begin{proof}
Let $M_1, M_2 \in \mathrm{Mod}^{(\varphi, \Gamma), \mathrm{ind-free}}_{/\mathbf{C}_{\Lambda}^{\dagger}}$. Let $r > 0$ be small enough such that $M_1$ and $M_2$ come from $(\varphi, \Gamma)$-modules over $\mathbf{C}_{\Lambda}^r$, which we still denote by $M_1$ and $M_2$ respectively, by an abuse of notation. Let $N \coloneqq \Hom_{\mathbf{C}_{\Lambda}^r} (M_1, M_2) \cong M_1^{\vee}\otimes_{\mathbf{C}_{\Lambda}^r} M_2$, which is a $(\varphi, \Gamma)$-module over $\fC_{\Lambda}^r$. Note that $N\otimes_{\fC_{\Lambda}^r} \widetilde{\fC}_X^{\dagger}$ is free over $\widetilde{\fC}_X^{\dagger}$. By Corollary~\ref{cor: free-(phi, Gamma)-modules}, after shrinking $r$ if necessary, we know that $N\otimes_{\fC_{\Lambda}^r} \fC_{\Lambda}^{[s, r]}$ is free over $\fC_{\Lambda}^{[s, r]}$ for all $0 < s \le r$. \\

\noindent \textbf{Step 1.} We first consider the induced finite free $(\varphi, \Gamma_0)$-module over $\fC_{\Lambda}^{[\frac{r'}{p}, r']}$ for some $r' > 0$ and some open normal subgroup $\Gamma_0 \subset \Gamma$.

Let $s = \frac{r}{p}$. Consider the induced $(\varphi, \Gamma)$-module $N\otimes_{\fC_{\Lambda}^r} \fC_{\Lambda}^{[s, r]}$, which is finite free over $\fC_{\Lambda}^{[s, r]}$. Choose a $\fC_{\Lambda}^{[s, r]}$-basis $e_1, \ldots, e_n$, and equip $N\otimes_{\fC_{\Lambda}^r} \fC_{\Lambda}^{[s, r]}$ with the finite-module norm given by the presentation
\[
N\otimes_{\fC_{\Lambda}^r} \fC_{\Lambda}^{[s, r]} \cong \bigoplus_{j=1}^n \fC_{\Lambda}^{[s, r]}\cdot e_j.
\]
Let $\Gamma_0 \subset \Gamma$ be an open normal subgroup of the form $\Gamma_0 \cong (p^a \Z_p)^d \rtimes (1+p^b\Z_p) \subset \Z_p^d \rtimes \Z_p^{\times}$ for some $a \ge 0$ and $b \ge 1$. We take $\Gamma_0$ small enough so that
\[
\lambda_t((\gamma-1)e_j) \le \frac{1}{2}
\]
for all $\gamma \in \Gamma_0$, $t \in [s, r]$, and $1 \le j \le n$.

As in the proof of Lemma~\ref{lem: equiv-phi-module-base-change-r-to-[s, r]}, for each integer $k \ge 0$, consider the induced $(\varphi, \Gamma_0)$-module
\[
N_k \coloneqq (N\otimes_{\mathbf{C}_{\Lambda}^r} \fC_{\Lambda}^{[s, r]})\otimes_{\fC_{\Lambda}^{[s, r]}, \varphi^k} \fC_{\Lambda}^{[sp^{-k}, rp^{-k}]}.
\]
Let $k_0 \ge 0$ be large enough so that Proposition~\ref{prop: C-complex-strict-exact} (with $\Gamma_0$ given above) is applicable for $[s', r']\coloneqq [sp^{-k_0}, rp^{-k_0}]$. Write $e_{j, k_0} \coloneqq e_j \otimes 1 \in N_{k_0}$. Note that for any $t' \in [s', r']$, we have
\[
\lambda_{t'}((\gamma-1)e_{j, k_0}) = \lambda_{p^{k_0}t'}((\gamma-1)e_j) \le \frac{1}{2}
\]
for each $\gamma \in \Gamma_0$, since $p^{k_0}t' \in [s, r]$.\\

\noindent \textbf{Step 2.} We now claim that
\[
\left(N_{k_0}\otimes_{\fC_{\Lambda}^{[s', r']}} \widetilde{\fC}_X^{[s', r']}\right)^{\Gamma_0 = 1} \subset N_{k_0}.
\]

Consider the quotient $\widetilde{\fC}_X^{[s', r']} / \fC_{\Lambda}^{[s', r']}$ and consider an element \[x = a_1e_{1, k_0}+\cdots +a_n e_{n, k_0} \in N_{k_0}\otimes_{\fC_{\Lambda}^{[s', r']}} \left(\widetilde{\fC}_X^{[s', r']} / \fC_{\Lambda}^{[s', r']}\right),\] with $a_j \in \widetilde{\fC}_X^{[s', r']} / \fC_{\Lambda}^{[s', r']}$, such that $\gamma(x) = x$ for all $\gamma \in \Gamma_0$. Write $\lambda_{s', r'}(\cdot) \coloneqq \max\{\lambda_s'(\cdot), \lambda_r'(\cdot)\}$. Pick an index $1 \le j_0 \le n$ such that $\lambda_{s', r'}(a_{j_0}) = \max_{1 \le j \le n}\{\lambda_{s', r'}(a_j) \}$. From the equalities
\[
\gamma(x) = \gamma(a_1)\gamma(e_{1, k_0})+\cdots + \gamma(a_n)\gamma(e_{n, k_0}) = x = a_1e_{1, k_0}+\cdots +a_n e_{n, k_0},
\]
we deduce that 
\[
\lambda_{s', r'}((\gamma-1)(a_{j_0})) \le \frac{1}{2}\lambda_{s', r'}(a_{j_0})
\]
for all $\gamma \in \Gamma_0$, since $\lambda_{s', r'}((\gamma-1)e_{j, k_0}) \le \frac{1}{2}$. Thus, 
\[
\lambda_{s', r'}(d^{(0)}(a_{j_0})) \le \frac{1}{2}\lambda_{s', r'}(a_{j_0}) \le \frac{3}{4}\lambda_{s', r'}(d^{(0)}(a_{j_0}))
\]
where the last inequality follows from Proposition~\ref{prop: C-complex-strict-exact}. This implies $d^{(0)}(a_{j_0}) = 0$, hence $a_{j_0} = 0$, again by Proposition~\ref{prop: C-complex-strict-exact}. Therefore $a_j = 0$ for all $1 \le j \le d$; namely, $x = 0$. This finishes the proof of the claim.\\

Finally, let $x \in (N\otimes_{\fC_{\Lambda}^r, \iota} \fC_{\Lambda}^{r'})\otimes_{\fC_{\Lambda}^{r'}} \widetilde{\fC}_X^{r'}$ such that $\varphi(x) = x$ and $\gamma(x) = x$ for all $\gamma \in \Gamma$. By Step 2 and Lemma~\ref{lem: equiv-phi-module-base-change-r-to-[s, r]} (note that $s' = \frac{r'}{p}$), we must have $x \in N\otimes_{\fC_{\Lambda}^r, \iota} \fC_{\Lambda}^{r'}$. This proves the desired full faithfulness.
\end{proof}

\subsection{Descending relative $(\varphi, \Gamma)$-modules, I\!I: the arrows (1), (5), (11)} \label{subsec: descent-relative-etale-(phi, Gamma)-modules} 

Recall from Proposition~\ref{prop:descending_phi_mod_along_3_and_9} and Remark~\ref{remark: descent of phi modules implies descent of phi Gamma modules} that the base change functor
\[
\mathrm{Mod}_{/\mathbf{A}_{\Lambda}^{\dagger}}^{(\varphi, \Gamma)} \lra \mathrm{Mod}_{/ \widetilde{\mathbf{A}}_X^{\dagger}}^{(\varphi, \Gamma)}
\]
is fully faithful. We show that this functor is actually an equivalence of categories.

\begin{theorem} \label{thm: decompletion-etale-(phi, Gamma)-modules}
The base change functor of $(\varphi, \Gamma)$-modules
\[
\mathrm{Mod}_{/\mathbf{A}_{\Lambda}^{\dagger}}^{(\varphi, \Gamma)} \lra \mathrm{Mod}_{/ \widetilde{\mathbf{A}}_X^{\dagger}}^{(\varphi, \Gamma)}
\]
along the map $\fA_{\Lambda}^{\dagger} \rightarrow \widetilde{\fA}_X^{\dagger}$ is an equivalence.     
\end{theorem}

\begin{proof}
The base change functor is fully faithful by Proposition~\ref{prop:descending_phi_mod_along_3_and_9}(3) and Lemma~\ref{lem: general-fact-equivalence-phi-Gamma-modules}(2). For essential surjectivity, let $\widetilde{M} \in \mathrm{Mod}_{/ \widetilde{\mathbf{A}}_X^{\dagger}}^{(\varphi, \Gamma)}$.\\

\noindent \textbf{Step 1.} We first show that when $\widetilde{M}$ is free over $\widetilde{\fA}_X^{\dagger}$, there exists a (unique) $\fA_{\Lambda}^{\dagger}$-submodule $M \subset \widetilde{M}$ which is a $(\varphi, \Gamma)$-module over $\fA_{\Lambda}^{\dagger}$ such that the natural map
\[
M\otimes_{\fA_{\Lambda}^{\dagger}} \widetilde{\fA}_X^{\dagger} \longrightarrow \widetilde{M}
\]
is an isomorphism of $(\varphi, \Gamma)$-modules.

The uniqueness holds by the full faithfulness of the base change functor, and we need to show the existence. Let $r > 0$ be sufficiently small such that we have 
$
\widetilde{M}^r\otimes_{\widetilde{\fA}_X^r} \widetilde{\fA}_X^{\dagger} \cong \widetilde{M}
$ 
as $(\varphi, \Gamma)$-modules, for some finite free $(\varphi, \Gamma)$-module $\widetilde{M}^r$ over $\widetilde{\fA}_X^r$. It suffices to show that $\widetilde{M}^r$ descends to a $(\varphi, \Gamma)$-module over $\fA_{\Lambda}^r$, after shrinking $r$ if necessary.

Choose an $\widetilde{\fA}_X^r$-basis $e_1, \ldots e_n$ for $\widetilde{M}^r$, and equip $\widetilde{M}^r$ with the finite-module norm given by the presentation
\[
\widetilde{M}^r \cong \bigoplus_{j=1}^n \widetilde{\fA}_X^r\cdot e_j.
\]
Let $\Gamma_0 \subset \Gamma$ be an open normal subgroup of finite index of the form $\Gamma_0 \cong (p^a\Z_p)^d \rtimes (1+p^b \Z_p) \subset \Z_p^d \rtimes \Z_p^{\times}$ for some integers $a 
\ge 0$ and $b \ge 1$. We choose $\Gamma_0$ small enough such that
\[
\lambda_r((\gamma-1)e_j) \le \frac{1}{6}
\]
for all $\gamma \in \Gamma_0$ and $1 \le j \le n$. For each integer $k \ge 0$, consider the induced $(\varphi, \Gamma_0)$-module
\[
\widetilde{M}^r_k\coloneqq \widetilde{M}^r\otimes_{\widetilde{\fA}_X^r, \varphi^k} \widetilde{\fA}_X^{rp^{-k}}.
\]
Let $r_1$ be as in Corollary~\ref{cor: complex-A-strict-exact} with respect to our choice of $\Gamma_0$. Pick $k_0$ large enough so that $r'\coloneqq rp^{-k_0} < r_1$. For each $e_{j, k_0}\coloneqq e_j\otimes 1 \in \widetilde{M}^r_{k_0}$, note that
\[
\lambda_{r'}((\gamma-1)e_{j, k_0}) = \lambda_r((\gamma-1)e_j) \le \frac{1}{6}
\]
for all $\gamma \in \Gamma_0$.

Before proceeding, we introduce some notation for convenience. For a matrix $T = (T_{jh})_{1 \le j, h \le n} \in \mathrm{Mat}_{n\times n}(\widetilde{\fA}_X^{r'})$, write $\overline{T} = (\overline{T}_{jh})$ for the $n \times n$ matrix with entries in $\widetilde{\fA}_X^{r'} / \fA_{\Lambda}^{r'}$, where $\overline{T}_{jh}$ being the image of $T_{jh}$ in $\widetilde{\fA}_X^{r'} / \fA_{\Lambda}^{r'}$ under the projection. Write $\lambda_{r'}(\overline{T})\coloneqq \max_{1 \le j, h \le n}\{\lambda_{r'}(\overline{T}_{jh})\}$. We use a similar notation for group cochains as well.\\

\noindent \textbf{Claim.} The $\Gamma_0$-module $\widetilde{M}_{k_0}^r$ descends to a $\Gamma_0$-stable free $\fA_{\Lambda}^{r'}$-submodule $M$, namely, the natural map
\[
M\otimes_{\fA_{\Lambda}^{r'}} \widetilde{\fA}_X^{r'} \longrightarrow \widetilde{M}_{k_0}^r
\]
is a $\Gamma_0$-equivariant isomorphism.\\

\noindent \emph{Proof of Claim.} Let $g_j \in \Hom_{\mathrm{cont}}(\Gamma_0, \widetilde{M}_{k_0}^r)$ denote the element given by $\gamma \mapsto \gamma(e_j)$ for $\gamma \in \Gamma_0$. We can write
\[
g_j(\gamma) = \sum_{1 \le h \le n, ~h \neq j} f_{jh}(\gamma)\cdot e_h+(1+f_{jj}(\gamma))\cdot e_j
\]
with $f_{jh}(\gamma) \in \widetilde{\fA}_X^{r'}$. Denote $f(\gamma)\coloneqq (f_{jh}(\gamma))_{1 \le j, h \le n}\in \mathrm{Mat}_{n\times n}(\widetilde{\fA}_X^{r'})$. It suffices to find a change-of-basis matrix $B \in \GL_n(\widetilde{\fA}_X^{r'})$ satisfying
\[
\gamma(B)(1+f(\gamma))B^{-1} \in \mathrm{Mat}_{n\times n}(\fA_{\Lambda}^{r'})
\]
for all $\gamma \in \Gamma_0$.

We will construct such a matrix $B$ inductively. Note that $\lambda_{r'}(f(\gamma)) \le \frac{1}{6}$ by our choice of $\Gamma_0$, and so $\lambda_{r'}(\overline{f}(\gamma)) \le \frac{1}{6}$. Suppose we already have $\lambda_{r'}(\overline{f}(\gamma)) \le \frac{1}{3\cdot 2^{\ell}}$ for some $\ell \ge 1$. The key observation is that $\overline{f}$ is close to a cocycle, so that the $2$-cochain $d^{(1)}(\overline{f})$ is close to $0$. We then use Corollary~\ref{cor: complex-A-strict-exact} to modify $1+f$ iteratively. More precisely, for any $\sigma, \beta \in \Gamma_0$, we have
\[
1+f(\sigma \beta) = (1+\sigma(f(\beta)))(1+f(\sigma)) = 1+f(\sigma)+\sigma(f(\beta))+\sigma(f(\beta))\cdot f(\sigma).
\]
So $d^{(1)}(f)(\sigma, \beta) = \sigma(f(\beta))\cdot f(\sigma)$. Thus, $\lambda_{r'}(d^{(1)}(\overline{f}_{jh})) \le (\frac{1}{3\cdot 2^{\ell}})^2$ for all $1 \le j, h \le n$. By Corollary~\ref{cor: complex-A-strict-exact}, there exists an element $\overline{b}_{jh} \in \widetilde{\fA}_X^{r'} / \fA_{\Lambda}^{r'}$ such that
\[
\lambda_{r'}(\overline{f}_{jh}(\gamma)+\gamma(\overline{b}_{jh})-\overline{b}_{jh}) \le \frac{3}{2}\lambda_{r'}(d^{(1)}(\overline{f}_{jh})) \le \frac{1}{3\cdot 2^{2\ell+1}}
\]
for each $\gamma \in \Gamma_0$. In particular, we have $\lambda_{r'}(d^{(0)}(\overline{b}_{jh})) \le \frac{1}{3\cdot 2^{\ell}}$, and so
\[
\lambda_{r'}(\overline{b}_{jh}) \le \frac{3}{2}\cdot \frac{1}{3\cdot 2^{\ell}} = \frac{1}{2^{\ell+1}}
\]
again by Corollary~\ref{cor: complex-A-strict-exact}.

For each $1 \le j, h \le n$, let $b_{jh} \in \widetilde{\fA}_X^{r'}$ be a preimage of $\overline{b}_{jh}$ such that $\lambda_{r'}(b_{jh}) = \lambda_{r'}(\overline{b}_{jh})$. Let $B^{(1)} = 1+(b_{jh})_{1 \le j, h \le n} \in \mathrm{Mat}_{n\times n}(\widetilde{\fA}_X^{r'})$. We have the following identity on matrix multiplications
\begin{align*}
& \gamma(1+b_{jh})\cdot (1+(f_{jh}(\gamma)))\cdot (1-(b_{jh})+o((3\cdot 2^{\ell+1})^{-1}) \\
=\,&  1+\gamma(b_{jh})+(f_{jh}(\gamma))-(b_{jh})+\gamma(b_{jh})
\cdot (f_{jh}(\gamma))\cdot (b_{jh}) - \gamma(b_{jh})\cdot (b_{jh})+o((3\cdot 2^{\ell+1})^{-1}),
\end{align*}
where $o((3\cdot 2^{\ell+1})^{-1})$ denotes some matrix $T \in \mathrm{Mat}_{n\times n}(\widetilde{\fA}_X^{r'})$ with $\lambda_{r'}(\overline{T}) \le (3\cdot 2^{\ell+1})^{-1}$. This implies
\[
\lambda_{r'}(\overline{\gamma(B^{(1)})(1+f(\gamma))(B^{(1)})^{-1} -1}) \le \frac{1}{3\cdot 2^{\ell+1}}. 
\]
Starting from the case $\ell = 1$ and iterating this step inductively, we obtain a sequence of matrices $B^{(1)}, B^{(2)}, \ldots \in \GL_n(\widetilde{\fA}_X^{r'})$ such that for each $\gamma \in \Gamma_0$, we have $\lambda_{r'}(B^{(\ell)}-1) \le \frac{1}{2^{\ell+1}}$ and
\[\lambda_{r'}(\overline{\gamma(B^{(\ell)})\cdots \gamma(B^{(1)})(1+f(\gamma)) (B^{(1)})^{-1}\cdots (B^{(\ell)})^{-1} - 1}) \le \frac{1}{3\cdot 2^{\ell+1}}.\]
Taking the limit $B\coloneqq \lim_{\ell \rightarrow \infty} B^{(\ell)}B^{(\ell-1)}\cdots B^{(1)}$, we obtain the desired matrix and prove the claim.\\

Now, recall that $\Gamma_0 \subset \Gamma$ is an open normal subgroup of finite index. Pick $g_1, \ldots, g_s \in \Gamma$ representing $\Gamma / \Gamma_0$, and consider $g_i\cdot M \subset \widetilde{M}_{k_0}^r$ as in the proof of Lemma~\ref{lem: general-fact-equivalence-phi-Gamma-modules}(3) for each $1 \le i \le s$. Note that $g_i\cdot M$ is $\Gamma_0$-stable and free over $\fA_{\Lambda}^{r'}$ with $(g_i\cdot M)\otimes_{\fA_{\Lambda}^{r'}} \widetilde{\fA}_X^{r'} \cong \widetilde{M}_{k_0}^r$. By applying a similar argument as in the proof of full faithfulness in Theorem~\ref{thm: fully-faithfulness-for-free-modules} (but using Corollary~\ref{cor: complex-A-strict-exact} instead of Proposition~\ref{prop: C-complex-strict-exact}) to $M$ and $g_i\cdot M$ for each $1 \le i \le s$, we deduce that for some $k_1 \ge k_0$, 
\[
N\coloneqq M\otimes_{\fA_{\Lambda}^{rp^{-k_0}}, \varphi^{k_1-k_0}} \fA_{\Lambda}^{rp^{-k_1}} \subset \widetilde{M}_{k_1}^r
\]
is stable under the $\Gamma$-action on $\widetilde{M}_{k_1}^r$. Furthermore, since the $\varphi$-action and $\Gamma$-action on $\widetilde{M}_{k_1}^r$ commute, we similarly deduce that the isomorphism
\[
\varphi_{\widetilde{M}_{k_1}^r}\colon \widetilde{M}_{k_1}^r\otimes_{\widetilde{\fA}_X^{rp^{-k_1}}, \varphi} \widetilde{\fA}_X^{rp^{-k_1-1}} \stackrel{\sim}{\longrightarrow} \widetilde{M}_{k_1}^r\otimes_{\widetilde{\fA}_X^{rp^{-k_1}}, \iota} \widetilde{\fA}_X^{rp^{-k_1-1}}
\]
descends to an isomorphism
\[
N\otimes_{\fA_{\Lambda}^{rp^{-k_1}}, \varphi} \fA_{\Lambda}^{rp^{-k_1-1}} \stackrel{\sim}{\longrightarrow} N\otimes_{\fA_{\Lambda}^{rp^{-k_1}}, \iota} \fA_{\Lambda}^{rp^{-k_1-1}}
\]
after increasing $k_1$ if necessary. This proves Step 1.\\

\noindent \textbf{Step 2.} We now consider the general situation where $\widetilde{M}$ is finite projective (but not necessarily free) over $\widetilde{\fA}_X^{\dagger}$.
By Corollary~\ref{cor: equiv-local-system-overconvergent-etale-(phi, Gamma)-mods}, there is a unique \'etale $\Z_p$-local system $\L$ on $X$ such that $\widetilde{D}^{\mathrm{int}, \dagger}(\L) \cong \widetilde{M}$. Choose a finite \'etale Galois cover $X'\rightarrow X$ over which $\L/p\L$ is trivialized and let $\widetilde{M}'\coloneqq \widetilde{M}\otimes_{\widetilde{\fA}_X^{\dagger}} \widetilde{\fA}_{X'}^{\dagger}$. From the construction of the functor $\widetilde{D}^{\mathrm{int}, \dagger}$, we deduce that $\widetilde{M}'/p$ is free over $\widetilde{\fA}_{X'}^{\dagger} / p = \widetilde{\fA}_{X'}/p$. Recall from the proof of Proposition~\ref{prop: perfect theory}(1) that $(\widetilde{\fA}_{X'}^{\dagger}, (p))$ is a henselian pair. In particular, $p$ lies in the Jacobson radical of $\widetilde{\fA}_{X'}^{\dagger}$. Thus, $\widetilde{M}'$ is finite free over $\widetilde{\fA}_{X'}^{\dagger}$ by Nakayama's lemma.

Let $\Lambda'$ denote the decompleting datum given by Construction~\ref{construction:wdc_finite_etale} applied to $X' \rightarrow X$ and $\Lambda$. By Step 1 for $\Lambda'$, $\widetilde{M}'$ descends uniquely to a $(\varphi, \Gamma)$-module $M'$ over $\fA_{\Lambda'}^{\dagger}$. Write $G$ for the Galois group of the cover $X'\rightarrow X$. By Propositions~\ref{prop: perfect theory} and \ref{prop:equivalence of categories finite etale} and Construction~\ref{construction:wdc_finite_etale}, the maps $\spec (\widetilde{\fA}_{X'}^{\dagger}) \rightarrow \spec (\widetilde{\fA}_X^{\dagger})$ and $\spec (\fA_{\Lambda'}^{\dagger}) \rightarrow \spec (\fA_{\Lambda}^{\dagger})$ are finite \'etale Galois covers with the same Galois group $G$ as $X' \rightarrow X$. Note that by the construction of the toric towers, the $G$-action on $\widetilde{M}'$ commutes with the $\varphi$ and $\Gamma$-actions.

By the uniqueness of descent, the $G$-action on $\widetilde{M}'$ induces an $\fA_{\Lambda'}^{\dagger}$-semilinear $G$-action on $M'$. Such an action corresponds to a descent datum for the cover $\spec (\fA_{\Lambda'}^{\dagger}) \rightarrow \spec (\fA_{\Lambda}^{\dagger})$. Hence, by faithfully flat descent, we obtain a $(\varphi, \Gamma)$-module over $\fA_{\Lambda}^{\dagger}$ which is the descent of $\widetilde{M}$.
\end{proof}

\begin{remark}
The strategy of reducing to the free case in the proof above is somewhat similar to argument in the proof of \cite[Theorem~4.35]{Andreatta-Brinon_Overconvergence}. However, note that the proof of \textit{loc. cit.} generalizes the method of \cite{Cherbonnier_Colmez} (valid when $X$ admits good reduction), while our argument is based on the results established in \S \ref{sec:decompleting_data} and hence is different.
\end{remark}

\begin{corollary} \label{cor: decompletion-etale-(phi, Gamma)-modules}
The base change functor
\[
\mathrm{Mod}_{/\fB_{\Lambda}^{\dagger}}^{(\varphi, \Gamma), \ett} \lra \mathrm{Mod}_{/ \widetilde{\fB}_X^{\dagger}}^{(\varphi, \Gamma), \ett} \quad (\textrm{resp.} ~~\mathrm{Mod}_{/\fC_{\Lambda}^{\dagger}}^{(\varphi, \Gamma), \ett} \lra \mathrm{Mod}_{/ \widetilde{\fC}_X^{\dagger}}^{(\varphi, \Gamma), \ett})
\]
along the map $\fB_{\Lambda}^{\dagger} \rightarrow \widetilde{\fB}_X^{\dagger}$ (resp. $\fC_{\Lambda}^{\dagger} \rightarrow \widetilde{\fC}_X^{\dagger}$) is an equivalence.
\end{corollary}

\begin{proof}
This follows from Theorem~\ref{thm: decompletion-etale-(phi, Gamma)-modules} and Proposition~\ref{prop:descending_phi_mod_along_3_and_9}(3) (see Remark~\ref{remark: descent of phi modules implies descent of phi Gamma modules}).    
\end{proof}

\subsection{From $p$-adic local systems to \'etale $(\varphi, \Gamma)$-modules} \label{ss:relative_D_rig} 

Combining Corollary \ref{cor: equiv-local-system-overconvergent-etale-(phi, Gamma)-mods} with the descent results established in this section, we obtain equivalences between the category of $p$-adic local systems and that of $(\varphi, \Gamma)$-modules over imperfect relative period rings.

\begin{theorem} \label{thm: equivalence-isogeny-local-system-overconvergent-etale-(phi, Gamma)-modules-C}
There are natural equivalences of $\otimes$-categories
\[
    D^{\mathrm{int}, \dagger} \colon \mathrm{Loc}_{\Z_p}(X) \isom  \mathrm{Mod}^{(\varphi, \Gamma)}_{/\fA^{\dagger}_{\Lambda}},
\]
\[
    D^{\dagger} \colon \mathrm{Loc}_{\Z_p}(X)_{\Q} \isom  \mathrm{Mod}^{(\varphi, \Gamma), \ett}_{/\fB^{\dagger}_{\Lambda}},
\]
and
\[
    D_{\mathrm{rig}}^{\dagger} \colon \mathrm{Loc}_{\Z_p}(X)_{\Q} \isom  \mathrm{Mod}^{(\varphi, \Gamma), \ett}_{/\fC^{\dagger}_{\Lambda} }.
\]

\end{theorem}

\begin{proof}
This follows from Corollary~\ref{cor: equiv-local-system-overconvergent-etale-(phi, Gamma)-mods}, Theorem~\ref{thm: decompletion-etale-(phi, Gamma)-modules}, and Corollary~\ref{cor: decompletion-etale-(phi, Gamma)-modules}. Note that the compatibility with tensor structures follows from Remark~\ref{rem: base-change-functor-compatible-with-tensor-dual}.  
\end{proof}

In particular, this finishes the proof of Theorem \ref{intro:theorem_equiv_of_categories} from the introduction.

\vspace{0.1in}
\section{Strongly decompleting data}
\label{sec:strongly_decompleting_data}

In this section, we develop a theory that allows us to construct a relative version of ``$\mathbf{N}_{\mathrm{dR}}(V)$'' of Berger in \cite{Berger-differential}. This  generalizes the work of Ohkubo \cite{Ohkubo-differential} (in the case of CDVF with imperfect residue fields) and Andreatta-Brinon \cite{Andreatta-Brinon} (in the relative case where $X$ is the rigid generic fiber of an \'etale algebra over $\mathcal{O}_K\langle T_1^{\pm}, \ldots, T_d^{\pm}\rangle$). Our construction will be a natural extension of Construction~\ref{construction:N_dR_for_all_HT_wts} in \S \ref{ss:NdR classical setting} and builds on the notion of \textit{strongly decompleting data}. 

\subsection{Strongly decompleting data} \label{subsec: definition-strongly-decompleting-data}

\subsubsection{Setup} 

Throughout, let $X=\spa(A_X, A^+_X)$ be an affinoid $K$-algebra together with a toric chart $h\colon X\rightarrow \mathbb{D}^d_K$. Recall that the tower $X_{\infty}=\varprojlim X_n\rightarrow X$ is an affinoid perfectoid object in the pro-\'etale site $X_{\proet}$, with associated affinoid perfectoid space $\widehat{X}_{\infty}$. Write $X_n=\spa(A_{X,n}, A^+_{X,n})$ for each $n\ge 1$, and write $\widehat{X}_{\infty}=\spa(\widehat{A}_{X, \infty}, \widehat{A}^+_{X, \infty})$. 
For $r$ sufficiently small\footnote{We shall assume that $r > 0$ is sufficiently small such that the statements in \S \ref{sec:decompleting_data} are applicable.} and for any $n\ge -\log_pr$, consider the map 
\[\iota_n\colon \widetilde{\mathbf{A}}^r_X\rightarrow \mathbb{B}_{\mathrm{dR}, X}^+(X_{\infty})\]
sending $\sum_{m=0}^{\infty}p^m[x_m]\mapsto \sum_{m=0}^{\infty}p^m[x_m^{\frac{1}{p^n}}]$, where $\mathbb{B}_{\mathrm{dR}, X}^+$ is the (horizontal) de Rham period sheaf on $X_{\proet}$ defined in \cite{Scholze_p_adic_Hodge}. This map is well-defined following the same computation as in \cite[Lemme~7.2]{Andreatta-Brinon}. Inverting $p$, we obtain $\iota_n\colon \widetilde{\mathbf{B}}^r_X\rightarrow \mathbb{B}_{\mathrm{dR}, X}^+(X_{\infty})$. 

\subsubsection{Topology on $\mathbb{B}_{\mathrm{dR}, X}^+(X_{\infty})$}

On $\mathbb{B}_{\mathrm{dR}, X}^+(X_{\infty})$, there are two useful topologies. First, the \emph{canonical topology} is defined to be the inverse limit topology on
\[\mathbb{B}_{\mathrm{dR}, X}^+(X_{\infty})\cong \varprojlim_{m} \mathbb{B}_{\mathrm{dR}, X}^+(X_{\infty})/t^m\mathbb{B}_{\mathrm{dR}, X}^+(X_{\infty})\]
where each $\mathbb{B}_{\mathrm{dR}, X}^+(X_{\infty})/t^m\mathbb{B}_{\mathrm{dR}, X}^+(X_{\infty})$ is equipped with the $p$-adic Banach space topology. On the other hand, one can endow $\mathbb{B}_{\mathrm{dR}, X}^+(X_{\infty})$ with the $t$-adic topology, which is finer than the canonical topology. Note that $\mathbb{B}_{\mathrm{dR}, X}^+(X_{\infty})$ is complete with respect to the $t$-adic topology. 

\begin{lemma} \label{lem:iota_n_extends_for_the_tilde_version}
Suppose that $p^{-n}\in [s,r]$. Then $\iota_n$ extends to a continuous map
\[
\iota_n\colon \widetilde{\fC}_X^{[s,r]} \rightarrow \mathbb{B}_{\mathrm{dR}, X}^+(X_{\infty})
\] 
where the source is equipped with the topology given by $\max\{\lambda_r, \lambda_s\}$-norm and the target is equipped with the canonical topology. 
\end{lemma}

\begin{proof}
Let $m \ge 1$, and consider the induced map
\[
\iota_n\colon \widetilde{\fB}_X^r \rightarrow  \mathbb{B}_{\mathrm{dR}, X}^+(X_{\infty})/t^m.
\]
By the computation as in the proof of \cite[Lemme~7.2]{Andreatta-Brinon}, there is a constant $C_{s, r, n, m} > 0$ (which may depend on $s, r, n, m$) such that for any $x \in \widetilde{\fB}_X^r$, we have
\[
|\iota_n(x)|_m \le C_{s, r, n, m} \cdot \max\{\lambda_r(x), \lambda_s(x)\}
\]
where $|\cdot|_m$ denotes the $p$-adic Banach sapce norm on $\mathbb{B}_{\mathrm{dR}, X}^+(X_{\infty})/t^m$. So the above map extends to a continuous map
\[
\iota_n\colon  \widetilde{\fC}_X^{[s,r]} \rightarrow \mathbb{B}_{\mathrm{dR}, X}^+(X_{\infty})/t^m.
\]
Passing to the inverse limit, we get the desired statement. 
\end{proof}

\subsubsection{$\mathcal{O}\mathbb{B}_{\mathrm{dR}, X}^+(X_{\infty})$} 

Let $\mathcal{O}\mathbb{B}_{\mathrm{dR}, X}^+$ be the de Rham period sheaf defined in \cite{Scholze_p_adic_Hodge}. Consider $\mathcal{O}\mathbb{B}_{\mathrm{dR}, X}^+(X_{\infty})$, which comes with a natural connection
\[\nabla\colon \mathcal{O}\mathbb{B}_{\mathrm{dR}, X}^+(X_{\infty})\rightarrow \mathcal{O}\mathbb{B}_{\mathrm{dR}, X}^+(X_{\infty})\otimes_{A_X}\Omega_{A_X/K}^1\]
and the theta map
\[\theta_X\colon \mathcal{O}\mathbb{B}_{\mathrm{dR}, X}^+(X_{\infty})\twoheadrightarrow \widehat{A}_{X,\infty}.\]
We equip $\mathcal{O}\mathbb{B}_{\mathrm{dR}, X}^+(X_{\infty})$ with the \emph{canonical topology}; that is, we equip 
\[\mathcal{O}\mathbb{B}_{\mathrm{dR}, X}^+(X_{\infty})\cong \varprojlim_m \mathcal{O}\mathbb{B}_{\mathrm{dR}, X}^+(X_{\infty})/(\mathrm{ker}\,\theta_X)^m\mathcal{O}\mathbb{B}_{\mathrm{dR}, X}^+(X_{\infty})\]
with the inverse limit topology, where each $\mathcal{O}\mathbb{B}_{\mathrm{dR}, X}^+(X_{\infty})/(\mathrm{ker}\,\theta_X)^m\mathcal{O}\mathbb{B}_{\mathrm{dR}, X}^+(X_{\infty})$ is equipped with the $p$-adic Banach space topology. 

Moreover, by \cite[Proposition 6.10]{Scholze_p_adic_Hodge}, there is an identification
\[\mathcal{O}\mathbb{B}_{\mathrm{dR}, X}^+(X_{\infty})\cong \mathbb{B}_{\mathrm{dR}, X}^+(X_{\infty})[\![u_1, \ldots, u_d]\!]\]
where $u_i=T_i-[T_i^{\flat}]$, and we recall that $\mathrm{ker}\,\theta_X=(t,u_1, \ldots, u_d)$.

\subsubsection{The ring $A_{X,n}[\![t, u_1, \ldots, u_d]\!]$}

For every integer $n\ge 0$, consider the subring
\[A_{X,n}[\![t, u_1, \ldots, u_d]\!]\subset \mathcal{O}\mathbb{B}_{\mathrm{dR}, X}^+(X_{\infty}).\]
This ring will play the role of $K_n [\![t]\!]$ in \S \ref{ss:NdR classical setting}. 
First, let us note that $A_{X,n}[\![t, u_1, \ldots, u_d]\!]$ is stable under the natural $\Gamma$-action on $\mathcal{O}\mathbb{B}_{\mathrm{dR}, X}^+(X_{\infty})$. We can endow $A_{X,n}[\![t, u_1, \ldots, u_d]\!]$ with the following two topologies:
\begin{enumerate}
\item[(i)] the subspace topology induced from the canonical topology on $\mathcal{O}\mathbb{B}_{\mathrm{dR}, X}^+(X_{\infty})$;
\item[(ii)] the inverse limit topology via 
\[A_{X,n}[\![t, u_1, \ldots, u_d]\!]=\varprojlim_{m} A_{X,n}[t, u_1, \ldots, u_d]/(t,u_1, \ldots, u_d)^m\]
where each $A_{X,n}[t, u_1, \ldots, u_d]/(t,u_1, \ldots, u_d)^m$, being a finite free $A_{X,n}$-module, is endowed with the natural $p$-adic Banach space topology.
\end{enumerate}

\begin{lemma}\label{lemma: two topologies} \hfill
\begin{enumerate}
\item The two topologies on $A_{X,n}[\![t, u_1, \ldots, u_d]\!]$ coincide. (This topology is referred to as the \emph{canonical topology}.)
\item $A_{X,n}[\![t, u_1, \ldots, u_d]\!]$ is a closed subspace of $\mathcal{O}\mathbb{B}_{\mathrm{dR}, X}^+(X_{\infty})$ (with respect the canonical topology). 
\end{enumerate}
\end{lemma}

\begin{proof}
This essentially follows from the definitions. It suffices to note that, for each $m$, there is a natural injection
\[A_{X,n}[t, u_1, \ldots, u_d]/(t,u_1, \ldots, u_d)^m\hookrightarrow \mathcal{O}\mathbb{B}_{\mathrm{dR}, X}^+(X_{\infty})/(\mathrm{ker}\,\theta_X)^m\mathcal{O}\mathbb{B}_{\mathrm{dR}, X}^+(X_{\infty})\]
identifying the source as a closed subspace of the target. Moreover, the $p$-adic Banach topology on the source coincides with the subspace topology induced from the $p$-adic Banach topology on the target. Taking inverse limit gives the desired claim. 
\end{proof}

\begin{remark}
One can also consider the $(t,u_1, \ldots, u_d)$-adic topology on $A_{X,n}[\![t, u_1, \ldots, u_d]\!]$, which is  finer than the canonical topology.
\end{remark}

\subsubsection{The ring $A_{X,n}[\![t, u_1, \ldots, u_d]\!]^{\nabla=0}$} \label{sss:rings_inside_OB_dR}

Taking the $\nabla=0$ part, we arrive at an inclusion
\[A_{X,n}[\![t, u_1, \ldots, u_d]\!]^{\nabla=0}\subset \mathcal{O}\mathbb{B}_{\mathrm{dR}, X}^+(X_{\infty})^{\nabla=0}=\mathbb{B}_{\mathrm{dR}, X}^+(X_{\infty}).\]
We can endow $A_{X,n}[\![t, u_1, \ldots, u_d]\!]^{\nabla=0}$ with the \emph{canonical topology} that is the subspace topology inherited from the canonical topology on $\mathbb{B}_{\mathrm{dR}, X}^+(X_{\infty})$. It is also the subspace topology inherited from \[A_{X,n}[\![t, u_1, \ldots, u_d]\!]^{\nabla=0}\subset A_{X,n}[\![t, u_1, \ldots, u_d]\!],\]
where the latter is equipped with the canonical topology. 

\begin{corollary}\label{cor: B_n closed subspace} \hfill
\begin{enumerate}
    \item 
$A_{X,n}[\![t, u_1, \ldots, u_d]\!]^{\nabla=0}$ is a closed subspace of $\mathbb{B}_{\mathrm{dR}, X}^+(X_{\infty})=\mathcal{O}\mathbb{B}_{\mathrm{dR}, X}^+(X_{\infty})^{\nabla=0}$ (with respect to the canonical topology).  
\item For each $m \ge 1$,  $A_{X,n}[\![t, u_1, \ldots, u_d]\!]^{\nabla=0}/t^m$ is a closed subspace of $\mathbb{B}_{\mathrm{dR}, X}^+(X_{\infty})/t^m$ where the target is endowed with the $p$-adic topology. 
\end{enumerate} 

\end{corollary}

\begin{proof}
Since $\nabla$ is continuous, $\mathbb{B}_{\mathrm{dR}, X}^+(X_{\infty})$ is a closed subspace of $\mathcal{O}\mathbb{B}_{\mathrm{dR}, X}^+(X_{\infty})$ under the canonical topology. For (1),  it suffices to notice  that 
\[ A_{X,n}[\![t, u_1, \ldots, u_d]\!]^{\nabla=0}=\mathbb{B}_{\mathrm{dR}, X}^+(X_{\infty})\cap A_{X,n}[\![t, u_1, \ldots, u_d]\!] \]  and then apply Lemma~\ref{lemma: two topologies}(2). Part (2) follows from a similar argument. 
\end{proof}

For later use, we show that, as a ring, $A_{X,n}[\![t, u_1, \ldots, u_d]\!]^{\nabla=0}$ can be identified with $A_{X,n}[\![t]\!]$. The construction follows verbatim as in \cite[Lemma~4.1.4]{Ohkubo-differential}. For completeness, we recall the construction. Consider the map
\[f\colon A_{X,n}[t, u_1, \ldots, u_d]\rightarrow A_{X,n}[\![t, u_1, \ldots, u_d]\!]\]
sending 
\[x \mapsto \sum_{n_1, \ldots, n_d\in \mathbb{Z}_{\ge 0}} \frac{(-1)^{n_1+\cdots+n_d}}{n_1!\cdots n_d!} u_1^{n_1}\cdots u_d^{n_d}\partial_1^{n_1}\cdot\circ \cdots \circ\partial_d^{n_d} (x).\] 
It is straightforward to verify that $f$ is a ring homomorphism such that
\begin{itemize}
\item $f(tx)=tf(x)$ for all $x\in A_{X,n}[t, u_1, \ldots, u_d]$;
\item $f(u_j)=0$ for all $j=1,\ldots, d$; 
\item $f(T_j) = T_j - u_j = [T_j^{\flat}]$ for all $j=1,\ldots, d$; 
\item the image of $f$ is contained in $A_{X,n}[\![t, u_1, \ldots, u_d]\!]^{\nabla=0}$.
\end{itemize}
Moreover, $f$ is continuous with respect to the $(t,u_1, \ldots, u_d)$-adic topology on the source and target, hence $f$ extends to a ring homomorphism
\[f\colon A_{X,n}[\![t, u_1, \ldots, u_d]\!]\rightarrow A_{X,n}[\![t, u_1, \ldots, u_d]\!]^{\nabla=0}.\]

\begin{lemma} \label{lemma:isomorphism_formal_power_series_with_horizontal_sections}
  The map $f$ induces an isomorphism 
\begin{equation}
\label{eq:isomorphism_formal_power_series_with_horizontal_sections}
\bar{f}\colon A_{X,n}[\![t]\!] \isom  A_{X,n}[\![t, u_1, \ldots, u_d]\!]^{\nabla=0}.
\end{equation}
\end{lemma}

\begin{proof}
Notice that $f$ is identity on the subring $A_{X,n}[\![t, u_1, \ldots, u_d]\!]^{\nabla=0}$. Hence $f$ is surjective. We obtain a surjection
\[\bar{f}\colon A_{X,n}[\![t]\!]\cong A_{X,n}[\![t, u_1, \ldots, u_d]\!]/(u_1, \ldots, u_d)\twoheadrightarrow A_{X,n}[\![t, u_1, \ldots, u_d]\!]^{\nabla=0}.\]
It remains to show $\bar{f}$ is injective. For this, it suffices to notice that $\bar{f}(t)=t$, thus the desired injectivity reduces to the injectivity of $\cl f$ on $A_{X, n}$. Since $\cl f$ is a ring homomorphism and both the source and target are reduced, it suffices to show that $\cl f$ is injective on $A_X$, 
which follows from the explicit formula of $f$ (in particular, the fact that $f (T_j) = T_j - u_j$).  
\end{proof} 

\begin{remark}
    As a byproduct, note that the $t$-adic topology on $A_{X,n}[\![t, u_1, \ldots, u_d]\!]^{\nabla=0}$ (via the isomorphism $\bar{f}\colon A_{X,n}[\![t]\!] \isom A_{X,n}[\![t, u_1, \ldots, u_d]\!]^{\nabla=0}$) is finer than the canonical topology (cf. \cite[Lemma~4.1.5]{Ohkubo-differential}).
\end{remark}

\begin{lemma} \label{lemma:sublemma_on_horizontal_sections}
The natural inclusion $A_{X,n}[\![t, u_1, \ldots, u_d]\!]^{\nabla=0}\subset A_{X,n}[\![t,u_1, \ldots, u_d]\!]$ induces an isomorphism
\[A_{X,n}[\![t, u_1, \ldots, u_d]\!]^{\nabla=0}[\![u_1, \ldots, u_d]\!]\cong A_{X,n}[\![t,u_1, \ldots, u_d]\!].\]
\end{lemma}

\begin{proof}
First, the inclusion $A_{X,n}[\![t, u_1, \ldots, u_d]\!]^{\nabla=0}\subset A_{X,n}[\![t,u_1, \ldots, u_d]\!]$ indeed extends to a ring homomorphism $g\colon A_{X,n}[\![t, u_1, \ldots, u_d]\!]^{\nabla=0}[\![u_1, \ldots, u_d]\!]\rightarrow A_{X,n}[\![t,u_1, \ldots, u_d]\!]$.  
To check that $g$ is an isomorphism, it suffices to check after modulo $(u_1, \ldots, u_d)$; we obtain 
\[ \bar{g}\colon A_{X,n}[\![t, u_1, \ldots, u_d]\!]^{\nabla=0}\rightarrow A_{X,n}[\![t]\!].
\] 
Composing with the isomorphism $\bar{f}\colon A_{X,n}[\![t]\!]\cong A_{X,n}[\![t, u_1, \ldots, u_d]\!]^{\nabla=0}$ in Lemma~\ref{lemma:isomorphism_formal_power_series_with_horizontal_sections}, one sees that the map 
\[\bar{f}\circ \bar{g}\colon A_{X,n}[\![t, u_1, \ldots, u_d]\!]^{\nabla=0}\rightarrow A_{X,n}[\![t, u_1, \ldots, u_d]\!]^{\nabla=0}\] is simply the identity; hence $\bar{g}$ is also an isomorphism.
\end{proof}

\subsubsection{Strongly decompleting data}

We are ready to define the notion of strongly decompleting data. Let $\Lambda=(R_{\Lambda}, \mathbf{A}_{\Lambda})$ be a decompleting datum (see Definition~\ref{definition:decompleting_datum}), with associated imperfect period rings $\mathbf{A}_{\Lambda}^r, \mathbf{B}_{\Lambda}^r, \mathbf{C}_{\Lambda}^{[s,r]}, \mathbf{C}_{\Lambda}^r$, etc. Also, recall the open subgroups $\Gamma_n\subset \Gamma$ defined in \S \ref{ss:toric_tower}, for integers $n\ge 0$.

\begin{definition}\label{defn: strongly decompleting data}
The decompleting datum $\Lambda=(R_{\Lambda}, \mathbf{A}_{\Lambda})$ is called \emph{strongly decompleting} if 
the following conditions hold:
\begin{enumerate}
\item[(a)] For $r$ sufficiently small, the image of the composition
\[\mathbf{A}^r_{\Lambda}\subset \widetilde{\mathbf{A}}^r_X\xrightarrow[]{\iota_n} \mathbb{B}_{\mathrm{dR}, X}^+(X_{\infty})\]
lies in $A_{X,n}[\![t, u_1, \ldots, u_d]\!]^{\nabla=0}$, for all $n\ge -\log_pr$.
\item[(b)] The composition \[ \mathbf{A}^r_{\Lambda}\xrightarrow[]{\iota_n} A_{X,n}[\![t, u_1, \ldots, u_d]\!]^{\nabla=0} \twoheadrightarrow A_{X, n}
\] is surjective for all $n\ge -\log_pr$.

\item[(c)] Let $r_n=p^{-n}$. There are an integer $n_0 \ge 0$ and a constant $C < p^{-\frac{1}{p-1}}$ such that for any $n \ge n_0$, we have
\[\lambda_{r_n}\left((\gamma-1)(x)\right)\le C\cdot \lambda_{r_n}(x)\]
for all $x\in \mathbf{C}_{\Lambda}^{[r_n, r_n]}$ and $\gamma\in \Gamma_n$.

\end{enumerate}
\end{definition}

\begin{remark} \hfill
\begin{enumerate}
\item As a consequence of condition (a), we have the continuous map \[
\iota_n\colon \mathbf{A}^r_{\Lambda}\rightarrow A_{X,n}[\![t, u_1, \ldots, u_d]\!]^{\nabla=0}.
\] 
By Corollary~\ref{cor: B_n closed subspace} and the proof of Lemma~\ref{lem:iota_n_extends_for_the_tilde_version}, this map further extends to a continuous map 
\begin{equation}
     \label{eq:iota_n_for_C_s,r}  \iota_n\colon \mathbf{C}_{\Lambda}^{[s,r]}\rightarrow A_{X,n}[\![t, u_1, \ldots, u_d]\!]^{\nabla=0}, 
\end{equation} whenever $p^{-n}\in [s,r]$, 
  with respect to the topology induced from the norm $\lambda_{s, r}\coloneqq \max\{\lambda_s, \lambda_r\}$ on the source and the canonical topology on the target.   Note that $\iota_n$ is compatible with the $\Gamma$-actions and satisfies $\iota_n (t) = p^{-n}t$. 
\item The composition in condition (b) is simply the restriction of $\theta\circ \varphi^{-n}\colon \widetilde{\mathbf{B}}^r_X\twoheadrightarrow \widehat{A}_{X, \infty}$ to $\mathbf{A}^r_{\Lambda}$.
\item Condition (c) is imposed so that a relative version of Lemma~\ref{lemma: estimation Gamma on C} holds for $\mathbf{C}_{\Lambda}^{[r_n, r_n]}$. This property is needed for our construction of the relative $\mathbf{N}_{\mathrm{dR}}$ functors in \S \ref{sec:relative_N_dR}.
\end{enumerate}
\end{remark}

\begin{lemma}\label{lemma: dense in B_n}
Suppose that $\Lambda=(R_{\Lambda}, \mathbf{A}_{\Lambda})$ is strongly decompleting. Let $[s,r]$ be an interval containing $p^{-n}$, where $r$ is sufficiently small. Then $\iota_n\colon \mathbf{C}_{\Lambda}^{[s,r]}\rightarrow A_{X,n}[\![t, u_1, \ldots, u_d]\!]^{\nabla=0}$ has dense image with respect to the $t$-adic topology, for all $n\ge -\log_pr$.
\end{lemma}

\begin{proof}
We need to show that, for any $m\ge 1$, the composition
\[\mathbf{C}_{\Lambda}^{[s,r]}\xrightarrow[]{\iota_n}A_{X,n}[\![t, u_1, \ldots, u_d]\!]^{\nabla=0} \twoheadrightarrow A_{X,n}[\![t, u_1, \ldots, u_d]\!]^{\nabla=0}/t^m\]
is surjective. Since $p^nt\in \mathbf{C}_{\Lambda}^{[s,r]}$ and $t=\iota_n(p^n t)$, thus by induction, we reduce to the case $m=1$. But this follows from Definition~\ref{defn: strongly decompleting data}(b). 
\end{proof}

\begin{corollary} \label{cor: t-completion-of-C}
Suppose that $\Lambda=(R_{\Lambda}, \mathbf{A}_{\Lambda})$ is strongly decompleting. Let $r > 0$ be sufficiently small (so that Corollary~\ref{cor: Frob decomp A uniform strict} is applicable, in particular). Suppose $\#\left([s,r]\cap p^{\mathbb{Z}}\right) = 1$ with $p^{-n} \in [s, r]$. Then the map 
\[
\varprojlim_{m \ge 1} \fC_{\Lambda}^{[s, r]} /t^m  \fC_{\Lambda}^{[s, r]} \rightarrow A_{X,n}[\![t, u_1, \ldots, u_d]\!]^{\nabla=0}
\]
induced by $\iota_n$ is an isomorphism.    
\end{corollary}

\begin{proof}
By Lemma~\ref{lemma: dense in B_n}, the map
\[
\iota_n\colon \fC_{\Lambda}^{[s, r]} / t^m \rightarrow A_{X,n}[\![t, u_1, \ldots, u_d]\!]^{\nabla=0} / t^m 
\]
is surjective for each $m \ge 1$. So we need to show that the above map is injective. Since $t$ is a non-zero divisor of $\fC_{\Lambda}^{[s, r]}$ and $A_{X,n}[\![t, u_1, \ldots, u_d]\!]^{\nabla=0}$, it suffices to show that the map is injective when $m = 1$.

By the assumption on $[s, r]$, the map $\iota_n\colon \widetilde{\fC}_X^{[s, r]}/t \rightarrow \mathbb{B}_{\mathrm{dR}, X}^+(X_{\infty})/t$ is injective. So it suffices to show that the map
\[
\fC_{\Lambda}^{[s, r]}/t \rightarrow \widetilde{\fC}_X^{[s, r]}/t
\]
is injective. By Corollary~\ref{cor: Frob decomp A uniform strict} and Lemma~\ref{lem: C breve dense in C tilde}, we have a $\fC_{\Lambda}^{[s, r]}$-linear projection $P\colon \widetilde{\fC}_X^{[s, r]} \rightarrow \fC_{\Lambda}^{[s, r]}$ which is a section to the inclusion $\fC_{\Lambda}^{[s, r]} \rightarrow \widetilde{\fC}_X^{[s, r]}$. Let $x \in \fC_{\Lambda}^{[s, r]}$ be an element such that $x = ty$ for some $y \in \widetilde{\fC}_X^{[s, r]}$. Then
\[
x = P(x) = P(ty) = tP(y)
\]
with $P(y) \in \fC_{\Lambda}^{[s, r]}$. This proves the statement.
\end{proof}

Let us also record the following lemma for later use. 

\begin{lemma} \label{lemma:fully_faithful_relative_B_dR} 
For all $n \in \Z_{ \ge 1}$,  the maps 
\begin{enumerate}
    \item $A_{X,n}[\![t, u_1, \ldots, u_d]\!]^{\nabla=0}\rightarrow \mathbb{B}_{\mathrm{dR}, X}^+(X_{\infty})$ 
    \item $A_{X,n}[\![t, u_1, \ldots, u_d]\!] \ra \mathcal{O}\mathbb{B}_{\mathrm{dR}, X}^+(X_{\infty})$ 
\end{enumerate}
 are faithfully flat. The same assertions hold when we replace $n$ by $\infty$. 
\end{lemma}

\begin{proof} Let us prove (1) first. 
Note that $A_{X,n}[\![t, u_1, \ldots, u_d]\!]^{\nabla=0} \cong A_{X, n}[\![t]\!]$ is noetherian. Since both rings are $t$-adically separated and complete and $t$ is a non-zero-divisor, it suffices to check after modulo $t$ by \cite[051C, 0912]{stacks-project}. After modulo $t$, the map becomes $A_{X,n}\rightarrow \widehat{A}_{X, \infty}$, which we claim is faithfully flat. For this, it suffices to show that both $A_{X,n}\rightarrow A_{X,n}\otimes_{K_n}\widehat{K}_{\infty}$ and $A_{X,n}\otimes_{K_n}\widehat{K}_{\infty}\rightarrow \widehat{A}_{X, \infty}$ are faithfully flat. The first arrow is clearly faithfully flat. For the second arrow, it suffices to show that $A_{X,n}^{\circ}\otimes_{\mathcal{O}_{K_n}}\mathcal{O}_{\widehat{K}_{\infty}}\rightarrow \widehat{A}_{X, \infty}^{\circ}$ is almost faithfully flat. By almost base change, we immediately reduce to the case of poly-annuli, which is well-known. The proof of (2) is similar --- in this case both rings are $(t, u_1, \ldots, u_d)$-adically separated and complete. 
\end{proof}

\subsection{Strongly decompleting data under finite \'etale extensions} \label{subsec: strongly-decompleting-finite-etale}

Let $h\colon X\rightarrow \mathbb{D}^d_K$ be a toric chart as before, and let $g:X'\rightarrow X$ be a finite \'etale morphism. Let $h'\colon X'\rightarrow \mathbb{D}^d_K$ be the composition of $h$ and $g$. The goal of this subsection is to show that, for any strongly decompleting datum for $(X,h)$, the induced datum for $(X',h')$ is also strongly decompleting.

Recall the following construction from \S \ref{ss:weakly decompleting finite etale}.  Let $\Lambda=(R_{\Lambda}, \mathbf{A}_{\Lambda})$ be a decompleting datum for $(X,h)$. The finite \'etale morphism $g:X'\rightarrow X$ induces a finite \'etale ring map $\widehat{A}_{X, \infty}\rightarrow \widehat{A}_{X', \infty}$.  We view $\widehat{A}_{X',\infty}$ as an object in $\left(\widehat{A}_{X, \infty}\right)_{\mathrm{f\acute{e}t}}$, and let $\fA'$ and $R'$ be the corresponding objects in $({\fA}_{\Lambda})_{\mathrm{f\acute{e}t}}$ and $(R_{\Lambda})_{\mathrm{f\acute{e}t}}$. By Theorem~\ref{thm: weakly decompleting data finite etale}, the pair $\Lambda'=(R', \fA')$ is a decompleting datum for $(X',h')$.

\begin{theorem} \label{thm: strongly decompleting under finite etale extensions}
If $\Lambda=(R_{\Lambda}, \mathbf{A}_{\Lambda})$ is strongly decompleting, so is $\Lambda'=(R', \fA')$.
\end{theorem}

\begin{proof}
For $r>0$ and for any $n\ge -\log_pr$, we have a commutative diagram
\begin{equation}
   \label{eq:commutative_diagram_is_pushout}
\begin{tikzcd}
\mathbf{A}^r_{\Lambda} \arrow[r, hookrightarrow]\arrow[d] & \widetilde{\mathbf{A}}^r_X \arrow[r,  "\theta \circ \iota_n"] \arrow[d]
& \widehat{A}_{X, \infty}\arrow[d] \\ 
\fA'^r \arrow[r, hookrightarrow] & \widetilde{\mathbf{A}}^r_{X'} \arrow[r, "\theta \circ \iota_n"] 
&\widehat{A}_{X', \infty}.  
\end{tikzcd}
\end{equation}

\noindent\textbf{Claim 1}: For $r$ sufficiently small, both squares in (\ref{eq:commutative_diagram_is_pushout}) are  pushout squares.\\

\noindent \textit{Proof of Claim 1.} The first square is a pushout square for $r$ sufficiently small by Construction~\ref{construction:wdc_finite_etale} and Proposition~\ref{prop:equivalence of categories finite etale}. For the second square, note that the top arrow $\theta \circ \iota_n$ factors as 
\[
\widetilde{\mathbf{A}}^r_X \xrightarrow{\varphi^{-n}} \widetilde{\mathbf{A}}^{p^n r}_X \xrightarrow{ \theta }  \widehat{A}_{X, \infty}
\]
and similarly for the bottom arrow. The desired claim thus follows from Proposition~\ref{prop: perfect theory}. \\

\noindent\textbf{Claim 2}: For $r$ sufficiently small and for all $n\ge -\log_pr$, the image of the composition
\[\fA'^r\subset \widetilde{\mathbf{A}}^r_{X'}\xrightarrow[]{\theta\circ \varphi^{-n}}\widehat{A}_{X', \infty}\]
lies in $A_{X', n}$, and the resulting map $\theta\circ\varphi^{-n}\colon \fA'^r\rightarrow A_{X', n}$ is surjective.\\

\noindent \textit{Proof of Claim 2.}
We start with a fixed sufficiently small $r$ and consider $n\ge -\log_pr$. Let $D_n\coloneqq   \fA'^r  \otimes_{\mathbf{A}^r_{\Lambda}, \theta \circ \varphi^{-n}} A_{X, n}$ and let  $D\coloneqq A_{X,\infty}\otimes_{A_{X,n}}D_n$. Then we obtain a base change diagram
\[
\begin{tikzcd}
\mathbf{A}^r_{\Lambda} \arrow[r, "\theta\circ\varphi^{-n}"] \arrow[d]& A_{X,n}\arrow[r, hookrightarrow]\arrow[d] & A_{X, \infty}\arrow[r, hookrightarrow]\arrow[d] & \widehat{A}_{X, \infty}\arrow[d ] \\ 
\fA'^r \arrow[r] & D_n \arrow[r, hookrightarrow] & D\arrow[r, hookrightarrow] & \widehat{A}_{X', \infty}
\end{tikzcd}.
\]
The map $\fA'^r \twoheadrightarrow D_n$ is surjective. Moreover, since $\fA'^r$ is finite \'etale over $\mathbf{A}^r_{\Lambda}$, the maps $D_n \ra D$ and $D \ra  \widehat{A}_{X', \infty}$ are injective. Thus, for every sufficiently large $n$, $D_n$ can be identified with the image of $\fA'^r$ under the map $\theta\circ \varphi^{-n}\colon \fA'^r \ra \widehat{A}_{X', \infty}$. 
Using the equivalence of categories between $ (\widehat{A}_{X,\infty})_{\mathrm{f\acute{e}t}}$ and $\left(A_{X,\infty}\right)_{\mathrm{f\acute{e}t}}\cong \varinjlim_n\left(A_{X,n}\right)_{\mathrm{f\acute{e}t}}$, we conclude that $D=A_{X', \infty}$ and thus $D_n=A_{X',n}$ for sufficiently large $n$ (say $n\ge n_0$). Therefore, it suffices to shrink $r$ so that $r\le p^{-n_0}$. 
This finishes the proof of Claim 2.\\

Let us return to the proof of the theorem. We first check condition (a) and (b) of Definition~\ref{defn: strongly decompleting data}. By Claim 2, it remains to show that for $r$ sufficiently small and for all $n\ge -\log_pr$, the image of the composition
\[\fA'^r\subset \widetilde{\mathbf{A}}^r_{X'}\xrightarrow[]{\iota_n} \mathbb{B}_{\mathrm{dR}, X'}^+(X'_{\infty})\]
lies in the subring 
\[A_{X', n}[\![t, u_1, \ldots, u_d]\!]\subset \mathcal{O}\mathbb{B}_{\mathrm{dR}, X'}^+(X'_{\infty}).\]
Note that $A_{X', n}[\![t, u_1, \ldots, u_d]\!]$ is $I$-adically complete with $I = (t, u_1, \ldots, u_d)$, hence the assertion follows  from the lifting property of finite \'etale maps along pro-infinitesimal thickenings, as shown in the following diagram:
\begin{equation*}
\begin{tikzcd}
\mathbf{A}^r_{\Lambda} \arrow[d]\arrow[r, "\iota_n"]& A_{X,n}[\![t, u_1, \ldots, u_d]\!]^{\nabla=0} \arrow[r]& A_{X', n}[\![t, u_1, \ldots, u_d]\!] \arrow[d, twoheadrightarrow] \\ 
\fA'^r \arrow[rr, "\theta\circ \varphi^{-n}"] \arrow[urr, dashrightarrow, "\exists !"] && A_{X',n}=A_{X', n}[\![t, u_1, \ldots, u_d]\!]/I  
\end{tikzcd}
\end{equation*}

Now, we check condition (c) of Definition~\ref{defn: strongly decompleting data}. Assume $\Lambda=(R_{\Lambda}, \mathbf{A}_{\Lambda})$ satisfies condition (c) for $n_0$ and $C$. We want to show that $\Lambda'=(R', \mathbf{A}')$ satisfies condition (c) for some integer $n'_0 \ge 0$ and a constant $C' < p^{-\frac{1}{p-1}}$.\\

\noindent{\textbf{Step 1.}} (residue norms on $R'$)\\
Pick $x_1, \ldots, x_d\in R'$ with $|x_i|\le 1$ such that
\[R_{\Lambda}\langle T_1, \ldots, T_d\rangle\rightarrow R'; \qquad T_i\mapsto x_i\]
is surjective. Let $|\cdot|_{\mathrm{res}, x_1, \ldots, x_d}$ denote the residue norm induced by this surjection. If the generators $x_1, \ldots, x_d$ is clear from the context, we simply write $|\cdot|_{\mathrm{res}}$. On the one hand, we always have $|\cdot|\le |\cdot|_{\mathrm{res}}$. On the other hand, since $|\cdot|$ and $|\cdot|_{\mathrm{res}}$ are equivalent, there exists a constant $C_*\ge 1$ such that $|\cdot|_{\mathrm{res}}\le C_* |\cdot|$.

Notice that for any integer $M\ge 0$, $\varphi^M(x_1), \ldots, \varphi^M(x_d)$ is also a set of generators of $R'$ over $R_{\Lambda}$; namely, the map
\[R_{\Lambda}\langle T_1, \ldots, T_d\rangle\rightarrow R'; \qquad T_i\mapsto \varphi^{M}(x_i)\]
is surjective. This follows from the fact that since $R'$ is finite \'etale over $R_{\Lambda}$, any $p$-basis of $R_{\Lambda}$ is also a $p$-basis of $R'$ (\cite[\S 1.1.2(ii)]{Berthelot-Messing}). For the rest of Step 1, we want to compare $|\cdot|$ with the residue norm $|\cdot|_{\mathrm{res}, \varphi^M(x_1), \ldots, \varphi^M(x_d)}$.

Pick $e_1, \ldots, e_{\ell}\in R'$ such that we have 
$
R_{\Lambda}=\bigoplus_{j=1}^{\ell} \varphi(R_{\Lambda})\cdot e_j$  
with $|e_j|\le 1$ for all $j$. Then we also have 
$R'=\bigoplus_{j=1}^{\ell} \varphi(R')\cdot e_j$.
Consider the sup-norm $|\cdot|_{e_1, \ldots, e_{\ell}}$ on $R'$ defined as 
\[|y|_{e_1, \ldots, e_{\ell}}\coloneqq \max_{j} |y_j|,\]
where $y=\sum_{j=1}^{\ell}y_j\cdot e_j$ with $y_j\in \varphi(R')$. Since $|\cdot|_{e_1, \ldots, e_{\ell}}$ is equivalent to $|\cdot|$, there exists $D\ge 1$ such that
\[|\cdot|\le |\cdot|_{e_1, \ldots, e_{\ell}} \le D |\cdot|.\]

\noindent\textbf{Claim 3:} For each integer $M\ge 0$, we have
\[|\cdot|\le |\cdot|_{\mathrm{res}, \varphi^M(x_1), \ldots, \varphi^M(x_d)} \le D^{1+p+p^2+\cdots+p^{M-1}}\cdot C_*^{p^M}|\cdot|\]
on $R'$.\\

\noindent\textit{Proof of Claim 3:} We only need to prove the second inequality. By induction, it suffices to prove for $M=1$. Let $y\in R'$. We have
$y=\sum_{j=1}^{\ell}\varphi(y_j)\cdot e_j$ 
with $y_j \in R'$. For each $1 \le j \le \ell$, write
\[
y_j = \sum_{i_1, \ldots, i_d \in \Z_{\ge 0}} a_{i_1, \ldots, i_d; j}\cdot x_1^{i_1}\cdots x_d^{i_d}
\]
with $a_{i_1, \ldots, i_d; j} \in R_{\Lambda}$, such that
$
|y_j|_{\mathrm{res}, x_1, \ldots, x_d} = \max_{i_1, \ldots, i_d}\{|a_{i_1, \ldots, i_d; j}|\}.
$ 
Then
\[
y = \sum_{i_1, \ldots, i_d \in \Z_{\ge 0}} (\sum_{j=1}^{\ell} a_{i_1, \ldots, i_d;j}^p\cdot e_j)\cdot \varphi(x_1)^{i_1} \cdots \varphi(x_d)^{i_d}.
\]
Consequently, 
\begin{align*}
|y| \ge & \,\,D^{-1}\cdot |y|_{e_1, \ldots, e_{\ell}} = D^{-1}\cdot \max_{j} \{|\varphi(y_j)|\}=D^{-1}\cdot (\max_{j}|y_j|)^p\\
\ge & \,\,D^{-1}\cdot \left(\max_j\{C_*^{-1}\cdot |y_j|_{\mathrm{res}, x_1, \ldots, x_d}\}\right)^p\\
= & \,\,D^{-1}\cdot C_*^{-p} \cdot \max_j \max_{i_1, \ldots, i_d} \{|a_{i_1, \ldots, i_d; j}|^p\} = D^{-1}\cdot C_*^{-p}\cdot \max_j \max_{i_1, \ldots, i_d} \{|a_{i_1, \ldots, i_d; j}^p|\}\\
\ge & \,\,D^{-1}\cdot C_*^{-p}\cdot |y|_{\mathrm{res}, \varphi(x_1), \ldots, \varphi(x_d)}.
\end{align*}
This proves Claim 3.\\

For later use, let $C_{**}\coloneqq D\cdot C_*$. Then we have
\[|\cdot|\le |\cdot|_{\mathrm{res}, \varphi^M(x_1), \ldots, \varphi^M(x_d)} \le  C_{**}^{p^M}|\cdot|\]
on $R'$, for every integer $M\ge 0$.\\

\noindent\textbf{Step 2.} (residue norms on ${\mathbf{A}'}^r$ and ${\mathbf{C}'}^{[r,r]}$)\\
Pick good lifts $\widetilde{x}_1, \ldots, \widetilde{x}_d\in {\mathbf{A}'}^{r_0}$ of $x_1, \ldots, x_d$ as in Proposition~\ref{prop: good lift weakly decompleting}. By shrinking $r_0$, we assume $C_{**}^{r_0}<p^{\frac{1}{2}}$ and the surjectivity as in Lemma~\ref{lem: Frob decomp A} holds for ${\fA'}^r$ for any $0 < r \le r_0$. By further shrinking $r_0$, we may assume that $\widetilde{x}_1, \ldots, \widetilde{x}_d$ generate ${\mathbf{A}'}^r$ over $\mathbf{A}_{\Lambda}^r$ for all $0<r\le r_0$; namely, the map
\[\mathbf{A}_{\Lambda}^r\langle T_1, \ldots, T_d\rangle \rightarrow {\mathbf{A}'}^r;\qquad T_i\mapsto \widetilde{X}_i\] 
is surjective. This can be deduced from that $({\fA'}^{\dagger}, (p))$ is a henselian pair (see the proof of Proposition~\ref{prop:equivalence of categories finite etale}). Consequently, the map
\[\mathbf{C}_{\Lambda}^{[r,r]}\langle T_1, \ldots, T_d\rangle \rightarrow {\mathbf{C}'}^{[r,r]};\qquad T_i\mapsto \widetilde{X}_i\] 
is also surjective, for all $0<r\le r_0$. Consider the residue norm $\lambda_{r,\mathrm{res}, \widetilde{x}_1, \ldots, \widetilde{x}_d}$ on ${\mathbf{C}'}^{[r,r]}$. When the generators $\widetilde{x}_1, \ldots, \widetilde{x}_d$ are clear from the context, we simply write $\lambda_{r,\mathrm{res}}$.\\

\noindent{\textbf{Claim 4:}} For any $0<r\le r_0$, we have 
\[\lambda_r\le \lambda_{r,\mathrm{res}, \widetilde{x}_1, \ldots, \widetilde{x}_d}\le C_{**}^r\cdot  \lambda_r\]
on ${\mathbf{C}'}^{[r,r]}$.\\

\noindent\textit{Proof of Claim 4:} For simplicity, we just write $\lambda_{r,\mathrm{res}}$. The first inequality is clear; we only need to prove the second one. For this, it suffices to prove $\lambda_{r,\mathrm{res}}(y)\le C_{**}^r\cdot  \lambda_r(y)$ for all $y\in {\mathbf{A}'}^r$.

Let $y\in {\mathbf{A}'}^r$ be a nonzero element. Express $y$ as a norm-compatible expansion $y=\sum_{n=0}^{\infty} p^ny_n$ as in Proposition~\ref{prop: good lift weakly decompleting}(2), where $y_n\in {\mathbf{A}'}^{r_0}$ such that
\begin{itemize}
\item $y_n$ is a good lift of $\overline{y}_n$ (in particular, $\lambda_s(y_n)=|\overline{y}_n|^s$ for all $0<s\le r_0$);
\item $\lambda_s(y)=\max_{n\ge 0} \{p^{-n}\lambda_s(y_n)\}$ for all $0<s\le r$.
\end{itemize}
Note that each $\overline{y}_n$ is an element in $R'$. Pick a presentation
\[\overline{y}_n=\sum_{i_1, \ldots, i_d\in \mathbb{Z}_{\ge 0}} \overline{a}^{(n)}_{i_1, \ldots, i_d}\cdot x_1^{i_1}\cdots x_d^{i_d}\]
with $\overline{a}^{(n)}_{i_1, \ldots, i_d}\in R_{\Lambda}$ such that
$
|\overline{y}_n|_{\mathrm{res}, x_1, \ldots, x_d}=\max_{i_1, \ldots, i_d}\{|\overline{a}^{(n)}_{i_1, \ldots, i_d}|\}.
$ 
Hence, $|\overline{a}^{(n)}_{i_1, \ldots, i_d}|\le C_{**}\cdot |\overline{y}_n|$ for all $i_1, \ldots, i_d\in \mathbb{Z}_{\ge 0}$. Consider
\[\widetilde{\overline{y}}_n\coloneqq \sum_{i_1, \ldots, i_d} \widetilde{\overline{a}}^{(n)}_{i_1, \ldots, i_d}\cdot \widetilde{x}_1^{i_1}\cdots \widetilde{x}_d^{i_d}\in {\mathbf{A}'}^r\]
where $\widetilde{\overline{a}}^{(n)}_{i_1, \ldots, i_d}$ is a good lift of $\overline{a}^{(n)}_{i_1, \ldots, i_d}$ to $\mathbf{A}^{r_0}_{\Lambda}$. Let $z_{(1)}\coloneqq \sum_{n=0}^{\infty}p^n\widetilde{\overline{y}}_n\in {\mathbf{A}'}^r$ and let $y_{(1)}\coloneqq y-z_{(1)}$. Note that 
\begin{align*}
\lambda_{r,\mathrm{res}}(z_{(1)})&\le \max_{n\ge 0}\max_{i_1, \ldots, i_d}\left\{p^{-n}\cdot \lambda_r\left(\widetilde{\overline{a}}^{(n)}_{i_1, \ldots, i_d}\right)\right\}=\max_{n\ge 0}\max_{i_1, \ldots, i_d}\left\{p^{-n}\cdot \left|\overline{a}^{(n)}_{i_1, \ldots, i_d}\right|^r\right\} \\
&\le C_{**}^r\cdot\max_{n\ge 0} \{p^{-n}|\overline{y}_n|^r\}=C_{**}^r\cdot\lambda_r(y).
\end{align*}
By the construction, $y_{(1)}=\sum_{n\ge 0} p^n\left(y_n-\widetilde{\overline{y}}_n\right)$. We have $\lambda_r(y_n-\widetilde{\overline{y}}_n)\le \max\{D_1, D_2, D_3\}$ where 
\[D_1=\lambda_r(y_n-[\overline{y}_n]),\]
\[D_2=\lambda_r\left([\overline{y}_n]-\sum_{i_1, \ldots, i_d}[\overline{a}^{(n)}_{i_1, \ldots, i_d}]\cdot [x_1]^{i_1}\cdots[x_d]^{i_d}\right),\]
and
\[D_3=\lambda_r\left(\sum_{i_1, \ldots, i_d}[\overline{a}^{(n)}_{i_1, \ldots, i_d}]\cdot [x_1]^{i_1}\cdots[x_d]^{i_d}-\widetilde{\overline{y}}_n\right).\]
Using the definition of good lifts, we have $D_1\le p^{-\frac{1}{2}}\cdot \lambda_r(y_n)$ and
\[D_3\le p^{-\frac{1}{2}}\cdot \max_{i_1, \ldots, i_d}\left|\overline{a}^{(n)}_{i_1, \ldots, i_d}\right|^r=p^{-\frac{1}{2}}\cdot |\overline{y}_n|^r_{\mathrm{res}, x_1, \ldots, x_d}\le p^{-\frac{1}{2}}\cdot C_{**}^r\cdot |\overline{y}_n|^r=p^{-\frac{1}{2}}\cdot C_{**}^r\cdot \lambda_r(y_n).\]
Moreover, an explicit computation of Witt coordinates implies
\[D_2\le p^{-1}\cdot \max_{i_1, \ldots, i_d}\left|\overline{a}^{(n)}_{i_1, \ldots, i_d}\right|^r\le p^{-1}\cdot C_{**}^r\cdot \lambda_r(y_n).\]
Putting these together, we conclude that $\lambda_r(y_{(1)})\le p^{-\frac{1}{2}}\cdot C_{**}^r\cdot \lambda_r(y)$. (Recall that $p^{-\frac{1}{2}}\cdot C_{**}^r\le p^{-\frac{1}{2}}\cdot C_{**}^{r_0}<1$.)

Applying the same construction as above to $y_{(1)}$ in place of $y$, we obtain $z_{(2)}\in {\mathbf{A}'}^r$ such that 
\[\lambda_{r, \mathrm{res}}(z_{(2)})\le C_{**}^r\cdot \lambda_r(y_{(1)})< C_{**}^r\cdot \lambda_r(y)\]
and such that $y_{(2)}\coloneqq y_{(1)}-z_{(2)}$ satisfies $\lambda_r(y_{(2)})\le p^{-\frac{1}{2}}\cdot C_{**}^r\cdot\lambda_r(y_{(1)})$. Repeating this, we obtain a converging sum $y=\sum_{j=1}^{\infty}z_{(j)}$. In particular, \[\lambda_{r,\mathrm{res}}(y)\le C_{**}^r\cdot \lambda_r(y).\]

\begin{remark} \label{rem: comparison-with-residue-norm-Frobenius-twists}
For every integer $M\ge 1$, we may consider the generators $\varphi^M(x_1), \ldots, \varphi^M(x_d)$ in place of $x_1, \ldots, x_d$. Notice that $\varphi^M(\widetilde{x}_i)$ is also a good lift of $\varphi^M(x_i)$ for each $i$. Since the surjectivity as in Lemma~\ref{lem: Frob decomp A} holds for ${\fA'}^r$ for any $0 < r \le r_0$, we have that $\varphi^M(\widetilde{x}_1), \ldots, \varphi^M(\widetilde{x}_d)$ generate ${\fA'}^r$ over $\fA_{\Lambda}^r$ for all $0 < r \le \frac{r_0}{p^M}$. Thus, using Step 1, the same argument as above implies that, for any $0<r\le \frac{r_0}{p^M}$, we have 
\[\lambda_r\le \lambda_{r,\mathrm{res}, \varphi^M(\widetilde{x}_1), \ldots, \varphi^M(\widetilde{x}_d})\le C_{**}^{p^Mr}\cdot  \lambda_r\]
on ${\mathbf{C}'}^{[r,r]}$.
\end{remark}

\noindent\textbf{Step 3.} (Finish the proof)\\
To proceed, we choose $n_1\ge n_0$ sufficiently large such that $r_{n_1}\le r_0$ and 
\[C_{**}^{r_{n_1}}\cdot C<p^{-\frac{1}{p-1}}.\]
Let $C'\coloneqq C_{**}^{r_{n_1}}\cdot C$. Consider the continuous $\Gamma$-action on ${\mathbf{C}'}^{[r_{n_1}, r_{n_1}]}$. There exists $n'_0\ge n_1$ such that
\[
\lambda_{r_{n_1}}\left((\gamma-1)(\widetilde{x}_i)\right)\le C
\]
for all $\gamma\in \Gamma_{n'_0}$ and all $i=1,\ldots, d$.
We claim that $\Lambda'=(R', \mathbf{A}')$ satisfies condition (c) for $n'_0$ and $C'$.

Let $n\ge n'_0$. Consider the surjection
\[
\mathbf{C}_{\Lambda}^{[r_n, r_n]}\langle T_1, \ldots, T_d\rangle\twoheadrightarrow {\mathbf{C}'}^{[r_n,r_n]}; \qquad T_i\mapsto \varphi^{n-n_1}(\widetilde{x}_i). 
\]
By Remark~\ref{rem: comparison-with-residue-norm-Frobenius-twists}, we obtain 
\[\lambda_{r_n}\le \lambda_{r_n, \mathrm{res}, \varphi^{n-n_1}(\widetilde{x}_1), \cdots, \varphi^{n-n_1}(\widetilde{x}_d)}\le \left(C_{**}^{p^{n-n_1}}\right)^{r_n}\cdot \lambda_{r_n}=C_{**}^{r_{n_1}}\cdot \lambda_{r_n}\]
on ${\mathbf{C}'}^{[r_n, r_n]}$. 
Let $y\in {\mathbf{C}'}^{[r_n, r_n]}$. We need to verify that 
$\lambda_{r_n}((\gamma-1)y)\le C'\cdot \lambda_{r_n}(y)$ 
for all $\gamma\in \Gamma_n$. Indeed, write
\[y=\sum_{i_1, \ldots, i_d\in \mathbb{Z}_{\ge 0}} a_{i_1, \cdots, i_d}\cdot \left(\varphi^{n-n_1}(\widetilde{x}_1)\right)^{i_1}\cdots \left(\varphi^{n-n_1}(\widetilde{x}_d)\right)^{i_d}\]
with $a_{i_1, \cdots, i_d}\in \mathbf{C}_{\Lambda}^{[r_n, r_n]}$ such that
$
\lambda_{r_n, \mathrm{res}, \varphi^{n-n_1}(\widetilde{x}_1), \cdots, \varphi^{n-n_1}(\widetilde{x}_d)}(y)=\max_{i_1, \ldots, i_d} \{\lambda_{r_n}(a_{i_1, \ldots, i_d})\}.
$
Then
\[\max_{i_1, \ldots, i_d} \{\lambda_{r_n}(a_{i_1, \ldots, i_d})\}\le C_{**}^{r_{n_1}}\cdot \lambda_{r_n}(y).\]
For $\gamma\in \Gamma_n$, we have $(\gamma-1)(y)=E_1+E_2$ where 
\[E_1=\sum_{i_1, \ldots, i_d} (\gamma(a_{i_1, \ldots, i_d})-a_{i_1, \ldots, i_d})\cdot \left(\varphi^{n-n_1}(\widetilde{x}_1)\right)^{i_1}\cdots \left(\varphi^{n-n_1}(\widetilde{x}_d)\right)^{i_d}\]
and
\[E_2=\sum_{i_1, \ldots, i_d}\gamma(a_{i_1, \ldots, i_d})\cdot \left(\gamma(\left(\varphi^{n-n_1}(\widetilde{x}_1)\right)^{i_1}\cdots \left(\varphi^{n-n_1}(\widetilde{x}_d)\right)^{i_d})-\left(\varphi^{n-n_1}(\widetilde{x}_1)\right)^{i_1}\cdots \left(\varphi^{n-n_1}(\widetilde{x}_d)\right)^{i_d}\right).\]
By assumption, 
\begin{align*}
\lambda_{r_n}(E_1)\le & \,\,\lambda_{r_n, \mathrm{res}, \varphi^{n-n_1}(\widetilde{x}_1), \cdots, \varphi^{n-n_1}(\widetilde{x}_d)}(E_1)\le C\cdot \max_{i_1, \ldots, i_d}\{\lambda_{r_n}(a_{i_1, \ldots, i_d})\} \\
=&\,\, C\cdot \lambda_{r_n, \mathrm{res}, \varphi^{n-n_1}(\widetilde{x}_1), \cdots, \varphi^{n-n_1}(\widetilde{x}_d)}(y)\le C\cdot C_{**}^{r_{n_1}}\cdot \lambda_{r_n}(y)=C'\cdot \lambda_{r_n}(y),
\end{align*}
and
\begin{align*}
\lambda_{r_n}(E_2) \le & \max_{i_1, \ldots, i_d}\left\{\lambda_{r_n}(a_{i_1, \ldots, i_d})\cdot \lambda_{r_n}\left(\gamma(\left(\varphi^{n-n_1}(\widetilde{x}_1)\right)^{i_1}\cdots\left(\varphi^{n-n_1}(\widetilde{x}_d)\right)^{i_d})-\left(\varphi^{n-n_1}(\widetilde{x}_1)\right)^{i_1}\cdots \left(\varphi^{n-n_1}(\widetilde{x}_d)\right)^{i_d}\right)\right\} \\
= & \max_{i_1, \ldots, i_d}\left\{\lambda_{r_n}(a_{i_1, \ldots, i_d})\cdot \lambda_{r_{n_1}}\left(\gamma(\widetilde{x}_1^{i_1}\cdots \widetilde{x}_d^{i_d})-\widetilde{x}_1^{i_1}\cdots \widetilde{x}_d^{i_d}\right)\right\} \\
\le & \,\, C\cdot \max_{i_1, \ldots, i_d}\left\{\lambda_{r_n}(a_{i_1, \ldots, i_d})\right\} =C\cdot \lambda_{r_n, \mathrm{res}, \varphi^{n-n_1}(\widetilde{x}_1), \cdots, \varphi^{n-n_1}(\widetilde{x}_d)}(y)
\le \,\, C\cdot C_{**}^{r_{n_1}}\cdot \lambda_{r_n}(y)=C'\cdot \lambda_{r_n}(y).
\end{align*}
This finishes the proof of the claim and Theorem~\ref{thm: strongly decompleting under finite etale extensions}.
\end{proof}

\vspace{0.1in}
\section{Relative Sen theory over $\mathbb{B}_{\mathrm{dR}, X}^+(X_{\infty})$}   \label{sec:relative_Tate_Sen}

We retain the notation from the previous section and consider the  following ($\Gamma$-equivariant) commutative diagram of rings from \S \ref{sss:rings_inside_OB_dR}.
\[
\begin{tikzcd}
    A_{X, \infty}[\![t, u_1, \ldots, u_d]\!]^{\nabla=0} \arrow[r] \arrow[d] 
    &  
    \mathbb{B}_{\mathrm{dR}, X}^+(X_{\infty})
    \arrow[d] \\
        A_{X, \infty}[\![t, u_1, \ldots, u_d]\!] \arrow[r]
    &
    \mO\mathbb{B}_{\mathrm{dR}, X}^+(X_{\infty}).
\end{tikzcd}
\]
The goal in this section is to prove the following decompletion results.

\begin{theorem} \label{thm:equivalence_relative_B_dR^+} \hfill
\begin{enumerate}
\item The base change functor  
$
\mathrm{Mod}_{/A_{X, \infty}[\![t, u_1, \ldots, u_d]\!]}^{\Gamma} \lra  \mathrm{Mod}_{/\mO\mathbb{B}_{\mathrm{dR}, X}^+(X_{\infty}) }^{\Gamma} 
$  
from finite projective $\Gamma$-modules over $A_{X, \infty}[\![t, u_1, \ldots, u_d]\!]$ to those over $\mO\mathbb{B}_{\mathrm{dR}, X}^+(X_{\infty})$ is an equivalence of categories. 

\item The base change functor  
$
\mathrm{Mod}_{/A_{X, \infty}[\![t, u_1, \ldots, u_d]\!]^{\nabla=0}}^{\Gamma} \lra  \mathrm{Mod}_{/\mathbb{B}_{\mathrm{dR}, X}^+(X_{\infty}) }^{\Gamma} 
$  
from finite projective $\Gamma$-modules over $A_{X, \infty}[\![t, u_1, \ldots, u_d]\!]^{\nabla=0}$ to those over $\mathbb{B}_{\mathrm{dR}, X}^+(X_{\infty}) $ is an equivalence of categories.
\end{enumerate} 
\end{theorem}

\begin{remark}
When $X$ has good reduction and the $\Gamma$-module in question is finite free over the base, this is essentially a result of \cite{Andreatta-Brinon} using inputs from \cite{Andreatta-Brinon_Overconvergence} on relative Tate--Sen formalism. Our proof generalizes that of \cite[Proposition 3.19]{Andreatta-Brinon} but also relies on some enhanced techniques from \cite{DLLZ2}. In particular, we do not make use of the notion of \emph{finite vectors} from \cite{Andreatta-Brinon} (cf. Lemma~3.13 of \textit{loc. cit..}). \end{remark}

Before the proof of Theorem~\ref{thm:equivalence_relative_B_dR^+}, let us record a corollary that will be used in \S \ref{sec:relative_N_dR}.

\begin{corollary} \label{cor:equivalence_relative_B_dR^+}
Let $0 < s \le r$ where $r$ is sufficiently small, and suppose that $p^{-n} \in [s, r]$ for some integer $n$. Let $\Lambda=(R_{\Lambda}, \mathbf{A}_{\Lambda})$  be a strongly decompleting datum. Let $M$ be a finite projective $\Gamma$-module over $\mathbf{C}_{\Lambda}^{[s,r]}$ equipped with a $\Gamma$-equivariant isomorphism  
\[ 
\sq f\colon M  \otimes_{\mathbf{C}_{\Lambda}^{[s,r]}, \iota_n} \mathbb{B}_{\mathrm{dR}, X}(X_{\infty}) \isom D \otimes_{A_{X, \infty}[\![t, u_1, \ldots, u_d]\!]^{\nabla=0}} \mathbb{B}_{\mathrm{dR}, X}(X_{\infty})
\]
of $\mathbb{B}_{\mathrm{dR}, X}(X_{\infty})$-modules, where $D$ is a finite projective $\Gamma$-module over $A_{X, \infty}[\![t, u_1, \ldots, u_d]\!]^{\nabla=0}$. Then $\sq f$ descends to a $\Gamma$-equivariant isomorphism 
\[ 
M  \otimes_{\mathbf{C}_{\Lambda}^{[s,r]}, \iota_n} A_{X, \infty}[\![t, u_1, \ldots, u_d]\!]^{\nabla=0} [\frac{1}{t}]  \isom D  [\frac{1}{t}] 
\]
of $A_{X, \infty}[\![t, u_1, \ldots, u_d]\!]^{\nabla=0}[\frac{1}{t}]$-modules.  
\end{corollary}

\begin{proof}
   Since $M$ is finite projective over  $\mathbf{C}_{\Lambda}^{[s,r]}$,   there exists a large enough integer $N$ such that 
  \[ \sq f(t^N \cdot M) \subset D \otimes_{A_{X, \infty}[\![t, u_1, \ldots, u_d]\!]^{\nabla=0}} \mathbb{B}_{\mathrm{dR}, X}^+(X_{\infty}). 
  \] 
Therefore, we have a $\Gamma$-equivariant map 
 \[
 \sq f\colon (t^N \cdot M) \otimes_{\mathbf{C}_{\Lambda}^{[s,r]}, \iota_n} \mathbb{B}_{\mathrm{dR}, X}^+(X_{\infty}) \lra  D \otimes_{A_{X, \infty}[\![t, u_1, \ldots, u_d]\!]^{\nabla=0}} \mathbb{B}_{\mathrm{dR}, X}^+(X_{\infty})
 \]
 between finite projective $\mathbb{B}_{\mathrm{dR}, X}^+(X_{\infty})$-modules. This further descends to a $\Gamma$-equivariant map 
 \[
 f\colon (t^N \cdot M) \otimes_{\mathbf{C}_{\Lambda}^{[s,r]}, \iota_n} A_{X, \infty}[\![t, u_1, \ldots, u_d]\!]^{\nabla=0}  \lra  D  
 \]
by the full faithfulness in Theorem~\ref{thm:equivalence_relative_B_dR^+}(2). Finally, $f$ becomes an isomorphism after inverting $t$ by faithfully flatness of the inclusion $A_{X,n}[\![t, u_1, \ldots, u_d]\!]^{\nabla=0} [\frac{1}{t}] \rightarrow \mathbb{B}_{\mathrm{dR}, X} (X_{\infty})$
(cf. Lemma~\ref{lemma:fully_faithful_relative_B_dR}). 
\end{proof}
 
The rest of this section contributes to the proof of Theorem \ref{thm:equivalence_relative_B_dR^+}. For simplicity of exposition, we will focus on the case $d= 1$ in what follows; the difference between the case of $d=1$ and that of general $d$ is entirely notational.

\subsection{Relative normalized traces} \label{ss:relative_TS_setup}

We need a theory of Tate's normalized traces in the relative setting, generalizing the one developed in \cite{Andreatta-Brinon_Overconvergence}. In this section, we adopt the following notation. 

\begin{notation}\label{notation:L_dR^+} \hfill
\begin{enumerate} 
\item  Write $\mathbf{L}_{\mathrm{dR}}^+ \coloneqq \mO\mathbb{B}_{\mathrm{dR}, X}^+(X_{\infty}) \cong \mathbb{B}_{\mathrm{dR}, X}^+(X_{\infty}) [\![u]\!]$ and $\mathbf{L}_{\mathrm{dR}}^{+, \nabla}  \cong \mathbb{B}_{\mathrm{dR}, X}^+(X_{\infty})$. This notation agrees with that in \cite{Andreatta-Brinon}.
\item Following the convention of \textit{loc. cit..}, we put $l_{\mathrm{dR}}^+ = A_{X, \infty} [\![t, u]\!]$ and $l_{\mathrm{dR}}^{+, \nabla} = A_{X, \infty} [\![t, u]\!]^{\nabla = 0}$.  
\item Let $I \subset \mathbf{L}_{\mathrm{dR}}^+$ denote the ideal $I = (t, u)$. By a slight abuse of notation, we also denote the ideal $(t, u) \subset  A_{X, \infty} [\![t, u]\!]$ by $I$. Note that both $\mathbf{L}_{\mathrm{dR}}^+$ and $A_{X, \infty} [\![t, u]\!]$ are $I$-adically complete and we have $\mathbf{L}_{\mathrm{dR}}^+ /I =   \widehat{A}_{X,\infty}$. 
\end{enumerate}
\end{notation}

For each $n\ge 1$, recall that  
$A_{X, n}=A_X\widehat{\otimes}_{K \gr{T}} K_n\langle T^{\frac{1}{p^n}}\rangle$ from Definition \ref{construction:A_psi_X_infty}. Consider the subrings 
\[B_{X, n} \coloneqq \Big(\bigcup_{m \ge 1} K_m \gr{ T^{\frac{1}{p^n}}}   \otimes_{K \gr{T}} A_{X} \Big)^{\wedge}\]
and
\[ C_{X, n} \coloneqq \Big(\bigcup_{m \ge 1} K_n \gr{ T^{\frac{1}{p^m}}}   \otimes_{K \gr{T}} A_{X} \Big)^{\wedge}\]
of $\widehat A_{X, \infty}$, where $()^{\wedge}$ stands for the $p$-adic completion. Note that $B_{X, n}=\widehat K_\infty \widehat\otimes_{K_n} A_{X, n}$ for each $n \ge 1$ and $C_{X, m} = C_{X, n} \otimes_{K_n} K_m$ for each $m\ge n$. Also note that 
\begin{equation} \label{eq:intersection_of_B_n_and_C_n}
    B_{X, n} \cap C_{X, n} = A_{X, n}
\end{equation}
for all $n \ge 1$, where the intersection is taken inside $\widehat{A}_{X, \infty}$.  For each $m \ge n$, the decomposition 
\[
\widehat{K}_{\infty}\langle T^{\frac{1}{p^m}} \rangle = \bigoplus_{i= 0, \ldots, p^{m-n}-1} \widehat{K}_{\infty}\langle T^{\frac{1}{p^n}}\rangle T^{\frac{i}{p^m}}
\]
induces a decomposition 
\[
B_{X, m} = \bigoplus_{i= 0, \ldots, p^{m-n}-1} B_{X, n} T^{\frac{i}{p^m}}. 
\]
This allows us to define the \textit{geometric normalized traces}
\[
P_{m,n} \coloneqq \frac{1}{p^{m-n}} \operatorname{Tr}_{B_{X, m}/B_{X, n}} \: \colon \: B_{X, m} \ra B_{X, n}. 
\]
Note that $P_{m,n}\left(\sum_{i =0, \ldots,   p^{m-n}-1}b_i \cdot T^{\frac{i}{p^m}}\right)
=b_0$.  These maps extend to $B_{X, n}$-linear maps 
\[
P_n\colon \widehat{A}_{X, \infty} \longrightarrow B_{X, n}.
\]
We refer to $P_n$ as the \textit{geometric trace maps}. On the other hand, we have the classical Tate normalized traces $R_{m, n} = [K_m : K_n]^{-1}  \mathrm{Tr}_{K_m/K_n}\colon K_m \ra K_n$ and $R_n\colon \widehat K_\infty\rightarrow K_n$, which extend to $C_{X, n}$-linear maps $R_{m, n}\colon C_{X, m} \ra C_{X, n}$ and $R_n\colon \widehat{A}_{X, \infty}\rightarrow C_{X, n}$. Note that there is $n_{\mathrm{st}} \ge 1$ such that for any $m \ge n \ge n_{\mathrm{st}}$, we have $[K_m : K_n] = p^{m-n}$. We refer to 
\[R_n\colon \widehat{A}_{X, \infty}\longrightarrow C_{X, n}\]
as the \textit{arithmetic trace maps}. 

By construction, for each $n \ge 1$, $B_{X, n}$ is stable under the action of $\Gamma$ (while $C_{X, n}$ is not) and the maps $P_n$ are $\Gamma$-equivariant. Moreover, we have $R_n \circ P_n = P_n \circ R_n$ and this composition 
\[
\tau_n \coloneqq R_n \circ P_n\colon \widehat{A}_{X, \infty} \ra A_{X, n}
\]
is $\Gamma$-equivariant.

\subsection{Decompletion at stage 0} \label{ss:relative_decompletion_degree_0}

The strategy for proving Theorem \ref{thm:equivalence_relative_B_dR^+} is to first prove an equivalence for the base change along $l^+_{\mathrm{dR}}/I^{r+1} \rightarrow \mathbf{L}^+_{\mathrm{dR}}/I^{r+1}$ (resp. $l^{+, \nabla}_{\mathrm{dR}}/t^{r+1} \rightarrow \mathbf{L}^{+, \nabla}_{\mathrm{dR}}/t^{r+1}$) for each $r$, and then take the inverse limit. The case for ``stage 0'' (i.e., $r=0$) is already dealt with in \cite{DLLZ2}. 

\begin{proposition}[Diao--Lan--Liu--Zhu] \label{prop:relative_Sen_theory_degree_0_whole_X} 
The base change functor  \[ \mathrm{Mod}_{/A_{X, \infty} }^{\Gamma } \lra  \mathrm{Mod}_{/\widehat{A}_{X, \infty} }^{\Gamma  } \] from the category of finite projective $\Gamma$-modules over $A_{X, \infty} $ to those over $\widehat{A}_{X, \infty}$ is an equivalence of categories. 
\end{proposition} 

\begin{proof}
This is \cite[Proposition~A.2.2.3]{DLLZ2} (also see \cite[Remark~A.1.3]{DLLZ2}). 
\end{proof}

\begin{proposition}[Diao--Lan--Liu--Zhu]\label{prop:relative_Sen_theory_degree_0} 
Let $U \subset X$ be any rational subset. Then there exists an open subgroup $\Gamma_{U} \subset \Gamma$ (of finite index) such that $U$ is stable under the action of $\Gamma_{U}$ and the following holds. 
 \begin{enumerate} 
    \item For all $d \in \Z_{\ge 1}$, the base change functor induces a natural bijection
 \[
 H^1_{\tu{cont}} (\Gamma_U, \GL_d(A_{U, \infty})) \isom  H^1_{\tu{cont}} (\Gamma_U , \GL_d(\widehat{A}_{U, \infty}))
 \]
     \item Let $M_0$ be a finite free module over  $A_{U, \infty}$ with a continuous semilinear $\Gamma_U$-action and let $M \coloneqq M_0 \otimes_{A_{U, \infty}} \widehat{A}_{U, \infty}$. Then there are natural isomorphisms 
      \[
 H^i_{\tu{cont}} (\Gamma_U,  M_0) \isom  H^i_{\tu{cont}} (\Gamma_U , M) 
 \] for $i = 0$ and $i = 1$. 
 \end{enumerate} 
 In particular, the base change functor  
 $\mathrm{Mod}_{/A_{U, \infty} }^{\Gamma_U, \tu{ f.f.}} \lra  \mathrm{Mod}_{/\widehat{A}_{U, \infty} }^{\Gamma_U, \tu{ f.f.}}$  from finite free $\Gamma$-modules over $A_{U, \infty} $ to those over $\widehat{A}_{U, \infty}$ is an equivalence of categories. 
\end{proposition} 

\begin{proof}
    This follows from  \cite[Proposition~A.2.2.1 \& A.2.1.1]{DLLZ2}. 
\end{proof}

\begin{remark}  
\begin{enumerate}
    \item There is a typo in the argument of \cite[Proposition~A.2.2.1]{DLLZ2}: instead of using the spectral norms on $A_{X, n}$ and $\widehat{A}_{X, \infty}$, one should work with the spectral norms on $A_{\mathbb D_{K}, n}$ and $\widehat{A}_{\mathbb{D}_K, \infty}$ in the case of the unit disk, and then take the completed tensor product norms on $A_{X, n} = A_{\mathbb D_{K}, n} \otimes_{K\gr{T}} A_{X}$ and $\widehat{A}_{X, \infty} = \widehat{A}_{\mathbb{D}_K, \infty} \widehat \otimes_{K\gr{T}} A_{X}$ (cf. \cite[Definition~2.1.10]{Kedlaya-LiuI}). The proof in \cite{DLLZ2} in fact refers to this product norm.  
    \item 
Part (1) of Proposition \ref{prop:relative_Sen_theory_degree_0} also follows from the relative Tate--Sen formalism of \cite{Andreatta-Brinon_Overconvergence}, using \cite[Theorem~2.4]{Andreatta-Brinon_Overconvergence} (with Proposition~2.6 of \textit{loc. cit..} corrected by \cite[Appendix~II]{Andreatta-Brinon}). More precisely, one applies the Tate--Sen formalism with 
\begin{enumerate}
    \item $\sq \Lambda = \widehat A_{X, \infty}$; $ \: G' = G = \Gamma$; $\: H = \{1\}$; $\: d = 1$; $m_0 = m_0 (G') = 0$. 
    
    \item $\Lambda_{n}^{(0)} = C_{X, n}; \: \Lambda_n^{(1)} = B_{X, n}$;    $\tau^{(0)}_{n} = R_n; \: \tau^{(1)}_{n} = P_n$ for $n \ge n_{\mathrm{st}}$.
    
\end{enumerate}  
\end{enumerate} 
\end{remark}

\subsection{Decompletion at stage $r$} \label{ss:all_degrees}

Now we analyze the base change functors along $l^+_{\mathrm{dR}}/I^{r+1} \rightarrow \mathbf{L}^+_{\mathrm{dR}}/I^{r+1}$ and $l^{+, \nabla}_{\mathrm{dR}}/t^{r+1} \rightarrow \mathbf{L}^{+, \nabla}_{\mathrm{dR}}/t^{r+1}$, for general $r$.

\begin{theorem} \label{thm:relative_Sen_theory:over_L_finite_level} 
 Let $r \ge 0$ be an integer. Recall that $I = (t, u)$. 
\begin{enumerate}
    \item  The base change functor 
    $\mathrm{Mod}_{/(l_{\mathrm{dR}}^+ /I^{r+1}) }^{\Gamma } \lra  \mathrm{Mod}_{/(\mathbf{L}_{\mathrm{dR}}^+ /I^{r+1})}^{\Gamma  }$  from the  category of finite projective $\Gamma$-modules over $l_{\mathrm{dR}}^+ /I^{r+1} = A_{X, \infty} [\![t, u]\!]/(t, u)^{r+1}$ to those over $\mathbf{L}_{\mathrm{dR}}^+ /I^{r+1}=\mathcal{O}\mathbb{B}^+_{\mathrm{dR}, X}(X_{\infty})/I^{r+1}$ is an equivalence of categories. 
    \item Similarly, the base change functor  $\mathrm{Mod}_{/(l^{+, \nabla}_{\mathrm{dR}}/t^{r+1})}^{\Gamma} \lra  \mathrm{Mod}_{/(\mathbf{L}^{+, \nabla}_{\mathrm{dR}}/t^{r+1})}^{\Gamma} $ from the category of finite projective $\Gamma$-modules over $l^{+, \nabla}_{\mathrm{dR}}/t^{r+1}$ to those over $\mathbf{L}^{+, \nabla}_{\mathrm{dR}}/t^{r+1}$ is an equivalence of categories. 
\end{enumerate}
\end{theorem} 

\begin{proof} 
We prove (1) by induction on $r$. Part (2) will follow from a similar argument. The case $r = 0$ is taken care of by \S \ref{ss:relative_decompletion_degree_0}. Suppose that the assertions hold up to $r -1$ for some $r \ge 1$. We prove the claim for $r$.  \\

\noindent \textbf{Step 1}. \\ 
Let us first prove the full faithfulness for finite free modules. First suppose that $M_0$ is a finite free $\Gamma$-module over $l_{\mathrm{dR}}^+ /I^{r+1}$ and let 
$
M \coloneqq M_0 \otimes_{l_{\mathrm{dR}}^+ /I^{r+1}} \mathbf{L}_{\mathrm{dR}}^+ /I^{r+1}.
$
Let $M_0' \coloneqq I^r M_0$ and $M_0'' \coloneqq M_0/M_0'$, and similarly define $M'$ and $M''$.  Consider the following commutative diagram of $\Gamma$-equivariant short exact sequences 
\[
\begin{tikzcd}
        0 \arrow[r]   &  M_0' \arrow[d]  \arrow[r] &  M_0  \arrow[r] \arrow[d]  &  M_0'' \arrow[d] \arrow[r] &  0 \\ 
    0 \arrow[r] &  M'  \arrow[r] &  M \arrow[r] &  M'' \arrow[r] &  0 
\end{tikzcd}
\]
Note that $M''_0$ is a finite free $l_{\mathrm{dR}}^+ /I^{r}$-module and $M'' = M''_0  \otimes_{l_{\mathrm{dR}}^+ /I^{r}} \mathbf{L}_{\mathrm{dR}}^+ /I^{r}$. Similarly, $M'_0$ is finite free over $A_{X, \infty}$ and $M' =M'_0 \otimes_{A_{X, \infty}} \widehat{A}_{X, \infty}$.\footnote{
Note that we have decompositions
\begin{align*}
    I^r l_{\mathrm{dR}}^+/ I^{r+1} l_{\mathrm{dR}}^+ = \bigoplus_{a + b = r} A_{X, \infty} t^a u^b; \quad 
    I^r \mathbf{L}_{\mathrm{dR}}^+/ I^{r+1} \mathbf{L}_{\mathrm{dR}}^+ = \bigoplus_{a + b = r} \widehat{A}_{X, \infty} t^a u^b
\end{align*}
as $A_{X, \infty}$- and $\widehat{A}_{X, \infty}$-modules, respectively.} 
Taking $\Gamma$-invariants, we obtain a commutative diagram 
\[
\begin{tikzcd}
        0 \arrow[r]   &  (M_0')^{\Gamma} \arrow[d, "h_1"]  \arrow[r] &  (M_0)^{\Gamma}  \arrow[r] \arrow[d, "h_2"]  &  (M_0'')^{\Gamma} \arrow[d, "h_3"] \arrow[r] &  H^1_{\tu{cont}} (\Gamma, M_0')  \arrow[d, "h_4"] \\ 
    0 \arrow[r] &  (M')^{\Gamma}  \arrow[r] &  (M)^{\Gamma} \arrow[r] &  (M'')^{\Gamma} \arrow[r] &   H^1_{\tu{cont}} (\Gamma, M') 
\end{tikzcd}
\]
with exact rows.  In the diagram above, the map $h_3$ is an isomorphism by induction hypothesis, $h_1$ and $h_4$ are isomorphisms by Proposition \ref{prop:relative_Sen_theory_degree_0}. This implies that $h_2$ is also an isomorphism. \\

\noindent \textbf{Step 2}. \\ 
Now we prove full faithfulness in general. Let $M_0$ be a finite projective $\Gamma$-module over $l_{\mathrm{dR}}^+ /I^{r+1} = A_{X, \infty} [\![t, u]\!]/(t, u)^{r+1}$ and let 
$ 
M\coloneqq M_0 \otimes_{l_{\mathrm{dR}}^+ /I^{r+1}} \mathbf{L}_{\mathrm{dR}}^+ /I^{r+1}.
$ 
It suffices to show that there exists an open subgroup $\Gamma' \subset \Gamma$ of finite index, such that $(M_0)^{\Gamma'} = M^{\Gamma'}$. Let $Q_0\coloneqq M_0/I M_0$ and $Q \coloneqq M/IM$. Since $Q_0$ is a finite projective $\Gamma$-module over $l_{\mathrm{dR}}^+ /I = A_{X, \infty}$ and $\Gamma$ is topologically finitely generated, there exists a sufficiently large integer $N$ such that, after replacing $\Gamma$ by an open subgroup $\Gamma' \subset  \Gamma$ of finite index, $Q_0$ is the base change of a finite projective $\Gamma'$-module $L_N$ over $A_{X, N}$ (as a $\Gamma'$-module over $A_{X, \infty}$). In particular, after further shrinking $\Gamma'$, there exists a finite covering of $X_N = \spa(A_{X, N}, A_{X, N}^\circ)$ by rational subsets $\{U_i\}_{i \in \Sigma}$ (where each $U_i = \spa(A_{U_i}, A_{U_i}^+)$ is stable under the action of $\Gamma'$), such that $L_N$ becomes a finite free $\Gamma'$-module over each $U_i$. Using the  toric tower  pulled back from $X_\infty \ra X_N$ via $U_i \ra X_N$, we obtain $A_{U_i, \infty}$ and $\widehat{A}_{U_i, \infty}$ as in Definition~\ref{construction:A_psi_X_infty}, where 
$A_{U_i, \infty} = A_{X, \infty} \otimes_{A_{X, N}} A_{U_i}$. For each $i\in \Sigma$, let 
\begin{itemize}
    \item $Q_{0, U_i} \coloneqq Q_{0} \otimes_{A_{X, \infty}} A_{U_i, \infty}$;   
    \item $M_{0, U_i} \coloneqq M_0  \otimes_{(A_{X, \infty} [\![t, u ]\!]/I^{r+1})} A_{U_i, \infty}  [\![t, u ]\!]/I^{r+1}$; 
    \item $M_{U_i} \coloneqq M_{0, U_i}  \otimes_{(A_{U_i, \infty} [\![t, u ]\!]/I^{r+1})} \mO\mathbb{B}_{\mathrm{dR}, X}^+(U_{i, \infty}) /I^{r+1}$
\end{itemize} 
where $U_{i, \infty}=U_i\times_{X_N}X_{\infty}\in X_{\proet}$ is an affinoid perfectoid object. By construction, each $Q_{0, U_i} = M_{0, U_i}/I$ is finite free over $A_{U_i, \infty}$. This in turn implies that each $M_{0, U_i}$ is finite free over $A_{U_i, \infty} [\![t, u ]\!]/I^{r+1}$, since $I$ is nilpotent. By Step 1, we know that $(M_{0, U_i})^{\Gamma'} = (M_{U_i})^{\Gamma'}$ for all $i$. To proceed, let us write $U_{ij}$ for the overlap $U_{ij} = U_i \cap U_j$, which is again a rational subset of $X$. We have $(M_{0, U_{ij}})^{\Gamma'} =  (M_{ U_{ij}})^{\Gamma'}$, again by Step 1. Note that $t, u$ form a regular sequence in $l_{\mathrm{dR}}^+$. Also note that for each $s \ge 1$ and each rational subset $U\in \{U_i\}_{i\in \Sigma}\cup\{U_{ij}\}_{i,j\in \Sigma}$, we have a short exact sequence 
\[
0 \ra (I^{s-1} M_0/I^s)_{U} \ra (M_0/I^s)_{U} \ra (M_0/I^{s-1})_U \ra 0
\]
when restricting to $U$, where $I^{s-1} M_0 / I^s \cong (M_0/I)^{\oplus s}$ by construction. Thus we have 
\[M_0 = \tu{Eq} \left(\prod_i M_{0, U_i} \rightrightarrows \prod_{i, j} M_{0, U_{ij}} \right)
\] 
by induction on $r$. This in turn implies that 
\[(M_0)^{\Gamma'}= \tu{Eq} \left(\prod_i (M_{0,  U_i})^{\Gamma'} \rightrightarrows \prod_{i, j} (M_{0,  U_{ij}})^{\Gamma'} \right)
\] 
after taking the $\Gamma'$-invariants. The same assertions hold for $M$ by faithfully flatness of  $l_{\mathrm{dR}}^+ \ra \mathbf{L}_{\mathrm{dR}}^+$. Consequently, we  have  $(M_0)^{\Gamma'} = M^{\Gamma'}$ as desired.  \\

\noindent \textbf{Step 3}. \\
Next we prove essential surjectivity. Let us first prove this in the finite free case. 
Let $X$ be a finite free $\Gamma$-module over $\mathbf{L}_{\mathrm{dR}}^+ /I^{r+1}$. Let $X' \coloneqq I^r X$ and $X'' \coloneqq X/X'$. Then we have a  $\Gamma$-equivariant short exact sequence 
\begin{equation} \label{eq:extension_of_X''_by_X'_1}
0 \ra X' \ra X \ra X'' \ra 0
\end{equation}
of $\mathbf{L}_{\mathrm{dR}}^+/I^{r+1}$-modules as in Step 1, where $X'$ is a finite free $\Gamma$-module over $\widehat{A}_{X, \infty}$ of rank $d (r+1)$ while $X''$ is finite free over $\mathbf{L}_{\mathrm{dR}}^+/I^{r}$ of rank $d$. By the induction hypothesis, we know that $X'$ and $X''$ descend to $\Gamma$-modules $X_0'$ and $X_0''$, over ${A}_{X, \infty}$ and $l_{\mathrm{dR}}^+/I^{r}$ respectively. Let $\cl{\mathfrak B} =  (\cl e_1, \ldots, \cl e_d)$ be a basis of $X''_0$ over $l_{\mathrm{dR}}^+/I^r$ (and thus a basis of $X''$ over $\mathbf{L}_{\mathrm{dR}}^+/I^r$). Let $e_i \in X$ be an element lifting $\cl e_i$ for each $i$. Then $(e_1, \ldots, e_d)$ form a basis of $X$ over $\mathbf{L}_{\mathrm{dR}}^+/I^{r+1}$. Write 
$ U\colon \Gamma \ra \GL_d (\mathbf{L}_{\mathrm{dR}}^+/I^{r+1})$ 
for the 1-cocycle describing the action of $\Gamma$ on $X$ relative to this basis. Our goal is to show that, up to a change of basis, the image of $U$ lies in $\GL_d (l_{\mathrm{dR}}^+/I^{r+1})$. In other words, we want to prove that the base change functor induces a natural bijection (of pointed sets) 
 \[
 H^1_{\tu{cont}} (\Gamma, \GL_d(l_{\mathrm{dR}}^+/ I^{r+1})) \isom  H^1_{\tu{cont}} (\Gamma, \GL_d( \mathbf{L}_{\mathrm{dR}}^+/ I^{r+1})).
 \]
 
First, let $\cl U\colon \Gamma \ra \GL_d (\mathbf{L}_{\mathrm{dR}}^+/I^{r})$ denote the reduction mod $I^r$, which has coefficients in $ l_{\mathrm{dR}}^+/I^{r}$. Using the isomorphism $A_{X, \infty} [t, u]^{\deg < r} \isom  l_{\mathrm{dR}}^+/I^{r} $ of $A_{X, \infty}$-modules, for each $\gamma \in \Gamma$,  we can write $U(\gamma)$ (uniquely) as
\[
U(\gamma ) = P(\gamma) + V(\gamma)
\]
with $P (\gamma) \in M_{d \times d} (A_{X, \infty} [t, u]^{\deg < r} )$ and $V(\gamma) \in M_{d \times d} ( I^r \mathbf{L}_{\mathrm{dR}}^+/ I^{r+1} \mathbf{L}_{\mathrm{dR}}^+)$.\footnote{This notation simply means that each entry of $V(\gamma)$ lives in $ I^r \mathbf{L}_{\mathrm{dR}}^+/ I^{r+1} \mathbf{L}_{\mathrm{dR}}^+ \subset  \mathbf{L}_{\mathrm{dR}}^+/ I^{r+1} \mathbf{L}_{\mathrm{dR}}^+$.} Moreover,  there exists a sufficiently large integer $N \ge n_{\mathrm{st}}$ such that $P(\gamma) \in M_{d \times d} (A_{X, N} [t, u]^{\deg < r} )$ for all $\gamma \in \Gamma$. 
 
Next, for $n \ge n_{\mathrm{st}}$, we extend the $\Gamma$-equivariant map $\tau_n\coloneqq  P_n \circ R_n\colon \widehat A_{X, \infty} \ra A_{X, n}$ (cf. \S \ref{ss:relative_TS_setup}) to an $A_{X, n}$-linear map  \[ \tau_n\colon  I^r \mathbf{L}_{\mathrm{dR}}^+/ I^{r+1} \mathbf{L}_{\mathrm{dR}}^+   = \bigoplus_{a + b = r} \widehat{A}_{X, \infty} t^a u^b  \ra  \bigoplus_{a + b = r} A_{X, n} t^a u^b  
 \]
by setting $\tau_n(t^a u^b x) = t^a u^b  (P_n \circ R_n(x) )$ for each $n \ge n_{\mathrm{st}}$. 
The extended map $\tau_n$ is again $\Gamma$-equivariant. 
For each $n \ge n_{\mathrm{st}}$, let us view $ \bigoplus_{a + b = r} A_{X, n} t^a u^b $ as a submodule of $A_{X, n} [\![t, u]\!]/I^{r+1}$. Using the map $\tau_N$ (with the  integer $N$ chosen earlier), we can modify $U$ into a new map $U^{(N)}\colon \Gamma \ra \GL_d( A_{X, N} [\![t, u]\!]/I^{r+1}) $ by setting 
\[
U^{(N)} (\gamma) \coloneqq P(\gamma) + \tau_N (V(\gamma)).
\]
Here, $\tau_N$ is defined entry-wise on the matrix $V(\gamma)$. It is straightforward to verify that  
\begin{equation}\label{eq: UN 1-cocycle}
U^{(N)} (\gamma \eta) = U^{(N)} (\gamma) \cdot \gamma (U^{(N)}(\eta))
\end{equation}
for all $\gamma, \eta \in \Gamma$; that is, $U^{(N)}$ is still a 1-cocycle (although $\tau_N$ is not a ring homomorphism). Now, using the same basis $(e_1,\ldots , e_d)$, the continuous $1$-cocycle $U^{(N)}$ corresponds to a finite free module $X_0$ over $l_{\mathrm{dR}}^+ /I^{r+1}$ with a $\Gamma$-action, such that its $\Gamma$-action descends to $A_{X, N} [\![t, u]\!]/I^{r+1}$. Moreover, if $x = \sum \alpha_i e_i \in X'$ is an element in $X'$ with $\alpha_i \in I^r$, then we have  
\[
\left( U(\gamma) - U^{(N)} (\gamma) \right) (x) = \sum_i \beta_i e_i = 0 
\]
since $\beta_i \in I^{2r} = 0$, thus $U(\gamma) (x) = U^{(N)} (\gamma) (x)$. This implies that $I^rX_0 = X' $ as a $\Gamma$-representation. In other words, we have constructed another extension 
\begin{equation}
     \label{eq:extension_of_X''_by_X'_2}
    0 \ra X' \ra X_0 \ra X'' \ra 0 
\end{equation}
of $X''$ by $X'$.

Let $c \in \tu{Ext}^1_{\widehat{A}_{X, \infty}, \Gamma} (X'', X')$ denote the extension class given by  $X$ in (\ref{eq:extension_of_X''_by_X'_1}) and let $c_0$ denote the extension class given by (\ref{eq:extension_of_X''_by_X'_2}). Let $\Delta(\gamma) \coloneqq U(\gamma) - U^{(N)} (\gamma) \in M_{d \times d} (I^r \mathbf{L}_{\mathrm{dR}}^+/I^{r+1}\mathbf{L}_{\mathrm{dR}}^+)$ denote the difference of the two $1$-cocycles $U$ and $U^{(N)}$. Then 
\[\Delta (\gamma \eta)  
= U^{(N)} (\gamma) \cdot \gamma(\Delta(\eta)) + \Delta(\gamma) \cdot \gamma(U^{(N)} (\eta)) \]
for all $\gamma, \eta \in \Gamma$. This satisfies the cocycle relation of an extension of $X''$ by $X'$ as a $\Gamma$-representation as in \cite[Remark~3.16 \& 3.18]{Andreatta-Brinon}, and in fact defines the extension class  $c - c_0$. We claim that there exists some matrix $N_\epsilon \in M_{d \times d} (I^r \mathbf{L}_{\mathrm{dR}}^+/I^{r+1}\mathbf{L}_{\mathrm{dR}}^+)$ such that 
\begin{equation}
    \label{eq:expression_of_U_epsilon}
U_\epsilon (\gamma) \coloneqq  
\Delta(\gamma) + U^{(N)}(\gamma) \cdot \gamma(N_\epsilon) - N_\epsilon \cdot U^{(N)} (\gamma)
 \end{equation}
has entries in $I^r l_{\mathrm{dR}}^+/I^{r+1}l_{\mathrm{dR}}^+$. 
To prove the claim, we view $M_{d \times d} (I^r\mathbf{L}_{\mathrm{dR}}^+/I^{r+1}\mathbf{L}_{\mathrm{dR}}^+)$ as a finite free $\widehat A_{X, \infty}$-module and consider a twisted action of $\Gamma$ on $M_{d \times d} (I^r\mathbf{L}_{\mathrm{dR}}^+/I^{r+1}\mathbf{L}_{\mathrm{dR}}^+)$ defined by
\begin{equation}\label{eq: twisted action}
\gamma \star M \coloneqq  U^{(N)}(\gamma) \cdot \gamma (M)\cdot (U^{(N)}(\gamma))^{-1}
\end{equation}
for $\gamma\in \Gamma$ and $M\in M_{d \times d} (I^r\mathbf{L}_{\mathrm{dR}}^+/I^{r+1}\mathbf{L}_{\mathrm{dR}}^+)$. Since $U^{(N)}$ defines a 1-cocycle, the twisted action \eqref{eq: twisted action} is indeed a $\Gamma$-action. Now  we consider a map  $\delta\colon \Gamma \ra M_{d \times d} (I^r\mathbf{L}_{\mathrm{dR}}^+/I^{r+1}\mathbf{L}_{\mathrm{dR}}^+)$ defined by 
\[ \delta (\gamma) \coloneqq \Delta(\gamma) \cdot (U^{(N)} (\gamma))^{-1} = U(\gamma)  \cdot (U^{(N)} (\gamma))^{-1} - 1.\] 
Using \eqref{eq: UN 1-cocycle}, one computes that
\[
\delta (\gamma \eta)   = \delta (\gamma) + U(\gamma) \cdot \gamma (\delta(\eta)) \cdot U^{(N)} (\gamma)^{-1}  =   \delta (\gamma) + U^{(N)}(\gamma) \cdot \gamma (\delta(\eta)) \cdot U^{(N)} (\gamma)^{-1} = \delta(\gamma) + \gamma \star \delta (\eta).
\]
where the second equality uses the fact that $U(\gamma) - U^{(N)}(\gamma)$ has entries in $I^r\mathbf{L}_{\mathrm{dR}}^+/I^{r+1}\mathbf{L}_{\mathrm{dR}}^+$. That is, $\delta$ defines a $1$-cocycle (with respect to the action \eqref{eq: twisted action}). By Proposition~\ref{prop:relative_Sen_theory_degree_0}(2), we know that $\delta$ is in fact cohomologous to a continuous $1$-cocycle $c\colon \Gamma \ra  M_{d \times d} (I^r l_{\mathrm{dR}}^+/I^{r+1}l_{\mathrm{dR}}^+)$. In other words, there exists a coboundary of the form $\gamma \mapsto \gamma \star  N_{\epsilon} - N_{\epsilon}$ for some $N_\epsilon \in M_{d \times d} (I^r \mathbf{L}_{\mathrm{dR}}^+/I^{r+1}\mathbf{L}_{\mathrm{dR}}^+)$ 
such that 
$
c(\gamma) = \delta (\gamma) + \gamma \star  N_{\epsilon} - N_{\epsilon}.
$ 
A direct computation yields
\[
\Delta(\gamma) + U^{(N)}(\gamma) \cdot \gamma(N_\epsilon) - N_\epsilon \cdot U^{(N)} (\gamma)=c(\gamma) \cdot U^{(N)} (\gamma)\in M_{d \times d}(I^r l_{\mathrm{dR}}^+/I^{r+1}l_{\mathrm{dR}}^+);
\]
hence, this $N_\epsilon$ satisfies (\ref{eq:expression_of_U_epsilon}) as desired. 

Finally, we modify the basis $(e_1, \ldots , e_d)$ using the change-of-basis matrix $1+ N_{\epsilon}$. Under this new basis, the action of $\Gamma$ is given by 
\[
   (1+N_{\epsilon})^{-1} \cdot U (\gamma)  \cdot \gamma (1+N_{\epsilon})
= (1- N_{\epsilon})  \cdot \left(U^{(N)}(\gamma)  + \Delta(\gamma) \right) \cdot \gamma (1+N_{\epsilon}) =   U^{(N)}(\gamma)  + U_{\epsilon}(\gamma),
\]
whose coefficients lie in $l_{\mathrm{dR}}^+/I^{r+1}$ by construction. This finishes the proof of essential surjectivity in the finite free case. \\

\noindent \textbf{Step 4}. \\
Lastly, we prove essential surjectivity in general. By Lemma~\ref{lem: general-fact-equivalence-phi-Gamma-modules}, it suffices to prove the statement after replacing $\Gamma$ by some open normal subgroup of finite index. Let $M$ be a finite projective $\Gamma$-module over $\mathbf{L}_{\mathrm{dR}}^+ /I^{r+1}$ and we put $Q \coloneqq M/I$. By the induction hypothesis (cf. Proposition~\ref{prop:relative_Sen_theory_degree_0_whole_X}), there exists a (unique) finite projective $\Gamma$-module $Q_0$  over $A_{X, \infty}$ such that $Q = Q_0 \otimes_{A_{X, \infty}} \widehat{A}_{X, \infty}$. Next,   as in Step 2, we pick a sufficiently large integer $N \ge n_{\mathrm{st}}$, an open normal subgroup $\Gamma' \subset \Gamma$ of finite index, and a finite covering $\{U_i\}_{i \in \Sigma}$ of $X_N$ (by rational subsets) stable under the action of $\Gamma'$, such that
\begin{itemize}
    \item $Q_0$ further descends to a finite projective $\Gamma'$-module $L_N$ over $A_{X, N}$, and 
    \item $L_N$ becomes finite free $\Gamma'$-module over each $U_i$.
\end{itemize}     
In particular, $Q \otimes_{ \widehat{A}_{X, \infty}}  \widehat{A}_{U_i, \infty} \cong L_N \otimes_{A_{X, N}} \widehat{A}_{U_i, \infty}$  is a finite free $\Gamma'$-module over $ \widehat{A}_{U_i, \infty}$ for each $i$. This further implies that $M$ becomes finite free (as a $\Gamma'$-module) over $\mO\mathbb{B}_{\mathrm{dR}, X}^+(U_{i, \infty}) /I^{r+1}$ for each $i$. Let $U_{ij}$ denote the overlap $U_i\cap U_j$ as in Step 2. For each $U\in \{U_i\}_{i\in \Sigma}\cup \{U_{ij}\}_{i,j\in \Sigma}$, let
\[M_{U} \coloneqq M \otimes_{(\mathbf{L}_{\mathrm{dR}}^+/I^{r+1})} \mO\mathbb{B}_{\mathrm{dR}, X}^+(U_{\infty}) /I^{r+1}.\]
By Step 1--3, $M_U$ descends to a (unique) finite free $\Gamma'$-module $M_{0, U}$ over $A_{U, \infty} [\![t, u]\!]/(t, u)^{r+1}$. Moreover, by full faithfulness (Step 1), this descent is functorial with respect to inclusions of the form $U_{ij} \ra U_i$ for each $i, j \in \Sigma$. 
This gives rise to a finite projective module $M_0$ defined as the equalizer 
\[M_0 \coloneqq \tu{Eq} \left(\prod_i M_{0, U_i} \rightrightarrows \prod_{i, j} M_{0, U_{ij}} \right),
\] 
which is stable under the $\Gamma'$-action (viewed as a submodule of $M$), since all the maps in the diagram are $\Gamma'$-equivariant. This finishes the proof of part (1) of the theorem. As mentioned earlier, part (2) follows from a similar argument.  
\end{proof}

\subsection{Decompletion over $\mathbb{B}_{\mathrm{dR}, X}^+(X_{\infty})$ and $ \mathcal{O}\mathbb{B}_{\mathrm{dR}, X}^+(X_{\infty})$} \label{subsec: decompletion-BdR+}

We are ready to prove Theorem~\ref{thm:equivalence_relative_B_dR^+}.

\begin{proof}[Proof of Theorem~\ref{thm:equivalence_relative_B_dR^+}] 
We only prove part (1), as the proof of part (2) is identical. First, the full faithfulness is  immediate from Theorem~\ref{thm:relative_Sen_theory:over_L_finite_level} as taking limits over $r$ commutes with taking $\Gamma$-invariants. It remains to prove essential surjectivity. Let $X$ be a finite projective $\Gamma$-module over $\mathbf{L}_{\mathrm{dR}}^{+}$.  Let $X_r \coloneqq X /I^{r+1} X$ and let $D_r$ denote the unique descent of $X_r$ over $l_{\mathrm{dR}}^+/I^{r+1}$. The natural projection $X_{r+1} \ra X_r$ induces a map $\pi_r\colon D_{r+1} \ra D_r$ by Theorem~\ref{thm:relative_Sen_theory:over_L_finite_level}.   The induced map 
\[h_r\colon D_{r+1} \otimes_{l_{\mathrm{dR}}^+/I^{r+2}} l_{\mathrm{dR}}^+/I^{r+1} \lra D_r \] 
is in fact a $\Gamma$-equivariant isomorphism, since it 
becomes an isomorphism after taking base change along $l_{\mathrm{dR}}^+/I^{r+1} \ra  \mathbf{L}_{\mathrm{dR}}^+/I^{r+1}$, which 
is faithfully flat (cf. Lemma~\ref{lemma:fully_faithful_relative_B_dR}). In particular, $\pi_r$ is surjective and $D_r = D_{r+1}/I^{r+1}$. Finally, taking inverse limit, we obtain the desired finite projective $\Gamma$-module $D = \varprojlim_{r} D_r$ over $l_{\mathrm{dR}}^+$. This finishes the proof of part (1). 
\end{proof}

\vspace{0.1in}
\section{Relative $p$-adic differential equations} \label{sec:relative_N_dR} 

The goal of this section is to construct relative versions of $\mathbf{N}_{\mathrm{dR}}$ and $\widetilde{\mathbf{N}}_{\mathrm{dR}}$ functors (see \S \ref{sss: NdR classical setting} and \S \ref{sss:tilde NdR classical setting}). In particular, for a de Rham local system $\L$, the resulting relative $(\varphi, \Gamma)$-module $\mathbf{N}_{\mathrm{dR}}(\L)$ can be viewed as a family version of $p$-adic differential equations.

Let $h\colon X\rightarrow \D^d_K$ be a toric chart and let $\Lambda=(R_{\Lambda}, \mathbf{A}_{\Lambda})$ be a strongly decompleting data with associated imperfect period rings \[\mathbf{A}^r_{\Lambda}, \mathbf{A}^{\dagger}_{\Lambda}, \mathbf{B}^r_{\Lambda},\mathbf{B}^{\dagger}_{\Lambda}, \mathbf{C}_{\Lambda}^{[s,r]}, \mathbf{C}^r_{\Lambda},\mathbf{C}^{\dagger}_{\Lambda}.\] We have seen that, from any \'etale $\mathbb{Z}_p$-local system $\mathbb{L}$ on $X$, one can associate an \'etale $(\varphi, \Gamma)$-module $\widetilde{M}^{\dagger} = \widetilde{D}_{\mathrm{rig}}^{\dagger} (\L)$ over $\widetilde{\mathbf{C}}^{\dagger}_X$ (see Corollary~\ref{cor: equiv-local-system-overconvergent-etale-(phi, Gamma)-mods}), which uniquely descends to an \'etale $(\varphi, \Gamma)$-module $M^{\dagger}=D^{\dagger}_{\mathrm{rig}}(\L)$ over $\mathbf{C}^{\dagger}_{\Lambda}$ by Corollary~\ref{cor: decompletion-etale-(phi, Gamma)-modules} and Theorem~\ref{thm: equivalence-isogeny-local-system-overconvergent-etale-(phi, Gamma)-modules-C}. In particular, for all $r$ sufficiently small and $s \in (0, r]$, we obtain a compatible family of $\widetilde{\mathbf{C}}^{[s, r]}_X$-modules $\widetilde{M}^{[s, r]}$ which descend to $\mathbf{C}^{[s, r]}_{\Lambda}$-modules $M^{[s, r]}$, compatible with $\varphi$- and $\Gamma$-structures. 

In \S \ref{subsec: tilde-NdR} (resp. \S \ref{subsec: NdR-preparation}--\ref{subsec: NdR-BL gluing}), we will construct a functor 
\[
\widetilde{\mathbf{N}}_{\mathrm{dR}}\colon 
\mathrm{Loc}^{\mathrm{dR}}_{\Z_p}(X) \lra \mathrm{Mod}^{(\varphi, \Gamma)}_{/\widetilde{\mathbf{C}}^{\dagger}_{X}}    
\]
\[(\,\textup{resp. }\mathbf{N}_{\mathrm{dR}}\colon 
\mathrm{Loc}^{\mathrm{dR}}_{\Z_p}(X) \lra \mathrm{Mod}^{(\varphi, \Gamma)}_{/\mathbf{C}^{\dagger}_{\Lambda}} \,) \]
from de Rham local systems to finite projective $(\varphi, \Gamma)$-modules over $\widetilde{\mathbf{C}}^{\dagger}_{X}$ (resp. over $\mathbf{C}^{\dagger}_{\Lambda}$) such that $\widetilde{\mathbf{N}}_{\mathrm{dR}} [\frac{1}{t}] \cong \sq M^\dagger [\frac{1}{t}]$ (resp. $\mathbf{N}_{\mathrm{dR}} [\frac{1}{t}] \cong M^{\dagger} [\frac{1}{t}]$). The constructions are similar to the ones in \S \ref{sss:tilde NdR classical setting} and \S \ref{sss: NdR classical setting}; namely, we modify $\widetilde{M}^{[s,r]}$ (resp. $M^{[s,r]}$) along the ``$t= 0$'' locus through Beauville--Laszlo gluing. We remark that although $\mathbf{N}_{\mathrm{dR}}(\L)$ is a descent of $\widetilde{\mathbf{N}}_{\mathrm{dR}}(\L)$ from $\widetilde{\mathbf{C}}^{\dagger}_{X}$ to $\mathbf{C}^{\dagger}_{\Lambda}$, we cannot construct $\mathbf{N}_{\mathrm{dR}}(\L)$ from $\widetilde{\mathbf{N}}_{\mathrm{dR}}(\L)$ using Corollary~\ref{cor: decompletion-etale-(phi, Gamma)-modules} because $
\widetilde{\mathbf{N}}_{\mathrm{dR}}(\L)$ is not \'etale in general. In \S \ref{subsec: NdR-pdE}, we explain how to view $\mathbf{N}_{\mathrm{dR}}(\L)$ as a ``relative $p$-adic differential equation''. Finaly in \S \ref{subsec: compatibility-NdR}, we study the compatibility of $\fN_{\mathrm{dR}}$ functor under certain extensions.

\subsection{Construction of $\sq{\mathbf{N}}_{\mathrm{dR}} (\L)$} \label{subsec: tilde-NdR}

In the rest of the section, let $\L$ be a de Rham $\Z_p$-local system on $X$ and let $D_{\mathrm{dR}}(\L)$ be the filtered vector bundle with integrable connection on $X$ constructed in \cite{Liu_Zhu}. Recall that $D_{\mathrm{dR}}(\L)$ is associated with the de Rham local system $\L$ in the sense of \cite[Definition 8.3]{Scholze_p_adic_Hodge}. We first construct $\widetilde{\mathbf{N}}_{\mathrm{dR}}(\L)$. The key ingredient is the following isomorphism.

\begin{lemma} \label{lemma:gluing_isomorphism_perfect_relative}
Suppose that $r$ is sufficiently small and $p^{-n}\in [s,r]$ for some integer $n$. Then we have a $\Gamma$-equivariant isomorphism 
    \[\widetilde{M}^{[s,r]}\otimes_{\widetilde{\mathbf{C}}^{[s,r]}_X, \iota_n} \mathbb{B}_{\mathrm{dR}, X}(X_{\infty})\cong \left(D_{\mathrm{dR}}(\mathbb{L})\otimes_{A_X}\mathcal{O}\mathbb{B}_{\mathrm{dR}, X}^+(X_{\infty})\right)^{\nabla=0}[\frac{1}{t}].\]
\end{lemma}

\begin{proof}
Fix an algebraic closure of $\mathrm{Frac}(A_X)$ containing $A_{X, \infty}$. Define $\overline{A}_{X}$ to be the union of all finite \'etale $A_X$-algebras contained in this algebraic closure, equipped with the spectral norm $|\cdot|$ extending that on $A_X$. Write $(\overline{A}_X)^{\wedge}$ for the completion of $\overline{A}_{X}$ with respect to $|\cdot|$.
Let $G_{X}\coloneqq \Gal\left(\overline{A}_{X}/A_X\right)$ and $H_{X}\coloneqq \Gal\left(\overline{A}_{X}/A_{X, \infty}\right) \cong \Gal\big((\overline{A}_X)^{\wedge} /\widehat{A}_{X, \infty} \big)$ so that $\Gamma=G_{X}/H_{X}$. Write 
$\mathbf{B}_{\mathrm{dR}, \cl X} \coloneqq \mathbb{B}_{\mathrm{dR},X}((\overline{A}_X)^{\wedge, \flat})$. Then we have 
$
\mathbb{B}_{\mathrm{dR}, X} (X_\infty )  = \big( \mathbf{B}_{\mathrm{dR}, \cl X}  \big)^{H_X}.
$ 
By the construction of $\widetilde{D}_{\mathrm{rig}}^{\dagger} (\L)$ (cf. \cite[Theorem~8.5.3 \& 8.5.5]{Kedlaya-LiuI}), we know that for $0 < s \le r$ with $r$ sufficiently small, we have $\Gamma$-equivariant isomorphism 
\[
\widetilde{M}^{[s,r]}\otimes_{\widetilde{\mathbf{C}}^{[s,r]}_X, \iota_n} \mathbf{B}_{\mathrm{dR}, \cl X}  \cong 
\widetilde{D}^{\mathrm{int}, r} (\L) \otimes_{\widetilde{\mathbf{A}}_X^r, \iota_n} \mathbf{B}_{\mathrm{dR}, \cl X}  \cong \L \otimes_{\Z_p} \mathbf{B}_{\mathrm{dR}, \cl X}.
\]
Here, $\L$ in the last term of the above isomorphism means the finite free $\Z_p$-module associated with the $\Z_p$-local system $\L$, by slight abuse of notation. In the second term, $\widetilde{D}^{\mathrm{int}, r} (\L) \coloneqq (\L\otimes_{\Z_p} \widetilde{\fA}^r)^{H_X}$ where $\widetilde{\fA} \coloneqq W((\overline{A}_X)^{\wedge, \flat})$.

On the other hand, by \cite[Theorem~7.6]{Scholze_p_adic_Hodge}, we have a $\Gamma$-equivariant isomorphism 
\[ 
 \L \otimes_{\Z_p} \mathbf{B}_{\mathrm{dR}, \cl X} \cong \left(D_{\mathrm{dR}}(\mathbb{L})\otimes_{A_X}\mathcal{O}\mathbb{B}_{\mathrm{dR}, X}(  (\overline{A}_X)^{\wedge, \flat})\right)^{\nabla=0}.  
\]
Therefore, we have a $\Gamma$-equivariant isomorphism 
\[
\widetilde{M}^{[s,r]}\otimes_{\widetilde{\mathbf{C}}^{[s,r]}_X, \iota_n} \mathbf{B}_{\mathrm{dR}, \cl X} \cong \left(D_{\mathrm{dR}}(\mathbb{L})\otimes_{A_X}\mathcal{O}\mathbb{B}_{\mathrm{dR}, X} (  (\overline{A}_X)^{\wedge, \flat})\right)^{\nabla=0}.  
\]
Taking $H_X$-invariants on both sides, we obtain the desired isomorphism as claimed in the lemma. 
\end{proof}

The rest of the construction is a direct generalization of Construction~\ref{construction:tilde NdR classical setting}.  This time, we consider the Beauville--Laszlo datum 
\[
\left(\widetilde{\mathbf{C}}^{[s,r]}_{X},\,\, \widetilde{\mathbf{C}}^{[s,r]}_{X}[\frac{1}{t}],\,\, \mathbb{B}^+_{\mathrm{dR}, X} (X_\infty)\right)
\]
when $\#\left([s,r]\cap p^{\mathbb{Z}}\right)=1$. Let us write
\[
\widetilde{N}^{\mathrm{dR}} \coloneqq (D_{\mathrm{dR}}(\mathbb{L})\otimes_{A_X}\mathcal{O}\mathbb{B}_{\mathrm{dR}, X}^+(X_{\infty}) )^{\nabla=0},
\]
which is a finite projective $\mathbb{B}^+_{\mathrm{dR}, X} (X_\infty)$-module by the proof of \cite[Theorem~7.6]{Scholze_p_adic_Hodge}. Then, by Lemma~\ref{lemma:gluing_isomorphism_perfect_relative} above and Proposition~\ref{prop:BL gluing}, the finite projective $\widetilde{\mathbf{C}}^{[s,r]}_X[\frac{1}{t}]$-module $\widetilde{M}^{[s,r]}[\frac{1}{t}]$ and the finite projective $\mathbb{B}^+_{\mathrm{dR},X}(X_\infty)$-module $\widetilde{N}^{\mathrm{dR}}$ glue to a finite projective $\widetilde{\mathbf{C}}^{[s,r]}_X$-module $\widetilde{N}^{[s,r]}$. By a similar argument as in Construction~\ref{construction:N_dR_for_all_HT_wts}, these $\widetilde{N}^{[s,r]}$ glue to a $(\varphi, \Gamma)$-module $\widetilde{N}^r$ (resp. $\widetilde{N}^{\dagger}$) over $\widetilde{\mathbf{C}}^r_{X}$ (resp. $\widetilde{\mathbf{C}}^{\dagger}_{X})$. This yields the desired functor 
\begin{equation}\label{eq: tilde NdR relative setting}
\widetilde{\mathbf{N}}_{\mathrm{dR}}\colon 
\mathrm{Loc}^{\mathrm{dR}}_{\Z_p}(X) \lra \mathrm{Mod}^{(\varphi, \Gamma)}_{/\widetilde{\mathbf{C}}^{\dagger}_{X}},
\end{equation}
sending $\L \mapsto \widetilde{N}^{\dagger}$.

\vspace{0.1in}
\subsection{Construction of $\mathbf{N}_{\mathrm{dR}}(\L)$, I: preparation} \label{subsec: NdR-preparation}
\noindent
\vspace{0.1in}

\noindent The rest of the section is devoted to the construction and properties of $\mathbf{N}_{\mathrm{dR}}(\L)$. Recall that $\mathbb{L}$ is a de Rham local system on $X$ and $M^{\dagger} \coloneqq D_{\mathrm{rig}}^{\dagger} (\L)$ is the \'etale $(\varphi, \Gamma)$-module over $\mathbf{C}^{\dagger}_{\Lambda}$ constructed in \S 
\ref{ss:relative_D_rig}.   
In particular, $M^{\dagger}$ corresponds to a compatible system of finite projective $\mathbf{C}^{[s,r]}_{\Lambda}$-modules $M^{[s,r]}$, equipped with $\Gamma$-actions and compatible Frobenius isomorphisms 
\begin{equation} \label{eq:frob_isom_relative}
 \psi =  \psi^{[s, r]}\colon   M^{[s, r]} \otimes_{\mathbf{C}^{[s,r]}_{\Lambda}, \varphi} \mathbf{C}^{[s/p, \,r/p]}_{\Lambda} \isom M^{[s/p, \,r/p]}, 
\end{equation}
for all $r$ sufficiently small and all $0<s\le r$.  As in the classical setting (see \S \ref{sss: NdR classical setting}), we will modify each $M^{[s,r]}$ along the ``$t=0$ locus'' to obtain another compatible system of $\mathbf{C}^{[s,r]}_{\Lambda}$-modules $N^{[s,r]}$ which glue to a $(\varphi, \Gamma)$-module $N^{\dagger}$ over $\mathbf{C}^{\dagger}_{\Lambda}$. (When $[s,r]\cap p^{\mathbb{Z}}=\emptyset$, we again put $N^{[s,r]}=M^{[s,r]}$ as $t$ is invertible on $\mathbf{C}^{[s,r]}_{\Lambda}$.)  

We first consider the situation where $\#\left([s,r]\cap p^{\mathbb{Z}}\right)=1$, say $p^{-n}\in [s,r]$ (where $r$ is sufficiently small). By Corollary~\ref{cor: t-completion-of-C}, $A_{X,n}[\![t, u_1, \ldots, u_d]\!]^{\nabla=0}$ is precisely the $t$-adic completion of $\mathbf{C}^{[s,r]}_{\Lambda}$. Thus \[(\mathbf{C}^{[s,r]}_{\Lambda}, \,\,\mathbf{C}^{[s,r]}_{\Lambda}[\frac{1}{t}],\,\, A_{X,n}[\![t, u_1, \ldots, u_d]\!]^{\nabla=0})\] is a Beauville--Laszlo datum (see Definition~\ref{definition: BL data}). 

For each $n \ge 1$, consider  
\[N^{\mathrm{dR}}_n\coloneqq \left(D_{\mathrm{dR}}(\mathbb{L})\otimes_{A_X} A_{X,n}[\![t,u_1,\ldots,u_d]\!]\right)^{\nabla=0}.\]
This is an $A_{X,n}[\![t, u_1, \ldots, u_d]\!]^{\nabla=0}$-module. Recall  $A_X[\![t, u_1, \ldots, u_d]\!]^{\nabla = 0}[\![u_1, \ldots, u_d]\!] \cong A_X[\![t, u_1, \ldots, u_d]\!]$ from Lemma~\ref{lemma:sublemma_on_horizontal_sections}. Write 
\[
\mathcal{D}^{\nabla =0} \coloneqq (D_{\mathrm{dR}}(\L) \otimes_{A_X} A_X[\![t, u_1, \ldots, u_d]\!])^{\nabla = 0}.
\]
By the proof of \cite[Proposition~8.9]{Katz_differential_equations} (applied to the $\Q$-algebra $A_X[\![t, u_1, \ldots, u_d]\!]^{\nabla = 0}$ as the base ring in place of $K$ in \textit{loc. cit.}), $\mathcal{D}^{\nabla =0}$ is a finite projective $A_X[\![t, u_1, \ldots, u_d]\!]^{\nabla = 0}$-module, and we have a natural isomorphism
\[
\mathcal{D}^{\nabla = 0} \otimes_{A_X[\![t, u_1, \ldots, u_d]\!]^{\nabla = 0}} A_X[\![t, u_1, \ldots, u_d]\!] \cong D_{\mathrm{dR}}(\L) \otimes_{A_X} A_X[\![t, u_1, \ldots, u_d]\!].
\]
Thus, 
\begin{align*}  
N_n^{\mathrm{dR}}  \: \cong \: &  (\mathcal{D}^{\nabla = 0}\otimes_{A_X[\![t, u_1, \ldots, u_d]\!]^{\nabla = 0}} A_{X,n}[\![t, u_1, \ldots, u_d]\!])^{\nabla = 0} \\ \cong \: &  \mathcal{D}^{\nabla = 0}\otimes_{A_X[\![t, u_1, \ldots, u_d]\!]^{\nabla = 0}} A_{X,n}[\![t, u_1, \ldots, u_d]\!]^{\nabla = 0},
\end{align*}
which is a finite projective $A_{X,n}[\![t, u_1, \ldots, u_d]\!]^{\nabla = 0}$-module of rank $\tu{rk}_{\Z_p} \L$. In particular, we have 
\[   N^{\mathrm{dR}}_{n+1} \cong 
N^{\mathrm{dR}}_n\otimes_{A_{X,n}[\![t, u_1, \ldots, u_d]\!]^{\nabla=0} }  A_{X,n+1}[\![t, u_1, \ldots, u_d]\!]^{\nabla=0}
\]
for each $n \ge 1$. 
Similarly, we have
\[
 \left(D_{\mathrm{dR}}(\mathbb{L})\otimes_{A_X}\mathcal{O}\mathbb{B}_{\mathrm{dR}, X}^+(X_{\infty})\right)^{\nabla=0} \cong \mathcal{D}^{\nabla = 0}\otimes_{A_X[\![t, u_1, \ldots, u_d]\!]^{\nabla = 0}} \mathbb{B}_{\mathrm{dR}, X}^+(X_{\infty}).
\]
The above isomorphisms immediately imply the following.

\begin{corollary} \label{cor:RHS_gluing_isom}
There is a $\Gamma$-equivariant isomorphism 
\[
N^{\mathrm{dR}}_n\otimes_{A_{X,n}[\![t, u_1, \ldots, u_d]\!]^{\nabla=0} }\mathbb{B}_{\mathrm{dR}, X}^+(X_{\infty})\cong \left(D_{\mathrm{dR}}(\mathbb{L})\otimes_{A_X}\mathcal{O}\mathbb{B}_{\mathrm{dR}, X}^+(X_{\infty})\right)^{\nabla=0}.
\] 
\end{corollary}

\vspace{0.1in}
\subsection{Construction of $\mathbf{N}_{\mathrm{dR}} (\L)$, I\!I: the gluing isomorphism} \label{subsec: NdR-gluing isom} 
\noindent
\vspace{0.1in}

\noindent Similar to the construction in \S \ref{sss: NdR classical setting}, we would like to apply Beauville--Laszlo gluing to the $\mathbf{C}^{[s,r]}_{\Lambda}[1/t]$-module $M^{[s,r]}[1/t]$ and the $A_{X,n}[\![t, u_1, \ldots, u_d]\!]^{\nabla=0}$-module $N^{\mathrm{dR}}_n$ over $A_{X,n}[\![t, u_1, \ldots, u_d]\!]^{\nabla=0}[\frac{1}{t}]$. For this, we need an isomorphism
\[M^{[s,r]}\otimes_{\mathbf{C}_{\Lambda}^{[s,r]}, \iota_n} A_{X,n}[\![t, u_1, \ldots, u_d]\!]^{\nabla=0}[\frac{1}{t}]\cong N^{\mathrm{dR}}_n[\frac{1}{t}].\] 
Parallel to our construction in the classical setting, we achieve the isomorphism in two steps. We first construct an isomorphism over $A_{X,\infty}[\![t, u_1, \ldots, u_d]\!]^{\nabla=0}[\frac{1}{t}]$ using the relative Sen theory developed in the previous section, then further descend it to $A_{X,n}[\![t, u_1, \ldots, u_d]\!]^{\nabla=0}[\frac{1}{t}]$ by studying the action of $\mathrm{Lie}\,\Gamma$.

\begin{proposition} \label{prop: iota_n-isomorphism-over-infty}
Suppose that $r$ is sufficiently small and $p^{-n}\in [s,r]$. There is a canonical $\Gamma$-equivariant isomorphism  
\begin{equation}
     \label{eq:map_iota_n_infty_relative} \quad 
  \iota_n\colon M^{[s,r]}\otimes_{\mathbf{C}^{[s,r]}_{\Lambda}, \iota_n} A_{X,\infty}[\![t, u_1, \ldots, u_d]\!]^{\nabla=0}[\frac{1}{t}]   \isom  N^{\mathrm{dR}}_n \otimes_{ A_{X,n}[\![t,u_1,\ldots,u_d]\!]^{\nabla=0} } A_{X,\infty}[\![t,u_1,\ldots,u_d]\!] ^{\nabla=0} [\frac{1}{t}].    
\end{equation} 
\end{proposition}
    
\begin{proof}
By Corollary~\ref{cor:equivalence_relative_B_dR^+}, it suffices to  check the existence of such an isomorphism after base change along 
$A_{X,\infty}[\![t, u_1, \ldots, u_d]\!]^{\nabla=0}\rightarrow \mathbb{B}_{\mathrm{dR}, X}^+(X_{\infty})$.  To this end, note that we have
\begin{multline*}
\left(M^{[s,r]}\otimes_{\mathbf{C}_{\Lambda}^{[s,r]}, \iota_n} A_{X,\infty}[\![t, u_1, \ldots, u_d]\!]^{\nabla=0}[\frac{1}{t}]\right)\otimes_{A_{X,\infty}[\![t, u_1, \ldots, u_d]\!]^{\nabla=0}}\mathbb{B}_{\mathrm{dR}, X}^+(X_{\infty}) \\ \cong \widetilde{M}^{[s,r]}\otimes_{\widetilde{\mathbf{C}}^{[s,r]}_X, \iota_n}\mathbb{B}_{\mathrm{dR}, X}(X_{\infty}).
\end{multline*}
On the other hand, by Corollary~\ref{cor:RHS_gluing_isom}, we have 
\[N^{\mathrm{dR}}_n\otimes_{A_{X,n}[\![t, u_1, \ldots, u_d]\!]^{\nabla=0} }\mathbb{B}_{\mathrm{dR}, X}^+(X_{\infty})\cong \left(D_{\mathrm{dR}}(\mathbb{L})\otimes_{A_X}\mathcal{O}\mathbb{B}_{\mathrm{dR}, X}^+(X_{\infty})\right)^{\nabla=0}.\] 
  The assertion thus follows from Lemma~\ref{lemma:gluing_isomorphism_perfect_relative}. 
\end{proof}

Next, we descend \eqref{eq:map_iota_n_infty_relative} to an isomorphism over $A_{X,n}[\![t, u_1, \ldots, u_d]\!]^{\nabla=0}[\frac{1}{t}]$.

\begin{proposition} \label{prop:image_of_iota_n_fixed_by_Gamma_n}
Let $n$ be a sufficiently large integer and suppose that $p^{-n}\in [s,r]$. Then the image of $M^{[s, r]}$ under the map $\iota_n$ in  \eqref{eq:map_iota_n_infty_relative} satisfies 
 $
 \iota_n (M^{[s,r]}) \subset N^{\mathrm{dR}}_n [\frac{1}{t}].
 $ 
 Consequently, $\iota_n$ in  \eqref{eq:map_iota_n_infty_relative} descends to a $\Gamma$-equivariant isomorphism 
\begin{equation} \label{eq:gluing_isom_relative_s,r}
  \iota_n\colon M^{[s,r]}\otimes_{\mathbf{C}^{[s,r]}_{\Lambda}, \iota_n} A_{X,n}[\![t, u_1, \ldots, u_d]\!]^{\nabla=0}[\frac{1}{t}]   \isom  N^{\mathrm{dR}}_n  [\frac{1}{t}] 
\end{equation} 
by faithfully flat descent along $A_{X,n}[\![t, u_1, \ldots, u_d]\!]^{\nabla=0}[\frac{1}{t}]\rightarrow A_{X,\infty}[\![t, u_1, \ldots, u_d]\!]^{\nabla=0}[\frac{1}{t}]$.
\end{proposition}

The proof of Proposition~\ref{prop:image_of_iota_n_fixed_by_Gamma_n} will occupy the rest of \S \ref{subsec: NdR-gluing isom}. We will follow a similar argument as in the proof of Proposition~\ref{lemma:image_of_iota_n_contained_in_level_n}. Fix a splitting $\Gamma \cong \Gamma_{\mathrm{geom}} \rtimes  \Gamma_{K}$ as in \S \ref{ss:toric_tower} and choose topological generators $\tau_1, \ldots, \tau_d \in \Gamma_{\mathrm{geom}}$ such that   
\begin{equation}\label{eq:generators_of_Gamma}
\tau_i(T_j^{\frac{1}{p^n}})= \begin{cases}  \zeta_{p^n}T_j^{\frac{1}{p^n}}, & \text{ if } i = j \\ 
T_j^{\frac{1}{p^n}} &\text{ otherwise.}
\end{cases} 
\end{equation}
This determines an isomorphism $\eta= (\eta_1, \ldots, \eta_d)\colon \Gamma_{\mathrm{geom}} \isom \Z_p^{d}$. We consider differential operators on $\mathbf{C}^{[s,r]}_{\Lambda}$ and $M^{[s,r]}$ induced from the action of $\mathrm{Lie}\, \Gamma$. More precisely, $\mathbf{C}^{[s,r]}_{\Lambda}$ is equipped with $(d+1)$ differential operators
$
\nabla_{0} \coloneqq \frac{\log \gamma_0}{\log \chi(\gamma_0)}
$ 
with $\gamma_0\in \Gamma_K$ sufficiently close to 1, and
$
\nabla_{i} \coloneqq \frac{\log (\tau_i^{p^m})}{\eta_i (\tau_i^{p^m})} = p^{-m}\cdot\log (\tau_i^{p^m})
$  
for $i=1, \ldots, d$, and for $m$ sufficiently large. Note that $\nabla_0$ (resp. $\nabla_i$) is independent of choice of $\gamma_0$ (resp. choice of $m$). Similarly, $M^{[s, r]}$ is also equipped with differential operators $\nabla_0, \nabla_1, \ldots, \nabla_d$ defined by the same formulae, making $(M^{[s, r]}, \{\nabla_i\}_{0 \le i \le d})$ into a differential module over $(\mathbf{C}_{\Lambda}^{[s, r]}, \{\nabla_i\}_{0 \le i \le d})$. 

Note that it suffices to prove the proposition for $[s,r]=[r_n, r_n]$ where $r_n=p^{-n}$. By Proposition~\ref{prop: iota_n-isomorphism-over-infty}, we have
\[
\iota_n(M^{[r_n, r_n]}) \subset (D_{\mathrm{dR}}(\L)\otimes_{A_X} A_{X, \infty}[\![t, u_1, \ldots, u_d]\!]))^{\nabla = 0}[\frac{1}{t}]. 
\]
So by the definition of $N_n^{\mathrm{dR}}$, it suffices to show that
\[
\iota_n(M^{[r_n, r_n]}) \subset D_{\mathrm{dR}}(\L)\otimes_{A_X} A_{X, n}[\![t, u_1, \ldots, u_d]\!][\frac{1}{t}]
\]
for any large enough $n$. Since the map $\iota_n$ in \eqref{eq:map_iota_n_infty_relative} is $\Gamma$-equivariant,  it is also compatible with the differential operators on both sides, where the right-hand side of $\iota_n$ can be viewed as a differential module over 
\[ \left(A_{X,\infty}[\![t,u_1,\ldots,u_d]\!] ^{\nabla=0} [\frac{1}{t}], \quad   t \frac{d}{dt},  \quad  t \cdot  [T_i^\flat]   \frac{d}{d [T_i^\flat]} \right).\]  In particular, we have $\iota_n \circ \nabla_i = \nabla_i \circ \iota_n$ for all $0 \le i \le d$. Inspired by the proof of Proposition~\ref{lemma:image_of_iota_n_contained_in_level_n}, we consider formal sums
\begin{equation}\label{eq: formal sum sigma}
\exp (\log \chi(\gamma_0) \cdot \nabla_0)\coloneqq \sum_{m=0}^{\infty}\frac{(-1)^m}{m!}\left(\sum_{k=1}^{\infty}\frac{(1-\gamma_0)^k}{k}\right)^m
\end{equation}
with $\gamma_0\in \Gamma_K$, and
\begin{equation}\label{eq: formal sum tau}
\exp ( p^n \cdot \nabla_i )\coloneqq \sum_{m=0}^{\infty}\frac{(-1)^m}{m!}\left(\sum_{k=1}^{\infty}\frac{(1-\tau_i^{p^n})^k}{k}\right)^m
\end{equation}
for all $i=1,\ldots, d$ and for $n$ sufficiently large.

\begin{lemma}\label{lemma: exp of log--relative setting}
There exists an integer $n_*$ such that for any $n\ge n_*$, the formal sums \eqref{eq: formal sum sigma} and \eqref{eq: formal sum tau} converge\footnote{Just like in Lemma~\ref{lemma: gamma as exp}, convergence of $\exp (\log \chi(\gamma_0) \cdot \nabla_0)$ means that both formal sums $\sum_{k=1}^{\infty}\frac{(1-\gamma_0)^k}{k}$ and $\sum_{m=0}^{\infty}\frac{(-1)^m}{m!}(\sum_{k=1}^{\infty}\frac{(1-\gamma_0)^k}{k})^m$ converge. Similarly for $\exp ( p^n \cdot \nabla_i )$.} as operators on $M^{[r_n, r_n]}$, and we have 
\[
\gamma_0  = \exp (\log \chi(\gamma_0) \cdot \nabla_0),\qquad \tau_i^{p^n}=\exp ( p^n \cdot \nabla_i )
\]
for all $\gamma_0\in \Gamma_{K_n}$ and for all $i=1, \ldots, d$.
\end{lemma}

\begin{proof}
The proof is similar to that of Lemma~\ref{lemma: gamma as exp}, except that Lemma~\ref{lemma: estimation Gamma on C} is replaced by Condition (c) in Definition~\ref{defn: strongly decompleting data}. We sketch the proof here for completeness.

Let $n_0$ and $C$ be as in Condition (c) of Definition~\ref{defn: strongly decompleting data}. Since $C<p^{-\frac{1}{p-1}}$, for the convergence conditions and identities in the statement of the lemma, it suffices to construct a Banach norm $\mu$ on $M^{[r_n, r_n]}$ which is equivalent to the finite-module norm over $\mathbf{C}_{\Lambda}^{[r_n, r_n]}$, such that
\[
\mu((\gamma-1)(x))\le C\cdot \mu(x)
\]
for all $x\in M^{[r_n, r_n]}$ and $\gamma\in \Gamma_n$. Choose any finite-module norm $\mu$ on the finite projective $\mathbf{C}_{\Lambda}^{[r_{n_0}, r_{n_0}]}$-module $M^{[r_{n_0}, r_{n_0}]}$. By continuity, there exists $n_*\ge n_0$ such that
\[
\mu((\gamma-1)(x))\le C\cdot \mu(x)
\]
for all $x\in M^{[r_{n_0}, r_{n_0}]}$ and $\gamma\in \Gamma_{n_*}$\footnote{Recall that $\Gamma_n$ stands for the open subgroup of $\Gamma$ corresponding to $\Gamma_{\mathrm{geom},n}\rtimes \Gamma_{K_n}\subset \Gamma_{\mathrm{geom}}\rtimes \Gamma_K$, where $\Gamma_{\mathrm{geom},n}=p^n\Gamma_{\mathrm{geom}}\subset \Gamma_{\mathrm{geom}}$.}. Then following the same argument as in Lemma~\ref{lemma: gamma as exp}, one sees that this $n_*$ works.
\end{proof}

We now return to the proof of Proposition~\ref{prop:image_of_iota_n_fixed_by_Gamma_n}. Let $n\ge n_*$ as in Lemma~\ref{lemma: exp of log--relative setting}. Recall that $D_{\mathrm{dR}}(\L)$ is finite projective over $A_X$. Let $\{U_f = \Spa(B_f, B_f^+)\}_{f \in \Sigma}$ be a finite covering of $A_X$ by rational subsets such that $D_{\mathrm{dR}}(\L)$ becomes free over each $U_f$. Note that $\Gamma$ acts trivially on each $B_f$. For each $m\ge 1$, write 
\[
\Spa(B_f, B_f^+)\times_{\D_K^d} \D_m=\Spa(B_{f, m}, B_{f, m}^+) 
\]
and let $B_{f, \infty}=\varinjlim_m B_{f, m}$.

Fix $f \in \Sigma$ and choose a $B_f$-basis $x_1, \ldots, x_{\ell}$ of $D_{\mathrm{dR}}(\L)\otimes_{A_X} B_f$. Let $x \in M^{[r_n, r_n]}$ be an arbitrary element and write $\iota_n(x)$ into the following (unique) expression 
\[ 
\iota_n (x) = \sum_{j=1}^{\ell} c_j \otimes x_j = \sum_{j=1}^{\ell} \sum_{i \in \Z} a_{i, j} t^i \otimes x_j, 
\]
where $c_j \in B_{f, \infty}[\![t, u_1, \ldots, u_d ]\!][\frac{1}{t}]$ and $a_{i, j} \in  B_{f,\infty}[\![u_1, \ldots, u_d]\!]$. 

For $\gamma_0 \in \Gamma_{K_n}$, the same argument as in the proof of Proposition~\ref{lemma:image_of_iota_n_contained_in_level_n} shows that 
\begin{align*}
    \iota_n (\gamma_0 \cdot x) \: & = 
\: \iota_n (\exp (\log \chi(\gamma_0) \cdot \nabla_0) (x)) \\ 
& = \: \sum_{m =0}^{\infty} \frac{1}{m!} (\log \chi (\gamma_0))^m (t \frac{d}{dt})^m  \Big( \sum_{j=1}^{\ell}\sum_{i \in \Z} a_{i, j} t^i \otimes x_j \Big) 
\\ & = \: 
\sum_{j=1}^{\ell}\sum_{i \in \Z } a_{i,j}  \cdot \chi(\gamma_0)^i t^i\otimes x_j.    
\end{align*}
On the other hand, we have
\[
\gamma_0 \cdot \iota_n (x)  =  \sum_{j=1}^{\ell}\sum_{i \in \Z} \gamma_0(a_{i, j} t^i) \otimes x_j  =  \sum_{j=1}^{\ell}\sum_{i \in \Z} \gamma_0(a_{i, j}) \cdot  \chi(\gamma_0)^it^i \otimes x_j. 
\]
This implies that $\gamma_0 (a_{i, j}) = a_{i, j}$ for all $i, j$. Since $\gamma_0$ acts trivially on each $u_j$, this further implies that $a_{i, j} \in (B_{f, \infty} )^{\Gamma_{K_n}}[\![u_1, \ldots, u_d]\!]$. 

Next, we consider the geometric components of $\Gamma$. Without loss of generality, let us consider an element $\gamma_1 \in \Gamma_{\mathrm{geom}, n}$ such that   $\eta_j (\gamma_1) = 0$ for all $j \ne 1$ but $\eta_1(\gamma_1) \ne 0$. Write $\delta_1 \coloneqq \eta_1(\gamma_1)$. We express each $a_{i, j}$ as a formal sum  
$a_{i, j} = \sum_{m = 0}^{\infty} a_{i, j, m} \cdot  u_1^m$ with $a_{i, j, m} \in B_{f, \infty} [\![u_2, \ldots, u_d ]\!]$ for each $i, j, m$. Using Lemma~\ref{lemma: exp of log--relative setting}, we compute
\begin{align*}
    \iota_n (\gamma_1 \cdot x) \: & = 
\: \iota_n (\exp ( \delta_1  \cdot \nabla_1) (x))  = \:    \iota_n \big(\sum_{k = 0}^{\infty} \frac{1}{k!}( \delta_1 )^k  ( \nabla_1)^k (x)\big)   \\ 
& = \: 
    \sum_{k =0}^{\infty} \frac{1}{k!}(\delta_1 )^k  ( \nabla_1)^k ( \iota_n (x)) 
  \\ 
& = \: \sum_{k =0}^{\infty} \frac{1}{k!} (\delta_1 t)^k ((u_1-T_1)\frac{d}{d u_1})^k  \Big( \sum_{j=1}^{\ell}\sum_{i \in \Z}\sum_{m=0}^{\infty} a_{i, j, m} \cdot u_1^m \cdot  t^i \otimes x_j \Big) \\
& = \: 
\sum_{j=1}^{\ell}  \sum_{i \in \Z}\sum_{m=0}^{\infty} a_{i,j, m}  \cdot  \Big( ([\epsilon]^{\delta_1}-1)(u_1-T_1)+u_1 \Big)^m \cdot t^{i}\otimes x_j,
\end{align*}
where we have used Lemma~\ref{lem: Gamma-action-identity} below for the last equality.
On the other hand, observe that 
\[
\gamma_1 (u_1) = T_1-[\epsilon]^{\delta_1}[T_1^{\flat}] = ([\epsilon]^{\delta_1}-1)(u_1-T_1)+u_1.
\]
Thus, we have 
\begin{align*}
\gamma_1 \cdot \iota_n (x)  \: & = \: \sum_{j=1}^{\ell}\sum_{i \in \Z}\sum_{m=0}^{\infty} \gamma_1(a_{i, j, m} \cdot  u_1^m \cdot t^i) \otimes x_j \\ 
& = \: \sum_{j=1}^{\ell}\sum_{i \in \Z}\sum_{m=0}^{\infty} \gamma_1(a_{i, j, m}) \cdot (([\epsilon]^{\delta_1}-1)(u_1-T_1)+u_1)^m \cdot  t^i \otimes x_j.
\end{align*}
Since $t$ divides $[\epsilon]^{\delta_1}-1$ in $A_{X, \infty}[\![t]\!]$, every element in $B_{f, \infty}[\![t, u_1, \ldots, u_d]\!]$ can be uniquely expressed as
\[
\sum_{i, m \ge 0} b_{i, m} (([\epsilon]^{\delta_1}-1)(u_1-T_1)+u_1)^m\cdot t^i
\]
with $b_{i, m} \in B_{f, \infty}[\![u_2, \ldots, u_d]\!]$. Hence, the computation above yields $\gamma_1(a_{i, j, m}) = a_{i, j, m}$ for all $i, j, m$. Since $\gamma_1$ acts trivially on $u_2, 
\ldots, u_d$, we conclude that the coefficients of $a_{i, j}$ 
(as power series in variables $u_1, \ldots, u_d$) lie in $(B_{f, \infty})^{\gamma_1 = 1}$. Putting everything together (and repeating the argument above for the other geometric components), we know that the coefficients of $a_{i, j}$ in fact lie in $(B_{f, \infty})^{\Gamma_n} = B_{f, n}$. In other words, each $c_j$ is an element in $B_{f, n}[\![t, u_1, \ldots, u_d]\!][\frac{1}{t}]$, so
\[
\iota_n(M^{[r_n, r_n]}) \subset D_{\mathrm{dR}}(\L)\otimes_{A_X} B_{f, n}[\![t, u_1, \ldots, u_d]\!][\frac{1}{t}].
\]
Since this holds for any $f \in \Sigma$ and $D_{\mathrm{dR}}(\L)$ is flat over $A_X$, we deduce
\[
\iota_n(M^{[r_n, r_n]}) \subset D_{\mathrm{dR}}(\L)\otimes_{A_X} A_{X, n}[\![t, u_1, \ldots, u_d]\!][\frac{1}{t}]
\]
as desired.
The following lemma is used in the proof above.

\begin{lemma} \label{lem: Gamma-action-identity}
For each $m \ge 0$, we have
\[
\sum_{k = 0}^{\infty} \frac{(\delta_1 t)^k}{k!}((u_1-T_1)\frac{d}{du_1})^k (u_1^m) = ((u_1-T_1)[\epsilon]^{\delta_1}+T_1)^m.
\]
\end{lemma}

\begin{proof}
Note that $[\epsilon]^{\delta_1} = \exp(\delta_1 t)$. So it suffices to show that for each $k, m \ge 0$, we have
\[
((u_1-T_1)\frac{d}{du_1})^k u_1^m = (\frac{d}{dx})^k \Big|_{x = 0} ((u_1-T_1)e^x+T_1)^m.
\]
When $k = 0$, both sides are equal to $u_1^m$. Suppose $k \ge 1$. From
\[
u_1^m = ((u_1-T_1)+T_1)^m = \sum_{i=0}^m \binom{m}{i} (u_1-T_1)^i \cdot T_1^{m-i},
\]
we deduce that 
\[
((u_1-T_1)\frac{d}{du_1})^k u_1^m = \sum_{i=0}^m i^k \binom{m}{i} (u_1-T_1)^i \cdot T_1^{m-i}.
\]
On the other hand, note that
\[
((u_1-T_1)e^x+T_1)^m = \sum_{i=0}^m \binom{m}{i} (u_1-T_1)^i e^{ix}\cdot  T^{m-i}.
\]
Therefore, we have 
\[
(\frac{d}{dx})^k \mid_{x = 0} ((u_1-T_1)e^x+T_1)^m = \sum_{i=0}^m i^k  \binom{m}{i} (u_1-T_1)^i e^{ix}\cdot T_1^{m-i} \mid_{x = 0} = \sum_{i=0}^m i^k \binom{m}{i} (u_1-T_1)^i \cdot T_1^{m-i},
\]
as desired.
\end{proof}

\subsection{Construction of $\mathbf{N}_{\mathrm{dR}} (\L)$, I\!I\!I: Beauville--Laszlo gluing} \label{subsec: NdR-BL gluing}
\noindent
\vspace{0.1in}

\noindent 
Now we are ready to finish the construction of $\mathbf{N}_{\mathrm{dR}}(\L)$. Given the gluing isomorphism \eqref{eq:gluing_isom_relative_s,r}, the rest of the construction goes \emph{verbatim} as in Construction \ref{construction:N_dR_for_all_HT_wts}; namely, $\mathbf{N}_{\mathrm{dR}}(\L)$ is obtained via a ``modification'' procedure along the ``$t=0$ locus'' using Beauville--Laszlo gluing. 

\begin{construction} \label{const: NdR}
First, when $ [s,r]\cap p^{\mathbb{Z}} = \{ p^{-n} \}$ with $r$ being sufficiently small, we apply Beauville--Laszlo gluing (cf. Proposition~\ref{prop:BL gluing}) to $M^{[s,r]}[\frac{1}{t}]$ and $N_{n}^{\mathrm{dR}}$ along the isomorphism \eqref{eq:gluing_isom_relative_s,r} to obtain a finite projective $\mathbf{C}^{[s,r]}_{\Lambda}$-module $N^{[s,r]}$. As in Construction~\ref{construction:N_dR_for_all_HT_wts}, $N^{[s, r]}$ is equipped with a natural continuous $\Gamma$-action. Moreover,  one checks that the construction above is compatible when we vary $r$ and $s$, and hence glue to $\mathbf{C}^{[s,r]}_{\Lambda}$-modules $N^{[s,r]}$ for arbitrary $[s,r]$ with $r$ sufficiently small, equipped with continuous $\Gamma$-actions. By the same argument as in Step 4 of Construction~\ref{construction:N_dR_for_all_HT_wts}, we also have a compatible family of Frobenius isomorphisms 
\[
\phi^{[s, r]}\colon  N^{[s, r]} \otimes_{\mathbf{C}^{[s,r]}_{\Lambda}, \varphi} \mathbf{C}^{[s/p, \,r/p]}_{\Lambda} \isom N^{[s/p, \,r/p]} 
\]
for $r$ sufficiently small and $0 < s \le r$. Consequently, the collection $\{N^{[s,r]}\}$ forms a uniformly finitely generated vector bundle over $\mathrm{Sp}\,\mathbf{C}^r_{\Lambda}$ (cf. Definition~\ref{def:uniform_finite_generation_over_C^r}), and hence induces a finite projective $(\varphi, \Gamma)$-module $N^r$ over $\mathbf{C}^{r}_{\Lambda}$ by Corollary~\ref{cor: fin-proj-vec-bun-Cr}. Finally, taking colimit as $r\rightarrow 0^+$, we obtain a $(\varphi, \Gamma)$-module $N^{\dagger}$ over $\mathbf{C}^{\dagger}_{\Lambda}$. Summing up, we obtain the desired functor 
\begin{equation}\label{eq: NdR relative setting}
\mathbf{N}_{\mathrm{dR}}\colon 
\mathrm{Loc}^{\mathrm{dR}}_{\Z_p}(X) \lra \mathrm{Mod}^{(\varphi, \Gamma)}_{/\mathbf{C}^{\dagger}_{\Lambda}},
\end{equation}
sending $\L \mapsto N^{\dagger}$.  
\end{construction}

\begin{remark} \label{remark: not unique descent}
By construction, $\mathbf{N}_{\mathrm{dR}}(\L)$ is \emph{a} descent of  $\widetilde{\mathbf{N}}_{\mathrm{dR}}(\L)$ as $(\varphi, \Gamma)$-modules. But we do not know whether it is the \emph{unique} $(\varphi, \Gamma)$-descent of $\widetilde{\mathbf{N}}_{\mathrm{dR}}(\L)$. For the same reason, one cannot use the theories of relative $B$-pairs or vector bundles on the relative Fargues-Fontaine curve to construct $\mathbf{N}_{\mathrm{dR}}(\L)$ due to the lack of such descent results. 
\end{remark}

\subsection{Relative $p$-adic differential xequations}\label{subsec: NdR-pdE}

To end this section, we explain how to regard $\mathbf{N}_{\mathrm{dR}} (\L)$ as a relative version of $p$-adic differential equation. This will finish the proof of Theorem \ref{intro:theorem_relative_Berger} from the introduction. We retain the notation $N^{[s,r]}$, $N^r$ and $N^{\dagger}$ from \S \ref{subsec: NdR-BL gluing}. Recall the topological generators $\tau_1, \ldots, \tau_d \in \Gamma_{\mathrm{geom}}$ satisfying \eqref{eq:generators_of_Gamma}, which induce an isomorphism $\eta= (\eta_1, \ldots, \eta_d)\colon \Gamma_{\mathrm{geom}} \isom \Z_p^{d}$ as in the proof of Proposition~\ref{prop:image_of_iota_n_fixed_by_Gamma_n}. 
For $r$ sufficiently small and $s \in (0, r]$, each $N^{[s, r]}$ (and thus $N^r$) is equipped with $(d+1)$ differential operators 
\[
\nabla_{0} \coloneqq \frac{\log (\gamma_0)}{\log \chi(\gamma_0)}, \qquad \nabla_{i} \coloneqq \frac{\log (\tau_i^{p^m})}{\eta_i (\tau_i^{p^m})} = p^{-m}\cdot \log (\tau_i^{p^m}),
\]    
which are well-defined for $\gamma_0\in \Gamma_K$ sufficiently close to 1 and for $m$ sufficiently large. Passing to the limit, we obtain differential operators $\nabla_{0}, \nabla_{1}, \ldots, \nabla_{d}$ on $N^{\dagger}=\mathbf{N}_{\mathrm{dR}}(\L)$.

\begin{proposition}  \label{prop:singularities_are_removed}
For any $\L \in \mathrm{Loc}^{\mathrm{dR}}_{\Z_p}(X)$, we have 
     \[ \nabla_i (\mathbf{N}_{\mathrm{dR}}(\L)) \subset t \cdot \mathbf{N}_{\mathrm{dR}}(\L)\]   
 for each  $i = 0, 1, \ldots, d$. 
\end{proposition}

\begin{proof}
    By the construction of $\mathbf{N}_{\mathrm{dR}}(\L)$, it suffices to show that $ \nabla_i (N^{[s, r]}) \subset t \cdot N^{[s, r]}$ for all $i$, when $r$ is sufficiently small and $[s, r] \cap p^{\Z} = p^{-n}$ for some $n$. For this, it suffices to show that we have the desired inclusion after taking the $t$-adic completion (that is, after taking base change along the map $\iota_n\colon \mathbf{C}^{[s,r]}_{\Lambda}\hookrightarrow A_{X,n}[\![t, u_1, \ldots, u_d]\!]^{\nabla=0}$). In other words, we need to show
    $
    \iota_n (\nabla_i (x)) = \nabla_i (\iota_n (x)) \in t \cdot N_{n}^{\mathrm{dR}}
    $ 
     for each $x \in N^{[s, r]}$. 
     Then it suffices to show that $\nabla_i (N_n^{\mathrm{dR}}) \subset t \cdot N_n^{\mathrm{dR}}$. Note that 
     \[ N_n^{\mathrm{dR}} \subset D_{\mathrm{dR}}(\mathbb{L})\otimes_{A_X}\mathcal{O}\mathbb{B}_{\mathrm{dR}, X}^+(X_{\infty})\] and  the action of $\Gamma$ on $D_{\mathrm{dR}}(\L)$ is trivial. Hence, we are reduced to check that \[\nabla_i (\mathcal{O}\mathbb{B}_{\mathrm{dR}, X}^+(X_{\infty})) \subset t \cdot \mathcal{O}\mathbb{B}_{\mathrm{dR}, X}^+(X_{\infty}).\] This is clear because
     $ \nabla_0 =  t \cdot (1+\pi) \cdot \frac{d}{d\pi}  =  t \frac{d}{dt}, $  and $ \nabla_i =  t \cdot  [T_i^\flat]   \frac{d}{d [T_i^\flat]}$ for $i = 1, \ldots, d$ on $\mathcal{O}\mathbb{B}_{\mathrm{dR}, X}^+(X_{\infty})$. 
\end{proof}

Thanks to Proposition~\ref{prop:singularities_are_removed}, we can introduce the following ``normalized'' differential operators. Let $u_0 = \frac{d \pi_K}{ d \pi}\in (\mathbf{C}^{\dagger}_K)^{\times}$ from Remark~\ref{remark:partial vs d/d pi_K}. 

\begin{definition}\label{defn: normalized differential operators} \hfill
\begin{enumerate}
\item On $\mathbf{C}^{\dagger}_{\Lambda}$, we define differential operators
\[ 
    d_{0} \coloneqq t^{-1}\cdot ((1+\pi)\cdot u_0)^{-1}\cdot \nabla_{0}
\quad \textup{ 
and  } 
\quad    d_{i} \coloneqq t^{-1}\cdot [T_i^{\flat}]^{-1}\cdot \nabla_{i} \] 
for $i=1, \ldots, d$.

\item On $\mathbf{N}_{\mathrm{dR}}(\L)$, we define differential operators 
\[
    D_{0} \coloneqq t^{-1}\cdot ((1+\pi)\cdot u_0)^{-1}\cdot \nabla_{0}
 \quad \textup{ and } \quad 
   D_{i} \coloneqq t^{-1}\cdot [T_i^{\flat}]^{-1}\cdot \nabla_{i} 
\]
for $i=1, \ldots, d$. 
\end{enumerate}
These make $(\mathbf{N}_{\mathrm{dR}}(\L), \{D_i\}_{0\le i\le d})$ into a differential module over $(\mathbf{C}^{\dagger}_{\Lambda}, \{d_i\}_{0\le i\le d})$. From this perspective, $\mathbf{N}_{\mathrm{dR}}(\L)$ can be viewed as a \emph{relative $p$-adic differential equation} over $\mathbf{C}^{\dagger}_{\Lambda}$. 
\end{definition}

In other words, Proposition~\ref{prop:singularities_are_removed} tells us that, while the operators $D_i$ \textit{a priori} have singularities along $t = 0$ on $D_{\mathrm{rig}}^{\dagger} (\L)$, the modification process in our construction of  $\mathbf{N}_{\mathrm{dR}}(\L)$ removes these singularities. 

\begin{remark}
As we shall see in \S \ref{subsec: normalized-differential-operators-closed-annuli}, in the case where $X$ is a closed poly-annulus over $K$, we have 
\[\nabla_0=t\cdot (1+\pi)\cdot \frac{d}{d\pi},\qquad \nabla_i=t\cdot [T^{\flat}_i]\cdot \frac{d}{d[T^{\flat}_i]}\]
for $i=1, \ldots, d$. Hence, we have $d_0=\frac{d}{d\pi_K}$ and $d_i=\frac{d}{d[T^{\flat}_i]}$. This justifies our normalization.
\end{remark}

\begin{remark}
     In the case of CDVFs with imperfect residue field, Ohkubo proposes a construction of $\mathbf{N}_{\mathrm{dR}}$ in \cite{Ohkubo-differential}. It is in fact a prototype of the relative $\mathbf{N}_{\mathrm{dR}}$ functor as it contains both arithmetic and geometric directions. However, in the proof of \cite[Proposition~4.2.12(i)]{Ohkubo-differential}, it is unclear to us how the faithful flatness of the inclusion $B_n \hookrightarrow \mathbb{B}_{\mathrm{dR}}^+$ produces the required isomorphism by constructing an isomorphism after taking base change along ${B_n} \hookrightarrow \mathbb{B}_{\mathrm{dR}}^+$. Hence, the construction in \emph{loc. cit.} seems to be incomplete.
\end{remark}

\subsection{Compatibility of $\fN_{\mathrm{dR}}$} \label{subsec: compatibility-NdR}

In this subsection, we study the compatibility of $\fN_{\mathrm{dR}}$ under certain extensions. Let $h_1\colon X_1 \rightarrow \D_K^d$ be a toric chart, and let $\Lambda_1 = (R_{\Lambda_1}, \fA_{\Lambda_1})$ be a strongly decompleting datum for $(X_1, h_1)$. Suppose we have a map $g\colon X_2 \rightarrow X_1$ of $K$-affinoid spaces such that the composite $h_2\coloneqq h_1\circ g\colon X_2 \rightarrow \D_K^d$ is a toric chart. Let $\Lambda_2 = (R_{\Lambda_2}, \fA_{\Lambda_2})$ be a strongly decompleting datum for $(X_2, h_2)$. Assume further that we have a map $R_{\Lambda_1} \rightarrow R_{\Lambda_2}$ of $\kappa(\!(\overline{\pi}_K)\!)$-algebras (resp. a map $\fA_{\Lambda_1} \rightarrow \fA_{\Lambda_2}$ of $\fA_{K}$-algebras) compatible with the map $\widetilde{R}_{X_1} \rightarrow \widetilde{R}_{X_2}$ (resp. $\widetilde{\fA}_{X_1} \rightarrow \widetilde{\fA}_{X_2}$). Note that this induces $(\varphi, \Gamma)$-equivariant maps $\fA_{\Lambda_1}^r \rightarrow \fA_{\Lambda_2}^r ~~\text{ and }~~ \fC_{\Lambda_1}^{[s, r]} \rightarrow \fC_{\Lambda_2}^{[s, r]}$, for all $r > 0$ sufficiently small and $0 < s \le r$. Consider the $\fN_{\mathrm{dR}}$ functor associated with $\Lambda_1$ (resp. $\Lambda_2$).

\begin{proposition} \label{prop: compatibility-NdR}
Let $\L \in \mathrm{Loc}^{\mathrm{dR}}_{\Z_p}(X_1)$ be a de Rham local system on $X_1$. Then there is natural isomorphism
\[
\mathbf{N}_{\mathrm{dR}}(\L)\otimes_{\fC_{\Lambda_1}^{\dagger}} \fC_{\Lambda_2}^{\dagger} \isom  \mathbf{N}_{\mathrm{dR}}(\L |_{X_2})
\]
of $(\varphi, \Gamma)$-modules.    
\end{proposition}

\begin{proof}
For $r > 0$ sufficiently small and $0 < s \le r$, let $M_1^{[s, r]}$ (resp. $M_2^{[s, r]}$) denote the finite projective $\fC_{\Lambda_1}^{[s, r]}$-module (resp. $\fC_{\Lambda_2}^{[s, r]}$-module) given by Theorem~\ref{thm: equivalence-isogeny-local-system-overconvergent-etale-(phi, Gamma)-modules-C} obtained from $\L$ and $\Lambda_1$ (resp. from $\L |_{X_2}$ and $\Lambda_2$). Write $N_1^{[s, r]}$ (resp. $N_2^{[s, r]}$) for the finite projective $\fC_{\Lambda_1}^{[s, r]}$-module (resp. $\fC_{\Lambda_2}^{[s, r]}$-module) in Construction~\ref{const: NdR} obtained from $\L$ and $\Lambda_1$ (resp. from $\L |_{X_2}$ and $\Lambda_2$). It suffices to produce compatible $(\varphi, \Gamma)$-equivariant isomorphisms
\[
N_1^{[s, r]} \otimes_{\fC_{\Lambda_1}^{[s, r]}} \fC_{\Lambda_2}^{[s, r]} \isom N_2^{[s, r]}
\]
for $0 < s \le r$ with $r$ sufficiently small. By construction, this amounts to constructing two isomorphisms after taking base changes along $\fC_{\Lambda_2}^{[s, r]} \rightarrow \fC_{\Lambda_2}^{[s, r]}[\frac{1}{t}]$ and $\iota_n\colon \fC_{\Lambda_2}^{[s, r]} \rightarrow A_{X_2, n}[\![t, u_1, \ldots, u_d]\!]^{\nabla = 0}$ respectively, such that the two isomorphisms agree over $A_{X_2, n}[\![t, u_1, \ldots, u_d]\!]^{\nabla = 0}[\frac{1}{t}]$. 

Note that, by construction, $\widetilde{D}^{\mathrm{int}, \dagger}$ in Corollary~\ref{cor: equiv-local-system-overconvergent-etale-(phi, Gamma)-mods} is compatible with the extension $\widetilde{\fA}_{X_1}^{\dagger} \rightarrow \widetilde{\fA}_{X_2}^{\dagger}$. So, by Corollary~\ref{cor: decompletion-etale-(phi, Gamma)-modules}, we have
\[
M_1^{[s, r]} \otimes_{\fC_{\Lambda_1}^{[s, r]}} \fC_{\Lambda_2}^{[s, r]} \cong M_2^{[s, r]}
\]
for $r$ sufficiently small. Thus, over $\fC_{\Lambda_2}^{[s, r]}[\frac{1}{t}]$, we have (for $r$ sufficiently small)
\[
N_{\Lambda_1}^{[s, r]} \otimes_{\fC_{\Lambda_1}^{[s, r]}} \fC_{\Lambda_2}^{[s, r]}[\frac{1}{t}] \cong M_{\Lambda_1}^{[s, r]} \otimes_{\fC_{\Lambda_1}^{[s, r]}} \fC_{\Lambda_2}^{[s, r]}[\frac{1}{t}] \cong M_{\Lambda_2}^{[s, r]}[\frac{1}{t}] \cong N_{\Lambda_2}^{[s, r]}[\frac{1}{t}]. 
\]

On the other hand, there is a natural isomorphism
\[
D_{\mathrm{dR}}(\L)\otimes_{A_{X_1}} A_{X_2} \cong D_{\mathrm{dR}}(\L |_{X_2})
\]
by the base change property of the $D_{\mathrm{dR}}$ functor. Write
\[
\mathcal{D}^{\nabla = 0}\coloneqq (D_{\mathrm{dR}}(\L)\otimes_{A_{X_1}} A_{X_1}[\![t, u_1, \ldots, u_d]\!])^{\nabla = 0}.
\]
Recall from \S \ref{subsec: NdR-preparation} the natural isomorphisms
\[
\mathcal{D}^{\nabla = 0} \otimes_{A_{X_1}[\![t, u_1, \ldots, u_d]\!]^{\nabla = 0}}  A_{X_1}[\![t, u_1, \ldots, u_d]\!] \cong D_{\mathrm{dR}}(\L)\otimes_{A_{X_1}} A_{X_1}[\![t, u_1, \ldots, u_d]\!]
\]
and
\[
(D_{\mathrm{dR}}(\L)\otimes_{A_{X_1}} A_{X_1, n}[\![t, u_1, \ldots, u_d]\!])^{\nabla = 0} \cong \mathcal{D}^{\nabla = 0} \otimes_{A_{X_1}[\![t, u_1, \ldots, u_d]\!]^{\nabla = 0}}  A_{X_1, n}[\![t, u_1, \ldots, u_d]\!]^{\nabla = 0}
\]
for each $n$. Hence, over $A_{X_2}[\![t, u_1, \ldots, u_d]\!]^{\nabla = 0}$, we obtain isomorphisms
\begin{align*}
 &  (D_{\mathrm{dR}}(\L)\otimes_{A_{X_1}} A_{X_1, n}[\![t, u_1, \ldots, u_d]\!])^{\nabla = 0} \otimes_{A_{X_1, n}[\![t, u_1, \ldots, u_d]\!]^{\nabla = 0}} A_{X_2, n}[\![t, u_1, \ldots, u_d]\!]^{\nabla = 0}\\
    \cong  \:\:  &  (\mathcal{D}^{\nabla = 0} \otimes_{A_{X_1}[\![t, u_1, \ldots, u_d]\!]^{\nabla = 0}}  A_{X_1, n}[\![t, u_1, \ldots, u_d]\!]^{\nabla = 0}) \otimes_{A_{X_1, n}[\![t, u_1, \ldots, u_d]\!]^{\nabla = 0}} A_{X_2, n}[\![t, u_1, \ldots, u_d]\!]^{\nabla = 0}\\
      \cong  \:\:  & \mathcal{D}^{\nabla = 0} \otimes_{A_{X_1}[\![t, u_1, \ldots, u_d]\!]^{\nabla = 0}}  A_{X_2, n}[\![t, u_1, \ldots, u_d]\!]^{\nabla = 0} \\
      \cong \:\:  & (\mathcal{D}^{\nabla = 0} \otimes_{A_{X_1}[\![t, u_1, \ldots, u_d]\!]^{\nabla = 0}}  A_{X_2, n}[\![t, u_1, \ldots, u_d]\!])^{\nabla = 0}\\
    \cong \:\:  &  (D_{\mathrm{dR}}(\L |_{X_2})\otimes_{A_{X_2}} A_{X_2, n}[\![t, u_1, \ldots, u_d]\!])^{\nabla = 0}.
\end{align*}

It remains to show that the two isomorphisms agree over $A_{X_2, n}[\![t, u_1, \ldots, u_d]\!]^{\nabla = 0}[\frac{1}{t}]$. Let $\iota_{n, 1}$ (resp. $\iota_{n, 2}$) be the isomorphism \eqref{eq:gluing_isom_relative_s,r} for $X_1$ (resp. for $X_2$). We need to prove that the following diagram commutes
\[
\begin{tikzcd}
    M_1^{[s,r]} \otimes_{\mathbf{C}_{\Lambda_1}^{[s,r]}, \iota_n } A_{X_1, n}[\![t, u_1, \ldots, u_d]\!]^{\nabla = 0}[\frac{1}{t}] \arrow[r, "\iota_{n, 1}"] \arrow[d]
    & (D_{\mathrm{dR}}(\L)\otimes_{A_{X_1}} A_{X_1, n}[\![t, u_1, \ldots, u_d]\!])^{\nabla = 0}[\frac{1}{t}] \arrow[d] \\
    M_2^{[s,r]} \otimes_{\mathbf{C}_{\Lambda_2}^{[s,r]}, \iota_n } A_{X_2, n}[\![t, u_1, \ldots, u_d]\!]^{\nabla = 0}[\frac{1}{t}]   \arrow[r, "\iota_{n, 2}"]
    & (D_{\mathrm{dR}}(\L|_{X_2})\otimes_{A_{X_2}} A_{X_2, n}[\![t, u_1, \ldots, u_d]\!])^{\nabla = 0}[\frac{1}{t}].
\end{tikzcd}
\] 
Consider $\overline{A}_{X_1}$ (resp. $\overline{A}_{X_2}$) as in the proof of Lemma~\ref{lemma:gluing_isomorphism_perfect_relative} for $X = X_1$ (resp. $X = X_2$). Choose a map $\overline{A}_{X_1} \rightarrow \overline{A}_{X_2}$ extending $A_{X_1, \infty} \rightarrow A_{X_2, \infty}$. This induces a map $\overline{A}_{X_1}^{\wedge, \flat} \rightarrow \overline{A}_{X_2}^{\wedge, \flat}$. For $i = 1, 2$, write
\[
\widetilde{\fA}_i \coloneqq W(\overline{A}_{X_i}^{\wedge, \flat}), \quad \mathbf{B}_{\mathrm{{dR}}, \overline{X}_i} \coloneqq \mathbb{B}_{\mathrm{dR}, X_i}(\overline{A}_{X_i}^{\wedge, \flat}).
\]
Denote by $H_{X_i}$ (for $i = 1, 2$) the Galois group $H_X$ in the proof of Lemma~\ref{lemma:gluing_isomorphism_perfect_relative} for $X = X_i$. By the construction of $\widetilde{D}^{\mathrm{int}, \dagger}$ in Corollary~\ref{cor: equiv-local-system-overconvergent-etale-(phi, Gamma)-mods}, there is $r > 0$ sufficiently small such that
\[
\widetilde{D}^{\mathrm{int}, \dagger}(\L) = (\L\otimes_{\Z_p} \widetilde{\fA}_1^r)^{H_{X_1}}\otimes_{\fA_{X_1}^r} \fA_{X_1}^{\dagger}
\]
and that
\[
\widetilde{D}^{\mathrm{int}, \dagger}(\L|_{X_2}) = (\L\otimes_{\Z_p} \widetilde{\fA}_2^r)^{H_{X_2}}\otimes_{\fA_{X_2}^r} \fA_{X_2}^{\dagger}.
\]
Hence, by the construction of $\iota_{n, 1}$ and $\iota_{n, 2}$ in \S \ref{sec:relative_N_dR} and Lemma~\ref{lemma:gluing_isomorphism_perfect_relative}, the desired commutativity follows from the commutativity of
\[
\begin{tikzcd}
    (\L\otimes_{\Z_p} \widetilde{\fA}_1^r)^{H_{X_1}} \arrow[r, "\iota_n"] \arrow[d]
    & (\L\otimes_{\Z_p} \mathbf{B}_{\mathrm{{dR}}, \overline{X}_1})^{H_{X_1}} \arrow[d] \\
    (\L\otimes_{\Z_p} \widetilde{\fA}_2^r)^{H_{X_2}}    \arrow[r, "\iota_n"]
    & (\L\otimes_{\Z_p} \mathbf{B}_{\mathrm{{dR}}, \overline{X}_2})^{H_{X_2}},
\end{tikzcd}
\]
where $\iota_n$ in top and bottom horizontal maps are given by $\iota_n\colon \widetilde{\fA}_i^r \rightarrow \mathbf{B}_{\mathrm{{dR}}, \overline{X}_i}$.

Lastly, by the construction of $(\varphi, \Gamma)$-structure in Construction~\ref{const: NdR}, it follows that the base change isomorphism is compatible with $\varphi$-actions and $\Gamma$-actions.    
\end{proof}

\vspace{0.1in}
\section{Examples of strongly decompleting data} \label{sec:example} 

In this section, we provide examples of strongly decompleting data for closed annuli (\S \ref{subsec: closed-annuli}), and those for ``punctured closed annuli'', that is, closed annuli of outer radii $1$ with finitely many open unit disks removed (\S \ref{subsec: punctured-closed-annuli}). For simplicity, we restrict our attention to the one-dimensional situation; the higher dimensional case of closed poly-annuli (with punctures) is similar, with differences that are purely notational. For the two examples, we also explicitly compute the normalized differential operators in \S \ref{subsec: normalized-differential-operators-closed-annuli}, and study the compatibility of the $\fN_{\mathrm{dR}}$ functor under base field extensions in \S \ref{subsec: compatibility-NdR-base-field-extension}. Lastly in \S \ref{subsec: relation-with-AB}, we briefly discuss how our notion of strongly decompleting data relates to the work of Andreatta--Brinon in \cite{Andreatta-Brinon_Overconvergence} and \cite{Andreatta-Brinon}.

Throughout the section, we retain the notation in \S \ref{ss: notation and convention} and Notation~\ref{notn: annuli}.

\subsection{Closed annuli} \label{subsec: closed-annuli}

We first consider the case where $X$ is a closed annulus. Let $X = X_{[\sigma_1, \sigma_2]} = A_K([\sigma_1, \sigma_2]; T)$ with $0 < \sigma_1 \le \sigma_2 \le 1$ (see Notation~\ref{notn: annuli}), and take the natural inclusion $h\colon X \hookrightarrow \D_K = \Spa(K\langle T \rangle, \mathcal{O}_K\langle T \rangle)$ (sending $T \mapsto T$) to be the toric chart. In this case, we have $\Gamma = \Gamma_{\mathrm{geom}} \rtimes \Gamma_K$ with $\Gamma_{\mathrm{geom}} \cong \Z_p$ and $\Gamma_K \subset \Z_p^{\times}$ is an open subgroup. We shall construct a natural example of a strongly decompleting datum for $(X, h)$. More generally, let $g\colon X' \rightarrow X$ be a finite \'etale morphism, and consider $h'\coloneqq h\circ g\colon X' \rightarrow \D_K$. Then the strongly decompleting datum for $(X,h)$ induces a strongly decompleting datum for $(X',h')$ by Theorem~\ref{thm: strongly decompleting under finite etale extensions}.

Note that $\widetilde{R}_X$ is the rational localization of $\widetilde{R}_{\D_K} = \widehat{K}_{\infty}^{\flat}\langle(T^{\flat})^{\frac{1}{p^{\infty}}}\rangle$ given by $\sigma_1 \le |T^{\flat}| \le \sigma_2$. Namely, $\widetilde{R}_X$ consists of elements of the form $\sum_{j\in \Z[p^{-1}]} a_j (T^{\flat})^j$ with $a_j \in \widehat{K}_{\infty}^{\flat}$ such that for any $\epsilon > 0$, there are only finitely many indices $j$ such that $\max\{|a_j| \sigma_1^j, |a_j| \sigma_2^j\} \ge \epsilon$. Consider the following subrings of $\widetilde{R}_X$ and $\widetilde{\fA}_X = W(\widetilde{R}_X)$.

\begin{definition} \label{def: closed-annuli-Lambda} \hfill
\begin{itemize}
\item Let $R_{[\sigma_1, \sigma_2]}$ be the ring of functions on the closed annulus $A_{\kappa(\!(\overline{\pi}_K)\!)}([\sigma_1, \sigma_2]; T^{\flat}) \subset \D_{\kappa(\!(\overline{\pi}_K)\!)}$. That is, 
\begin{align*}
R_{[\sigma_1, \sigma_2]} & \coloneqq \Gamma_{\kappa(\!(\overline{\pi}_K)\!)}([\sigma_1, \sigma_2]; T^{\flat}) \\
    & = \left\{\sum_{i\in \Z} a_i(T^{\flat})^i ~|~ a_i \in \kappa(\!(\overline{\pi}_K)\!), ~~|a_i|\sigma^i \rightarrow 0 \textrm{ as } i \rightarrow \pm \infty \textrm{ for all } \sigma_1 \le \sigma \le \sigma_2 \right\} \subset \widetilde{R}_X.    
\end{align*}
 
\item For each $0 < r < r_K$, let
\[
\fA_{[\sigma_1, \sigma_2]}^r\coloneqq \left\{\sum_{i\in \Z} b_i [T^{\flat}]^i ~|~ b_i\in \fA_K^r, ~~ \lambda_r(b_i)\sigma^{ri} \rightarrow 0 \textrm{ as } i \rightarrow \pm \infty \textrm{ for all } \sigma_1 \le \sigma \le \sigma_2 \right\}.
\]
There is a natural inclusion $\fA_{[\sigma_1, \sigma_2]}^r \subset \widetilde{\fA}_X^r$.

\item Let 
\[
\fA_{[\sigma_1, \sigma_2]}^{\dagger}\coloneqq \varinjlim_{0 < r < r_K} \fA_{[\sigma_1, \sigma_2]}^r,
\]
and let $\fA_{[\sigma_1, \sigma_2]}$ be the $p$-adic completion of $\fA_{[\sigma_1, \sigma_2]}^{\dagger}$. Note that $\fA_{[\sigma_1, \sigma_2]}$ is a subring of $\widetilde{\fA}_X$, and every element in $\fA_{[\sigma_1, \sigma_2]}$ can be written uniquely as $\sum_{i \in \Z} b_i [T^{\flat}]^i$ with $b_i \in \fA_K$. 
\end{itemize}    
\end{definition}

\noindent \textbf{Step 1.} We first show that $\Lambda_{[\sigma_1, \sigma_2]}\coloneqq (R_{[\sigma_1, \sigma_2]}, \fA_{[\sigma_1, \sigma_2]})$ is weakly decompleting.\\
We start with some preliminary facts on the $\lambda_r$ norm. For any $0 < r < r_K$, consider the norm $\lambda'_r$ on $\fA_{[\sigma_1, \sigma_2]}^r$ defined by
\[
\lambda'_r\left(\sum_{i \in \Z} b_i [T^{\flat}]^i\right) \coloneq \sup_{\sigma_1 \le \sigma \le \sigma_2} \sup_{i \in \Z} \{\lambda_r(b_i)\sigma^{ri} \}.
\] 

\begin{lemma} \label{lem: closed-annuli-lambda'-equals-lambda}
For each $0 < r < r_K$, we have $\lambda'_r = \lambda_r$ on $\fA_{[\sigma_1, \sigma_2]}^r$.    
\end{lemma}

To prove this, first note that $\lambda'_r \ge \lambda_r$ by the definition. Let $s\colon \mathbf{E}_K \rightarrow \fA_K$ be the section of the mod $p$ reduction map in Construction~\ref{construction:section_s_from_E_K}, which satisfies the properties in Lemma~\ref{lemma:section}. For each $y = \sum_{i\in \Z} y_i(T^{\flat})^i \in R_{[\sigma_1, \sigma_2]}$, set
$
s(y)\coloneqq \sum_{i\in \Z} s(y_i)[T^{\flat}]^i.
$ 
Then for any $0 < r < r_K$, we have
\[
\lambda_r(s(y_i))\sigma^{ri} = |y_i|^r \sigma^{ri} \rightarrow 0
\]
as $i \rightarrow \pm \infty$ for all $\sigma_1 \le \sigma \le \sigma_2$, and thus $s(y) \in \fA_{[\sigma_1, \sigma_2]}^r$. This gives a section $s\colon R_{[\sigma_1, \sigma_2]} \rightarrow \fA_{[\sigma_1, \sigma_2]}^r$ of the mod $p$ reduction map. In particular, $\fA_{[\sigma_1, \sigma_2]}^r \rightarrow R_{[\sigma_1, \sigma_2]}$ is surjective. Note that
\[
\lambda'_r(s(y)) = \sup_{\sigma_1 \le \sigma \le \sigma_2} \sup_{i \in \Z} \{\lambda_r(s(y_i))\sigma^{ri} \} = \sup_{\sigma_1 \le \sigma \le \sigma_2} \sup_{i \in \Z} \{|y_i|^r\sigma^{ri} \} = |y|^r.
\]
Since $\lambda_r' \ge \lambda_r$, we must have
\begin{equation} \label{eq: lambda'-equals-lambda-section}
\lambda_r'(s(y)) = \lambda_r(s(y)) = |y|^r.    
\end{equation}

\begin{lemma} \label{lemma: 3.5 weakly decompleting}
Let $x\in \fA_{[\sigma_1, \sigma_2]}$ and let $\overline{x}\in R_{[\sigma_1, \sigma_2]}$ be the reduction of $x$ mod $p$. Then for each $k\in \mathbb{Z}_{\ge 0}$, we have
\[
w_k\left(\frac{x-s(\overline{x})}{p}\right)\ge \inf\left\{w_{k+1}(x), w_0(x)-\frac{k+1}{r_K}\right\}.
\]
\end{lemma}

\begin{proof}
By Lemma~\ref{lemma: w_k}, we have 
\[
w_k\left(\frac{x-s(\overline{x})}{p}\right)\ge \inf\left\{w_{k+1}(x), w_0(x)-\frac{k+1}{r}\right\}
\]
for all $0<r<r_K$, so the statement holds.
\end{proof}

Now for each $x \in \fA_{[\sigma_1, \sigma_2]}$, we construct a sequence $x_0, x_1, \ldots \in \fA_{[\sigma_1, \sigma_2]}$ inductively by setting $x_0=x$ and $x_{n+1}=\frac{1}{p}(x_n-s(\overline{x}_n))$, where $\overline{x}_n\in R_{[\sigma_1, \sigma_2]}$ stands for $x_n$ mod $p$. Then we have $x = \sum_{n=0}^{\infty}p^n \cdot s(\overline{x}_n)$.

\begin{lemma} \label{lemma: 3.6 weakly decompleting}
Let $x\in \mathbf{A}_{[\sigma_1, \sigma_2]}$ and write $x=\sum_{n=0}^{\infty}p^n \cdot s(\overline{x}_n)$ as above. Then for each $n\in \mathbb{Z}_{\ge 0}$, we have 
\[
v(\overline{x}_n)\ge \inf_{i\le n}\left\{w_i(x)-\frac{n-i}{r_K}\right\}.
\]
\end{lemma}

\begin{proof}
The proof is similar to \cite[Lemma~7.4]{Colmez08}. By Lemma~\ref{lemma: 3.5 weakly decompleting}, we have 
\[
w_k(x_{n+1})\ge \inf\left\{w_{k+1}(x_n), w_0(x_n)-\frac{k+1}{r_K}\right\}.
\]
By induction on $n$, we deduce that
\[
w_k(x_n)\ge \inf\left\{w_{k+n}(x), \inf_{0\le i\le n-1}\left\{w_i(x)-\frac{k+n-i}{r_K}\right\}\right\}.
\]
Taking $k=0$ yields the statement.
\end{proof}

\begin{proof}[Proof of Lemma~\ref{lem: closed-annuli-lambda'-equals-lambda}]
For $x\in \mathbf{A}^r_{[\sigma_1, \sigma_2]}$, set $v'_r(x)\coloneqq -\frac{1}{r}\log_p(\lambda_r(x))$. We need to show $v'_r(x)\ge v_r(x)$. Write $x=\sum_{n=0}^{\infty}p^n \cdot s(\overline{x}_n)$ as above. This is a norm compatible expansion by the identity \eqref{eq: lambda'-equals-lambda-section}, so
\[
\lambda_r(x) = \max_{n \ge 0} \{p^{-n}\lambda_r(s(\overline{x}_n))\}
\]
by Proposition~\ref{prop: tech lemma norm compatible expansion}.

It suffices to show that $v'_r(p^n\cdot s(\overline{x}_n))\ge v_r(x)$ for all $n$. By Lemma~\ref{lemma: 3.6 weakly decompleting}, we have
\begin{align*}
&v'_r(p^n\cdot s(\overline{x}_n))= \frac{n}{r}+v(\overline{x}_n)
\ge  \frac{n}{r}+\inf_{0\le i\le n}\left\{w_i(x)-\frac{n-i}{r_K}\right\}\\
=& \inf_{0\le i\le n} \left\{(w_i(x)+\frac{i}{r})+(n-i)(\frac{1}{r}-\frac{1}{r_K})\right\}
\ge  \inf_{0\le i\le n} \left\{w_i(x)+\frac{i}{r}\right\}
\ge   v_r(x)
\end{align*}
as desired.
\end{proof}

\begin{corollary} \label{cor: closed-annuli-A^r-C} \hfill
\begin{enumerate}
\item For each $0 < r < r_K$, $\fA_{[\sigma_1, \sigma_2]}^r$ is complete with respect to $\lambda_r$, and
\[
\fA_{[\sigma_1, \sigma_2]}^r = \widetilde{\fA}_X^r \cap \fA_{[\sigma_1, \sigma_2]}
\]
as subrings of $\widetilde{\fA}_X$.

\item For each $0 < s \le r < r_K$, define $\fB_{[\sigma_1, \sigma_2]}^r$ and $\fC_{[\sigma_1, \sigma_2]}^{[s, r]}$ as in Definition~\ref{def: imperfect-period-rings-B-C}. We have
\[
\fC_{[\sigma_1, \sigma_2]}^{[s, r]} = \left\{\sum_{i\in \Z} b_i[T^{\flat}]^i ~|~ b_i \in \fC_K^{[s, r]}, ~~\lambda_t(b_i)\sigma^{ti} \rightarrow 0 \textrm{ as } i \rightarrow \pm \infty \textrm{ for all } t \in [s, r] \textrm{ and } \sigma_1 \le \sigma \le \sigma_2 \right\}.
\]
\end{enumerate}    
\end{corollary}

\begin{proof}
It directly follows from Lemma~\ref{lem: closed-annuli-lambda'-equals-lambda} that $\fA_{[\sigma_1, \sigma_2]}^r$ is complete with respect to $\lambda_r$. By the definition, we have $\fA_{[\sigma_1, \sigma_2]}^r \subset \widetilde{\fA}_X^r \cap \fA_{[\sigma_1, \sigma_2]}$. For the converse inclusion, let $x\in \widetilde{\fA}_X^r \cap \fA_{[\sigma_1, \sigma_2]}$ and write $x=\sum_{n=0}^{\infty}p^n\cdot s(\bar{x}_n)$ as above. Since each $p^n\cdot s(\bar{x}_n)\in \mathbf{A}^r_{[\sigma_1, \sigma_2]}$, it suffices to show that $v'_r(p^n\cdot s(\bar{x}_n))\rightarrow +\infty$ as $n\rightarrow +\infty$. Recall from the proof of Lemma~\ref{lem: closed-annuli-lambda'-equals-lambda} that
\[
v'_r(p^n\cdot s(\overline{x}_n))\ge  \inf_{0\le i\le n} \left\{(w_i(x)+\frac{i}{r})+(n-i)(\frac{1}{r}-\frac{1}{r_K})\right\}.
\]
The last term goes to $+\infty$ as $n\rightarrow +\infty$, by the same argument as in the proof of Proposition~\ref{prop: good lift weakly decompleting}(2). This proves (1). By the definition of $\fC_{[\sigma_1, \sigma_2]}^{[s, r]}$, part (2) also follows from Lemma~\ref{lem: closed-annuli-lambda'-equals-lambda}. 
\end{proof}

\begin{proposition} \label{prop: closed-annuli-weakly-decompleting}
$\Lambda_{[\sigma_1, \sigma_2]} = (R_{[\sigma_1, \sigma_2]}, \fA_{[\sigma_1, \sigma_2]})$ is a weakly decompleting datum for $(X,h)$.  
\end{proposition}

\begin{proof}
Recall from Corollary~\ref{cor: closed-annuli-A^r-C}(1) that $\fA_{[\sigma_1, \sigma_2]}^r = \widetilde{\fA}_X^r \cap \fA_{[\sigma_1, \sigma_2]}$ for each $0 < r< r_K$. We need to verify the conditions (1)--(8) in Definition~\ref{defn: weakly decompleting datum} (for some sufficiently small $r_* > 0$). We will assume $0 < r < r_K$ below.

\begin{enumerate}
\item Since the extension $\kappa(\!(\overline{\pi}_K)\!) / \kappa(\!(\pi)\!)$ is separable, $\{\overline{\pi}^{\frac{i}{p}} ~|~ 0 \le i \le p-1\}$ gives a $p$-basis of $\kappa(\!(\overline{\pi}_K)\!)$. We claim that $R_{[\sigma_1, \sigma_2]}$ admits a $p$-basis given by $\{\overline{\pi}^{\frac{i}{p}}(T^{\flat})^{\frac{j}{p}} ~|~ 0 \le i, j \le p-1\}$. Indeed, it follows from the explicit description of $R_{[\sigma_1, \sigma_2]}$ that these elements generate $\varphi^{-1}(R_{[\sigma_1, \sigma_2]})$ as an $R_{[\sigma_1, \sigma_2]}$-module. Moreover, these elements do not have any non-trivial linear relation over $R_{[\sigma_1, \sigma_2]}$ because they provide a $p$-basis of $\kappa(\!(\overline{\pi}_K)\!)[\![(T^{\flat})^{\pm 1}]\!]$ by the proof of \cite[Lemma~1.1.3]{deJong-dieudonnemodule} and $R_{\sigma_1, \sigma_2} \subset \kappa(\!(\overline{\pi}_K)\!)[\![(T^{\flat})^{\pm 1}]\!]$.

\item Note that $\kappa(\!(\overline{\pi}_K)\!)^{\mathrm{perf}}$ is dense in $\widehat{K}_{\infty}^{\flat}$. So it follows the definition that $(R_{[\sigma_1, \sigma_2]})^{\mathrm{perf}}$ is dense in $\widetilde{R}_{[\sigma_1, \sigma_2]}$.

\item It is straightforward to check that $\fA_{[\sigma_1, \sigma_2]}^r$ is stable under the $\Gamma$-action. Thus, $\fA_{[\sigma_1, \sigma_2]}^{\dagger}$ is stable under the $\Gamma$-action, and so is its $p$-adic completion $\fA_{[\sigma_1, \sigma_2]}$.

\item This is noted in Corollary~\ref{cor: closed-annuli-A^r-C}(1).

\item We already observed that $\fA_{[\sigma_1, \sigma_2]}^r/p \rightarrow R_{[\sigma_1, \sigma_2]}$ is surjective, so it suffices to show that the map is injective. Let $x, y \in \fA_{[\sigma_1, \sigma_2]}^r$, and write
\[
x = \sum_{n=0}^{\infty} p^n s(\overline{x}_n), \quad y = \sum_{n=0}^{\infty} p^n s(\overline{y}_n)
\]
as above. Suppose $\overline{x} = \overline{y} \in R_{[\sigma_1, \sigma_2]}$. Then $s(\overline{x}_0) = s(\overline{y}_0)$, so 
\[
x-y = \sum_{n=1}^{\infty} p^n(s(\overline{x}_n) - s(\overline{y}_n)) \in p\fA_{[\sigma_1, \sigma_2]}^r,
\]
as desired.

\item For any $y \in R_{[\sigma_1, \sigma_2]}$, we have seen that $s(y) \in \fA_{[\sigma_1, \sigma_2]}^r$ and $\lambda_r(s(y)) = |y|^r$. On the other hand, let $x \in \fA_{[\sigma_1, \sigma_2]}^r$ and write
$ 
x = \sum_{n=0}^{\infty} p^n s(\overline{x}_n).
$ 
Since this is a norm compatible expansion, by Proposition~\ref{prop: tech lemma norm compatible expansion}, we have
$
\lambda_r(x) = \max_{n \ge 0} \{p^{-n}\lambda_r(s(\overline{x}_n))\}.
$ 
In particular, $\lambda_r(x) \ge \lambda_r(s(\overline{x}_0)) = |\overline{x}|^r$. Thus, the quotient norm on $\fA_{[\sigma_1, \sigma_2]}^r/p \cong R_{[\sigma_1, \sigma_2]}$ induced from $\lambda_r$ on $\fA_{[\sigma_1, \sigma_2]}^r$ agrees with $|\cdot|^r$ on $R_{[\sigma_1, \sigma_2]}$.

\item It is clear from the definition that $R_{[\sigma_1, \sigma_2]}$ is an affinoid $\kappa(\!(\overline{\pi}_K)\!)$-algebra.

\item Let $0 < s \le r < r_K$ with $s, r \in \Q$. By Proposition~\ref{prop:identification_C_with_annulus}, there exists $r_* > 0$ small enough such that for any $0 < s \le r \le r_*$, we have a natural identification 
\[
\fC_{K}^{[s, r]} \cong \Gamma_{K'_0}([\alpha, \beta]; \pi_K)
\]
where $\alpha = p^{- \frac{p}{p-1} \cdot \frac{r}{e_K}}$ and $\beta = p^{-\frac{p}{p-1} \cdot \frac{s}{e_K}}$. Thus, by Corollary~\ref{cor: closed-annuli-A^r-C}(2), $\fC_{[\sigma_1, \sigma_2]}^{[s, r]}$ is a reduced affinoid $K_0'$-algebra for each $0 < s \le r \le r_*$.  
\end{enumerate} 
\end{proof}

We also consider finite \'etale extensions of $X$.

\begin{definition} \label{def: closed-annuli-finite-etale}
Let $g\colon X' \rightarrow X$ be a finite \'etale morphism and let $h'=h\circ g:X'\rightarrow \D_K$. Let $\Lambda_{[\sigma_1, \sigma_2]}'\coloneqq (R_{[\sigma_1, \sigma_2]}', \fA_{[\sigma_1, \sigma_2]}')$ be the pair obtained by applying Construction~\ref{construction:wdc_finite_etale} to $\Lambda_{[\sigma_1, \sigma_2]} =(R_{[\sigma_1, \sigma_2]}, \mathbf{A}_{[\sigma_1, \sigma_2]})$ and $g\colon X' \rightarrow X$. 
\end{definition}

\begin{corollary}
Let $(X', h')$ be as in Definition~\ref{def: closed-annuli-finite-etale}. Then $\Lambda_{[\sigma_1, \sigma_2]}' = (R'_{[\sigma_1, \sigma_2]}, \fA_{[\sigma_1, \sigma_2]}')$ is a weakly decompleting datum for $(X',h')$.  
\end{corollary}

\begin{proof}
This follows from Proposition~\ref{prop: closed-annuli-weakly-decompleting} and Theorem~\ref{thm: weakly decompleting data finite etale}.     
\end{proof}

\hfill

\noindent \textbf{Step 2.} We now prove that $\Lambda_{[\sigma_1, \sigma_2]} = (R_{[\sigma_1, \sigma_2]}, \fA_{[\sigma_1, \sigma_2]})$ is decompleting.\\
We follow a similar computation as in \cite[\S 5]{Kedlaya_Hochschild-Serre-p-analy-gp}. In particular, we use the following result from \cite{Kedlaya_Hochschild-Serre-p-analy-gp} and refer the reader to Section~3 of \textit{loc. cit.} for relevant notions. 

\begin{theorem}[{\cite[Theorem~4.1]{Kedlaya_Hochschild-Serre-p-analy-gp}}] \label{thm: hoschild-serre-p-analy-gp}
Let $\Gamma_1$ be a profinite $p$-analytic group, and $\Gamma_2 \subset \Gamma_1$ a pro-$p$ procyclic subgroup. Let $M$ be an analytic $\Gamma_1$-module which is a Banach space over a non-archimedean field of characteristic $p$ with a non-trivial norm. If $H^i_{\mathrm{cont}}(\Gamma_2, M) = 0$ for all $i \ge 0$, then $H^i_{\mathrm{cont}}(\Gamma_1, M) = 0$ for all $i \ge 0$. 
\end{theorem}

Since $R_{[\sigma_1, \sigma_2]}'$ is finite \'etale over $R_{[\sigma_1, \sigma_2]}$, any $p$-basis of $R_{[\sigma_1, \sigma_2]}$ is also a $p$-basis of $R_{[\sigma_1, \sigma_2]}'$ by \cite[\S 1.1.2(ii)]{Berthelot-Messing}. So we have an isomorphism
\[
\varphi^{-1}(R_{[\sigma_1, \sigma_2]}') / R_{[\sigma_1, \sigma_2]}' \cong \bigoplus_{0 \le i, j \le p-1, ~(i, j)\neq (0, 0)} (1+\overline{\pi})^{\frac{i}{p}} (T^{\flat})^{\frac{j}{p}}\cdot R_{[\sigma_1, \sigma_2]}',
\]
which is strict by Theorem~\ref{thm: open mapping thm} since $\pi \in R_{[\sigma_1, \sigma_2]}'$. Note that the $\Gamma$-action on $\varphi^{-1}(R_{[\sigma_1, \sigma_2]}') / R_{[\sigma_1, \sigma_2]}'$ is analytic by \cite[Proposition~5.3]{Kedlaya_Hochschild-Serre-p-analy-gp}. Write $\overline{y}_{i, j}\coloneqq (1+\overline{\pi})^{\frac{i}{p}}(T^{\flat})^{\frac{j}{p}}$.

\begin{proposition} \label{prop: closed-annuli-decompleting}
$\Lambda_{[\sigma_1, \sigma_2]} = (R_{[\sigma_1, \sigma_2]}, \fA_{[\sigma_1, \sigma_2]})$ is a decompleting datum for $(X,h)$.    
\end{proposition}

\begin{proof}
Let $\Gamma_0 \subset \Gamma$ be an open normal subgroup. It suffices to show that the cochain complex
\[
\mathcal{C}^{(\bullet)}(\Gamma_0, \varphi^{-1}(R_{[\sigma_1, \sigma_2]}')/R_{[\sigma_1, \sigma_2]}')
\]
is strict exact at $d^{(\ell)}$ for each $\ell \ge 0$. 

We will first show the exactness. By Theorem~\ref{thm: hoschild-serre-p-analy-gp}, it suffices to consider the case where $\Gamma_0$ is of the form
\[
\Gamma_0 \cong p^a\Z_p \rtimes (1+p^b \Z_p) \subset \Z_p\rtimes \Z_p^{\times}
\]
for some integers $a \ge 0$ and $b \ge 1$. 

Let $\gamma \in p^a\Z_p \subset \Gamma_0$ be an element such that
\[
\gamma(1+\overline{\pi}) = 1+\overline{\pi}, \quad \gamma(T^{\flat}) = (1+\overline{\pi})^{p^m}\cdot T^{\flat} 
\]
for some $m \ge 1$ (e.g. we have such $\gamma$ for $m = a$). Let $\lambda \in 1+p^b\Z_p \subset \Gamma_0$ be an element such that 
\[
\lambda(1+\overline{\pi}) = (1+\overline{\pi})^{p^{m'}}(1+\overline{\pi}), \quad \lambda(T^{\flat}) = T^{\flat}
\]
for some $m' \ge 2$.

We check that the following inequalities hold
\begin{equation} \label{eq: closed-annuli-inequality-for-Gamma-action}
|(\gamma^{p^n}-1)\cdot \overline{x}| \le p^{-\frac{p^{n+m+1}}{p-1}}|\overline{x}|,  \quad |(\lambda^{p^n}-1)\cdot \overline{x}| \le |\overline{\pi}|\cdot p^{-\frac{p^{n+m'+1}}{p-1}}|\overline{x}|    
\end{equation}
for all $\overline{x} \in R_{[\sigma_1, \sigma_2]}$ and all integers $n \ge 1$. Indeed, recall that
\[
R_{[\sigma_1, \sigma_2]} = \left\{\sum_{i\in \Z} a_i(T^{\flat})^i ~|~ a_i \in \kappa(\!(\overline{\pi}_K)\!), ~~|a_i|\sigma^i \rightarrow 0 \textrm{ as } i \rightarrow \pm \infty \textrm{ for all } \sigma_1 \le \sigma \le \sigma_2 \right\},
\]
and for $\overline{x} = \sum_{i\in \Z} a_i(T^{\flat})^i \in R_{[\sigma_1, \sigma_2]}$, we have
$
|\overline{x}| = \sup_{\sigma_1 \le \sigma \le \sigma_2} \sup_{i \in \Z} \{|a_i|\sigma^i \}.
$ 
Thus, it suffices to check the inequalities~\eqref{eq: closed-annuli-inequality-for-Gamma-action} for $\overline{x} = a_i({T^{\flat}})^i$, which follows from a direct computation: for example, we have
\[
(\gamma^{p^n}-1)\cdot (T^{\flat})^i = ((1+\overline{\pi})^{i\cdot p^{n+m}}-1)\cdot (T^{\flat})^i,
\]
and
\[
|(1+\overline{\pi})^{i\cdot p^{n+m}}-1| \le |\overline{\pi}^{p^{n+m}}| = p^{-\frac{p^{n+m+1}}{p-1}}.
\]

\noindent \textbf{Claim 1.} For each $0 \le i, j \le p-1$ with $j \neq 0$, $(\gamma^{p^n}-1)$ acts bijectively on $R_{[\sigma_1, \sigma_2]}'\cdot \overline{y}_{i, j}$ for all $n \gg 0$ sufficiently large.\\

\noindent \emph{Proof of Claim 1.} 
For each $n \ge 1$, we have
$(\gamma^{p^n}-1)\cdot \overline{y}_{i, j} = ((1+\overline{\pi})^{jp^{n+m-1}}-1)\cdot \overline{y}_{i, j}$. Write 
\[
a_{j, n} = (1+\overline{\pi})^{jp^{n+m-1}}-1 = (1+\overline{\pi}^{p^{n+m-1}})^j-1.
\]
Since $1 \le j \le p-1$, the term $j\overline{\pi}^{p^{n+m-1}}$ dominates the other terms of $a_{j, n}$ in the binomial expansion, and $|a_{j, n}| = p^{-\frac{p^{n+m}}{p-1}}$. 

On the other hand, we have $|(\gamma^{p^n}-1)\cdot \overline{x}| \le p^{-\frac{p^{n+m+1}}{p-1}}|\overline{x}|$ for $\overline{x} \in R_{[\sigma_1, \sigma_2]}$ by \eqref{eq: closed-annuli-inequality-for-Gamma-action}. Let $x_1, \ldots, x_d \in R_{[\sigma_1, \sigma_2]}'$ be generators of $R_{[\sigma_1, \sigma_2]}'$ as an $R_{[\sigma_1, \sigma_2]}$-module. Then there exists $n_0 \ge 0$ such that
\[
|(\gamma^{p^{n_0}}-1)(x_k)| \le \frac{1}{2}|x_k|
\]
for all $1 \le k \le d$. Let $\ell \ge 0$ be an integer sufficiently large so that $(\frac{1}{2})^{p^{\ell}} \le p^{-\frac{p^{n_0+m+1}}{p-1}}$. Writing $z_k\coloneqq x_k^{p^{\ell}}$ for each $k$, we have
\[
|(\gamma^{p^{n_0}}-1)(z_k)| = |((\gamma^{p^{n_0}}-1)(x_k))^{p^{\ell}}| \le (\frac{1}{2})^{p^\ell}|x_k|^{p^{\ell}} = (\frac{1}{2})^{p^\ell}|z_k| \le p^{-\frac{p^{n_0+m+1}}{p-1}}|z_k|.
\]
Here in the first equality, we used that $R'_{[\sigma_1, \sigma_2]}$ has characteristic $p$. Thus,
\[
|(\gamma^{p^n}-1)(z_k)| \le p^{-\frac{p^{n+m+1}}{p-1}}|z_k|
\]
for all $n \ge n_0$.

Since any finite $p$-basis of $R_{[\sigma_1, \sigma_2]}$ is a $p$-basis of $R_{[\sigma_1, \sigma_2]}'$, $z_1, \ldots, z_d$ also generate $R_{[\sigma_1, \sigma_2]}'$ as an $R_{[\sigma_1, \sigma_2]}$-module. Consider the strict surjection
$
\bigoplus_{k=1}^d R_{[\sigma_1, \sigma_2]}\cdot e_k \twoheadrightarrow R_{[\sigma_1, \sigma_2]}'
$ 
given by $e_k \mapsto z_k$. Let $c > 0$ be a constant such that the residue norm on $R_{[\sigma_1, \sigma_2]}'$ given by the finite-module norm on the left-hand side is bounded above by $c$ times the norm on $R_{[\sigma_1, \sigma_2]}'$. It follows from the above inequalities that
\begin{equation} \label{ineq: gamma-action}
|(\gamma^{p^n}-1)\cdot \overline{x}| \le c\cdot p^{-\frac{p^{n+m+1}}{p-1}}|\overline{x}|    
\end{equation}
for all $\overline{x} \in R_{[\sigma_1, \sigma_2]}'$.
Since
\[
a_{j, n}^{-1}(\gamma^{p^n}-1)(\overline{x}\cdot \overline{y}_{i, j}) = (\overline{x}+a_{j, n}^{-1}(1+a_{j, n})(\gamma^{p^n}-1)(\overline{x}))\cdot \overline{y}_{i, j},
\]
we deduce that $a_{j, n}^{-1}\cdot (\gamma^{p^n}-1)$ acts on $R_{[\sigma_1, \sigma_2]}'\cdot \overline{y}_{i, j}$ as the identity plus an operator of norm strictly less than $1$, for all $n \gg 0$ sufficiently large. Thus, Claim 1 holds.\\

\noindent \textbf{Claim 2.} For each $1 \le i \le p-1$, $(\lambda^{p^n}-1)$ acts bijectively on $R_{[\sigma_1, \sigma_2]}'\cdot \overline{y}_{i, 0}$ for all $n \gg 0$ sufficiently large.\\

\noindent \emph{Proof of Claim 2.} We have
\[
(\lambda^{p^n}-1)\cdot \overline{y}_{i, 0} = ((1+\overline{\pi})^{i((1+p^{m'})^{p^n}-1)/p}-1)\cdot \overline{y}_{i, 0}
\]
for each $n \ge 1$. Since $m' \ge 2$\footnote{This condition is for the case $p=2$. When $p \ge 3$, we can take $m' \ge 1$.} and $1 \le i \le p-1$, we have
\[
v_p(i((1+p^{m'})^{p^n}-1)/p) = m'+n-1.
\]
So $|((1+\overline{\pi})^{i((1+p^{m'})^{p^n}-1)/p}-1)| = p^{-\frac{p^{n+m'}}{p-1}}$. Thus, Claim 2 follows from \eqref{eq: closed-annuli-inequality-for-Gamma-action} and a similar computation as in the proof of Claim 1, by showing that for some constant $c > 0$, we have
\begin{equation} \label{ineq: sigma-action}
|(\lambda^{p^n}-1)\cdot \overline{x}| \le c\cdot p^{-\frac{p^{n+m'+1}}{p-1}}|\overline{x}|   
\end{equation}
for all $\overline{x} \in R_{[\sigma_1, \sigma_2]}'$ and $n \gg 0$.\\

Now by Theorem~\ref{thm: hoschild-serre-p-analy-gp}, we deduce from Claim 1 and Claim 2 that $\mathcal{C}^{(\bullet)}(\Gamma_0, R_{[\sigma_1, \sigma_2]}'\cdot \overline{y}_{i, j})$ is exact at $d^{(\ell)}$ for each $\ell \ge 0$ and $0 \le i, j \le p-1$ with $(i, j) \neq (0, 0)$.

For strictness, we consider an affinoid $\kappa(\!(\overline{\pi}_K)\!)$-algebra
$
S \coloneqq R_{[\sigma_1, \sigma_2]}' \langle \frac{\overline{U}}{\overline{\pi}}, \frac{\overline{\pi}}{\overline{U}} \rangle
$ 
where $\overline{U}$ is an auxiliary free variable. Note that $S$ is equipped with the spectral norm which extends the norm $|\cdot|$ on $R_{[\sigma_1, \sigma_2]}'$, and $R_{[\sigma_1, \sigma_2]}'$ is closed in $S$. Furthermore, setting $\Gamma$ to act trivially on $\overline{U}$, $S$ has the induced $\Gamma$-action extending that on $R_{[\sigma_1, \sigma_2]}$.
Consider the cochain complex
\[
\mathcal{C}_{S}^{(\bullet)}\coloneq \mathcal{C}^{(\bullet)}(\Gamma_0, \bigoplus_{0 \le i, j \le p-1, ~(i, j)\neq (0, 0)} S\cdot \overline{y}_{i, j}).    
\]
By a similar computation as in the proof of Claims 1 and 2 (together with Theorem~\ref{thm: hoschild-serre-p-analy-gp}), the above complex is exact at $d^{(\ell)}$ for each $\ell \ge 0$. Note that $d^{(\ell)}$ is $\kappa(\!(\overline{U})\!)$-linear and $\overline{U} \in S$ is a topologically nilpotent unit, since $|\overline{U}| = |\overline{\pi}|$. Thus, by Theorem~\ref{thm: open mapping thm}, we conclude that $d^{(\ell)}$ in $\mathcal{C}_{S}^{(\bullet)}$ is strict for each $\ell \ge 0$. 

Finally, note that each $R_{[\sigma_1, \sigma_2]}'\cdot \overline{y}_{i, j}$ is closed in $S\cdot \overline{y}_{i, j}$ with compatible norms, and that the differential maps of the complex $\mathcal{C}_{S}^{(\bullet)}$ are compatible with those of $\mathcal{C}_{R_{[\sigma_1, \sigma_2]}'}^{(\bullet)}\coloneq \mathcal{C}^{(\bullet)}(\Gamma_0, \varphi^{-1}(R_{[\sigma_1, \sigma_2]}')/R_{[\sigma_1, \sigma_2]}')$. So for each $\ell \ge 0$, the quotient $\mathcal{C}_{R_{[\sigma_1, \sigma_2]}}^{(\ell)}/\ker(d^{(\ell)})$ is a closed subspace of $\mathcal{C}_{S}^{(\ell)}/\ker(d^{(\ell)})$ with compatible norms, and $\mathrm{image}(d^{(\ell)})$ for $\mathcal{C}_{R_{[\sigma_1, \sigma_2]}'}^{(\ell)}$ is a closed subspace of $\mathrm{image}(d^{(\ell)})$ for $\mathcal{C}_{S}^{(\ell)}$ with compatible norms. Therefore, the strictness of $\mathcal{C}_{S}^{(\bullet)}$ implies the strictness of $\mathcal{C}_{R_{[\sigma_1, \sigma_2]}'}^{(\bullet)}$ at $d^{(\ell)}$ for each $\ell \ge 0$.
\end{proof}

\hfill

\noindent \textbf{Step 3.} Lastly, we prove that $\Lambda_{[\sigma_1, \sigma_2]} = (R_{[\sigma_1, \sigma_2]}, \fA_{[\sigma_1, \sigma_2]})$ is strongly decompleting.

\begin{proposition}\label{prop: closed-annuli-strongly-decompleting}
$\Lambda_{[\sigma_1, \sigma_2]} = (R_{[\sigma_1, \sigma_2]}, \fA_{[\sigma_1, \sigma_2]})$ is a strongly decompleting datum for $(X,h)$.     
\end{proposition}

\begin{proof}
We need to check the conditions (a)-(c) of Definition~\ref{defn: strongly decompleting data} for $\Lambda_{[\sigma_1, \sigma_2]} = (R_{[\sigma_1, \sigma_2]}, \fA_{[\sigma_1, \sigma_2]})$.
By Proposition~\ref{prop:identificaiton_for_mathbf_A^r} and the definition of $\fA_{[\sigma_1, \sigma_2]}^r$, we have
\[
\fA_{[\sigma_1, \sigma_2]}^r = \Big\{\sum_{i, j \in \Z} c_{i, j}\pi_K^i[T^{\flat}]^j ~|~ c_{i, j} \in \mathcal{O}_{K_0'}, ~~|c_{i, j}| \cdot p^{-\frac{pr}{p-1}\cdot \frac{i}{e_K}}\cdot \sigma^{rj}\rightarrow 0 \text{ as } |i|+|j| \rightarrow \infty \textrm{ for all\ } \sigma_1 \le \sigma \le \sigma_2 \Big\}.
\]

We first check the condition (a), i.e., the image of the composition
\[
\mathbf{A}^r_{[\sigma_1, \sigma_2]} \subset \widetilde{\mathbf{A}}^r_X\xrightarrow[]{\iota_n} \mathbb{B}_{\mathrm{dR}, X}^+(X_{\infty})
\]
lies in $A_{X, n}[\![t,u]\!]^{\nabla=0}$ where $u\coloneqq T-[T^{\flat}]$, for all $n\ge -\log_pr$. It suffices to show that the image lies in $A_{X, n}[\![t,u]\!]$. To this end, we compute $\iota_n([T^{\flat}]^{\pm 1})$ and $\iota_n(\pi_K^{\pm 1})$. Indeed, we have
\[\iota_n([T^{\flat}])=[(T^{\flat})^{p^{-n}}]=(T-u)^{p^{-n}}=T^{p^{-n}}(1-\frac{u}{T})^{p^{-n}}.\]
If we write $u'\coloneqq \iota_n([T^{\flat}])-T^{p^{-n}}$, then $A_{X,n}[\![u']\!]\subset A_{X,n}[\![u]\!]$. Taking inverse, we can write $\iota_n([T^{\flat}]^{-1})=T^{-p^{-n}}+u''$ with $A_{X,n}[\![u'']\!]\subset A_{X,n}[\![u]\!]$. To compute $\iota_n(\pi_K)$, recall from Construction~\ref{construction:pi_K} that the polynomial $P_K\in \mathbf{A}^+_{K'_0}[X]$ lifts the defining polynomial $\overline{P}_K\in \mathbf{E}_{K'_0}[\overline{X}]$ of $\overline{\pi}_K$, and that $\pi_K$ is a simple root of $P_K$ lifting $\overline{\pi}_K$. By \cite[Proposition~6.6]{Colmez08}, $\theta\circ\varphi^{-n}(\pi_K)$ is a uniformizer $\varpi_n$ of $K_n$, and it is a simple root of $\theta\circ\varphi^{-n}(P_K)$. By Hensel's Lemma, $\varphi^{-n}(P_K)$ has a unique root in $K_n[\![t]\!]$ lifting $\theta\circ\varphi^{-n}(\pi_K)$. Thus, we have $\iota_n(\pi_K)=\varpi_n+t'$ such that $A_{X,n}[\![t']\!] \subset A_{X,n}[\![t]\!]$. Taking the inverse, we obtain $\iota_n(\pi_K^{-1}) = \varpi_n^{-1}+t''$ satisfying $A_{X,n}[\![t'']\!] \subset A_{X,n}[\![t]\!]$.

Now, let $x = \sum_{i,j\in \mathbb{Z}} c_{i,j} \pi_K^i [T^{\flat}]^j \in \mathbf{A}^r_{[\sigma_1, \sigma_2]}$. We need to show that $\iota_n(x)\in A_{X,n}[\![t,u]\!]$. We may assume $r = p^{-n}$. Write 
$
x=x_{+,+}+x_{+,-}+x_{-,+}+x_{-,-},
$ 
where 
\[
x_{+,+}=\sum_{i,j\ge 0} c_{i,j}\pi_K^i[T^{\flat}]^j, \quad x_{+,-}=\sum_{i\ge 0, j<0} c_{i,j}\pi_K^i[T^{\flat}]^j,  
\]
\[
x_{-,+}=\sum_{i<0,j\ge 0} c_{i,j}\pi_K^i[T^{\flat}]^j, ~~\textrm{ and } ~~ x_{-,-}=\sum_{i,j<0} c_{i,j}\pi_K^i[T^{\flat}]^j.
\]
It suffices to show that $\iota_n(x_{+,+})\in A_{X, n}[\![t', u']\!]$ (resp. $\iota_n(x_{+,-})\in A_{X, n}[\![t', u'']\!]$, $\iota_n(x_{-,+})\in A_{X, n}[\![t'', u']\!]$, $\iota_n(x_{-,-})\in A_{X, n}[\![t'', u'']\!]$). We will only compute for $x_{+,+}$, as the other terms follow from similar computations. We have
\[
\iota_n(x_{+,+}) = \sum_{i,j\ge 0} \iota_n(c_{i,j})\left(\varpi_n+t'\right)^i \left(T^{p^{-n}}+u'\right)^j.
\]
For integers $a,b \ge 0$, the $(t')^a(u')^b$-term has coefficient
\[
w_{a,b}=\sum_{i,j\ge 0}\iota_n(c_{i,j})\varpi_n^{i-a}\cdot T^{p^{-n}(j-b)}\cdot \binom{i}{a}\binom{j}{b}.
\]
For each $j \ge 0$, let 
\[
c_j \coloneqq \sum_{i \ge 0}\iota_n(c_{i,j})\varpi_n^i\cdot\binom{i}{a}\binom{j}{b}.
\]
Then $w_{a,b}=\varpi_n^{-a}\cdot T^{-p^{-n}b}\cdot \sum_{j\ge 0} c_j \cdot T^{p^{-n}j}$. To show $w_{a,b}\in A_{X,n}$, it remains to check that $|c_j|\cdot \sigma^{p^{-n}j}\rightarrow 0$ as $j \rightarrow \infty$, for all $\sigma\in [\sigma_1, \sigma_2]$. Indeed, 
\[
|c_j|\le \max_{i\ge 0}\{|\iota_n(c_{i,j})|\cdot |\varpi_n|^i\}= \max_{i\ge 0}\{|c_{i,j}|\cdot p^{-\frac{i}{e_K \cdot p^{n-1}(p-1)}}\}.
\]
By the assumption, we have
\[
|c_{i,j}|\cdot p^{-\frac{i}{e_K \cdot p^{n-1}(p-1)}}\cdot \sigma^{p^{-n}j} = |c_{i,j}|\cdot p^{-\frac{pr}{(p-1)}\cdot \frac{i}{e_K}}\cdot \sigma^{rj}\rightarrow 0
\]
as $|i|+|j|\rightarrow +\infty$, as desired.

The condition (b) of Definition~\ref{defn: strongly decompleting data} follows directly from the explicit description.

Lastly, we check the condition (c). It suffices to check that there exists an integer $n_0 \ge 0$ such that for any $n \ge n_0$, $b \in \fC_K^{[r_n, r_n]}$ and $i \in \Z$, we have
\[
\lambda_{r_n}((\gamma-1)(b\cdot [T^{\flat}]^i)) \le p^{-\frac{p+1}{p}}\cdot \lambda_{r_n}(b\cdot [T^{\flat}]^i)
\]
for all $\gamma \in \Gamma_n$.

Since $\Gamma = \Gamma_{\mathrm{geom}} \rtimes \Gamma_K$ with $\Gamma_{\mathrm{geom}}$ acting trivially on $\fC_K^{[r_n, r_n]}$, by Lemma~\ref{lemma: estimation Gamma on C}, there exists an integer $n_0 \ge 0$ such that
\[
\lambda_{r_n}((\gamma-1)(b)) \le p^{-\frac{p+1}{p}}\cdot \lambda_{r_n}(b)
\]
for any $b \in \fC_K^{[r_n, r_n]}$ and $\gamma \in \Gamma_n$ with $n \ge n_0$. On the other hand, $\Gamma_K$ acts trivially on $[T^{\flat}]$, and for any $\gamma \in \Gamma_{\mathrm{geom}}$ and $i \in \Z$, we have
\[
(\gamma-1)([T^{\flat}]^i) = ([\epsilon]^{i\cdot \eta(\gamma)}-1)\cdot[T^{\flat}]^i = ((1+\pi)^{i\cdot \eta(\gamma)}-1)\cdot [T^{\flat}]^i
\]
where $\eta \colon \Gamma_{\mathrm{geom}} \isom \Z_p$ is as in the proof of Proposition~\ref{prop:image_of_iota_n_fixed_by_Gamma_n}. So for any $\gamma \in \Gamma_{\mathrm{geom}}$ with $\eta(\gamma) \in p^n \Z_p$,
\[
\lambda_{r_n}((\gamma-1)([T^{\flat}]^i)) = \lambda_{r_n}((1+\pi)^{i\cdot \eta(\gamma)}-1))\cdot \lambda_{r_n}([T^{\flat}]^i) \le p^{-\frac{p}{p-1}}\cdot \lambda_{r_n}([T^{\flat}]^i).
\]
Here, we used in the first equality that $\lambda_{r_n}(b\cdot [T^{\flat}]^i) = \lambda_{r_n}(b)\cdot \lambda_{r_n}([T^{\flat}]^i)$ for any $b \in \fC_K^{[r_n, r_n]}$. The second inequality follows from the computation of $\lambda_{r_n}((1+\pi)^{i\cdot \eta(\gamma)}-1))$ in the proof of Lemma~\ref{lemma: estimation Gamma on C}, using that $i\cdot \eta(\gamma) \in p^n \Z_p$. Therefore, for any $n \ge n_0$, we deduce from the equality 
\[
(\gamma-1)(b\cdot [T^{\flat}]^i) = (\gamma-1)(b)\cdot [T^{\flat}]^i+\gamma(b)\cdot (\gamma-1)([T^{\flat}]^i)
\]
that for all $\gamma \in \Gamma_n$,
\[
\lambda_{r_n}((\gamma-1)(b\cdot [T^{\flat}]^i)) \le p^{-\frac{p+1}{p}}\cdot \lambda_{r_n}(b)\cdot \lambda_{r_n}([T^{\flat}]^i) = p^{-\frac{p+1}{p}}\cdot \lambda_{r_n}(b\cdot [T^{\flat}]^i).
\]
This concludes the proof. 
\end{proof}

As an immediate consequence, we obtain a strongly decompleting datum for $(X',h')$ as well.

\begin{corollary}
Let $(X', h')$ be as in Definition~\ref{def: closed-annuli-finite-etale}. Then 
$\Lambda_{[\sigma_1, \sigma_2]}' = (R'_{[\sigma_1, \sigma_2]}, \fA_{[\sigma_1, \sigma_2]}')$ is a strongly decompleting datum for $(X',h')$.  
\end{corollary}

\begin{proof}
This follows from Proposition~\ref{prop: closed-annuli-strongly-decompleting} and Theorem~\ref{thm: strongly decompleting under finite etale extensions}.
\end{proof}

\subsection{Closed annuli with punctures} \label{subsec: punctured-closed-annuli}

In this subsection, we consider closed annuli $A_K([\delta, 1];T)$ for some $\delta\in (0,1]\cap p^{\Q}$ with finitely many open unit disks removed from the boundary circle $\{|T|=1\}$; they are ``closed annuli with punctures''. The goal is to explicitly construct natural strongly decompleting data for these objects.

Let $q(T) \in \mathcal{O}_{K_0}[T]$ be a monic polynomial of degree $h \ge 1$ of the form \[q(T)=\prod_{i=1}^h (T-[\beta_i])\]
for some $\beta_1,\ldots, \beta_h\in k_K^{\times}$. In particular, $q(0) \in \mathcal{O}_{K_0}^{\times}$. For any $\delta\in (0,1]\cap p^{\Q}$, let $X_{[\delta, 1]}$ denote the closed annulus $A_K([\delta, 1]; T)$ as before, and consider the rational subset
\[
X=X_{\delta, q}\coloneqq X_{[\delta, 1]} \cap \{|q(T)| = 1\},
\]
equipped with the toric chart $h\colon X \hookrightarrow \D_K$ sending $T \mapsto T$. Roughly speaking, $X$ is the closed annulus $X_{[\delta, 1]}$ minus open unit disks centered at $[\beta_i]$'s.

Let $\overline{q}\in k_K[T]$ denote the reduction of $q$ modulo $p$. Then $\widetilde{R}_X$ is the rational localization of $\widetilde{R}_{[\delta, 1]}$ given by $|\overline{q}(T^{\flat})| = 1$. Consider the following subrings of $\widetilde{R}_X$ and $\widetilde{\fA}_X = W(\widetilde{R}_X)$.

\begin{definition} \label{def: punctured-closed-annuli-datum} \hfill
\begin{itemize}
\item Let $R_{\delta, q}$ be a subring of $\widetilde{R}_X$ defined by
\[
R_{\delta, q}\coloneqq \left\{\sum_{i \in \Z} a_i(T^\flat)^i+ \sum_{j=1}^{\infty}\sum_{k=0}^{h-1} b_{jk}\frac{(T^{\flat})^k}{\overline{q}(T^{\flat})^j} ~\bigg|~ \begin{aligned} &
a_i, b_{jk} \in \kappa(\!(\overline{\pi}_K)\!), ~~|b_{jk}| \rightarrow 0 \text{ as } j \rightarrow +\infty,\\
&|a_i|\sigma^i \rightarrow 0 \text{ as } i \rightarrow \pm \infty \text{ for all } \delta \le \sigma \le 1 \end{aligned} \right\}.
\]

\item For each $0 < r < r_K$, let $\fA_{\delta, q}^r$ be a subring of $\widetilde{\fA}_X$ given by
\[
\fA_{\delta, q}^r\coloneqq \left\{\sum_{i\in \Z} a_i[T^{\flat}]^i+\sum_{j=1}^{\infty}\sum_{k=0}^{h-1} \frac{b_{jk} [T^{\flat}]^k}{q([T^{\flat}])^j} ~\bigg|~ \begin{aligned} & a_i, b_{jk} \in \fA_K^r,~~\lambda_r(b_{jk}) \rightarrow 0 \text{ as } j \rightarrow +\infty,\\ &\lambda_r(a_i)\sigma^{ri} \rightarrow 0 \text{ as } i \rightarrow \pm \infty \text{ for all } \delta \le \sigma \le 1\end{aligned}  \right\}.
\]

\item Let 
\[
\fA_{\delta, q}^{\dagger} \coloneqq \varinjlim_{0 < r < r_K}  \fA_{\delta, q}^r,
\]
and let $\fA_{\delta, q}$ be the $p$-adic completion of $\fA_{\delta, q}^{\dagger}$. Note that $\fA_{\delta, q}$ is a subring of $\widetilde{\fA}_X$.
\end{itemize}    
\end{definition}

\noindent \textbf{Step 1}. We first show that $\Lambda_{\delta, q}\coloneqq (R_{\delta, q}, \fA_{\delta, q})$ is weakly decompleting.\\
We will follow a similar argument as in Step 1 of \S \ref{subsec: closed-annuli}. For any $0 < r < r_K$, consider the norm $\lambda'_r$ on $\fA_{\delta, q}^r$ given by
\[
\lambda'_r\left(\sum_{i\in \Z} a_i[T^{\flat}]^i+\sum_{j=1}^{\infty}\sum_{k=0}^{h-1}b_{jk}\frac{[T^{\flat}]^k}{q([T^{\flat}])^j}\right)\coloneqq \max\left\{ \sup_{\delta \le \sigma \le 1} \sup_{i\in \Z} \{\lambda_r(a_i)\sigma^{ri}\},\,\, \sup_{j \ge 1}\sup_{0 \le k \le h-1}\{\lambda_r(b_{jk})\} \right\}. 
\]

\begin{lemma} \label{lem: lambda-equals-lambda'-puncture}
For each $0 < r < r_K$, we have $\lambda'_r = \lambda_r$ on $\fA_{\delta, q}^r$.    
\end{lemma}

\begin{proof}
We have $\lambda'_r \ge \lambda_r$ by the definition. Let $s\colon \mathbf{E}_K \rightarrow \fA_K$ be the section of the mod $p$ reduction map in Construction~\ref{construction:section_s_from_E_K}, which satisfies the properties in Lemma~\ref{lemma:section}. For each 
\[
y = \sum_{i \in \Z} a_i(T^\flat)^i+ \sum_{j=1}^{\infty}\sum_{k=0}^{h-1} b_{jk}\frac{(T^{\flat})^k}{\overline{q}(T^{\flat})^j} \in R_{\delta, q},
\]
set
\[
s(y)\coloneqq \sum_{i \in \Z} s(a_i)[T^\flat]^i+ \sum_{j=1}^{\infty}\sum_{k=0}^{h-1} s(b_{jk})\frac{[T^{\flat}]^k}{q([T^{\flat}])^j}.
\]
Then for any $0 < r < r_K$, we have
\[
\lambda_r(s(a_i))\sigma^{ri} = |a_i|^r \sigma^{ri} \rightarrow 0
\]
as $i \rightarrow \pm \infty$ for all $\delta \le \sigma \le 1$, and
\[
\lambda_r(s(b_{jk})) = |b_{jk}|^r \rightarrow 0
\]
as $j \rightarrow +\infty$. So $s(y) \in \fA_{\delta, q}^r$. This yields a section $s\colon R_{\delta, q} \rightarrow \fA_{\delta, q}^r$ of the mod $p$ reduction map. In particular, $\fA_{\delta, q}^r \rightarrow R_{\delta, q}$ is surjective. Note that
\begin{align*}
\lambda'_r(s(y)) &= \max\left\{ \sup_{\delta \le \sigma \le 1} \sup_{i\in \Z} \{\lambda_r(s(a_i))\sigma^{ri}\},\,\, \sup_{j \ge 1}\sup_{0 \le k \le h-1}\{\lambda_r(s(b_{jk}))\} \right\} \\
    & = \max\left\{ \sup_{\delta \le \sigma \le 1} \sup_{i\in \Z} \{|a_i|^r\sigma^{ri}\}, \,\,\sup_{j \ge 1}\sup_{0 \le k \le h-1}\{|b_{jk}|^r\} \right\} = |y|^r.
\end{align*}
Since $\lambda'_r \ge \lambda_r$, this implies $\lambda'_r(s(y)) = \lambda_r(s(y)) = |y|^r$.

Now for each $x \in \fA_{\delta, q}$, we construct a sequence $x_0, x_1, \ldots\in \fA_{\delta, q}$ inductively by setting $x_0=x$ and $x_{n+1}=\frac{1}{p}(x_n-s(\overline{x}_n))$, where $\overline{x}_n\in R_{\delta, q}$ stands for $x_n$ mod $p$. Then we have $x = \sum_{n=0}^{\infty}p^n \cdot s(\overline{x}_n)$, and the statements analogous to Lemmas~\ref{lemma: 3.5 weakly decompleting} and \ref{lemma: 3.6 weakly decompleting} hold. Consequently, we can proceed in the same way as in the proof of Lemma~\ref{lem: closed-annuli-lambda'-equals-lambda} to conclude.
\end{proof}

\begin{corollary} \label{cor: puncture-A^r-C} \hfill
\begin{enumerate}
\item For each $0 < r < r_K$, $\fA_{\delta, q}^r$ is complete with respect to $\lambda_r$, and
\[
\fA_{\delta, q}^r = \widetilde{\fA}_X^r \cap \fA_{\delta, q}
\]
as subrings of $\widetilde{\fA}_X$.

\item For each $0 < s \le r < r_K$, define $\fB_{\delta, q}^r$ and $\fC_{\delta, q}^{[s, r]}$ as in Definition~\ref{def: imperfect-period-rings-B-C}. Then $\fC_{\delta, q}^{[s, r]}$ consists of elements of the form
\[
\sum_{i\in \Z} a_i[T^{\flat}]^i+\sum_{j=1}^{\infty}\sum_{k=0}^{h-1}b_{jk}\frac{[T^{\flat}]^k}{q([T^{\flat}])^j}, \quad a_i, b_{jk} \in \fC_K^{[s, r]}
\]
satisfying the convergence conditions $\lambda_t(a_i)\sigma^{ti} \rightarrow 0$ as $i \rightarrow \pm \infty$ for all $t \in [s, r]$ and $\delta\le \sigma \le 1$, and $\lambda_t(b_{jk}) \rightarrow 0$ as $j \rightarrow +\infty$ for all $t \in [s, r]$.
\end{enumerate} 
\end{corollary}

\begin{proof}
Given Lemma~\ref{lem: lambda-equals-lambda'-puncture}, this follows from a similar argument as in the proof of Corollary~\ref{cor: closed-annuli-A^r-C}.    
\end{proof}

\begin{proposition} \label{prop: punctured-closed-annuli-strongly-decompleting}
$\Lambda_{\delta, q} = (R_{\delta, q}, \fA_{\delta, q})$ is a weakly decompleting datum for $(X, h)$.    
\end{proposition}

\begin{proof}
Note first that $R_{\delta, q}$ has a $p$-basis given by $\{\overline{\pi}^{\frac{i}{p}}(T^{\flat})^{\frac{j}{p}} ~|~ 0 \le i, j \le p-1\}$: it follows from direct computation that these elements generate $\varphi^{-1}(R_{\delta, q})$ as an $R_{\delta, q}$-module. Furthermore, these elements do not have any non-trivial linear relation over $R_{\delta, q}$ since they provide a $p$-basis of $\kappa(\!(\overline{\pi}_K)\!)[\![(T^{\flat})^{\pm 1}]\!]$ by the proof of \cite[Lemma~1.1.3]{deJong-dieudonnemodule} and $R_{\delta, q} \subset \kappa(\!(\overline{\pi}_K)\!)[\![(T^{\flat})^{\pm 1}]\!]$. 

Now, given Corollary~\ref{cor: puncture-A^r-C}, the statement follows from the same argument as in the proof of Proposition~\ref{prop: closed-annuli-weakly-decompleting}.
\end{proof}

We also consider finite \'etale extensions of $X$.

\begin{definition}
Let $g\colon X' \rightarrow X$ be a finite \'etale morphism and let $h'=h\circ g:X'\rightarrow \D_K$. Let $\Lambda_{\delta, q}'\coloneqq (R_{\delta, q}', \fA_{\delta, q}')$ be the pair obtained by applying Construction~\ref{construction:wdc_finite_etale} to $\Lambda_{\delta, q} = (R_{\delta, q}, \fA_{\delta, q})$ and $g\colon X' \rightarrow X$. 
\end{definition}

\noindent \textbf{Step 2}. We now prove that $\Lambda_{\delta, q} = (R_{\delta, q}, \fA_{\delta, q})$ is decompleting.\\
We will follow a similar argument as in Step 2 of \S \ref{subsec: closed-annuli}. Since $R_{\delta, q}'$ is finite \'etale over $R_{\delta, q}$, any $p$-basis of $R_{\delta, q}$ is also a $p$-basis of $R_{\delta, q}'$ by \cite[\S 1.1.2(ii)]{Berthelot-Messing}. So we have an isomorphism
\[
\varphi^{-1}(R_{\delta, q}') / R_{\delta, q}' \cong \bigoplus_{0 \le i, j \le p-1, ~(i, j)\neq (0, 0)} (1+\overline{\pi})^{\frac{i}{p}} (T^{\flat})^{\frac{j}{p}}\cdot R_{\delta, q}',
\]
which is strict by Theorem~\ref{thm: open mapping thm} since $\pi \in R_{\delta, q}'$. Write $\overline{y}_{i, j}\coloneqq (1+\overline{\pi})^{\frac{i}{p}}(T^{\flat})^{\frac{j}{p}}$.

\begin{proposition}
$\Lambda_{\delta, q} = (R_{\delta, q}, \fA_{\delta, q})$ is a decompleting datum for $(X, h)$.    
\end{proposition}

\begin{proof}
We follow a similar argument as in the proof of Proposition~\ref{prop: closed-annuli-decompleting}. Let $\Gamma_0 \subset \Gamma$ be an open normal subgroup. It suffices to show that the cochain complex
\[
\mathcal{C}^{(\bullet)}(\Gamma_0, \varphi^{-1}(R_{\delta, q}')/R_{\delta, q}')
\]
is strict exact at $d^{(\ell)}$ for each $\ell \ge 0$. 

We first show the exactness. By Theorem~\ref{thm: hoschild-serre-p-analy-gp}, it suffices to consider the case where $\Gamma_0$ is of the form
\[
\Gamma_0 \cong p^a\Z_p \rtimes (1+p^b \Z_p) \subset \Z_p\rtimes \Z_p^{\times}
\]
for some integers $a \ge 0$ and $b \ge 1$. Let $\gamma, \lambda \in \Gamma_0$ as in the proof of Proposition~\ref{prop: closed-annuli-decompleting}, so that
\[
\gamma(1+\overline{\pi}) = 1+\overline{\pi}, \quad \gamma(T^{\flat}) = (1+\overline{\pi})^{p^m}\cdot T^{\flat} 
\]
and
\[
\lambda(1+\overline{\pi}) = (1+\overline{\pi})^{p^{m'}}(1+\overline{\pi}), \quad \lambda(T^{\flat}) = T^{\flat}
\]
for some $m, m' \ge 1$. Applying a similar argument as in the proof of Proposition~\ref{prop: closed-annuli-decompleting} (based on Theorem~\ref{thm: hoschild-serre-p-analy-gp}), we reduce to showing that there is a constant $c > 0$ satisfying
\begin{equation} \label{eq: puncture-ineq-Gamma-action}
|(\gamma^{p^n}-1)\cdot \overline{x}| \le c\cdot p^{-\frac{p^{n+m+1}}{p-1}}|\overline{x}|,  \quad |(\lambda^{p^n}-1)\cdot \overline{x}| \le c\cdot p^{-\frac{p^{n+m'+1}}{p-1}}|\overline{x}|    
\end{equation}
for all $\overline{x} \in R_{\delta, q}$ and all integers $n \ge 1$. Recall that
\[
R_{\delta, q}\coloneqq \left\{\sum_{i \in \Z} a_i(T^\flat)^i+ \sum_{j=1}^{\infty}\sum_{k=0}^{h-1} b_{jk}\frac{(T^{\flat})^k}{\overline{q}(T^{\flat})^j} ~\bigg|~ \begin{aligned} &
a_i, b_{jk} \in \kappa(\!(\overline{\pi}_K)\!), ~~|b_{jk}| \rightarrow 0 \text{ as } j \rightarrow +\infty,\\
&|a_i|\sigma^i \rightarrow 0 \text{ as } i \rightarrow \pm \infty \text{ for all } \delta \le \sigma \le 1 \end{aligned} \right\}.
\]
and for $\overline{x} = \sum_{i \in \Z} a_i(T^\flat)^i+ \sum_{j=1}^{\infty}\sum_{k=0}^{h-1} b_{jk}\frac{(T^{\flat})^k}{\overline{q}(T^{\flat})^j} \in R_{\delta, q}$, we have
\[
|\overline{x}| = \max\left\{ \sup_{\delta \le \sigma \le 1} \sup_{i\in \Z} \{|a_i|\sigma^i\}, \,\, \sup_{j \ge 1}\sup_{0 \le k \le h-1}\{|b_{jk}|\} \right\}.
\]
Thus, it suffices to check the inequalities~\eqref{eq: puncture-ineq-Gamma-action} for $a_i(T^\flat)^i$ and $b_{jk}\frac{(T^{\flat})^k}{\overline{q}(T^{\flat})^j}$. The inequalities for $a_i(T^\flat)^i$ is already checked in the proof of Proposition~\ref{prop: closed-annuli-decompleting}. 

Note that for any $g \in \Gamma$, we have
\[
(g-1)(\overline{q}(T^{\flat})^{-1}) = -(g-1)\overline{q}(T^{\flat})\cdot g(\overline{q}(T^{\flat})^{-1})\cdot \overline{q}(T^{\flat})^{-1},
\]
so
\[
|(g-1)(\overline{q}(T^{\flat})^{-1})| \le |(g-1)\overline{q}(T^{\flat})|.
\]
Using the identity
\[
(g-1)(\overline{x}\cdot \overline{y}) = (g-1)(\overline{x})\cdot \overline{y}+g(\overline{x})\cdot (g-1)\overline{y}
\]
for any $\overline{x}, \overline{y} \in R_{\delta, q}$, we deduce that the inequalities \eqref{eq: puncture-ineq-Gamma-action} hold for $b_{jk}\frac{(T^{\flat})^k}{\overline{q}(T^{\flat})^j}$ as well. 

The strictness at $d^{(\ell)}$ of the complex $\mathcal{C}^{(\bullet)}(\Gamma_0, \varphi^{-1}(R_{\delta, q}')/R_{\delta, q}')$ also follows from a similar argument as in the proof of Proposition~\ref{prop: closed-annuli-decompleting}, by considering $R_{\delta, q}'\langle \frac{\overline{U}}{\overline{\pi}}, \frac{\overline{\pi}}{\overline{U}} \rangle$.   
\end{proof}

\noindent \textbf{Step 3}. Lastly, we show that $\Lambda_{\delta, q} = (R_{\delta, q}, \fA_{\delta, q})$ is strongly decompleting.

\begin{proposition}
$\Lambda_{\delta, q} = (R_{\delta, q}, \fA_{\delta, q})$ is a strongly decompleting datum for $(X, h)$.    
\end{proposition}

\begin{proof}
We check the conditions (a)-(c) of Definition~\ref{defn: strongly decompleting data} for $\Lambda_{\delta, q} = (R_{\delta, q}, \fA_{\delta, q})$. Recall that
\[
\fA_K^r = \left\{\sum_{i\in \Z} c_i\pi_K^i ~\Big|~ c_i\in \mathcal{O}_{K_0'}, ~~|c_i|\cdot p^{-\frac{pr}{p-1}\cdot \frac{i}{e_K}} \rightarrow 0 \text{ as } i \rightarrow -\infty \right\}
\]
and
\[
\fA_{\delta, q}^r\coloneqq \left\{\sum_{i\in \Z} a_i[T^{\flat}]^i+\sum_{j=1}^{\infty}\sum_{k=0}^{h-1} \frac{b_{jk} [T^{\flat}]^k}{q([T^{\flat}])^j} ~\bigg|~ \begin{aligned} & a_i, b_{jk} \in \fA_K^r,~~\lambda_r(b_{jk}) \rightarrow 0 \text{ as } j \rightarrow +\infty,\\ &\lambda_r(a_i)\sigma^{ri} \rightarrow 0 \text{ as } i \rightarrow \pm \infty \text{ for all } \delta \le \sigma \le 1\end{aligned}  \right\}.
\]

\noindent \textbf{Condition (a)}. We need to check that the image of the composite
\[
\mathbf{A}^r_{\delta, q} \subset \widetilde{\mathbf{A}}^r_X\xrightarrow[]{\iota_n} \mathbb{B}_{\mathrm{dR}, X}^+(X_{\infty})
\]
lies in $A_{X, n}[\![t,u]\!]^{\nabla=0}$ where $u\coloneqq T-[T^{\flat}]$, for all $n\ge -\log_p r$. As in the proof of Proposition~\ref{prop: closed-annuli-strongly-decompleting}, write $\iota_n([T^{\flat}]) = T^{p^{-n}}+u'$. By \textit{loc. cit.}, it suffices to show that
\[
\iota_n(q([T^{\flat}])^{-j}) \in A_{X, n}[\![t, u']\!] ~~(\subset A_{X, n}[\![t, u]\!])
\]
for each $j \ge 1$, and that the norms of its coefficients in $A_{X,n}$ are bounded uniformly in $j$.

Let $q_n(T) \in \mathcal{O}_{K_0}[T]$ be the monic polynomial whose coefficients are given by applying $\varphi^{-n}$ on the coefficients of $q(T)$. Then
\begin{equation} \label{eq: q_n-q}
(q_n(T^{p^{-n}}))^{p^n} \equiv q(T) \mod p\mathcal{O}_{K_0}[T^{p^{-n}}].    
\end{equation}
Thus, by considering the Taylor expansion, we have
\begin{equation} \label{eq: iota_n-q}
\iota_n(q([T^{\flat}])^{-j}) = q_n(T^{p^{-n}}+u')^{-j} \in A_{X, n}[\![u']\!]    
\end{equation}
whose coefficients in $A_{X, n}$ have norm at most $1$ for any $j$. This shows condition (a).

\noindent \textbf{Condition (b)}. We need to check that the composite
\[
\fA_{\delta, q}^r \stackrel{\iota_n}{\rightarrow} A_{X, n}[\![t, u]\!]^{\nabla = 0} \rightarrow A_{X, n}
\]
is surjective for all $n \ge -\log_p r$.\\ 

Note that for the closed annulus $X_{[\delta, 1]}$ given by $\delta \le \sigma \le 1$, the composite
\[
\fA_{[\delta, 1]}^r \stackrel{\iota_n}{\rightarrow} A_{X_{[\delta, 1]}, n}[\![t, u]\!]^{\nabla = 0} \rightarrow A_{X_{[\delta, 1]}, n}
\]
is surjective for any $n \ge -\log_p r$ by Proposition~\ref{prop: closed-annuli-strongly-decompleting}. So condition (b) follows from the above computations \eqref{eq: q_n-q} and \eqref{eq: iota_n-q}.\\

\noindent \textbf{Condition (c)}. Let $r_n = p^{-n}$. It suffices to show the claim that there is an integer $n_0 \ge 0$ such that for any $n \ge n_0$, we have
\[
\lambda_{r_n}((\gamma-1)(x)) \le p^{-1}\cdot \lambda_{r_n}(x)
\]
for all $x \in \fC_{\delta, q}^{[r_n, r_n]}$ and $\gamma \in \Gamma_n$.

Take $n_0$ to be an integer $n_1$ as in Lemma~\ref{lemma: estimation Gamma on C}. Using the identity
\[
(\gamma-1)(x\cdot y) = (\gamma-1)(x)\cdot y+\gamma(x)\cdot (\gamma-1)y
\]
for any $x, y \in \fC_{\delta, q}^{[r_n, r_n]}$ and $\gamma \in \Gamma$, we deduce from the proof of Proposition~\ref{prop: closed-annuli-strongly-decompleting} that it suffices to show
\[
\lambda_{r_n}((\gamma-1)(q([T^{\flat}]))^{-j}) \le p^{-\frac{p+1}{p}}
\]
for any $n \ge n_0$, $\gamma \in \Gamma_n$, and $j \ge 1$. Note that
\[
(\gamma-1)(q([T^{\flat}]))^{-j}) = -(\gamma-1)(q([T^{\flat}]))^j)\cdot \gamma(q([T^{\flat}]))^{-j})\cdot q([T^{\flat}]))^{-j}.
\]
Thus, we have
\[
\lambda_{r_n}((\gamma-1)(q([T^{\flat}]))^{-j}) \le \lambda_{r_n}((\gamma-1)(q([T^{\flat}]))^j)) \le p^{-\frac{p+1}{p}}\lambda_{r_n}(q([T^{\flat}]))^j)) = p^{-\frac{p+1}{p}}
\]
where the second inequality follows from the proof of Proposition~\ref{prop: closed-annuli-strongly-decompleting}. This concludes the proof. 
\end{proof}

\subsection{Normalized differential operators for closed annuli (with punctures)} \label{subsec: normalized-differential-operators-closed-annuli}

Keep the notation as in \S \ref{subsec: closed-annuli} and \S \ref{subsec: punctured-closed-annuli}. We explicitly compute the normalized differential operators $d_0$ and $d_1$ on $\mathbf{C}^{\dagger}_{[\sigma_1, \sigma_2]}$ and $\fC_{\delta, q}^{\dagger}$ defined in Definition~\ref{defn: normalized differential operators}. Again, although we only compute in the one-dimensional situation, the general case for poly-annuli (with punctures) is similar.

\begin{lemma} \label{lem: normalized-differential-operators}
On $\mathbf{C}^{\dagger}_{[\sigma_1, \sigma_2]}$ and $\fC_{\delta, q}^{\dagger}$, we have
\[\nabla_0=t\cdot (1+\pi)\cdot \frac{d}{d\pi},\qquad \nabla_1=t\cdot [T^{\flat}]\cdot \frac{d}{d[T^{\flat}]}.\]
Consequently, we have $d_0=\frac{d}{d\pi_K}$ and $d_1=\frac{d}{d[T^{\flat}]}$.
\end{lemma}

\begin{proof}
We first consider the case of $\mathbf{C}^{\dagger}_{[\sigma_1, \sigma_2]}$. Recall that elements in $\mathbf{C}^{\dagger}_{[\sigma_1, \sigma_2]}$ are two-variable formal power series subject to certain convergence conditions. Since $\Gamma_K$ acts trivially on $[T^{\flat}]$, the computation of $\nabla_0$ follows from \cite[Lemme~4.2]{Berger-differential}. To compute $\nabla_1$, since $\Gamma_{\mathrm{geom}}$ acts trivially on $\mathbf{C}^{\dagger}_K$, it suffices to show $\nabla_1([T^{\flat}]^n)=nt\cdot [T^{\flat}]^n$ for all $n\in \Z$. Consider the topological generator $\tau=\tau_1$ of $\Gamma_{\mathrm{geom}}$ as in the proof of Proposition~\ref{prop:image_of_iota_n_fixed_by_Gamma_n} so that $\tau([T^{\flat}])=[\varepsilon]\cdot[T^{\flat}]$. A direct computation yields (for any integer $m\ge 0$)
\begin{align*}
\nabla_1([T^{\flat}]^n)=&\: p^{-m}\cdot \log(\tau^{p^m})([T^{\flat}]^n) = -p^{-m}\cdot \sum_{k=1}^{\infty}\frac{1}{k}(1-\tau^{p^m})^k([T^{\flat}]^n)\\
=&\: -p^{-m}\cdot \left(\sum_{k=1}^{\infty}\frac{1}{k}\left(1-[\varepsilon]^{p^mn}\right)^k\right)\cdot [T^{\flat}]^n=p^{-m}\cdot \log([\varepsilon]^{p^mn})\cdot [T^{\flat}]^n =nt\cdot [T^{\flat}]^n.
\end{align*}

For $\fC_{\delta, q}^{\dagger}$, we need to compute the effect of $\nabla_1$ on terms of the form $\frac{[T^{\flat}]^k}{q([T^{\flat}])^j}$. But this follows from previous computation and the Leibniz rule. 
\end{proof}

\subsection{Compatibility of $\mathbf{N}_{\mathrm{dR}}$ for closed annuli (with punctures) under base field extensions} \label{subsec: compatibility-NdR-base-field-extension}

Keep the notation as in \S \ref{subsec: closed-annuli} and \S \ref{subsec: punctured-closed-annuli}. In both the cases of closed annuli and closed annuli with punctures, we study the compatibility of the $\fN_{\mathrm{dR}}$ functor (see Construction~\ref{const: NdR}) under base field extensions.

Let $L/K$ be a finite field extension and let $0 < \sigma_1 \le \sigma_2 \le 1$. Consider the closed annulus
\[
X_K \coloneqq X_{K, [\sigma_1, \sigma_2]} = A_K([\sigma_1, \sigma_2]; T)
\]
with the toric chart $h_K\colon X_K \hookrightarrow \D_K$ sending $T \mapsto T$, and consider its base change along $K \rightarrow L$
\[
X_L\coloneqq X_{L, [\sigma_1, \sigma_2]} = A_L([\sigma_1, \sigma_2]; T),
\]
with the toric chart $h_L\colon X_L \hookrightarrow \D_L$. Let 
\[\Gamma_{(L)} \cong \Gamma_{\mathrm{geom}}\rtimes \Gamma_L\] be the Galois group for the toric tower of $\D_L$.  

By Proposition~\ref{prop: closed-annuli-strongly-decompleting}, we have a strongly decompleting datum $\Lambda_{K, [\sigma_1, \sigma_2]} = (R_{K, [\sigma_1, \sigma_2]}, \fA_{K, [\sigma_1, \sigma_2]})$ (resp. $\Lambda_{L, [\sigma_1, \sigma_2]} = (R_{L, [\sigma_1, \sigma_2]}, \fA_{L, [\sigma_1, \sigma_2]})$) for $(X_{K, [\sigma_1, \sigma_2]}, h_K)$ (resp. for $(X_{L, [\sigma_1, \sigma_2]}, h_L)$) given by Definition~\ref{def: closed-annuli-Lambda}. By construction, there is a natural $(\varphi, \Gamma_{(L)})$-equivariant injection
\[
\fA_{K, [\sigma_1, \sigma_2]}^{r} \hookrightarrow \fA_{L, [\sigma_1, \sigma_2]}^{r}
\]
for all $r > 0$ sufficiently small, which induces a natural $(\varphi, \Gamma_{(L)})$-equivariant injection
\[
\fC_{K, [\sigma_1, \sigma_2]}^{[s,r]} \hookrightarrow \fC_{L, [\sigma_1, \sigma_2]}^{[s,r]}
\]
for $0<s\le r$. 

\begin{proposition} \label{prop: closed-annuli-NdR-base-field-change-compatibility}
Let $\L \in \mathrm{Loc}^{\mathrm{dR}}_{\Z_p}(X_K)$ be a de Rham local system on $X_{K, [\sigma_1, \sigma_2]}$. Then there is natural isomorphism
\[
\mathbf{N}_{\mathrm{dR}}(\L)\otimes_{\fC_{K, [\sigma_1, \sigma_2]}^{\dagger}} \fC_{L, [\sigma_1, \sigma_2]}^{\dagger} \isom  \mathbf{N}_{\mathrm{dR}}(\L |_{X_L})
\]
of $(\varphi, \Gamma_{(L)})$-modules.
\end{proposition}

\begin{proof}
This follows from a similar argument as in the proof of Proposition~\ref{prop:compatibility_with_L/K} and \ref{prop: compatibility-NdR}. For completeness, we explain some details. Let $N_{K, [\sigma_1, \sigma_2]}^{[s, r]}$ (resp. $N_{L, [\sigma_1, \sigma_2]}^{[s, r]}$) denote the finite projective $\fC_{K, [\sigma_1, \sigma_2]}^{[s, r]}$-module (resp. $\fC_{L, [\sigma_1, \sigma_2]}^{[s, r]}$-module) in Construction~\ref{const: NdR} obtained from $\L$ (resp. from $\L |_{X_L}$). It suffices to produce compatible $(\varphi, \Gamma_{(L)})$-equivariant isomorphisms
\[
N_{K, [\sigma_1, \sigma_2]}^{[s, r]} \otimes_{\fC_{K, [\sigma_1, \sigma_2]}^{[s, r]}} \fC_{L, [\sigma_1, \sigma_2]}^{[s, r]} \isom N_{L, [\sigma_1, \sigma_2]}^{[s, r]}
\]
for $0 < s \le r$ with $r$ sufficiently small. By construction, this amounts to constructing two isomorphisms after taking base changes along $\fC_{L, [\sigma_1, \sigma_2]}^{[s, r]} \rightarrow \fC_{L, [\sigma_1, \sigma_2]}^{[s, r]}[\frac{1}{t}]$ and $\iota_n\colon \fC_{L, [\sigma_1, \sigma_2]}^{[s, r]} \rightarrow A_{X_L, n}[\![t, u]\!]^{\nabla = 0}$ respectively, such that the two isomorphisms agree over $A_{X_L, n}[\![t, u]\!]^{\nabla = 0}[\frac{1}{t}]$. Here, $u = T - [T^{\flat}]$. 

By construction, $\widetilde{D}^{\mathrm{int}, \dagger}$ in Corollary~\ref{cor: equiv-local-system-overconvergent-etale-(phi, Gamma)-mods} is compatible with extensions of base field. Therefore, by Corollary~\ref{cor: decompletion-etale-(phi, Gamma)-modules}, we have
\[
M_{K, [\sigma_1, \sigma_2]}^{[s, r]} \otimes_{\fC_{K, [\sigma_1, \sigma_2]}^{[s, r]}} \fC_{L, [\sigma_1, \sigma_2]}^{[s, r]} \cong M_{L, [\sigma_1, \sigma_2]}^{[s, r]}
\]
for $r$ sufficiently small. Inverting $t$, we obtain
\[
N_{K, [\sigma_1, \sigma_2]}^{[s, r]} \otimes_{\fC_{K, [\sigma_1, \sigma_2]}^{[s, r]}} \fC_{L, [\sigma_1, \sigma_2]}^{[s, r]}[\frac{1}{t}] \cong M_{K, [\sigma_1, \sigma_2]}^{[s, r]} \otimes_{\fC_{K, [\sigma_1, \sigma_2]}^{[s, r]}} \fC_{L, [\sigma_1, \sigma_2]}^{[s, r]}[\frac{1}{t}] \cong M_{L, [\sigma_1, \sigma_2]}^{[s, r]}[\frac{1}{t}] \cong N_{L, [\sigma_1, \sigma_2]}^{[s, r]}[\frac{1}{t}],
\]
for $r$ sufficiently small. On the other hand, there is a natural isomorphism
\[
D_{\mathrm{dR}}(\L)\otimes_{A_{X_K}} A_{X_L} \cong D_{\mathrm{dR}}(\L |_{X_L}) 
\]
by the base change property of the $D_{\mathrm{dR}}$ functor. By a similar argument as in the proof of Proposition~\ref{prop: compatibility-NdR}, we obtain a natural isomorphism
\[
(D_{\mathrm{dR}}(\L)\otimes_{A_{X_K}} A_{X_K, n}[\![t, u]\!])^{\nabla = 0} \otimes_{A_{X_K, n}[\![t, u]\!]^{\nabla = 0}} A_{X_L, n}[\![t, u]\!]^{\nabla = 0} \cong (D_{\mathrm{dR}}(\L |_{X_L})\otimes_{A_{X_L}} A_{X_L, n}[\![t, u]\!])^{\nabla = 0}.
\]

It remains to show that the two isomorphisms agree over $A_{X_L, n}[\![t, u]\!]^{\nabla = 0}[\frac{1}{t}]$. Let $\iota_{n, K}$ (resp. $\iota_{n, L}$) be the isomorphism \eqref{eq:gluing_isom_relative_s,r} for $X_K$ (resp. for $X_L$). We need to prove that the following diagram commutes
\[
\begin{tikzcd}
    M_{K, [\sigma_1, \sigma_2]}^{[s,r]} \otimes_{\mathbf{C}_{K, [\sigma_1, \sigma_2]}^{[s,r]}, \iota_n } A_{X_K, n}[\![t, u]\!]^{\nabla = 0}[\frac{1}{t}] \arrow[r, "\iota_{n, K}"] \arrow[d]
    & (D_{\mathrm{dR}}(\L)\otimes_{A_{X_K}} A_{X_K, n}[\![t, u]\!])^{\nabla = 0}[\frac{1}{t}] \arrow[d] \\
    M_{L, [\sigma_1, \sigma_2]}^{[s,r]} \otimes_{\mathbf{C}_{L, [\sigma_1, \sigma_2]}^{[s,r]}, \iota_n } A_{X_L, n}[\![t, u]\!]^{\nabla = 0}[\frac{1}{t}]   \arrow[r, "\iota_{n, L}"]
    & (D_{\mathrm{dR}}(\L|_{X_L})\otimes_{A_{X_L}} A_{X_L, n}[\![t, u]\!])^{\nabla = 0}[\frac{1}{t}].
\end{tikzcd}
\] 
Consider $\overline{A}_X$ as in the proof of Lemma~\ref{lemma:gluing_isomorphism_perfect_relative} for $X = X_K$, and let
\[
\widetilde{\fA}\coloneqq W(\overline{A}_X^{\wedge, \flat}), \quad \mathbf{B}_{\mathrm{{dR}}, \overline{X}} \coloneqq \mathbb{B}_{\mathrm{dR}, X_K}(\overline{A}_X^{\wedge, \flat}).
\]
Denote by $H_{X_K}$ (resp. $H_{X_L}$) the Galois group $H_X$ in the proof of Lemma~\ref{lemma:gluing_isomorphism_perfect_relative} for $X = X_K$ (resp. for $X = X_L$). By the construction of $\widetilde{D}^{\mathrm{int}, \dagger}$ in Corollary~\ref{cor: equiv-local-system-overconvergent-etale-(phi, Gamma)-mods}, there is $r > 0$ sufficiently small such that
\[
\widetilde{D}^{\mathrm{int}, \dagger}(\L) = (\L\otimes_{\Z_p} \widetilde{\fA}^r)^{H_{X_K}}\otimes_{\fA_{X_K}^r} \fA_{X_K}^{\dagger}
\]
and that
\[
\widetilde{D}^{\mathrm{int}, \dagger}(\L|_{X_L}) = (\L\otimes_{\Z_p} \widetilde{\fA}^r)^{H_{X_L}}\otimes_{\fA_{X_L}^r} \fA_{X_L}^{\dagger}.
\]
Hence, by the construction of $\iota_{n, K}$ and $\iota_{n, L}$ in \S \ref{sec:relative_N_dR} and Lemma~\ref{lemma:gluing_isomorphism_perfect_relative}, the desired commutativity follows from the commutativity of
\[
\begin{tikzcd}
    (\L\otimes_{\Z_p} \widetilde{\fA}^r)^{H_{X_K}} \arrow[r, "\iota_n"] \arrow[d]
    & (\L\otimes_{\Z_p} \mathbf{B}_{\mathrm{{dR}}, \overline{X}})^{H_{X_K}} \arrow[d] \\
    (\L\otimes_{\Z_p} \widetilde{\fA}^r)^{H_{X_L}}    \arrow[r, "\iota_n"]
    & (\L\otimes_{\Z_p} \mathbf{B}_{\mathrm{{dR}}, \overline{X}})^{H_{X_L}},
\end{tikzcd}
\]
where $\iota_n$ in both top and bottom horizontal maps are given by $\iota_n\colon \widetilde{\fA}^r \rightarrow \mathbf{B}_{\mathrm{{dR}}, \overline{X}}$.

Lastly, by the construction of $(\varphi, \Gamma_{(L)})$-structure in Construction~\ref{const: NdR}, it follows that the base change isomorphism is compatible with $\varphi$-actions and $\Gamma_{(L)}$-actions.
\end{proof}

Similarly, let $0 < \delta \le 1$ with $\delta \in p^{\Q}$, and let $q(T) \in \mathcal{O}_{K_0}[T]$ be a monic polynomial of the form $q(T)=\prod_{i=1}^m (T-[\beta_i])$
with $\beta_i\in k_K^{\times}$. Consider the rational subset
\[
X_{K, \delta, q} \coloneqq X_{K, [\delta, 1]} \cap \{|q(T)| = 1\}\subset  X_{K, [\delta, 1]}
\]
with the toric chart $h_K\colon X_{K, \delta, q} \hookrightarrow \D_K$ sending $T \mapsto T$. Let $X_{L, \delta, q} = X_{K, \delta, q}\times_K L = X_{L, [\delta, 1]} \cap \{|q(T)| = 1\}$ with the toric chart $h_L\colon X_{L, \delta, q} \hookrightarrow \D_L$. By Proposition~\ref{prop: punctured-closed-annuli-strongly-decompleting}, we have a strongly decompleting datum $\Lambda_{K, \delta, q} = (R_{K, \delta, q}, \fA_{K, \delta, q})$ (resp. $\Lambda_{L, \delta, q} = (R_{L, \delta, q}, \fA_{L, \delta, q})$) for $(X_{K, \delta, q}, h_K)$ (resp. for $(X_{L, \delta, q}, h_L)$) given by Definition~\ref{def: punctured-closed-annuli-datum}.

\begin{proposition}
Let $\L \in \mathrm{Loc}^{\mathrm{dR}}_{\Z_p}(X_{K, \delta, q})$ be a de Rham local system on $X_{K, \delta, q}$. Then there is natural isomorphism
\[
\mathbf{N}_{\mathrm{dR}}(\L)\otimes_{\fC_{K, \delta, q}^{\dagger}} \fC_{L, \delta, q}^{\dagger} \isom  \mathbf{N}_{\mathrm{dR}}(\L |_{X_{L, \delta, q}})
\]
of $(\varphi, \Gamma_{(L)})$-modules.    
\end{proposition}

\begin{proof}
This follows from essentially the same argument as in the proof of Proposition~\ref{prop: closed-annuli-NdR-base-field-change-compatibility}.    
\end{proof}

\subsection{Relation with the work of Andreatta--Brinon} \label{subsec: relation-with-AB}

In this subsection, we briefly explain the relation between the notion of strongly decompleting data and the imperfect period rings studied in \cite{Andreatta-Brinon_Overconvergence} and \cite{Andreatta-Brinon}. Let $A_0$ be an integral domain obtained from $W(k)\langle T_1^{\pm 1}, \ldots, T_d^{\pm 1}\rangle$ by a finite number of iterations of the following operations:
\begin{itemize}
\item[(F\'Et)] the $p$-adic completion of a finite \'etale extension;
\item[(Loc)] the $p$-adic completion of a localization.
\end{itemize}
Consider the situation where $X$ has an integral model given by $A\coloneq A_0\otimes_{W(k)} \mathcal{O}_K$. In this case, there is a toric chart $h\colon X \rightarrow \D_K^d$ given by the natural map $\mathcal{O}_K\langle T_1, \ldots, T_d\rangle \rightarrow A$.  

\begin{remark}
In \cite{Andreatta-Brinon_Overconvergence} and \cite{Andreatta-Brinon}, a larger class of rings for $A_0$ is studied; the operations also include the $p$-adic completion of any \'etale (not necessarily finite) extension and the completion with respect to any ideal including $p$ (see \cite[\S 2.1]{Andreatta-Brinon} and \cite[Remark~2.4]{du-liu-moon-shimizu-completed-prismatic-F-crystal-loc-system}). To be compatible with our definition of toric charts and for simplicity of exposition, we restrict our attention to the setup above.
\end{remark}

Consider $R_{\Lambda} \subset \widetilde{R}_X$ given by the ring denoted as $\mathbf{E}_S$ in \cite[\S 4.1]{Andreatta-Brinon_Overconvergence}. By Lemma~4.1 and its proof in \textit{loc. cit.}, we have
$ 
R_{\Lambda} = \kappa(\!(\overline{\pi}_K)\!)\langle (T_1^{\flat})^{\pm 1}, \ldots, (T_d^{\flat})^{\pm 1} \rangle
$ 
when $A = \mathcal{O}_K\langle T_1^{\pm}, \ldots, T_d^{\pm}\rangle$, and in general, $R_{\Lambda}$ is obtained from $\kappa(\!(\overline{\pi}_K)\!)\langle (T_1^{\flat})^{\pm 1}, \ldots, (T_d^{\flat})^{\pm 1} \rangle$ by a finite number of operations of finite \'etale extensions and $\overline{\pi}$-adic completions of localizations. In particular, $R_{\Lambda}$ is an affinoid $\kappa(\!(\overline{\pi}_K)\!)$-algebra, and 
\[\{\overline{\pi}_K^{\frac{i_0}{p}} (T_1^{\flat})^{\frac{i_1}{p}} \cdots (T_d^{\flat})^{\frac{i_d}{p}}, ~~0 \le i_0, \ldots, i_d \le p-1 \}\] gives a $p$-basis of $R_{\Lambda}$.

Let $\fA_{\Lambda} \subset \widetilde{\fA}_X$ be the ring denoted by $\fA_S$ in \cite[\S 4.3]{Andreatta-Brinon_Overconvergence}. Note that $\widetilde{\fA}_X^r$ (resp. $\fA_{\Lambda}^r$) is denoted therein by $\widetilde{\fA}_X^{(0, r]}$ (resp. $\fA_S^{(0, r]}$). Then one can check using the results in \S 4 of \textit{loc. cit.} that $\Lambda = (R_{\Lambda}, \fA_{\Lambda})$ is a weakly decompleting datum for $(X,h)$; for example, $R_{\Lambda}^{\mathrm{perf}}$ is dense in $\widetilde{R}_X$ by \cite[Proposition~4.9(iii)]{Andreatta-Brinon_Overconvergence}, and the map $\fA_{\Lambda}^r/p \rightarrow R_{\Lambda}$ is an isomorphism for $r$ sufficiently small by \cite[Proposition~4.28(e)]{Andreatta-Brinon_Overconvergence}. Furthermore, using the explicit description of $\fA_{\Lambda}$ in the proof of Proposition~4.42 of \textit{loc. cit.} and a similar argument as in Step 1 of \S \ref{subsec: closed-annuli}, one can deduce that the map $\fA_{\Lambda}^r/p \rightarrow R_{\Lambda}$ is strict for $r$ sufficiently small. We leave the details to the reader.

Next, one can check that $\Lambda = (R_{\Lambda}, \fA_{\Lambda})$ is a decompleting datum by applying a similar argument as in Step 2 of \S \ref{subsec: closed-annuli}. Indeed, since $R_{\Lambda}$ has a $p$-basis 
\[\{(1+\overline{\pi})^{\frac{i_0}{p}} (T_1^{\flat})^{\frac{i_1}{p}} \cdots (T_d^{\flat})^{\frac{i_d}{p}}, ~~0 \le i_0, \ldots, i_d \le p-1 \},\] it reduces to showing the inequalities as in \eqref{ineq: gamma-action} and \eqref{ineq: sigma-action} for some constant $c > 0$. These hold when $R_{\Lambda} = \kappa(\!(\overline{\pi}_K)\!)\langle (T_1^{\flat})^{\pm 1}, \ldots, (T_d^{\flat})^{\pm 1} \rangle$ by direct computation, and are preserved under (F\'Et) (after enlarging $c$ if necessary). For (Loc), note that if $f \in R_{\Lambda}$ is invertible, we have
\[
(\gamma-1)(f^{-1}) = -(\gamma-1)f\cdot \gamma(f^{-1})\cdot f^{-1}
\]
for all $\gamma \in \Gamma$. Using the identity
\[
(\gamma-1)(\overline{x}\cdot \overline{y}) = (\gamma-1)(\overline{x})\cdot \overline{y}+\gamma(\overline{x})\cdot (\gamma-1)\overline{y}
\]
for any $\overline{x}, \overline{y} \in R_{\Lambda}$, one can deduce that the inequalities are also preserved under (Loc) (after enlarging $c$). Again, we leave the details to the reader.

Finally, one can check that the conditions (a)-(c) of Definition~\ref{defn: strongly decompleting data} hold for $\Lambda = (R_{\Lambda}, \fA_{\Lambda})$, so that it is a strongly decompleting datum. Indeed, Condition (a) holds by \cite[Proposition~8.6]{Andreatta-Brinon}. Condition (b) can be verified using the computation in the proof of \cite[Lemme~8.2]{Andreatta-Brinon}. When $A = \mathcal{O}_K\langle T_1^{\pm}, \ldots, T_d^{\pm}\rangle$, Condition (c) can be verified by a direct computation. In general, Condition (c) is stable under (F\'Et) by Theorem~\ref{thm: strongly decompleting under finite etale extensions}, and the stability under (Loc) also follows from a similar argument as in the proof of Theorem~\ref{thm: strongly decompleting under finite etale extensions}.

\vspace{0.1in}
\appendix
\section{Finitely generated modules over $\fC_{\Lambda}^r$} \label{section: finitely generated modules}

The goal of this appendix is to study finitely generated (resp. finite projective) modules over certain imperfect period rings associated to a weakly decompleting data and investigate their relations with coherent sheaves (resp. vector bundles) on the associated quasi-Stein spaces. This is used in the discussion of $\varphi$-modules over imperfect period rings of type $\mathbf{C}$ (for example, $\fC_{\Lambda}^r$) in \S \ref{sec: decompletion-relative-etale-phiGamma-modules}.

\subsection{Setup}

Let $h\colon X \ra \D_K^d$ be a toric chart and let $\Lambda = (R_{\Lambda}, \fA_{\Lambda})$ be a weakly decompleting datum for $(X, h)$ as define in \S \ref{sec:weakly_decompleting_data}. Let 
\[ \mathbf{A}^r_{\Lambda}, \mathbf{A}^{\dagger}_{\Lambda}, \mathbf{B}^r_{\Lambda}, \mathbf{B}^{\dagger}_{\Lambda}, \mathbf{C}^{[s,r]}_{\Lambda}, 
\mathbf{C}^r_{\Lambda}, \mathbf{C}^{\dagger}_{\Lambda}
\] 
be the associated imperfect period rings considered in \S \ref{subsection: weakly decompleting data}. 
Let $r_0$ be as in Proposition~\ref{prop: good lift weakly decompleting}. We may assume $r_0\in \mathbb{Q}$. Recall from Notation~\ref{notn: rs-rational} that in practice, we require that all superscripts in the period rings are in $\mathbb{Q}$. We also recall that  
$
\fC_{\Lambda}^r \isom \varprojlim_{0 <s \le r}\mathbf{C}_{\Lambda}^{[s,r]}
$ 
is an isomorphism of topological rings by Lemma~\ref{lem: C breve dense in C tilde}(1).

\subsection{Coherent sheaves on $\mathrm{Sp}\,\mathbf{C}_{\Lambda}^{[s,r]}$} 

We first show that, for  $0< s <s' \le r < r_0$ (resp. for $0 < s \le r' < r < r_0$), the inclusions $\mathbf{C}_{\Lambda}^{[s,r]}\rightarrow \fC_{\Lambda}^{[s',r]}$  (resp. $\mathbf{C}_{\Lambda}^{[s,r]}\rightarrow \fC_{\Lambda}^{[s,r']}$) are rational localizations of affinoid $K'_0$-algebras (see Conditions (7), (8) in Definition~\ref{defn: weakly decompleting datum} and Remark~\ref{rem: imperfect-period-rings-affinoid}). 

\begin{proposition} \label{prop: rational localizations C^r} \hfill
\begin{enumerate}
\item Let $0<s<s'\le r<r_0$ with $s,s',r\in \mathbb{Q}$. Choose a positive integer $d \ge 1$ such that $d'\coloneqq \frac{p}{p-1}\cdot \frac{1}{e_K}\cdot s'd\in \mathbb{Z}$. There is a natural isomorphism \footnote{Note that the rational localization $\mathbf{C}_{\Lambda}^{[s,r]}\left\langle \frac{\pi_K^d}{p^{d'}} \right\rangle$ is independent of the choice of $d$.}
\[
\fC_{\Lambda}^{[s',r]}\cong \mathbf{C}_{\Lambda}^{[s,r]}\left\langle \frac{\pi_K^d}{p^{d'}} \right\rangle
\]
identifying $\fC_{\Lambda}^{[s',r]}$ as a rational localization of $\mathbf{C}_{\Lambda}^{[s,r]}$, where the convergence condition is with respect to the $\mathrm{max}\{\lambda_r, \lambda_s\}$-norm on $\mathbf{C}_{\Lambda}^{[s,r]}$.

\item Let $0<s\le r'< r<r_0$ with $s,r',r\in \mathbb{Q}$. Choose a positive integer $d \ge 1$ such that $d'\coloneqq \frac{p}{p-1}\cdot \frac{1}{e_K}\cdot r'd\in \mathbb{Z}$. There is a natural isomorphism \footnote{Likewise, the rational localization $\mathbf{C}_{\Lambda}^{[s,r]}\left\langle \frac{p^{d'}}{\pi_K^d} \right\rangle$ is independent of the choice of $d$.}
\[
\fC_{\Lambda}^{[s,r']}\cong \mathbf{C}_{\Lambda}^{[s,r]}\left\langle \frac{p^{d'}}{\pi_K^d} \right\rangle
\]
identifying $\fC_{\Lambda}^{[s,r']}$ as a rational localization of $\mathbf{C}_{\Lambda}^{[s,r]}$, where the convergence is with respect to the $\mathrm{max}\{\lambda_r, \lambda_s\}$-norm on $\mathbf{C}_{\Lambda}^{[s,r]}$.
\end{enumerate}
   
\end{proposition}

\begin{proof}
We only prove (1), as the proof of (2) is similar.
By definition, $\mathbf{C}_{\Lambda}^{[s,r]}\left\langle \frac{\pi_K^d}{p^{d'}} \right\rangle$ consists of sums 
$
x=\sum_{i=0}^{\infty} x_i\left(\frac{\pi_K^d}{p^{d'}}\right)^i 
$ 
with $x_i\in \mathbf{C}_{\Lambda}^{[s,r]}$ such that, for each $t\in [s,r]$, $\lambda_t(x_i)\rightarrow 0$ as $i\rightarrow +\infty$. For such an $x$, note that $\lambda_t\left(\frac{\pi_K^d}{p^{d'}}\right)\le 1$ for all $t \in [s',r]$. So we have
\[
\lambda_t\left(x_i\left(\frac{\pi_K^d}{p^{d'}}\right)^i\right)\le \lambda_t(x_i)\rightarrow 0
\]
as $i\rightarrow +\infty$, for any $t\in [s', r]$. Thus, $x$ defines an element of $\fC_{\Lambda}^{[s',r]}$. This implies 
\[
\mathbf{C}_{\Lambda}^{[s,r]}\left\langle \frac{\pi_K^d}{p^{d'}} \right\rangle\subset \fC_{\Lambda}^{[s',r]}.
\]

To show the converse inclusion, we need to express each $x\in \fC_{\Lambda}^{[s',r]}$ as a sum $x=\sum_{i=0}^{\infty} x_i\left(\frac{\pi_K^d}{p^{d'}}\right)^i$ with $x_i\in \mathbf{C}_{\Lambda}^{[s,r]}$ such that, for any $t\in [s,r]$, $\lambda_t(x_i)\rightarrow 0$ as $i\rightarrow +\infty$.

By definition, we can write $x = \sum_{n=0}^{\infty} y_n$ with $y_n\in \mathbf{B}^r_{\Lambda}$ such that for some $0 < \varepsilon < 1$, we have $\lambda_t(y_n) < \varepsilon^n$ for all $n\ge 0$ and for all $t\in [s',r]$. The idea is to write each $y_n$ in terms of a norm compatible expansion (cf. Definition~\ref{defn: norm compatible expansion}). By Proposition~\ref{prop: good lift weakly decompleting}(2), for each $n\ge 0$, we can write $y_n= \sum_{m\in \mathbb{Z}} p^m y_n^{(m)}$ where
\begin{itemize}
\item $y_n^{(m)}\in \fA_{\Lambda}^{r_0}$ for all $m$;
\item for some $i_n\in \mathbb{Z}$, $y_n^{(m)}=0$ for all $m<i_n$;
\item $\lambda_{r_0}(y_n^{(m)}) = \left|\overline{y_n^{(m)}}\right|^{r_0}$ for all $m$ (by Remark~\ref{rmk: equivalent condition}, this is equivalent to 
\[
\lambda_t(y_n^{(m)})=\left|\overline{y_n^{(m)}}\right|^t
\]
for all $t\in (0,r_0]$);
\item $\lambda_t(y_n)=\mathrm{max}_{m\in\mathbb{Z}}\{p^{-m}\cdot \lambda_t(y_n^{(m)})\}$ for all $t\in (0, r]$.
\end{itemize}
By scaling $x$ by some $p$-power, we may assume $i_0=0$. We put 
$
z_0\coloneqq \sum_{n=0}^{\infty}\sum_{m=0}^{\infty}p^my_n^{(m)},
$ 
and for each $j\ge 1$, set 
$
z_j\coloneqq p^{-j}\sum_{n=0}^{\infty}y_n^{(-j)}.
$ \\

\noindent{\textbf{Claim 1:}} $z_j\in \mathbf{B}^r_{\Lambda}$ for all $j\ge 0$.\\

\noindent \textit{Proof of Claim 1}. For $z_0$, we have $\sum_{m=0}^{\infty} p^m y_n^{(m)}\in \mathbf{A}_{\Lambda}^r$ for all $n$ and 
\[
\lambda_r\left(\sum_{m=0}^{\infty}p^my_n^{(m)}\right)\le \lambda_r(y_n)\rightarrow 0
\]
as $n\rightarrow +\infty$. For $j\ge 1$, it suffices to show $\sum_{n=0}^{\infty} y_n^{(-j)} \in \mathbf{A}_{\Lambda}^r$. Again, we have $y_n^{(-j)}\in \mathbf{A}_{\Lambda}^r$ and $\lambda_r\left(y_n^{(-j)}\right)\le \lambda_r(y_n)\rightarrow 0$ as $n\rightarrow +\infty$. This proves Claim 1.\\

\noindent{\textbf{Claim 2:}} The sum $\sum_{j=0}^{\infty} z_j$ converges to $x$ with respect to $\mathrm{max}_{t\in [s',r]}\{\lambda_t\}$.\\

\noindent \textit{Proof of Claim 2}. 
Pick an unbounded increasing sequence $0 \le n_1\le n_2 \le \cdots$ of integers such that, for each $j\ge 1$, we have $i_0, i_1, \ldots, i_{n_j}>-j$.\footnote{Such a sequence exists. If the sequence $\{i_n\}$ is bounded below, say $i_n>-M$ for all $n$, then we simply pick $n_1 = n_2=\cdots = n_{M-1}=0$ and $n_j=j$ for all $j\ge M$. Otherwise, we have $i_n\rightarrow-\infty$ as $n\rightarrow +\infty$. Then we take $n_j$ to be the largest integer such that $i_0, i_1, \ldots, i_{n_j}>-j$.} Then 
\[
z_j = p^{-j} \sum_{n=0}^{\infty} y_n^{(-j)} = p^{-j} \sum_{n>n_j} y_n^{(-j)}.
\]
So for any $t\in [s',r]$, we have
\[
\lambda_t(z_j)\le \mathrm{max}_{n>n_j}\{\lambda_t(y_n)\}\rightarrow 0
\]
as $j\rightarrow +\infty$. This proves Claim 2. \\

It remains to show that each $z_j$ can be written as $z'_j\cdot \left(\frac{\pi_K^d}{p^{d'}}\right)^{\alpha_j}$ for some $\alpha_j\in \mathbb{Z}_{\ge 0}$ satisfying
\begin{itemize}
\item $\alpha_j\rightarrow +\infty$ as $j\rightarrow +\infty$;
\item $\lambda_t(z'_j)\rightarrow 0$ as $j\rightarrow +\infty$ for all $t\in [s,r]$.
\end{itemize}
To this end, let $n_1\le n_2\le \cdots$ be the unbounded increasing sequence of integers as in Claim 2, so that $z_j=\sum_{n>n_j} p^{-j}y_n^{(-j)}$. For each $j$, by the assumption of norm compatible expansions, we have
\[
\left|\overline{y_n^{(-j)}}\right|^t=\lambda_t\left(y_n^{(-j)}\right)\le p^{-j}\varepsilon^n
\]
for all $n>n_j$ and all $t\in [s',r]$. In particular,
$ 
\left|\overline{y_n^{(-j)}}\right|\le \left(p^{-j}\varepsilon^n\right)^{\frac{1}{s'}}
$ 
for all $n>n_j$. Take 
\[
\alpha_j\coloneqq \lfloor\frac{j}{d'}\rfloor=\lfloor\frac{j}{ds'}\cdot e_K\cdot \frac{p-1}{p}\rfloor,
\]
and let 
\[
z'_j\coloneqq z_j\cdot \pi_K^{-d\alpha_j}\cdot p^{d'\alpha_j} = p^{d'\alpha_j-j}\sum_{n>n_j} y_n^{(-j)}\pi_K^{-d\alpha_j}.
\]
We need to check that $\lambda_t(z'_j) \rightarrow 0$ as $j\rightarrow +\infty$ for all $t\in [s,r]$. We estimate $\lambda_t\left(y_n^{(-j)}\pi_K^{-d\alpha_j}\right)$ when $n > n_j$. By assumption, both $y_n^{(-j)}$ and $\pi_K^{-1}$ satisfy the condition in Lemma~\ref{lem: norm-product-of-reductions}; namely, 
\[
\lambda_{r_0}\left(y_n^{(-j)}\right) = \left|\overline{y_n^{(-j)}}\right|^{r_0}, \quad \lambda_{r_0}(\pi_K^{-1}) = \left|\overline{\pi}_K^{-1}\right|^{r_0}, \quad \text{and} \quad \left|\overline{y_n^{(-j)}} \cdot \overline{\pi}_K^{-1}\right| = \left|\overline{y_n^{(-j)}}\right| \cdot \left|\overline{\pi}_K^{-1}\right|.
\]
Applying Lemma~\ref{lem: norm-product-of-reductions}, we obtain 
\[
\lambda_t\left(y_n^{(-j)}\cdot \pi_K^{-d\alpha_j}\right)=\left|\overline{y_n^{(-j)}}\right|^t\cdot \left|\overline{\pi}_K\right|^{-d\alpha_j\cdot t}\le \left(p^{-j}\varepsilon^n\right)^{\frac{t}{s'}}\cdot p^{\frac{p}{p-1}\cdot\frac{1}{e_K}\cdot d\alpha_j\cdot t}
\]
for any $t\in [s,r]$ and $n > n_j$. Consequently, 
\[
\lambda_t(z'_j)\le \max_{n>n_j}\left\{p^{-(d'\alpha_j-j)}\cdot \left(p^{-j}\varepsilon^n\right)^{\frac{t}{s'}}\cdot p^{\frac{p}{p-1}\cdot\frac{1}{e_K}\cdot d\alpha_j\cdot t}\right\} = \max_{n>n_j}\left\{p^{(j-d'\alpha_j)(1-\frac{t}{s'})}\cdot \varepsilon^{n\cdot \frac{t}{s'}}\right\}.
\]
Note that $0\le j-d'\alpha_j \le d'-1$ and $1-\frac{t}{s'}\le 1-\frac{s}{s'}$, and that $\varepsilon^{n\cdot \frac{t}{s'}}<\varepsilon^{n_j\cdot \frac{t}{s'}}<\varepsilon^{n_j\cdot \frac{s}{s'}}\rightarrow 0$ as $j\rightarrow +\infty$. Hence, for any $t\in [s,r]$, $\lambda_t(z'_j)\rightarrow 0$ as $j\rightarrow +\infty$.
\end{proof}

\begin{corollary} \label{cor: C^I-flatness} \hfill
\begin{enumerate}
\item Let $I' \subset I\subset (0,r_0)$ be two closed intervals with all endpoints in $\mathbb{Q}$. Then the inclusion $\fC_{\Lambda}^I \rightarrow \fC_{\Lambda}^{I'}$ is flat.

\item Let $0<s\le r<r_0$ with $s,r\in \mathbb{Q}$. Then the map $\fC_{\Lambda}^r \rightarrow \mathbf{C}_{\Lambda}^{[s,r]}$ is flat.
\end{enumerate}
\end{corollary}

\begin{proof}
For (1), write $I = [s,r]$ and $I' = [s', r']$. The map  $\fC_{\Lambda}^I \rightarrow \fC_{\Lambda}^{I'}$ factors as $\mathbf{C}_{\Lambda}^{[s, r]} \rightarrow \fC_{\Lambda}^{[s, r']} \rightarrow \fC_{\Lambda}^{[s', r']}$. Both maps are flat by Proposition~\ref{prop: rational localizations C^r} and \cite[\S 7.3.2, Corollary~6]{BGR}. Part (2) follows from part (1), Lemma~\ref{lem: C breve dense in C tilde}(1), and \cite[Remark~3.2]{ST}.
\end{proof}

Next, we put everything in the language of rigid analytic geometry. For any closed interval $I=[s,r]\subset (0,r_0)$ with $s,r\in \mathbb{Q}$, we consider the rigid analytic variety $\mathrm{Sp}\, \fC_{\Lambda}^I$ over $K'_0$. For any two such intervals $I$, $I'$ with $I\cap I'\neq \emptyset$, we have 
\[
\mathrm{Sp}\, \fC_{\Lambda}^I\cap \mathrm{Sp}\, \fC_{\Lambda}^{I'}=\mathrm{Sp}\, \fC_{\Lambda}^{I\cap I'}
\]
and 
\[
\mathrm{Sp}\, \fC_{\Lambda}^I\cup \mathrm{Sp}\, \fC_{\Lambda}^{I'}=\mathrm{Sp}\, \fC_{\Lambda}^{I\cup I'}.
\]
Furthermore, if $I=\bigcup_{j=1}^m I_j$ where $I_1, \ldots, I_m$ are closed intervals with all endpoints in $\mathbb{Q}$, then $\left\{\mathrm{Sp}\, \fC_{\Lambda}^{I_j}\right\}_{j=1}^m$ is an admissible covering of $\mathrm{Sp}\, \fC_{\Lambda}^I$.

\begin{lemma} \label{lemma: C^I faithfully flat}
Let $I=\bigcup_{j=1}^m I_j$ be as above. Then the map $\fC_{\Lambda}^I \rightarrow \bigoplus_{j=1}^m \fC_{\Lambda}^{I_j}$ is faithfully flat.
\end{lemma}

\begin{proof}
The map is flat by Corollary~\ref{cor: C^I-flatness}, so it suffices to check that the image of 
\[
\spec\left(\bigoplus_{j=1}^m \fC_{\Lambda}^{I_j}\right)\rightarrow \spec \fC_{\Lambda}^I
\]
contains $\mathrm{MaxSpec}\, \fC_{\Lambda}^I=\mathrm{Sp}\, \fC_{\Lambda}^I$. This is clear since the image contains $\bigcup_{j=1}^m \mathrm{Sp}\, \fC_{\Lambda}^{I_j}=\mathrm{Sp}\, \fC_{\Lambda}^I$.
\end{proof}

By Kiehl's result on coherent sheaves (cf. \cite[\S 9.4]{BGR}), there is an equivalence of categories
\[
\{\textrm{finitely generated }\fC_{\Lambda}^I\textrm{-modules}\}\xrightarrow[]{\sim}\{\textrm{coherent sheaves on }\mathrm{Sp}\, \fC_{\Lambda}^I\}
\]
sending $M$ to $\widetilde{M}\coloneqq M\otimes_{\fC_{\Lambda}^I} \mathcal{O}_{\mathrm{Sp}\, \fC_{\Lambda}^I}$. It restricts to an equivalence of categories
\[
\{\textrm{finite projective }\fC_{\Lambda}^I\textrm{-modules}\}\xrightarrow[]{\sim}\{\textrm{vector bundles on }\mathrm{Sp}\, \fC_{\Lambda}^I\}.
\]

\subsection{Quasi-Stein spaces}

Let $r\in (0,r_0)$ and pick a decreasing sequence $s_0>s_1>\cdots$ of rational numbers in $(0,r]$ such that $s_i\rightarrow 0$ as $i\rightarrow \infty$. This gives an increasing sequence of closed intervals 
$ 
I_0=[s_0, r] \subset I_1=[s_1, r]\subset\cdots
$ 
whose union is $(0,r]$. Consider the rigid analytic space obtained as the union of the increasing sequence of rigid analytic spaces
\[
\mathrm{Sp}\, \fC_{\Lambda}^{I_0}\subset \mathrm{Sp}\, \fC_{\Lambda}^{I_1}\subset \cdots
\]
which is a \emph{quasi-Stein space} in the sense of \cite[Definition~2.1.4]{Bellovin}. By abusing the notation, we denote this union rigid analytic space as $\mathrm{Sp}\, \fC_{\Lambda}^r$. Note that $\mathrm{Sp}\, \fC_{\Lambda}^r$ does not depend on the choice of the $s_i$'s. 

In what follows, we consider coherent sheaves on the quasi-Stein space $\mathrm{Sp}\, \fC_{\Lambda}^r$. By Kiehl's result, giving a coherent sheaf (resp. vector bundle) $\mathcal{M}$ on $\mathrm{Sp}\, \fC_{\Lambda}^r$ is equivalent to giving a collection of finitely generated (resp. finite projective) $\fC_{\Lambda}^{I_i}$-modules $M_i$ together with isomorphisms
\[
M_j\otimes_{\fC_{\Lambda}^{I_j}} \fC_{\Lambda}^{I_i}\xrightarrow[]{\sim} M_i
\]
for each $i < j$, satisfying the cocycle conditions. We refer to $\{M_i\}_{i\ge 0}$ as a \emph{presentation} of the coherent sheaf (resp. vector bundle) $\mathcal{M}$. Note that the module of global sections $M=\Gamma(\mathrm{Sp}\, \fC_{\Lambda}^r, \mathcal{M})$ coincides with the $\fC_{\Lambda}^r$-module $\varprojlim_{i\ge 0} M_i$. We equip $M$ with the inverse limit topology from the isomorphism $M\cong \varprojlim_{i\ge 0} M_i$, where each finitely generated $\fC_{\Lambda}^{I_i}$-module $M_i$ is equipped with the finite-module topology. This topology on $M$ is independent of the presentation $\{M_i\}_{i\ge 0}$.

\begin{remark}
Since $\fC_{\Lambda}^{I_j} \rightarrow \fC_{\Lambda}^{I_i}$ is flat with dense image for each $i<j$, the inverse limit $\fC_{\Lambda}^r\cong \varprojlim_{i} \fC_{\Lambda}^{I_i}$ is a \emph{Fr\'echet-Stein algebra} in the sense of \cite[\S 3]{ST}. In the language of \emph{loc. cit.}, the collection $\{M_i\}_{i\ge 0}$ as above is called a ``coherent sheaf'' over $\fC_{\Lambda}^r$, and the $\fC_{\Lambda}^r$-module $M\coloneqq \varprojlim_{i\ge 0} M_i$ is called the ``module of global sections'' of $\{M_i\}_{i\ge 0}$. Moreover, a $\fC_{\Lambda}^r$-module that arises in this way is called \emph{coadmissible}; we will adopt this terminology in the statement of Proposition~\ref{prop: uniformly finitely generated coherent sheaves} below.
\end{remark}

The following facts are proved in \cite[\S 3]{ST}.

\begin{lemma} \label{lemma: tech lemma quasi-Stein}
Let $\mathcal{M}$ be a coherent sheaf over $\mathrm{Sp}\, \fC_{\Lambda}^r$. Let $\{M_i\}_{i\ge 0}$ be a presentation of $\mathcal{M}$, and let $M = \varprojlim_{i\ge 0} M_i$.
\begin{enumerate}
\item For each $i$, if we equip the finitely generated $\fC_{\Lambda}^{I_i}$-module $M_i$ with the finite-module topology, then the map $M\rightarrow M_i$ has dense image.

\item $R^1\varprojlim_{i\ge 0} M_i=0$.

\item For each $i$, the map $M\otimes_{\fC_{\Lambda}^r} \fC_{\Lambda}^{I_i}\rightarrow M_i$ is an isomorphism of $\fC_{\Lambda}^{I_i}$-modules.
\end{enumerate}
\end{lemma}

\begin{proof}
See \cite[\S 3 Theorem, Corollary~3.1]{ST}.    
\end{proof}

\subsection{Uniform finite generation}

In general, the module of global sections is not necessarily finitely generated. We establish a criterion that ensures finite generation.

\begin{definition} \label{def:uniform_finite_generation_over_C^r}
A coherent sheaf $\mathcal{M}$ over $\mathrm{Sp}\, \fC_{\Lambda}^r$ is called \emph{uniformly finitely generated} if there exists a positive integer $N$ and a covering $(0,r] = \bigcup_{t=1}^{\infty} J_t$ by closed intervals $J_t$ (with rational endpoints) such that each $\fC_{\Lambda}^{J_t}$-module $\Gamma(\mathrm{Sp}\, \fC_{\Lambda}^{J_t}, \mathcal{M})$ is generated by at most $N$ elements.
\end{definition}

\begin{proposition} \label{prop: uniformly finitely generated coherent sheaves}
Taking the module of global sections induces an equivalence of categories
\[\left\{\textrm{uniformly finitely generated coherent sheaves over }\mathrm{Sp}\, \fC_{\Lambda}^r\right\}\xrightarrow[]{\sim}\left\{\textrm{finitely generated coadmissible }\fC_{\Lambda}^r\textrm{-modules}\right\}\]
with the inverse functor given by $M\mapsto M\otimes_{\fC_{\Lambda}^r} \mathcal{O}_{\mathrm{Sp}\, \fC_{\Lambda}^r}$.
\end{proposition}

\begin{proof}
Our argument is partially inspired by the proof of \cite[Proposition~2.1.13]{KPX}. Let $\mathcal{M}$ be a uniformly finitely generated coherent sheaf over $\mathrm{Sp}\, \fC_{\Lambda}^r$ with module of global sections $M$. By Lemma~\ref{lemma: tech lemma quasi-Stein}(3), it suffices to show that $M$ is a finitely generated $\fC_{\Lambda}^r$-module. By re-ordering and shrinking the closed intervals $J_t$'s if necessary, we may assume there is a positive integer $N$ and a covering $(0,r] = \bigcup_{t=1}^{\infty} J_t$ where
\begin{itemize}
\item $J_t=[s_t, r_t]$ with $0<s_t\le r_t\le r$ and $s_t, r_t\in \mathbb{Q}$;
\item $r_1=r$;
\item $s_{t+1}\le r_{t+2}<s_t\le r_{t+1}$ for all $t\ge 1$;
\item $r_t, s_t\rightarrow 0$ as $t\rightarrow \infty$;
\item each $\fC_{\Lambda}^{J_t}$-module $M_{(t)}\coloneqq \Gamma(\mathrm{Sp}\, \fC_{\Lambda}^{J_t}, \mathcal{M})$ is generated by at most $N$ elements (by Lemma~\ref{lemma: tech lemma quasi-Stein}(3), we have $M_{(t)}\cong M\otimes_{\fC_{\Lambda}^r} \fC_{\Lambda}^{J_t}$).
\end{itemize}
In particular, the intervals $J_1, J_3, J_5,\ldots$ (resp., $J_2, J_4, J_6,\ldots$) are pairwise disjoint.

We claim that there exist global sections $f_1^{\mathrm{even}}, \ldots, f_N^{\mathrm{even}}\in M$ (resp., $f_1^{\mathrm{odd}}, \ldots, f_N^{\mathrm{odd}}\in M$) which generate $M_{(t)}$ for all even $t$ (resp., odd $t$). 

For each $t$, fix a finite-module norm $|\cdot|_t$ on $M_{(t)}$. Since $\fC_{\Lambda}^{J_t}$ is noetherian, $|\cdot|_t$ is a Banach norm. By assumption, for each $t$, the $\fC_{\Lambda}^{J_t}$-module $M_{(t)}$ is generated by some $f_{t,1}, \ldots, f_{t,N}\in M_{(t)}$. Since $M$ is dense in $M_{(t)}$ by Lemma~\ref{lemma: tech lemma quasi-Stein}(1), applying \cite[Lemma~2.1.12]{KPX}, we may assume that $f_{t,1}, \ldots, f_{t,N}\in M$ for all $t$. Also by \textit{loc. cit.}, for each $t$, there exists $\varepsilon_t > 0$ such that the following holds: for any $f'_1, \ldots, f'_N\in M_{(t)}$, if $|f'_n-f_{t,n}|_t<\varepsilon_t$ for all $n=1, \ldots, N$, then $f'_1, \ldots, f'_N$ also generate $M_{(t)}$.

We prove the claim in the even case; the proof for the odd case is identical. We will inductively construct $\alpha_{2t, 1}, \ldots, \alpha_{2t, N}\in \mathbb{Q}_p$ and $m_{2t, 1}, \ldots, m_{2t, N}\in \mathbb{Z}_{\ge 0}$ such that
\begin{itemize}
\item for each $n=1, \ldots, N$, the sum $\displaystyle f_n^{\mathrm{even}}\coloneqq \sum_{t=1}^{\infty} \alpha_{2t, n}\pi_K^{m_{2t, n}}f_{2t, n}$ converges in $M$;

\item $f_1^{\mathrm{even}}, \ldots, f_N^{\mathrm{even}}$ generate $M_{(2t)}$ for all $t\ge 1$.
\end{itemize}
For $t=1$, we put $\alpha_{2,1}=\ldots=\alpha_{2,N}=1$ and $m_{2, 1}=\ldots, m_{2,N}=0$. Then choose $\alpha_{2t, n}$'s and $m_{2t, n}$'s satisfying
\begin{enumerate}
\item[(i)] $m_{2t, n}\ge m_{2t-2, n}$ for all $n=1, \ldots, N$;
\item[(ii)] $\displaystyle \left|\frac{\alpha_{2t, n}\pi_K^{m_{2t, n}}}{\alpha_{2t', n}\pi_K^{m_{2t', n}}}\cdot f_{2t,n}\right|_{2t'}< \varepsilon_{2t'}$ for all $t'<t$ and $n=1,\ldots, N$;
\item[(iii)] $\displaystyle \left|\frac{\alpha_{2t', n}\pi_K^{m_{2t', n}}}{\alpha_{2t, n}\pi_K^{m_{2t, n}}}\cdot f_{2t',n}\right|_{2t}< \varepsilon_{2t}$ for all $t'<t$ and $n=1,\ldots, N$;
\item[(iv)] $\left|\alpha_{2t, n}\pi_K^{m_{2t, n}}f_{2t, n}\right|_{t''}\le p^{-t}$ for all $t''<2t-1$ and $n=1,\ldots, N$.
\end{enumerate}
Indeed, these conditions can be achieved simultaneously. Let $\omega\coloneqq p^{-\frac{p}{p-1}\cdot\frac{1}{e_K}}$. We need 
\begin{itemize}
\item for all $t'<t$ and $1 \le n \le N$, 
\[
|\alpha_{2t, n}|\cdot \omega^{s_{2t'}m_{2t, n}}< \omega^{s_{2t'}m_{2t',n}}\left|\frac{f_{2t, n}}{\alpha_{2t',n}}\right|^{-1}_{2t'}\cdot \varepsilon_{2t'};
\]
\item for all $t'<t$ and $1 \le n \le N$, 
\[
|\alpha_{2t, n}|\cdot \omega^{r_{2t}m_{2t, n}}> \omega^{r_{2t}m_{2t', n}}|\alpha_{2t',n}f_{2t',n}|_{2t}\cdot \varepsilon_{2t}^{-1};
\]
\item for all $t''<2t-1$ and $1 \le n \le N$, 
\[
|\alpha_{2t, n}|\cdot \omega^{s_{t''}m_{2t, n}}\le |f_{2t, n}|^{-1}_{t''}\cdot p^{-t}.
\]
\end{itemize}
Note that $r_{2t} < s_{2t'}$ for all $t'<t$ and $r_{2t}<s_{t''}$ for all $t''<2t-1$. Hence, by making $m_{2t, n}-m_{2t',n}$ sufficiently large, we can always find $\alpha_{2t, n}\in \mathbb{Q}_p$ satisfying all the conditions above. 

The conditions (ii)-(iv) guarantee that $f_n^{\mathrm{even}}=\sum_{t=1}^{\infty} \alpha_{2t, n}\pi_K^{m_{2t, n}}f_{2t, n}$ converges in $M$, and that 
\[
\left|f_{2t, n}-\frac{1}{\alpha_{2t, n}\pi_K^{m_{2t, n}}}\cdot f_n^{\mathrm{even}}\right|_{2t}<\varepsilon_{2t}
\]
for all $n=1, \ldots, N$. Thus, $f_1^{\mathrm{even}}, \ldots, f_N^{\mathrm{even}}$ generate $M_{(2t)}$ for all $t$, and the claim holds.

Now, the $2N$ global sections $f_1^{\mathrm{even}}, \ldots, f_N^{\mathrm{even}}, f_1^{\mathrm{odd}}, \ldots, f_N^{\mathrm{odd}}$ generate $M_{(t)}$ for all $t\ge 1$. For each $i\ge 1$, let $I_i\coloneqq \bigcup_{t=1}^{i}J_t=[s_i, r]$ and $M_i\coloneqq \Gamma(\mathrm{Sp}\, \fC_{\Lambda}^{I_i}, \mathcal{M})$. By Lemma~\ref{lemma: C^I faithfully flat}, $f_1^{\mathrm{even}}, \ldots, f_N^{\mathrm{even}}, f_1^{\mathrm{odd}}, \ldots, f_N^{\mathrm{odd}}$ also generate the $\fC_{\Lambda}^{I_i}$-module $M_i$. Let $K_i$ be the kernel of the induced surjection $(\fC_{\Lambda}^{I_i})^{2N} \twoheadrightarrow M_i$. Then $\{K_i\}_{i \ge 1}$ gives a coherent sheaf $\mathcal{K}$ over $\mathrm{Sp}\, \fC_{\Lambda}^r$. Applying Lemma~\ref{lemma: tech lemma quasi-Stein}(2) to $\mathcal{K}$, we deduce that the induced map $(\fC_{\Lambda}^r)^{2N} \rightarrow M$ is surjective.
\end{proof}

\begin{corollary} \label{cor: fin-proj-vec-bun-Cr}
Taking the module of global sections induces an equivalence of categories
\[\left\{ \begin{matrix} \textrm{uniformly finitely generated} \\  \textrm{vector bundles over }\mathrm{Sp}\, \fC_{\Lambda}^r \end{matrix}  \right\} \isom \left\{\textrm{finite projective }\fC_{\Lambda}^r\textrm{-modules}\right\}\]
with the inverse functor given by $M\mapsto M\otimes_{\fC_{\Lambda}^r} \mathcal{O}_{\mathrm{Sp}\, \fC_{\Lambda}^r}$.
\end{corollary}

\begin{proof}
We follow a similar argument as in the proof of \cite[Corollary~2.6.8]{Kedlaya-LiuII}, and provide here more details. One can see easily that any finite projective $\fC_{\Lambda}^r$-module is coadmissible (by considering an idempotent matrix). So it suffices to show that, for any uniformly finitely generated vector bundle $\mathcal{M}$ over $\mathrm{Sp}\, \fC_{\Lambda}^r$, its module of global sections $M$ is a finite projective $\fC_{\Lambda}^r$-module. We already know from Proposition~\ref{prop: uniformly finitely generated coherent sheaves} that $M$ is finitely generated.

To show projectivity, pick a surjection $\delta\colon F\coloneqq (\fC_{\Lambda}^r)^{\oplus m}\twoheadrightarrow M$ of $\fC_{\Lambda}^r$-modules, which gives rise to a morphism of vector bundles $\delta\colon \mathcal{F}\rightarrow \mathcal{M}$ over $\mathrm{Sp}\, \fC_{\Lambda}^r$. Choose an admissible covering $\mathrm{Sp}\, \fC_{\Lambda}^r = \bigcup_{i=0}^{\infty}\mathrm{Sp}\, \fC_{\Lambda}^{I_i}$ as before, and let $\{F_i\}_{i\ge 0}$ and $\{M_i\}_{i\ge 0}$ be the corresponding presentations of $\mathcal{F}$ and $\mathcal{M}$. By Lemma~\ref{lemma: tech lemma quasi-Stein}(3), $\delta(\mathrm{Sp}\, \fC_{\Lambda}^{I_i})\colon F_i\rightarrow M_i$ is surjective, and so $\delta\colon \mathcal{F}\rightarrow \mathcal{M}$ is surjective. Let $\mathcal{N}$ (resp. $N$, $N_i$) be the kernel of $\mathcal{F}\twoheadrightarrow \mathcal{M}$ (resp. $F\twoheadrightarrow M$, $F_i\twoheadrightarrow M_i$). Since $M_i$'s are finite projective by assumption, each exact sequence $0\rightarrow N_i\rightarrow F_i\rightarrow M_i\rightarrow 0$ splits. In particular, $\mathcal{N}$ is a vector bundle over $\mathrm{Sp}\, \fC_{\Lambda}^r$ with presentation $\{N_i\}_{i\ge 0}$, and its module of global sections is equal to $N\cong \varprojlim_{i\ge 0}N_i$.

We want to show that the exact sequence $0\rightarrow \mathcal{N}\rightarrow \mathcal{F}\rightarrow \mathcal{M}\rightarrow 0$ is split (so $0\rightarrow N\rightarrow F\rightarrow M\rightarrow 0$ splits). The splittings of $0\rightarrow N_i\rightarrow F_i\rightarrow M_i\rightarrow 0$ define a class in $H^1(\mathrm{Sp}\, \fC_{\Lambda}^r, \mathcal{M}^{\vee}\otimes \mathcal{N})$, and it suffices to show $H^1(\mathrm{Sp}\, \fC_{\Lambda}^r, \mathcal{M}^{\vee}\otimes \mathcal{N})=0$.

In fact, we claim that $H^1(\mathrm{Sp}\, \fC_{\Lambda}^r, \mathcal{G})=0$ for all coherent sheaves $\mathcal{G}$ over $\mathrm{Sp}\, \fC_{\Lambda}^r$. Write $\mathfrak{U}$ for the admissible covering $\mathrm{Sp}\, \fC_{\Lambda}^r = \bigcup_{i=0}^{\infty}\mathrm{Sp}\, \fC_{\Lambda}^{I_i}$, and let $G_i\coloneqq \mathcal{G}(\mathrm{Sp}\, \fC_{\Lambda}^{I_i})$. Since $\mathcal{G}$ is acyclic on each $\mathrm{Sp}\, \fC_{\Lambda}^{I_i}$, the cohomology $H^1(\mathrm{Sp}\, \fC_{\Lambda}^r, \mathcal{G})$ is computed by the \v{C}ech cohomology group $\check{H}^1(\mathfrak{U}, \mathcal{G})$. We have $\check{H}^1(\mathfrak{U}, \mathcal{G})=0$, as $R^1\varprojlim_{i\ge 0} G_i=0$ by Lemma~\ref{lemma: tech lemma quasi-Stein}(2).
\end{proof}

\vspace{0.3in}

\bibliographystyle{amsalpha}
\bibliography{library}

\vspace{0.1in}

\begin{tabular}{l}
Hansheng Diao\\
Yau Mathematical Sciences Center \& Department of Mathematical Sciences, 
Tsinghua University   \\
Beijing, China\\
\textit{E-mail address: }\texttt{hdiao@mail.tsinghua.edu.cn }\\
\\ 
Yong Suk Moon\\
Beijing Institute of Mathematical Sciences and Applications\\
Beijing, China\\
\textit{E-mail address: }\texttt{ysmoon@bimsa.cn }\\
\\
Zijian Yao\\
Department of Mathematics, University of California Santa Barbara\\
Santa Barbara, California, USA\\
\textit{E-mail address: }\texttt{yao@math.ucsb.edu }
\end{tabular}

\end{document}